\documentclass[11pt]{article}
\usepackage{amsfonts,mathrsfs}
\usepackage{graphicx}
\usepackage{amsmath,mathabx}
\usepackage{amsthm}
\usepackage{mathtools}
\usepackage{amssymb}
\usepackage{xparse}
\usepackage{xcolor}
\usepackage{float}
\usepackage{tikz}

\usepackage{amsfonts}
\usepackage{dsfont}

\usepackage{tabulary}
\usepackage{booktabs}

\definecolor{taborange}{RGB}{150,55,0}

\definecolor{tabviolet}{RGB}{148,103,189}

\usepackage[colorlinks,citecolor=blue]{hyperref}

\newtheorem{theorem}{Theorem}

\newtheorem{lemma}[theorem]{Lemma}
\newtheorem{proposition}[theorem]{Proposition}
\newtheorem{corollary}[theorem]{Corollary}

\newtheorem{claim}[theorem]{Claim}

\theoremstyle{definition} 

\newtheorem{remark}[theorem]{Remark}

\numberwithin{theorem}{section}
\numberwithin{equation}{section}

\newcommand{\T}{{\mathbb T}}
\newcommand{\R}{{\mathbb R}}
\newcommand{\Q}{{\mathbb Q}}
\newcommand{\C}{{\mathbb C}}
\newcommand{\N}{{\mathbb N}}
\newcommand{\D}{{\mathbb D}}
\newcommand{\bj}{{\boldsymbol{{\mathsf j}}}}

\newcommand{\bn}{{\boldsymbol{\mathsf n}}}
\newcommand{\bk}{{\boldsymbol{\mathsf k}}}
\newcommand{\bl}{{\boldsymbol{\mathsf l}}}
\newcommand{\bq}{{\boldsymbol{\mathsf q}}}
\newcommand{\bt}{{\boldsymbol{\mathsf t}}}
\newcommand{\bu}{{\boldsymbol{\mathsf u}}}
\newcommand{\bv}{{\boldsymbol{\mathsf v}}}
\newcommand{\bm}{{\boldsymbol{\mathsf m}}}
\newcommand{\bT}{{\boldsymbol{\mathsf T}}}

\newcommand{\bs}{{\boldsymbol{\mathsf s}}}

\newcommand{\dif}{\mathrm{d}}

\newcommand{\Exp}{\mathbb{E}}
\newcommand{\E}{\mathbb{E}}
\newcommand{\prob}{\mathbb{P}}
\newcommand{\filt}{\mathscr{F}}

\newcommand{\Prob}{\prob}
\renewcommand{\Pr}{\prob}
\DeclareDocumentCommand \one { o }
{%
\IfNoValueTF {#1}
{\mathbf{1}  }
{\mathbf{1}\left\{ {#1} \right\} }%
}

\newcommand{\lawof}{\mathscr{L}}

\newcommand{\Beta}{\operatorname{Beta}}
\newcommand{\Exponential}{\operatorname{Exp}}

\newcommand{\Gam}{\operatorname{Gamma}}

\newcommand{\Unif}{\operatorname{Unif}}
\newcommand{\As}{\operatorname{a.s.}}

\newcommand{\lawequals}{\overset{\mathscr{L}}{=}}

\DeclareDocumentCommand{\Prto} {o O{\infty} } {
\IfNoValueTF {#1}
{\overset{\Pr}{\longrightarrow}}
{ \xrightarrow[ #1 \to #2]{\Pr }}
}
\DeclareDocumentCommand{\Asto} {o O{\infty}} {
\IfNoValueTF {#1}
{\overset{\operatorname{a.s.}}{\longrightarrow}}
{
\xrightarrow[ #1 \to #2]{\operatorname{a.s.} }
}
}

\DeclareDocumentCommand{\Lpto} {o O{\infty} G{2}} {
\IfNoValueTF {#1}
{\overset{\operatorname{L}^{#3}}{\longrightarrow}}
{
\xrightarrow[ #1 \to #2]{\operatorname{L}^{#3} }
}
}

\DeclareDocumentCommand{\Mgfto} {o O{\infty}} {
\IfNoValueTF {#1}
{\overset{\operatorname{mgf}}{\longrightarrow}}
{ \xrightarrow[ #1 \to #2]{\operatorname{mgf} }}
}

\DeclareDocumentCommand{\Wkto} {o O{\infty}} {
\IfNoValueTF {#1}
{\overset{(d)}{\longrightarrow}}
{ \xrightarrow[ #1 \to #2]{(d) }}
}

\newcommand{\lgr}{\lg_r}

\author{Gaultier Lambert\thanks{Department of Mathematics, KTH Stockholm, Sweden.  {\tt email:glambert@kth.se}}, 
Joseph Najnudel\thanks{Department of Mathematics, University of Bristol, UK. \tt{joseph.najnudel@bristol.ac.uk}}, 
Elliot Paquette\thanks{Department of Mathematics, McGill University, Montreal, Canada.  {\tt email:elliot.paquette@mcgill.ca}.}\; and Ofer Zeitouni\thanks{Faculty of Mathematics and Computer Science, Weizmann Institute of Sciences, Rehovot, Israel.  {\tt email:ofer.zeitouni@weizmann.ac.il}.}}
\title{The mass of the critical C$\beta$E chaos}

\begin{document}
\maketitle
\begin{abstract}
  We study the maximum of the logarithm of the modulus of the characteristic polynomial of the circular beta ensemble C$\beta$E, taken over the unit circle. In earlier work, it was shown that this maximum, after subtracting an explicit deterministic centering, converges in distribution. The limit is the sum of three terms: a deterministic constant, a scaled Gumbel random variable, and an independent random variable arising from a critical multiplicative chaos (a derivative martingale) built from the characteristic polynomial. In this paper we identify the law of this critical chaos: its logarithm is a Gumbel random variable plus an explicit constant. We do this by showing that the critical chaos coming from the C$\beta$E is a deterministic multiple of the critical Gaussian multiplicative chaos on the circle, whose total mass has a known law. Consequently, the limiting law of the recentered maximum is a shifted and scaled sum of two independent Gumbel variables. Similar results hold for the imaginary part of the logarithm of the characteristic polynomial. This completes the proof of the random matrix part of the Fyodorov-Hiary-Keating conjecture.

  \end{abstract}
\section{Introduction}
The Circular-$\beta$ ensemble (C$\beta$E) is a distribution on $n$ points $(e^{i\omega_1},e^{i\omega_2}, \dots, e^{i\omega_n})$ on the unit circle with a joint density given by
\begin{equation}
\label{eq-density}
\frac{1}{Z_{n,\beta}} \prod_{1 \leq j < k \leq n} | e^{i \omega_j} - e^{i \omega_k} |^{\beta} \dif \omega_1 \cdots \dif \omega_n.
\end{equation}
In the special case of $\beta=2$ this is the joint distribution of eigenvalues of a Haar-distributed unitary random matrix.  The characteristic polynomial $X_n(z) \coloneqq \prod_{j=1}^n (1-e^{i\omega_j} z)$ of the C$\beta$E has attracted considerable interest, for its connections to the theories of logarithmically--correlated fields and (when $\beta=2$) analytic number theory.

A particular quantity of interest is 
\begin{equation}
\label{eq-280126}
M_n\coloneqq\max_{|z|=1}
\log |X_n(z)|.
\end{equation} Let
\begin{equation}
\label{eq-mn}
m_n = \log n - \tfrac{3}{4}\log\log n.
\end{equation}
The random matrix
part of the Fyodorov--Hiary--Keating conjecture \cite{FHK12}
states that in the special case that $\beta=2$, $M_n-m_n$ converges in distribution 
towards a limiting random variable ${\bf R}_2$, with
\begin{equation}
\label{eq-R2}
P({\bf R}_2\in dx)=4e^{2x} K_0(2e^{x})dx,
\end{equation}
where $K_0$ is the modified Bessel function of the second kind.
It was later observed in \cite{SZ} that the probability density in \eqref{eq-R2} is the law of the sum of two independent Gumbel random variables.

For general $\beta>0$, an important step forward was obtained by
\cite{CMN}, who proved that $M_n-\sqrt{2/\beta}m_n$ is tight. 
In \cite{PaquetteZeitouni03}
this was strengthened to convergence, with a description of the limiting law.
\begin{theorem}[Theorem 1.1 of \cite{PaquetteZeitouni03}]
\label{old-main}
The sequence of random variables $M_n-\sqrt{2/\beta}m_n$ converges in distribution
to a random variable ${\bf R}_\beta$. Further,
\begin{equation}
\label{eq-convolution}
{\bf R}_\beta=
C_\beta+G_\beta+\frac{1}{\sqrt{2\beta}}\log (\mathscr{B}_\infty(\beta)) ,
\end{equation}
where $C_\beta$ is an (implicit) constant, $G_\beta$ is Gumbel distributed with
parameter $1/\sqrt{2\beta}$, and the random variable $\mathscr{B}_\infty(\beta)$ is  independent
of $G_\beta$.
\end{theorem}
The random variable $\mathscr{B}_\infty(\beta)$ is constructed in \cite{PaquetteZeitouni03} as a certain \textit{derivative martingale} (which is a \textit{critical chaos}), and we review the construction in Section \ref{sec-overview} below. 
In this article we identify its law
and hence resolve the Fyodorov--Hiary--Keating conjecture for random matrices \cite{FHK12}.
\begin{theorem}\label{thm:main}
The random variable $\log (\mathscr{B}_\infty(\beta))$ has the law
\[
\log (\mathscr{B}_\infty(\beta)) = \operatorname{Gumbel}(1) + \bar C_\beta
\]
where $\bar C_\beta=\log \sqrt{2}+\mathfrak{g}_\beta $ is an explicit constant, with  $\mathfrak{g}_\beta$ given in \cite[(9.6)]{PaquetteZeitouni03}.
\end{theorem}

\subsection{OPUC preliminaries}
A major advance in the study of $M_n$ was achieved in \cite{CMN}, who used
the Orthogonal Polynomials on the Unit Circle (OPUC) representation of the  C$\beta$E measure due to \cite{KillipNenciu}; we refer to \cite{SimonBook} for an encyclopedic account of the OPUC theory.
Let $\left\{ \gamma_k \right\}$ be independent, complex, rotationally invariant random variables for which $|\gamma_k|^2 = \Beta(1, \beta(k+1)/2),$ that is with density on $[0,1]$ proportional to $(1-x)^{\beta(k+1)/2-1}$. 
As in \cite{CMN} it is convenient to introduce a latent-variable coupling as 
\begin{equation}
\gamma_j = \sqrt{
\frac{E_j}{E_j+\Gamma_j}
} e^{i\Theta_j},
\quad
\text{where}
\quad
\begin{aligned}
&E_j \sim \Exponential(1), \\
&\Gamma_j \sim \Gam(\tfrac{\beta}{2}(j+1)), \\
&\Theta_j \sim \Unif([0,2\pi])
\end{aligned}
\label{eq:betagamma}
\end{equation}
with all variables independent.  We define the variables
\begin{equation}
X_j = \sqrt{E_j} \Re e^{i\Theta_j},
\quad
Y_j = \sqrt{E_j} \Im e^{i\Theta_j},
\quad
\text{and}
\quad
Z_j = X_j+iY_j.
\label{eq:gaussians}
\end{equation}
Then $\left\{ \{X_j,Y_j\}\!:\! j \geq 0 \right\}$ are jointly iid $N(0,\tfrac 12)$ and independent of $\left\{ \Gamma_j\! :\! j \geq 0\right\}.$   Set 
\begin{equation}
\label{eq-betaj}
\beta_j \coloneqq  \sqrt{\frac{\beta}{2}(j+1)},
\end{equation}
 noting that $\Exp \Gamma_j = \beta_j^2.$ Then, we have that $\gamma_j = Z_j/\sqrt{\Gamma_j+|Z_j|^2} \approx \sqrt{\tfrac{2}{\beta (j+1)}} Z_j$.

The \emph{Szeg\H{o}} recurrence is, for all $k \geq 0,$
\begin{equation}
\begin{pmatrix}
\Phi_{k+1}(z) \\
\Phi^*_{k+1}(z)
\end{pmatrix}
\coloneqq
\begin{pmatrix}
z & -\overline{\gamma_k} \\
-\gamma_k z & 1
\end{pmatrix}
\begin{pmatrix}
\Phi_{k}(z) \\
\Phi^*_{k}(z)
\end{pmatrix},
\quad
\biggl\{
\begin{aligned}
  &\Phi_0(z) \equiv 1, \\
  &\Phi_k^*(z) = z^k \overline{ \Phi_k(1/\overline{z})}.
\end{aligned}
\biggr\}
\label{eq:szego}
\end{equation}
where $\Phi_{k}^{*}$ and $\Phi_{k}$ are polynomials of degree at most $k$.
We introduce the factor $Q_k(z) \coloneqq z\tfrac{\Phi_k(z)}{\Phi_k^*(z)}$, which is a function of $\gamma_0,\ldots,\gamma_{k-1}$ and
can be checked to satisfy the autonomous recurrence
\begin{equation}\label{eq:Q}
Q_{k+1}(z) \coloneqq z \frac{Q_k(z) - \overline{\gamma}_{k}}{1-\gamma_{k}Q_k(z)}.
\end{equation}
Then the polynomials $\left\{ \Phi_k^* \right\}$ satisfy the recurrence
\begin{equation}\label{eq:Phi}
\log \Phi^*_{k+1}(z)=\log\Phi^*_{k}(z) +\log(1-\gamma_k Q_k(z))
,\quad \Phi_0^*(z) = 1,
\end{equation}
When $|z|=1$, the factor $Q_k$ is a complex phase, and $\theta\mapsto\arg\bigl(e^{i\theta}\Phi_k(e^{i\theta})/\Phi_k^*(e^{i\theta})\bigr)$ has a continuous version $\Psi_k(\theta)$, the \emph{Pr\"ufer phase} (see \cite[Lemma 2.3]{CMN}).  It satisfies
\begin{equation}\label{eq:Psi}
\Psi_{k+1}(\theta)=\Psi_{k}(\theta) + \theta - 2\Im\left( \log(1-\gamma_ke^{i\Psi_k(\theta)})\right)
,\quad \Psi_0(\theta) = \theta.
\end{equation}
However, in this paper we shall also need the behavior of $Q_k(z)$ for $|z| <  1$.

We recall the following basic properties of $Q_k(z)$ (see \cite[Proposition 5.1]{CN}).
\begin{lemma}\label{lem:Q}
The following holds for $Q_j(z)$, for all  $z \in \D$.
\begin{enumerate}
\item For all $j\geq 0$, $|Q_j(z)|\leq 1.$
\item If $|z|=r$, then the following recurrence holds
  \begin{equation}\label{eq:Qmod}
	|Q_{j+1}(z)|^2 = 
	r^2-r^2\frac{(1-|\gamma_{j}|^2)(1-|Q_j(z)|^2)}{|1-\gamma_{j}Q_j(z)|^2}.
  \end{equation}
\item For all $n,j\geq 0$, the following bound holds on the conditional expectation
  \[
	r^{2j} |Q_n(z)|^2 \le \Exp[ |Q_{n+j}(z)|^2 \mid \gamma_0,\ldots,\gamma_{n-1}] \leq 1.
  \]
\end{enumerate}
\end{lemma}

These polynomials $\{\Phi_{k}^{*}\}$ can be used to give an effective representation for the characteristic polynomial $X_n$ by setting $\alpha$ to be a uniformly distributed element of the unit circle, independent of $\left\{ \gamma_k : k \geq 0 \right\},$ and setting for any $\theta \in \R,$
\begin{equation}\label{eq:charpoly}
X_n(e^{i\theta})
\coloneqq \Phi_{n-1}^*(e^{i\theta}) - \alpha e^{i\theta} \Phi_{n-1}(e^{i\theta})
=\Phi_{n-1}^*(e^{i\theta})\bigl( 1-\alpha e^{i \Psi_{n-1}(\theta)}\bigr).
\end{equation}
See \cite[(2.2)]{CMN} for details of the demonstration that this indeed has the law of the characteristic polynomial of C$\beta$E.

We also recall the \emph{relative Pr\"ufer phase} \cite[Lemma 2.4]{CMN} given by the recurrence
\begin{equation}\label{eq:prufer}
\psi_{k+1}(\theta)=\psi_{k}(\theta) + \theta - 2\Im\left( \log(1-\gamma_ke^{i\psi_k(\theta)}) - \log(1-\gamma_k) \right)
,\quad \psi_0(\theta) = \theta.
\end{equation}
In law $\{\psi_{k}(\theta) : k \in \N, \theta \in [0,2\pi]  \}$ is equal to $\{\Psi_{k}(\theta)-\Psi_{k}(0) : k \in \N,\theta \in [0,2\pi] \}.$

We study the real and imaginary parts of $\log\Phi_k^*(re^{i\theta})$ through the field $\varphi_{k}(r,\theta) \coloneqq 2\Re(\sigma\log\Phi_k^*(re^{i\theta}))$, where $\sigma\in\{1,i\}$.  By \eqref{eq:Phi}, it satisfies the recurrence
\begin{equation}\label{eq:upsilon}
\varphi_{k+1}(r,\theta)=\varphi_{k}(r,\theta) + 2\Re \big\{ \sigma \left( \log(1-\gamma_k Q_k(r,\theta)) \right) \big\}
,\quad \varphi_0(r,\theta) = 0,
\end{equation}
The entire analysis of $\varphi$ occurs for a fixed $\sigma,$ and so we shall not display the dependence on $\sigma.$
We note that for fixed $k$, the map $(\theta,r)\mapsto \varphi_{k}(r,\theta)$
is a.s. continuous on $[0,2\pi]\times (0,1]$.

It is convenient to  introduce a continuous time process $(\zeta_t^\C : t \geq 0),$ which we will use for some comparisons.  As $\left\{ Z_k \right\}$ are independent standard complex Gaussians with $\Exp |Z_k|^2 = 1,$ we may enlarge the probability space to include a complex Brownian motion $(B_t : t \geq 0)$ such that
\begin{equation}\label{eq:W}
[ B_t, B_t] =0, \quad
[ B_t, \overline{B_t}] =2t
\quad
\text{for all }
t \geq 0,
\quad\text{and}\quad
\sqrt{2}Z_k = B_{k+1}-B_k
\quad
\text{for all }
k \geq 0,
\end{equation}
In terms of $B$ we define  for $t \geq 2^{\bn}$
\begin{equation}\label{eq:Z}
\zeta_t^\C(\bn;r,\theta) \coloneqq 
-
\sqrt{\frac{4}{\beta}}
\int_{2^{\bn}}^t
\frac{Q_{[s]}(r,\theta)}{\sqrt{[ s ]+1}}
dB_s
+
2\log \Phi^*_{2^{\bn}}(re^{i\theta})
\end{equation}
We also let for $\sigma \in \{1,i\}$
\begin{equation}
\label{eq-zetareal}
\zeta_t(\bn;r,\theta) \coloneqq \Re\big( \sigma \zeta_t^{\C}(\bn;r,\theta)\big) . 
\end{equation}
Then, for $k\geq2^{\bn}$, $\zeta_k(\bn;r,\theta)$ is a good approximation of $\varphi_{k}(r,\theta)$, since $\tfrac12\zeta_k^\C(\bn;r,\theta)$ is close to $\log\Phi_k^*(re^{i\theta})$; see Lemma \ref{lem:zetaphitail} and Proposition \ref{prop:zetaphitail}.
Set also
\begin{equation}\label{eq:filtration}
 \filt_k:=\sigma\!\left(\Gamma_0,\ldots,\Gamma_{k-1},
                 B_s:0\leq s\leq k\right),\qquad k\in\N.
\end{equation}
Thus $\gamma_0,\ldots,\gamma_{k-1}$ and $Q_k$ are $\filt_k$-measurable, while
$\Gamma_k$ and the path $(B_s-B_k)_{s\in[k,k+1]}$ are independent of $\filt_k$.  All conditional expectations below are taken with respect to this filtration.


\subsection{Overview of strategy and re-formulation of main result}
\label{sec-overview}
Our proof builds on ideas contained in \cite{lambert2024subcritical}, who identified the total mass of the \textit{subcritical}  multiplicative chaos associated with the C$\beta$E. We describe the connection below, after introducing some terminology.

We will need the so-called derivative martingale.
Set 
\begin{equation}
\lgr(k) \coloneqq \log (k \wedge \tfrac{1}{1-r^2}).\label{lgr}
\end{equation}
Define the random measure and its total mass
\begin{equation}\label{eq:Bk}
\mathscr{D}_k(r,\theta)
\coloneqq
\tfrac{1}{2\pi}
e^{ \sqrt{ \tfrac{\beta}{2}} \varphi_k(r,\theta) - \lgr k}
\bigl(\sqrt{2}\lgr k - \sqrt{\tfrac{\beta}{4}}\varphi_k(r,\theta) \bigr)_+
\quad
\text{and}
\quad
\mathscr{B}_k(r) \coloneqq
\int_0^{2\pi}
\mathscr{D}_k(r,\theta)
d\theta.
\end{equation}
For simplicity, we suppress the $r$ in the case that $r=1$.
In \cite{PaquetteZeitouni03} it is shown this converges to a non-degenerate random measure when $r=1$.

\begin{theorem}[Theorem 1.7 of \cite{PaquetteZeitouni03}]\label{thm:Bn}
For any $\sigma \in \left\{1, i \right\}$ and any $\beta > 0,$
there is an almost surely nonnegative, finite random variable $\mathscr{B}_\infty$ and an almost surely finite, nonatomic random measure $\mathscr{D}_\infty$ so that
\[
\mathscr{D}_{2^{\bn}}(1,\theta) \dif \theta
\Asto[\bn] \mathscr{D}_\infty(\dif \theta)
\quad
\text{and}
\quad
\mathscr{B}_{2^{\bn}}
\Asto[\bn] \mathscr{B}_\infty.
\]
Furthermore for any $\epsilon > 0$ there is a compact $K =K_\epsilon\Subset (0,\infty)$ so that, with
\[
\chi(\theta) = \one[{ (\sqrt{2}\log k - \sqrt{\tfrac{\beta}{4}}\varphi_k(\theta))/\sqrt{\log k} \not\in K }],
\]
then for any $k \in \N$
\begin{equation}\label{eq:Btight}
\Pr
\biggl(\int_0^{2\pi}
e^{ \sqrt{ \tfrac{\beta}{2}} \varphi_k(\theta) - \log k}\bigl|\sqrt{2}\log k - \sqrt{\tfrac{\beta}{4}}\varphi_k(\theta) \bigr| \chi(\theta) d\theta > \epsilon
\biggr)
<\epsilon.
\end{equation}
\end{theorem}
The appearance of the term $\sqrt{2}\log k - \sqrt{\tfrac{\beta}{4}}\varphi_k(\theta)$ in \eqref{eq:Btight} is what makes $\mathscr{D}_\infty$ into a derivative (or critical) martingale.

The overarching strategy, as pioneered in \cite{CN} and implemented in \cite{lambert2024subcritical} in the context of the subcritical chaos associated with the C$\beta$E, is to consider an  approximation to
$\mathscr{D}_\infty$ using  $r<1$ and then taking $r\to 1$. To perform this, our proof integrates additional ideas from \cite{CMN} and \cite{PaquetteZeitouni03}.  
The reason why the approximation with $r<1$ is useful is that
 if we fix $r<1$ and send $k\to\infty$, then we have from \cite[Proposition 3.1]{CN} that the process $\varphi_k(r,\theta)$ converges to a Gaussian process.
\begin{theorem}[{\cite[Proposition 3.1]{CN} or corollary of  \cite[Theorem 1]{JM15}}]\label{thm:fieldconverges}
Fix $r < 1$. Then, on sending $k\to\infty$, we have
\[
\varphi_k(r,\theta)
\Asto[k]
\varphi_\infty(r,\theta),
\]
a stationary (in $\theta$) and continuous Gaussian field having covariance 
\[
\Exp \varphi_\infty(r,\theta)\varphi_\infty(r',0)
=-\frac{4}{\beta} \log |1-rr'e^{i\theta}|.
\]
The a.s. convergence is uniform in $\theta$.
\end{theorem}
In particular, we may take $k\to\infty$ for fixed $r<1$ to arrive at
\[
\mathscr{D}_k(r,\theta)
\Asto[k]
\mathscr{D}_\infty(r,\theta)
\coloneqq
\tfrac{1}{2\pi}
e^{ \sqrt{ \tfrac{\beta}{2}} \varphi_\infty(r,\theta) + \log(1-r^2)}
\bigl(\sqrt{2}\log\tfrac{1}{1-r^2} - \sqrt{\tfrac{\beta}{4}}\varphi_\infty(r,\theta) \bigr)_+
\]
As this is a Gaussian process, existing theory is well suited to showing that this converges as a measure when $r\to 1$, and moreover converges to a specific 
(critical) Gaussian Multiplicative Chaos (CGMC), whose law (originally conjectured in \cite{FBo}) is known, see \cite{Remy} and \cite{CN}. 
\begin{remark}
\label{rem-almostgaussian}
We will see in Section \ref{sec-prelim} that the convergence (for a fixed $r$) as $k\to\infty$ can be quantified, and in fact that
$\varphi_\infty(r,\theta)-\varphi_{k}(r,\theta)$ with $k=\lceil 1/(1-r^2)\rceil$ has tails bounded uniformly in $r,\theta$ (see Corollary~\ref{cor:ezfieldincrement}),
and also that the difference $\varphi_t(r,\theta)-\varphi_t(1,\theta)$ is small when $t$ is much smaller than  $1/(1-r^2)$ (see Corollary~\ref{cor:radialincrement}).
\end{remark}

The main task is therefore, in some sense, to be able to interchange the order of limits, that is show that the known  limit when one first takes $k\to\infty$ and then $r\to 1$ is the same (as it turns out, up to a multiplicative constant) as the one where one takes $r\to 1$ followed by $k\to\infty$. The main difficulty is that the variable $\mathscr{B}_\infty$ does not possess bounded mean, and so appropriate truncations have to be introduced before one can apply standard (first/second moment) methods. 

The truncations alluded to above are given in terms of an event, which involves all radii smaller than a given threshold. Toward this end,
define
\begin{equation}
\label{eq-rk}
r_{\bk}^2 \coloneqq 1-2^{-\bk} \quad \bk \in \N,  
\end{equation}
and for  $L > 0$, $a \in (-\infty,0]$ and $\theta \in [0,2\pi]$, define the event 
\begin{equation}\label{eq:Abarrier}  
\mathscr{A}_{\bk}^{L,a}(\theta) 
\coloneqq 
\left\{
\sqrt{ \tfrac{\beta}{8}} \varphi_\infty(r_{\bj},\theta)
\leq \log \tfrac{1}{1-r_{\bj}^2}
+ a\left(\log \tfrac{1}{1-r_{\bj}^2}\right)^{1/8}
+ L
:
\,
\forall
\,
\bj \leq \bk
\right\}.
\end{equation}
When we take $a=0$, we shall suppress it from the notation.  We note that by Remark \ref{rem-almostgaussian},
the event $\mathscr{A}_{\bk}^{L,a}(\theta)$ can be thought of as a proxy for the  \textit{barrier} event 
\[\left\{
\sqrt{ \tfrac{\beta}{8}} \varphi_{2^{\bj}} (1,\theta)
\leq \log 2^{\bj}
+ a\left(\log 2^{\bj}
\right)^{1/8}
+ L
:
\,
\forall 
\,
\bj \leq \bk
\right\}.\]
This point of view is at the heart of the arguments in Section \ref{sec:proof-main-comparison} below.

The key reason for introducing the event $\mathscr{A}_{\bk}^{L,a}$ is that it allows for certain Gaussian $L^2$ computations. This is contained in the following proposition.

\begin{proposition}\label{prop:criticalchaos}
On taking $r\to 1$, we have in-probability weak convergence of measures
\[
\mathscr{D}_\infty(r,\theta)d\theta
\Prto[r][1]
\operatorname{CGMC}(\dif\theta),
\]
where $\operatorname{CGMC}(\dif\theta)$ is a critical Gaussian multiplicative chaos which is independent of $\beta$
and satisfies
\begin{equation}\label{eq:Remysmagiccarpet}
\int_0^{2\pi}\operatorname{CGMC}(\dif\theta) \lawequals \frac{\sqrt{2}}{\operatorname{Exp}(1)}.
\end{equation}
Moreover, we have that for any continuous $f : \T \to [0,1]$ and any $a \in (-\infty,0)$ and $p\in (0,2)$, 
\begin{equation}
\label{eq-140924a}
\int_0^{2\pi} 
\one_{  \mathscr{A}_{\bk}^{L,a}(\theta)  }f(\theta)
\mathscr{D}_\infty(r_{\bk},\theta) 
\dif \theta 
\Lpto[\bk]{p}
\int_0^{2\pi} 
\one_{  \mathscr{A}_\infty^{L,a}(\theta) }
f(\theta)
\operatorname{CGMC}(\dif \theta)
\Prto[L]
\int_0^{2\pi}  
f(\theta)
\operatorname{CGMC}(\dif \theta).
\end{equation}
The estimate \eqref{eq-140924a} holds also with $a=0$ and $p=2$.
\end{proposition}
\begin{proof}
The first assertion is known.  The convergence of the derivative martingale is
\cite[Theorem 4]{DRSV} (noting that the presence of the $(\cdot)_+$ in the definition of $\mathscr{D}_k(r,\theta)$  is inessential since with probability approaching $1$ as $k\to\infty$, the term enclosed in the parentheses is non-negative).  That the limit does not depend on the approximation is
\cite[Theorem 1.1]{JunnilaSaksman}, and the case of the circle, with the field approximated by its values $\varphi_\infty(r,\cdot)$ at radius $r$ (the Poisson regularization), is
\cite[Theorem 1.3 and Example 6.1]{JunnilaSaksman}.
  The equality  \eqref{eq:Remysmagiccarpet} follows from \cite[Theorem 1.1]{Remy} following the computation around \cite[(1.12)]{Remy}, together with the identification of the critical chaos as the limit of renormalized subcritical chaos \cite{APS19,PowellReview}; beware the constant in \cite[(1.12)]{Remy}.

To prove \eqref{eq-140924a}, we
 begin with the case $a=0$. In that case, by \cite[Proposition 3.1]{Lacoin}
(see also \cite[Proposition 14]{DRSV}) we have that
$\int_0^{2\pi} 
\one_{  \mathscr{A}_{\bk}^{L,0}(\theta)  }f(\theta)
\mathscr{D}_\infty(r_{\bk},\theta) 
\dif \theta $ converges in $L^2$ to an $L^2$  random variable, denoted $\mathscr{T}_f^L$.
We next show that 
\begin{equation}
\label{eq-14012026oa}\mathscr{T}_f^L=\int_0^{2\pi} 
\one_{  \mathscr{A}_\infty^{L,0}(\theta) }
f(\theta)
\operatorname{CGMC}(\dif \theta).
\end{equation}
 Indeed,
for fixed $\bk'$, set
\[
 \mathscr{A}^{\bk',L,0,\complement}
 :=\left\{\exists\,\theta\in[0,2\pi],\ \exists\,\bk>\bk':
 \mathscr A_{\bk'}^{L,0}(\theta)\text{ holds but }
 \mathscr A_{\bk}^{L,0}(\theta)\text{ fails}\right\}.
\]
This event is measurable because the fields are continuous in $\theta$.
Due to \cite[Proposition 2.4]{Lacoin}
(compare also the random-matrix upper-envelope estimate
\cite[Lemma 9.3]{PaquetteZeitouni03}),
for each $\delta>0$ there is $\bk_\delta$ so that 
$\mathbb{P}\big(\mathscr{A}^{\bk_\delta,L,0,\complement}\big)<\delta$.
On the other hand, let $\mathscr F^{\mathrm G}_{\bk'}$ denote the
sigma-field generated by the Gaussian field $\varphi_\infty(r_\bj,\cdot)$,
$\bj\leq\bk'$.  The convergence in probability
\[
\int_0^{2\pi}f(\theta)\mathscr{D}_\infty(r_{\bk},\theta)\dif\theta
\Prto[\bk]
\int_0^{2\pi}f(\theta)\operatorname{CGMC}(\dif\theta)
\]
also holds when $f$ is a bounded random function, Borel in $\theta$ and
$\mathscr F^{\mathrm G}_{\bk'}$-measurable for a fixed $\bk'$; see
\cite[Proposition B.1]{Lacoin} and \cite[Theorem 4.4]{JunnilaSaksman}.
The event $\mathscr A_{\bk_\delta}^{L,0}(\theta)$ is
$\mathscr F^{\mathrm G}_{\bk_\delta}$-measurable, so we may take $f(\theta)\one_{\mathscr A_{\bk_\delta}^{L,0}(\theta)}$ as the test function.  It follows that,
in the double limit with $\bk\to\infty$ followed by $\delta\to0$,
\[
\int_0^{2\pi} 
\one_{  \mathscr{A}_{\bk_\delta}^{L,0}(\theta)  }f(\theta)
\mathscr{D}_\infty(r_{\bk},\theta) 
\dif \theta \Prto
\int_0^{2\pi} 
\one_{  \mathscr{A}_{\infty}^{L,0}(\theta) }
f(\theta)
\operatorname{CGMC}(\dif \theta).
\]
This shows that indeed \eqref{eq-14012026oa} holds, completing the proof of the first limit in \eqref{eq-140924a} when $a=0$. The proof of the second limit in that case follows again from \cite[Proposition 2.4]{Lacoin}.

To handle the case $a<0$, note that the (positive) left side of \eqref{eq-140924a} with $a<0$ is bounded above by the same with $a=0$, and therefore the claimed $L^p$ convergence follows from convergence in probability and the uniform integrability of the $p$ power.
The main difference between the case $a=0$ and $a<0$ is that for $a=0$, \cite[Proposition 2.4]{Lacoin}  shows (by controlling the supremum of $\varphi_\infty(r,\theta)$ over $r>r_{\bk'}$ and all $\theta$)  that the set
$\bigcup_\theta\big(\bigcup_{\bk>\bk'} (\mathscr{A}_{\bk}^{L,0}(\theta))^\complement\big)$ has probability tending to $0$ as $\bk'\to \infty$.
Unfortunately, this is not true when $0$ is replaced by $a<0$ (it would be true if we used the barrier $a\log |\log (1-r)|$ with $a>-1$ instead of $a|\log(1-r)|^{1/8}$).

To overcome this, we claim that 
\begin{equation}
\label{eq-15012026b}
\limsup_{\bk'\to\infty}\limsup_{\bk>\bk'}\mathbb{E} \int_0^{2\pi} f(\theta) \one_ {\mathscr{A}_{\bk'}^{L,a}(\theta)  } \one_ {\mathscr{A}_{\bk}^{L,0}(\theta) \cap (\mathscr{A}_{\bk}^{L,a}(\theta))^\complement  }\mathscr{D}_\infty(r_{\bk},\theta) 
\dif \theta =0
\end{equation}
Indeed, this is a Gaussian computation using  the union bound
\begin{equation}
\label{eq-15012026a}
\begin{aligned}
&
2\pi\mathbb{E} \int_0^{2\pi} f(\theta) \one_ {\mathscr{A}_{\bk'}^{L,a}(\theta)  } \one_ {\mathscr{A}_{\bk}^{L,0}(\theta) \cap (\mathscr{A}_{\bk}^{L,a}(\theta))^\complement  }
\mathscr{D}_\infty(r_{\bk},\theta)
\dif \theta\\
&\leq 
\sum_{\bj=\bk'}^{\bk}  \mathbb{E} \int_0^{2\pi} \one_ {\mathscr{A}_{\bk'}^{L,a}(\theta)  } \one_ {\mathscr{A}_{\bk}^{L,0}(\theta)  }
\mathscr{D}_\infty(r_{\bk},\theta)\one_{\left\{
\sqrt{ \tfrac{\beta}{8}} \varphi_\infty(r_{\bj},\theta)
\geq \log \tfrac{1}{1-r_{\bj}^2}
+ a\left(\log \tfrac{1}{1-r_{\bj}^2}\right)^{1/8}
+ L\right\}}=:\sum_{\bj=\bk'}^\bk \mathcal{U}_{\bj,\bk,\bk'}\end{aligned}
\end{equation}
By the martingale property of $\mathscr{D}_\infty(r_{\bk},\theta)$ in $\bk$, a Girsanov change of measure turns $\mathcal{U}_{\bj,\bk,\bk'}$ into a ballot estimate for a Brownian motion, as in
\cite[Section 9, especially (9.3)]{PaquetteZeitouni03}.  The ballot density at time $\bj\log2$ carries the factor $\bj^{-3/2}$, and so
\[ \mathcal{U}_{\bj,\bk,\bk'}\leq C_1L\,\bj^{-3/2} \int_L^{|a|\bj^{1/8}+L}(L+h)^2 e^{-h^2
/2\bj}\,\dif h\leq C_2\, \bj^{3/8}\,\bj^{-3/2}= C_2\,\bj^{-9/8}
\]
for $\bj$ large, where $C_2$ depends only on $L$ and $a$.
In particular,
combining this estimate with \eqref{eq-15012026a} yields \eqref{eq-15012026b}. But then, to get the first limit in \eqref{eq-140924a} we can fix $\bk'$, replace the left expression in \eqref{eq-140924a}
by $\int_0^{2\pi} 
\one_{  \mathscr{A}_{\bk}^{L,0}(\theta)  }\one_{  \mathscr{A}_{\bk'}^{L,a}(\theta)  }f(\theta)
\mathscr{D}_\infty(r_{\bk},\theta) 
\dif \theta$ and obtain that this converges in $L^1$  (as first $\bk\to\infty$ and then $\bk'\to\infty$) to the middle  term in \eqref{eq-140924a}.
To show the second convergence in \eqref{eq-140924a}, we argue similarly, showing that 
\[\limsup_{L\to \infty}
\limsup_{\bk\to\infty}
\mathbb{E}
\int_0^{2\pi}  
\one_{ ( \mathscr{A}_{\bk}^{L,a}(\theta) )^\complement \cap \mathscr{A}_{\bk}^{L,0}(\theta)}
\mathscr{D}_\infty(r_{\bk},\theta) 
\dif \theta 
=0,\]
and then use the result for $a=0$.
\end{proof}

Let $f : \T \to [0,1]$ be a continuous function, and let 
\[
\mathscr{D}^{L,a}_{2^{\bn}}(r_{\bk},f) \coloneqq 
\int_0^{2\pi}  
\one_{  \mathscr{A}_{\bk}^{L,a}(\theta) }
f(\theta)
\mathscr{D}_{2^{\bn}}(r_{\bk},\theta)
\,\dif \theta.
\]
We note that by continuity and Theorem \ref{thm:fieldconverges}, we may take $n\to\infty$ on both sides for any fixed $r_{\bk} < 1$, and so we have 
\[
\mathscr{D}^{L,a}_{\infty}(r_{\bk},f) \coloneqq 
\int_0^{2\pi}  
\one_{  \mathscr{A}_{\bk}^{L,a}(\theta) }
f(\theta)
\mathscr{D}_{\infty}(r_{\bk},\theta)
\,\dif \theta.
\]
Set 
\begin{equation}
\label{eq-En}
\Exp_{\bn}(\cdot) \coloneqq \Exp\left( \cdot ~ \vert ~ \filt_{2^{\bn}}\right), \quad \mathbb P_{\bn}(\cdot) \coloneqq \mathbb P\left( \cdot ~ \vert ~ \filt_{2^{\bn}}\right), 
\end{equation}
  Then by the uniform integrability provided by \eqref{eq-140924a},
\[
\Exp_{\bn}\left( \mathscr{D}^{L,a}_{\infty}(r_{\bk},f) \right)
\Prto[k]
\Exp_{\bn}\left( 
\int_0^{2\pi}
\one_{  \mathscr{A}_\infty^{L,a}(\theta) }
f(\theta)
\operatorname{CGMC}(\dif \theta)
\right)
\eqqcolon
\Exp_{\bn}\left( \operatorname{CGMC}^{L,a}(f) \right)
\]

The main comparison result we show is the following.
\begin{proposition}\label{prop:main}
For any continuous $f : \T \to [0,1]$ and $a<0$,
\[
\mathscr{D}_{2^{\bn}}(1,f) - e^{\mathfrak{g}_\beta}\Exp_{\bn}\left( \operatorname{CGMC}^{L,a}(f) \right)
\Prto[L,\bn] 0,
\]
where $\mathfrak{g}_\beta$ is a deterministic constant given in \cite[(9.6)]{PaquetteZeitouni03}.
\end{proposition}
\noindent Here, for an array $X^L_\bn$ of random variables, $X^L_\bn \Prto[L,\bn] 0$ means that for all $\epsilon > 0$
\[
\lim_{L \to \infty} \limsup_{\bn \to \infty } \Pr( |X^L_\bn| > \epsilon)  = 0.
\]
As a corollary, we conclude the main theorem of this paper.
\begin{theorem}\label{thm:maincor}
Almost surely,
\[
\mathscr{D}_\infty(\dif \theta) = 
{e^{\mathfrak{g}_\beta}}\operatorname{CGMC}(\dif \theta).
\]
\end{theorem}
\begin{proof}
Fix a continuous $f : \T \to [0,1]$.
From the convergence in Theorem \ref{thm:Bn}, we have
\[
\mathscr{D}_{2^{\bn}}(1,f)
\Asto[\bn]\mathscr{D}_\infty(f).
\]
From Proposition \ref{prop:main}, we conclude that for any  $a<0$,
\[
{e^{\mathfrak{g}_\beta}}\Exp_{\bn}\left( \operatorname{CGMC}^{L,a}(f) \right)
=
\mathscr{D}_{2^{\bn}}(1,f)
+({e^{\mathfrak{g}_\beta}} \Exp_{\bn}\left( \operatorname{CGMC}^{L,a}(f) \right)-    \mathscr{D}_{2^{\bn}}(1,f)
)
\Prto[L,\bn]  \mathscr{D}_\infty(f).
\]
From martingale convergence, we have 
\[
\Exp_{\bn}\left( \operatorname{CGMC}^{L,a}(f) \right)
\Asto[n]\operatorname{CGMC}^{L,a}(f).
\]
Hence by \eqref{eq-140924a} of
Proposition \ref{prop:criticalchaos}, we conclude 
\[
\Exp_{\bn}\left( \operatorname{CGMC}^{L,a}(f) \right)
\Prto[L,\bn] \operatorname{CGMC}(f).
\]
Thus
\(\mathscr D_\infty(f)=e^{\mathfrak g_\beta}\operatorname{CGMC}(f)\)
almost surely for each fixed continuous \(f:\mathbb T\to[0,1]\).  Apply
this identity on one countable convergence-determining subset of
\(C(\mathbb T)\), intersect the resulting probability-one events, and use
linearity and uniform approximation.  The two finite Borel measures then
agree on that common probability-one event, which proves the theorem.
\end{proof}

Hence the remainder of the paper concerns the proof of Proposition \ref{prop:main}, which we develop in a sequence of reductions.

\subsection{Discussion and related work}
\label{sec:related}

\emph{The random matrix side of the Fyodorov--Hiary--Keating conjecture.}
The conjecture of \cite{FHK12} (see \cite{FK} for a detailed account) was derived from the freezing scenario for random energy models with logarithmically correlated potentials \cite{FBo,FLDR09}, and numerical evidence for it was given in \cite{FGK18}.
On the rigorous side, the leading order $M_n=(1+o(1))\log n$ for $\beta=2$ was proved in \cite{ABB}, and the second order correction $-\tfrac34\log\log n$ for $\beta=2$ was identified in \cite{PaquetteZeitouni02}.
For general $\beta>0$, tightness of $M_n-\sqrt{2/\beta}\,m_n$ was proved in \cite{CMN}, and convergence in distribution, in the form of Theorem \ref{old-main}, in \cite{PaquetteZeitouni03}.
Theorem \ref{thm:main} identifies the remaining random variable $\mathscr{B}_\infty(\beta)$, and so the limit law ${\bf R}_\beta$ is the law of an explicit constant plus $(2\beta)^{-1/2}$ times the sum of two independent Gumbel variables of parameter $1$; for $\beta=2$ this agrees, up to a deterministic shift, with the law \eqref{eq-R2}, in view of the observation of \cite{SZ}.
The freezing scenario also produces predictions for other quantities and other ensembles.
For instance, \cite{FyodorovSimm16} predicted the law of the extreme values of $\log|\det(x-H)|-\E \log |\det(x-H)|$ for a GUE matrix $H$, with $x$ ranging over an interval of the bulk of the spectrum, and \cite{FLD} predicted the moments of the position of the maximum, both for GUE characteristic polynomials and for logarithmically correlated Gaussian processes.

\emph{Other random matrix ensembles.}
Extreme values of characteristic polynomials have been studied for several other ensembles.
For the Gaussian $\beta$-ensembles, a central limit theorem for the logarithm of the characteristic polynomial was proved in \cite{ABZ}, 
see also \cite{BMP}, and strong approximations of the characteristic polynomial through its transfer matrix recursion were developed in \cite{LambertPaquetteHyperbolic,LambertPaquetteEdge}.
The leading order of the maximum for a class of almost Gaussian logarithmically correlated fields coming from random matrices was obtained in \cite{LambertPaquette01}.
Recently, \cite{AugeriZeitouni25} obtained the second order asymptotics of the maximum of the characteristic polynomial of random Jacobi matrices, which partially confirms the prediction of \cite{FyodorovSimm16} for the Gaussian $\beta$-ensembles.
The maximum of the characteristic polynomial has also been studied for Wigner matrices and $\beta$ ensembles \cite{BourgadeLopattoZeitouni}, for the Ginibre ensemble \cite{Lambert20Ginibre} and its i.i.d. entries counterpart
\cite{CL25}, for two dimensional log-gases \cite{LLZ}, 
and for random permutation matrices \cite{CookZeitouni20}. The closely related question of the maximal fluctuation of the eigenvalue counting function of Hermitian ensembles was treated in \cite{CFLW}.

\emph{The Riemann zeta function.}
The conjecture of \cite{FHK12,FK} also concerns the Riemann zeta function: if $\tau$ is uniformly distributed on $[T,2T]$, then $\max_{|h|\leq 1}\log|\zeta(\tfrac12+i\tau+ih)|$ should equal $\log\log T-\tfrac34\log\log\log T$ plus a random variable which converges in distribution, with a limit law described in \cite{FHK12}.
The leading order $\log\log T$ was proved in \cite{Najnudel} under the Riemann hypothesis and unconditionally in \cite{ABBRS}; see \cite{HarperBourbaki} for an exposition.
For a randomized model of the zeta function, built from independent uniform phases indexed by the primes, the second order term $-\tfrac34\log\log\log T$ was obtained in \cite{ABH17}, which made the connection with branching random walks explicit.
For the zeta function itself, the upper bound with the correct second order term was proved in \cite{Harper19}, an upper bound with an exponentially decaying upper tail was proved in \cite{ABRI}, and the matching lower bound was proved in \cite{ABRII}, so that the recentred maximum is tight.
Sharp upper bounds on intervals of length $(\log T)^{\theta}$ were obtained in \cite{ArguinHamdan}, and the moments of the zeta function on short intervals were studied in \cite{AOR21}.
To our knowledge, the convergence in distribution of the recentred maximum on the critical line remains open.
Slightly off the critical line, namely for $\sigma>\tfrac12$ with $\sigma-\tfrac12$ of order at most a power of $\log\log T$ divided by $\log T$, and for maxima over expanding intervals of length $(\log T)^{\theta}$ with $\theta>0$, the convergence in distribution of the recentred maximum of $\log|\zeta(\sigma+i\tau+ih)|$ to a Gumbel law was recently proved in \cite{OuimetMax26}; see also \cite{OuimetPD26,OuimetHaar26} for the associated Gibbs measures.
Finally, on the multiplicative chaos side, \cite{SaksmanWebb20} proved that the randomly shifted zeta function on the critical line converges in distribution, as a random generalized function, to the product of a smooth random function and a complex Gaussian multiplicative chaos.

\emph{Random matrices and multiplicative chaos.}
The convergence of powers of the modulus of the characteristic polynomial, normalized by their expectation, to Gaussian multiplicative chaos was proved for the CUE in the $L^2$ phase in \cite{Webb15} and throughout the subcritical ($L^1$) phase in \cite{NSW20}.
Analogous results were proved for a class of unitary invariant Hermitian ensembles in \cite{BWW18}, for the other classical compact groups in \cite{ForkelKeating21}, and for the Gaussian orthogonal and symplectic ensembles in \cite{Kivimae24}; for regularized eigenvalue counting statistics see \cite{LOS18}.
See also \cite{BourgadeFalconet} for a two-dimensional Liouville quantum gravity measure arising from unitary Brownian motion, \cite{NPS23} for the holomorphic multiplicative chaos arising from the CUE, and \cite{KeatingWong22} for the connection with moments of moments.
For the C$\beta$E with general $\beta>0$, the convergence to subcritical chaos was proved in \cite{lambert2024subcritical}, whose strategy we follow, as explained in Section \ref{sec-overview}.

The work \cite{CN} is closely related, but concerns a different object.
For $n\geq 1$, the Bernstein--Szeg\H{o} approximation of order $n$ associated with the Verblunsky coefficients $\gamma_0,\ldots,\gamma_{n-1}$ is the probability measure on the circle proportional to $|\Phi_n^*(e^{i\theta})|^{-2}\dif\theta$; its first $n$ Verblunsky coefficients are $\gamma_0,\ldots,\gamma_{n-1}$.
In \cite{CN} it is shown that, as $n\to\infty$, these measures converge to a random probability measure, which is identified with the normalized Gaussian multiplicative chaos $\operatorname{GMC}^{\gamma}$ with $\gamma=\sqrt{2/\beta}$, whenever $\beta\geq 2$.
In particular, the coupling between $\beta$ and the chaos parameter is fixed, and the critical chaos appears only at $\beta=2$.
In contrast, Theorem \ref{thm:maincor} concerns the derivative martingale built from $e^{\sqrt{\beta/2}\,\varphi_k}$, which for $\sigma=1$ is $|\Phi_k^*|^{\sqrt{2\beta}}$, and which is at the critical exponent for every $\beta>0$; the limit is, up to the deterministic factor $e^{\mathfrak{g}_\beta}$, the same critical chaos for all $\beta$.
We also mention the recent work \cite{CaiFangGuo26}, which claims convergence to critical Gaussian multiplicative chaos for measures built from the characteristic polynomial of complex i.i.d.\ matrices, on disks $\{|z|\leq r\}$ with $r<1$, with a covariance that is corrected by the fourth cumulant of the entries.

The convergence to Gaussian chaos is also related to the evaluation of Fisher--Hartwig asymptotics, see \cite{BDHK26}, \cite{Kel26} and references therein.

\emph{Critical chaos.}
Critical Gaussian multiplicative chaos was constructed via the derivative martingale in \cite{DRSV} and via a Seneta--Heyde renormalization in \cite{DRSV2}.
The independence of the limit from the approximation scheme was proved in \cite{JunnilaSaksman,Powell18}, and the critical chaos was identified as the limit of renormalized subcritical chaos in \cite{APS19} (see \cite{PowellReview} for general dimension); see \cite{PowellReview,BerestyckiPowell} for surveys, and \cite{Lacoin}, whose results we use in the proof of Proposition \ref{prop:criticalchaos}.
In the setting of branching processes, the corresponding objects are the derivative martingales of branching Brownian motion \cite{LS} and of branching random walks, together with their Seneta--Heyde normalizations \cite{BK97,AidekonShi14}; for Mandelbrot cascades see \cite{BKNSW14,BKNSW15}.
In all of these settings, the limit of the derivative martingale gives the random shift in the limit law of the recentred maximum: see \cite{LS,Aidekon,Madaule17} for branching processes, \cite{Madaule,DRZ17} for logarithmically correlated Gaussian fields, and \cite{BDZ,BL} for the two-dimensional discrete Gaussian free field.
Theorem \ref{old-main} is the analogous statement for the C$\beta$E.

The law of the total mass of subcritical Gaussian multiplicative chaos on the circle was predicted in \cite{FBo} and proved in \cite{Remy}, with an alternative proof in \cite{CN}; the corresponding law on the unit interval, predicted in \cite{FLDR09}, was proved in \cite{RemyZhu}.
The law \eqref{eq:Remysmagiccarpet} of the total mass of the critical chaos on the circle follows from the subcritical formula of \cite{Remy} by the renormalized limit of \cite{APS19,PowellReview}.

Theorem \ref{thm:maincor} can be viewed as a universality statement for critical chaos: the field $\varphi_k$ is not Gaussian, and its derivative martingale converges to a deterministic multiple of the critical Gaussian multiplicative chaos, with a limit that does not depend on $\beta$.
For subcritical chaos, such universality with respect to non-Gaussian approximations is well developed; see, for example, \cite{Junnila20} and the random matrix results cited above.
Critical chaos arising from non-Gaussian fields has been obtained in several other settings, including critical Brownian multiplicative chaos built from planar Brownian motion \cite{Jego21}, critical chaos arising from random multiplicative functions \cite{GW25b,HKK26} (see also \cite{Harper20,GuZhang26} for the related low moment problems).

\subsection*{Acknowledgements}
Research by E. Paquette was supported by a Discovery Grant from the
Natural Science and Engineering Research Council (NSERC),
RGPIN-2025-04643, an NSERC--FRQNT
NOVA Grant, a CIFAR AI Catalyst Grant, US Air Force Office of Scientific Research
(AFOSR) Award FA9550261B232, and a gift from the Google Canada corporation.
Research by O. Zeitouni was partially supported by the Israel Science
Foundation grant number 615/24.

\subsection*{Statement on the use of AI tools}
An early version of this document was written entirely without AI
tools.  The current version has been modified using
large language models (LLMs).  These were used to migrate the early
version of the proof, which contained a flaw, to a different barrier
setup.  They were also used for large reorganizations of the text, to
assist with writing, and to find errors in the written proof.  The
document has been revised repeatedly by the authors for style and correctness.

\subsection{Notation and conventions}
We use the Landau notation $O_n(\cdot),o_n(\cdot)$; when $n$ is clear from the context, we eliminate it from the notation.
If $c_n$ is a complex sequence, we write $O(c_k)$, $o(c_k)$ to mean $O(|c_k|)$ and $o(|c_k|)$. Throughout the paper, we use $\bk,\bn$ etc.
when discussing parameters in a dyadic exponent, so that typically $k=2^\bk$.
We write $\N =\{1,2,\dots\}$ and $\overline{\N} = \N \cup \{\infty\}$.
For convenience and easy reference, we record in
Tables \ref{tab:notation2}--\ref{tab:notation4} a glossary of notation.

\begin{table}[H]
\centering
\begin{tabular}{ r p{1.8cm} p{11cm}}
\toprule
Symbol & Ref. & Meaning \\
\midrule
$X_n(z)$ & \eqref{eq:charpoly} & The characteristic polynomial of the C$\beta$E, $X_n(z) = \prod_{j=1}^n (1-e^{i\omega_j}z)$. \\
$M_n$ &\eqref{eq-280126}& The maximum $M_n = \max_{|z|=1} \log|X_n(z)|$. \\
$m_n$ & \eqref{eq-mn} & The centering sequence $m_n = \log n - \tfrac{3}{4}\log\log n$. \\
${\bf R}_\beta$ & \eqref{eq-convolution} & The limiting random variable for $M_n - \sqrt{2/\beta}\, m_n$. \\
$\mathscr{B}_\infty$ & Thm.~\ref{old-main} & A nonnegative random variable appearing in the decomposition of ${\bf R}_\beta$, whose law is identified in Theorem~\ref{thm:main}. \\
$\gamma_k$ & \eqref{eq:betagamma} & The Verblunsky coefficients; complex, rotationally invariant with $|\gamma_k|^2 \sim \Beta(1, \frac{\beta(k+1)}{2})$. \\
$\filt_k$ & \eqref{eq:filtration} & The enlarged latent-variable filtration generated by $\Gamma_0,\ldots,\Gamma_{k-1}$ and the Brownian path through time $k$; it contains $\gamma_0,\ldots,\gamma_{k-1}$.\\
$\Phi_k, \Phi_k^*$ & \eqref{eq:szego} & The OPUC and reversed OPUC polynomials from the Szeg\H{o} recurrence. \\
$Q_k(z)$ & \eqref{eq:Q} & The ratio $Q_k(z) = z\Phi_k(z)/\Phi_k^*(z)$, which is $\filt_{k}$-adapted. For $|z|=1$, $|Q_k(z)|=1$; for $|z|<1$, $|Q_k(z)|<1$. \\
$\Psi_k(\theta)$ & \eqref{eq:Psi} & The Pr\"ufer phase, a continuous version of $\arg(e^{i\theta}\Phi_k(e^{i\theta})/\Phi_k^*(e^{i\theta}))$. \\
$\psi_k(\theta)$ & \eqref{eq:prufer} & The relative Pr\"ufer phase; in law equals $\Psi_k(\theta) - \Psi_k(0)$. \\

$\sigma$ & \eqref{eq:upsilon} & Either $1$ or $i$, selecting real or imaginary part of the log-correlated field. \\
$\varphi_k(r,\theta)$ & \eqref{eq:upsilon} & The log-correlated field $\varphi_k(r,\theta) = 2\Re(\sigma\log\Phi_k^*(re^{i\theta}))$. \\
$\lgr(k)$ & \eqref{lgr} & The cutoff logarithm $\lgr(k) = \log(k \wedge \tfrac{1}{1-r^2})$. \\
$\mathscr{D}_k(r,\theta)$ & \eqref{eq:Bk} & The derivative martingale density. \\
$\mathscr{B}_k(r)$ & \eqref{eq:Bk} & The total mass $\mathscr{B}_k(r) = \int_0^{2\pi} \mathscr{D}_k(r,\theta)\,d\theta$. \\
$r_{\bk}$ & \eqref{eq-rk} & The dyadic radii $r_{\bk}^2 = 1 - 2^{-\bk}$. \\
$\mathscr{A}_{\bk}^{L,a}(\theta)$  & \eqref{eq:Abarrier} & The barrier event controlling $\varphi_\infty(r_{\bj},\theta)$ for $\bj \leq \bk$. \\
$\mathscr{A}_{\bk,\bj}^{L,a}(\theta)$  & \eqref{eq:ABbarrier} &
The fixed-time barrier controlling
$\varphi_{2^{\bj}}(r_{\bl},\theta)$ over
$1\leq\bl\leq\bk\wedge\bj$. \\
$\mathscr{B}_{\bk,\bl}^{L,\underline{a}}(\theta)$&Lem. \ref{lem-1.20}& The intersection of the events $\mathscr{A}_{\bk,\bj}^{L,a_\bj}(\theta)$ 
for $\bj\geq \bl$.\\
$\mathscr{E}_{\bk,\bj}^{L,a}$ & \eqref{eq:physical-interface-barrier} &
The fixed-radius barrier controlling
$\varphi_{2^{\bl}}(r_{\bk},\theta)$ for $\bl\geq\bj$. \\
$\operatorname{CGMC}$ & Prop.~\ref{prop:criticalchaos} & Critical Gaussian multiplicative chaos.\\
$\mathfrak{g}_\beta$ & Prop.~\ref{prop:main} & A deterministic constant relating $\mathscr{D}_\infty$ to CGMC. \\
$\Exp_\bn$ & \eqref{eq-En} & The conditional expectation $\Exp(\cdot \mid \filt_{2^{\bn}})$. \\
\bottomrule
\end{tabular}
\caption{Table of notation: basic objects, derivative martingales, and barrier events.}
\label{tab:notation2}
\end{table}

\begin{table}[H]
\centering
\begin{tabular}{ r p{1.8cm} p{11cm}}
\toprule
Symbol & Ref. & Meaning \\
\midrule
$E_j, \Gamma_j, \Theta_j$ & \eqref{eq:betagamma} & Components of the Verblunsky coefficients: $E_j \sim \Exponential(1)$, $\Gamma_j \sim \Gam(\frac{\beta}{2}(j+1))$, $\Theta_j \sim \Unif([0,2\pi])$. \\
$Z_j$ & \eqref{eq:gaussians} & The complex Gaussian $Z_j = \sqrt{E_j}e^{i\Theta_j}$, with $\Exp|Z_j|^2 = 1$. \\
$\beta_j$ &\eqref{eq-betaj}  & The rescaled parameter $\beta_j = \sqrt{\frac{\beta}{2}(j+1)}$. \\
$B_t$ & \eqref{eq:W} & A complex Brownian motion with $[B_t, \overline{B_t}] = 2t$ and $\sqrt{2}Z_k = B_{k+1} - B_k$. \\
$\zeta_t^\C(\bn;r,\theta)$ & \eqref{eq:Z} & A Gaussian approximation to $2\log\Phi^*_t(re^{i\theta})$, defined for $t \geq 2^{\bn}$. \\
$\zeta_t(\bn;r,\theta)$ &\eqref{eq-zetareal}  & The real projection $\zeta_t = \Re(\sigma\zeta_t^\C)$. \\
$\operatorname{Tru}_N$ & \eqref{eq:good} & The truncation event on which the Verblunsky coefficients satisfy $|\gamma_j| \ll \sqrt{\log j/j}$ for $j \geq N$.\\
$\operatorname{Tru}_{N,J}$& \eqref{eq:goodtruncated} & The restriction of the truncation event to the coefficient interval $[N,J)$.\\
$\mathscr{G}_{\bn}$ & \eqref{eq:globalgood} & The global good event: $\operatorname{Tru}_{2^{\bn}}$ intersected with uniform control of $|\zeta^\C - 2\log\Phi^*|$, a bound on $\operatorname{Osc}_{2^\bn}$, and derivative bounds. \\
$\widehat{\mathscr{G}}_{\bn}$ & \eqref{eq:globalgooda} & Another global good event: $\operatorname{Tru}_{2^{\bn}}$ intersected with a bound on $\operatorname{Osc}_{2^\bn}$ and a derivative bound. \\
$\operatorname{Osc}_k$ & \eqref{eq:osc-def} & The oscillation of $\log|\Phi_k^*|$ on the unit circle. \\
$\mathscr N_\bu$ & \eqref{eq:scale-failure-event} & The event that the good event fails at scale $2^\bu$, through the derivative bound or the linearization increment. \\
$S_{\bk,\bn}^{\bt}$ & \eqref{eq:StLt} & The deviation of the quadratic variation $\langle \zeta \rangle$ from its typical value. \\
$L_{\bk,\bn}^{\bt}$ & \eqref{eq:StLt} & The linearization error $\sqrt{\frac{\beta}{2}}(\varphi_{2^{\bt}} - \zeta_{2^{\bt}})$. \\
$\mathscr{M}_{\bk,\bt}$ & \eqref{eq:Mn} & The exponential martingale $\frac{1}{2\pi}\exp(\sqrt{\frac{\beta}{2}}\varphi_{2^{\bt}}(r_{\bk},\theta) - \lgr(2^{\bt}))$. \\
$\Q_{\bk,\bn}^{\bt}$ & \eqref{eq:Qt} & A martingale change of measure tilting towards the derivative martingale. \\
$\Delta_{\bk,\bm}^{[\bn]}$&
\eqref{eq:buffered-gap-definition}&
The deficit (distance to the barrier)
$\sqrt{2}\lgr(2^{\bm})-\sqrt{\frac{\beta}{4}}
\zeta_{2^{\bm}}(\bn;r_{\bk},\theta)$ with fixed initial exponent $\bn$.\\
\bottomrule
\end{tabular}
\caption{Table of notation: Gaussian approximation and technical events.}
\label{tab:notation3}
\end{table}

\begin{table}[H]
\centering
\begin{tabular}{ r p{2.2cm} p{10.6cm}}
\toprule
Symbol & Ref. & Meaning \\
\midrule
$R_j(r_{\bk},\theta)$ & \eqref{eq-Rjnew} & The time dilation
$(1-|Q_j(r_{\bk},\theta)|^2)/(1-r_{\bk}^2)$. \\
$q_{\bk,\bl}(s)$ & \eqref{eq:brownian-clock} & The Brownian clock on the block starting at $\bl$. \\
$\mathscr V_{\bk,\bl}^{\nu,[\bn]}$ & \eqref{eq:buffer-event} & The event that
the deficit at the left endpoint, with fixed initial exponent $\bn$, has
square-root scale. \\
$\mathscr W_{\bu}^{B}$ & \eqref{eq:step2-radial-window} & The event $R_{2^{\bu}}\leq C_B2^{\bu}\bu^B$ of typical time dilation. \\
$\widetilde{\mathscr E}_{\bk,\bl}^{\bu,L,a,[\bn]}$ &
\eqref{eq:clock-adapted-barrier} & The Brownian barrier on the block
$[\bl,\bu]$, with initial exponent $\bn$. \\
$\widetilde{\mathscr E}_{\bk,\bn}^{\mathcal P,L,a}$ &
\eqref{eq:brownian-clock-barrier} & The Brownian barrier:
the intersection of the block barriers over the partition and the
terminal block. \\
$Z_{\bk,\bn}$ & \eqref{eq:averaged-FK-factor} & The ratio
$\Exp_{\bn}(\mathscr M_{\bk,\infty})/\mathscr M_{\bk,\bn}$, written via
\eqref{eq:com} as
$\mathbb Q_{\bk,\bn}^{\infty}(e^{S_{\bk,\bn}^{\infty}+L_{\bk,\bn}^{\infty}})$. \\
$c_{\bk}$ & Lemma~\ref{lem:step3-late-forgetting} & The constant with
$\mathbb Q_{\bk,\bt_\bk}^{\infty}(\operatorname{Tru}_{2^{\bt_\bk}}e^{S_{\bk,\bt_\bk}^{\infty}})=c_{\bk}+o_{\bk}(1)$;
identified in Proposition~\ref{lem:averaged-FK-transfer}. \\
\bottomrule
\end{tabular}
\caption{Table of notation: restoration, mixing of the clock, and the Brownian barrier.}
\label{tab:notation4}
\end{table}

\section{Preliminaries}
\label{sec-prelim}
\subsection{Marginally Gaussian approximation}

We shall need good single-$\theta$ estimates for $\varphi_k(r,\theta)$ which capture some of its random walk properties in $k$ (in particular we shall need a form of the ballot theorem). 
To do this, we shall use a partial Gaussian approximation that exploits that the Verblunsky coefficients are already nearly Gaussian, as follows from the latent-variables representation  \eqref{eq:betagamma}--\eqref{eq-betaj} and \eqref{eq:W}.

We will need to truncate the variables $\{\gamma_j\}$. 
Introduce the $\filt_J$-measurable events 
\begin{equation}
\operatorname{Tru}_{N,J} \coloneqq \left\{ 
\max_{t \in [j,j+1]} | B_t - B_{j}| \leq 4\sqrt{\log j}
\text{ and }
|  \tfrac2{\beta} \Gamma_j - j| \leq 4\sqrt{2/\beta} \sqrt{ j \log j}
\text{ for all }  j \in [N,J)\right\},
\label{eq:goodtruncated}
\end{equation}
and set, 
for $N\in\N$
\begin{equation}\label{eq:good}
\operatorname{Tru}_N \coloneqq
\bigcap_{J\geq N} \operatorname{Tru}_{J,J+1}=\bigcap_{J\geq N} \operatorname{Tru}_{N,J}.
\end{equation}
(On $\operatorname{Tru}_N $ we have that $|Z_j|\ll \sqrt{\log j}$ and $|\gamma_j|\ll \sqrt{\log j/j}$ for all $j\geq N$.)
 Note that
the random variables  $\sqrt{\beta/2} \big\{ \frac{\frac2{\beta} \Gamma_j -  j}{\sqrt j}\big\}$ are approximately standard Gaussians, so 
\begin{equation} \label{Truest}
\mathbb P[\operatorname{Tru}_N^\complement ] \ll \sum_{j\ge N} \exp(-8 \log j ) \ll N^{-7}.
\end{equation}

%

The next lemma controls the error in considering the process $ \tfrac12 \zeta_{k}^\C$ instead of $\log \Phi^*_{k}(re^{i\theta})$.
\begin{lemma} \label{lem:zetaphitail}
For any $p \ge 0$,  there is a constant $C_p$ such that for all $r \in [0,1]$,  all $\theta \in [0,2\pi]$ and for all $n$
\[
\Exp\!\left[
 \one\{\operatorname{Tru}_{2^{\bn}}\}
 \sup_{k \geq 2^{\bn}} \exp \Big( p \bn^{-1}2^{\bn/2}
 \big| \tfrac12 \zeta_{k}^\C(\bn;r,\theta)
       - \log \Phi^*_{k}(re^{i\theta})\big| \Big)
 \mathrel{\Big|}\filt_{2^{\bn}}\right] \le C_p, \quad\As .
\]
\end{lemma}

\begin{proof}
By \eqref{eq:Phi}, with $z=re^{i\theta}$ and $k\ge 2^{\bn}$, 
\[
\log \Phi^*_{k}(z)-\log\Phi^*_{2^{\bn}}(z) = \sum_{2^{\bn} \leq j<k}\log(1-\gamma_j Q_j(z))  .
\]
The random variables $\gamma_j$ are independent, rotation-invariant, so the RHS is a (complex-valued) $\filt_k$ martingale, even when additionally conditioning on $|\gamma_j|, j\geq 1$. Then, using that $|Q_j(z)|\le 1$ and $|\gamma_j|^2\ll \tfrac{\log j}{j}$ for $j\ge 2^{\bn}$  on $\operatorname{Tru}_{2^{\bn}}$,
the quadratic variation of the martingale  
\[
\sum_{2^{\bn} \leq j<k}\log(1-\gamma_j Q_j(z))  + \gamma_j Q_j(z)
\]
is controlled on
$ \operatorname{Tru}_{2^{\bn}} $ by 
\[
\sum_{j\ge 2^{\bn}} |\gamma_j|^4  \ll  \bn^2 2^{-\bn}\]
and its increments are bounded by
\[
\sup_{j\ge 2^{\bn}}\big| \log(1-\gamma_j Q_j(z))  + \gamma_j Q_j(z) \big| \ll \bn 2^{-\bn}
\quad\text{on } \operatorname{Tru}_{2^{\bn}} .
\]
Conditioning on the event $\operatorname{Tru}_{2^{\bn}}$ does not destroy the martingale property.  Indeed, that event constrains only the moduli of the Brownian increments and the $\Gamma_j$, and the phases of the $Z_j$ remain independent and uniform.
{Freedman's inequality (Lemma \ref{lem:freedman}) therefore gives, for any $\lambda>0$,}
\begin{equation} \label{logtb1}
\mathbb P_\bn \bigg[\operatorname{Tru}_{2^{\bn}}\cap
 \bigg\{ \sup_{k\ge 2^{\bn}} \bigg| \log \Phi^*_{k}(z)-\log\Phi^*_{2^{\bn}}(z) + \sum_{2^{\bn}\le j<k}  \gamma_j Q_j(z) \bigg| \ge \lambda\bigg\}\bigg]
 \ll \exp\bigg(- c \bigg(\frac{\lambda^2}{\bn^2 2^{-\bn}}\wedge \frac{\lambda}{\bn 2^{-\bn}} \bigg)\bigg) . 
\end{equation}

Moreover, by \eqref{eq:Z}, 
for $k\ge 2^{\bn}$ integer, 
\begin{equation}\label{eq:Z2}
\tfrac12 \zeta_{k}^\C(\bn;r,\theta) =  -  
\sum_{2^{\bn}\le j < k}
\frac{Q_j(z)}{\beta_j}Z_j  + \log \Phi^*_{2^{\bn}}(z) .
\end{equation}

Then, as $\gamma_j = Z_j/\sqrt{|Z_j|^2+\Gamma_j}$,
\begin{equation}
\label{eq-250826a}
\tfrac12 \zeta_{k}^\C(\bn;r,\theta)  -\log\Phi^*_{2^{\bn}}(z) + \sum_{2^{\bn}\le j<k}  \gamma_j Q_j(z)=  \sum_{2^{\bn}\le j <k}Q_j(z) Z_j
\left(\frac{1}{\sqrt{|Z_j|^2+\Gamma_j}}-\frac{1}{\beta_j}\right) .
\end{equation}
The RHS is also a martingale (even when conditioned on $\operatorname{Tru}_{2^{\bn}}$)  and, on the latter event,   its quadratic variation is bounded by 
\[
\sum_{j\ge 2^{\bn}} \frac{ (\log j)^2}{j^2}  \ll  \bn^2 2^{-\bn}\]
while its increments are bounded by 
\[
\sup_{j\ge 2^{\bn}} |Z_j| \left|\frac{1}{\sqrt{|Z_j|^2+\Gamma_j}}-\frac{1}{\beta_j} 
\right| \ll \bn 2^{-\bn}.
\] 
This yields, again by an application of Lemma \ref{lem:freedman}, the same tail bound as \eqref{logtb1}, and thus concludes the proof.
\end{proof}

Furthermore, we can upgrade Lemma \ref{lem:zetaphitail} to a bound which holds uniformly in $\theta$.  
\begin{proposition} \label{prop:zetaphitail}
We have 
\[
  \sup_{r \in [0,1]}
  \sup_{\theta \in [0,2\pi]}\sup_{\bl \geq \bn}  
  |\zeta^{\C}_{2^{\bl}}(\bn;r,\theta)-2\log \Phi_{2^{\bl}}^*(r,\theta)|
\Prto[\bn] 0.
\]
\end{proposition}
\begin{proof}
When $r=1$, this is \cite[Proposition 3.1]{CMN}, where it is stated as tightness for fixed $\bn$; the same argument gives convergence in probability as $\bn\to\infty$.

Now for controlling $r < 1$, we note that $Q_k(z)$ are analytic in the unit disk, and so is $\log \Phi^*_{k}$ for all $k$.  Hence, so is $re^{i\theta}\mapsto\zeta^{\C}_{k}(\bn;r,\theta)$.  Thus by the maximum modulus principle, we have 
\[
\sup_{r}
\sup_{\theta}\sup_{\bj \geq \bn}  
|\zeta^{\C}_{2^{\bj}}(\bn;r,\theta)-2\log \Phi_{2^{\bj}}^*(r,\theta)|
=
\sup_{\theta}\sup_{\bj \geq \bn}
|\zeta^{\C}_{2^{\bj}}(\bn;1,\theta)-2\log \Phi_{2^{\bj}}^*(1,\theta)|.
\]
\end{proof}
As a consequence, there is a sequence of $\epsilon_\bn \to 0$ so that the event 
\[
  \{
    \sup_{r}
    \sup_{\theta}
    \sup_{\bj \geq \bn}
    |\zeta^{\C}_{2^\bj}(\bn;r,\theta)-2\log \Phi_{2^\bj}^*(r,\theta)|
    \leq \epsilon_\bn
    \}
\]
holds with probability tending to $1$ as $\bn\to\infty$.  
Enlarging this sequence if necessary, we assume from now on that
\begin{equation}\label{eq:global-linearization-tolerance}
 \epsilon_\bn\geq
 C\sum_{\bu\geq\bn}\bu^2 2^{-\bu/4}\to_{\bn\to\infty} 0.
\end{equation}

For \(k\geq0\), let
\begin{equation}\label{eq:osc-def}
 \operatorname{Osc}_k:=\max_{|z|=1}\log|\Phi_k^*(z)|-\min_{|z|=1}\log|\Phi_k^*(z)|
\end{equation}
be the oscillation of \(\log|\Phi_k^*|\) on the unit circle.  It is
finite, since \(\Phi_k^*\) has no zeros on the closed unit disk (see \cite{SimonBook}).
Next we
introduce the combined good event
\begin{align}\label{eq:globalgood}
{\mathscr{G}_{\bn}}
\coloneqq  \operatorname{Tru}_{2^{\bn}}
&\bigcap 
\{
\sup_{r}
\sup_{\theta}
\sup_{\bj \geq \bn}
|\zeta^{\C}_{2^\bj}(\bn;r,\theta)-2\log \Phi_{2^\bj}^*(r,\theta)|
\leq \epsilon_\bn
\}
\bigcap
\left\{\operatorname{Osc}_{2^\bn}
       \leq C_{\rm in}\bn\log2\right\} \\
&\bigcap 
\{ \max_{|z|\leq 1}|Q_{j}'(z)| \leq j^{C'},\, \forall j\geq 2^{\bn}\}.\nonumber
\end{align}
The exponents \(C_{\rm in}\) and \(C'\) are fixed in
\eqref{eq:derivative-constant-hierarchy}, after
Lemma~\ref{lem:derivativecontroltwo}, and with that choice
\(\mathscr G_\bn\) has probability tending to one; see the paragraph
following \eqref{eq:derivative-constant-hierarchy}.
For \(\bm\geq\bn\), we also define the restriction of this event to
times up to \(2^{\bm}\),
\begin{align}\label{eq:globalgood-finite}
 \mathscr G_{\bn,\bm}
 \coloneqq \operatorname{Tru}_{2^{\bn},2^{\bm}}
 &\bigcap
 \Big\{
 \sup_{r}\sup_{\theta}\sup_{\bn\leq\bj\leq\bm}
 |\zeta^{\C}_{2^\bj}(\bn;r,\theta)-2\log \Phi_{2^\bj}^*(r,\theta)|
 \leq \epsilon_\bn
 \Big\}
 \bigcap
 \left\{\operatorname{Osc}_{2^\bn}
       \leq C_{\rm in}\bn\log2\right\} \notag\\
 &\bigcap
 \{ \max_{|z|\leq 1}|Q_{j}'(z)| \leq j^{C'},\
 2^{\bn}\leq j\leq 2^{\bm}\},
\end{align}
which is \(\filt_{2^{\bm}}\)-measurable, decreasing in \(\bm\), and
satisfies \(\mathscr G_{\bn,\infty}=\mathscr G_\bn\).
We sometimes need a slightly less restrictive event $\widehat{\mathscr{G}}_{\bn}$, defined as
\begin{equation}\label{eq:globalgooda}
{\widehat{\mathscr{G}}_{\bn}}
\coloneqq \operatorname{Tru}_{2^{\bn}}
\cap 
\left\{\operatorname{Osc}_{2^\bn}\leq C_{\rm in}\bn\log2\right\}\cap
\{ \max_{|z|\leq 1}|Q_{j}'(z)| \leq j^{C'},\, \forall j\geq 2^{\bn}\}.
\end{equation}

We define \(\widehat{\mathscr{G}}_{\bn,\bl}\) from
\eqref{eq:globalgooda} in the same way as \(\mathscr G_{\bn,\bl}\) is
obtained from \eqref{eq:globalgood}.
 Since $\mathscr{G}_{\bn}\subset
\mathscr{G}_{\bn,\bl}$, the latter also holds with probability tending to $1$ as $\bn\to\infty$, uniformly in $\bl$.
Note also  that on the event ${\widehat{\mathscr{G}}_{\bn}}$ (and hence, also on ${\mathscr{G}}_{\bn}$) we have that 
\begin{equation}
\label{eq-040426a}
|Q_j(z)|\geq 1-(1-|z|^2)j^{C'}, \quad \forall j\geq 2^{\bn}.
\end{equation}

\subsection{The quadratic variation process}

By construction, \eqref{eq:Z} and $\zeta_t$ are continuous time martingales for $t \geq 2^{\bn}.$  Its bracket process is given by, for $u > 2^{\bn}$,
\[
\frac{\beta}{4}
\langle {\zeta}_u(\bn;r,\theta)\rangle = 
\int_{2^{\bn}}^{u} 
\frac{|Q_{[s]}(r,\theta)|^2}
{[s]+1}\,\dif s,
\]
so that for $u \in \N$ we have
\begin{equation}
\label{eq-301024}
\frac{\beta}{4}
\langle \zeta_{u}(\bn;r,\theta) \rangle = \sum_{j=2^{\bn}}^{u-1} \frac{|Q_j(r,\theta)|^2}{j+1}.
\end{equation}

Recall that  $\varphi_{k}(r,\theta) = 2\Re(\sigma\log\Phi^*_k(re^{i\theta}))$, $r_{\bk}^2=1-2^{-\bk}$ and \eqref{lgr}, 
we define for $\bk, \bn , \bt \in\N$ with $\bt\ge \bn$ and $\theta\in[0,2\pi]$, 
\begin{equation}\label{eq:StLt}
\begin{aligned}
&
S_{\bk,\bn}^{\bt} 
=S_{\bk,\bn}^{\bt}(\theta) \coloneqq \tfrac{\beta}{4}\langle \zeta_{2^{\bt}}(\bn;r_{\bk},\theta) \rangle - \lg_{r_{\bk}}(2^{\bt})+
\lg_{r_{\bk}}(2^{\bn}), \\
&
L_{\bk,\bn}^{\bt} 
=
L_{\bk,\bn}^{\bt}(\theta) 
\coloneqq \sqrt{\tfrac{\beta}{2}}\left(\varphi_{2^{\bt}}(r_{\bk},\theta)-\zeta_{2^{\bt}}(\bn;r_{\bk},\theta)\right).\\
\end{aligned}
\end{equation}
The linearization errors  $L_{\bk,\bn}^{\bt}$ are controlled uniformly by Proposition~\ref{prop:zetaphitail}.
The quantity
\[
 \tfrac{\beta}{4}\langle\zeta_{2^{\bt}}(\bn;r_{\bk},\theta)\rangle
\]
is the \textbf{Brownian time} elapsed between $2^{\bn}$ and $2^{\bt}$: by the
Dambis--Dubins--Schwarz theorem, $\zeta$ is a Brownian motion when run in
it.  Its typical value is $\lg_{r_{\bk}}(2^{\bt})-\lg_{r_{\bk}}(2^{\bn})$,
and $S_{\bk,\bn}^{\bt}$, the \textbf{clock deviation}, is the difference
between the two.

Define, for any $\bk\in \overline{\N}\coloneqq \N \cup \{\infty\}$ and $\bt\in\N$,
\begin{equation}\label{eq:Mn}
\mathscr{M}_{\bk,\bt}=\mathscr{M}_{\bk,\bt}(\theta) \coloneqq
\tfrac{1}{2\pi}
\exp\big( \sqrt{\tfrac{\beta}{2}}\varphi_{2^{\bt}}(r_{\bk},\theta) -  \lg_{r_{\bk}}(2^{\bt})\big).
\end{equation}

We also introduce the martingale $\filt_{2^{\bn}}$-conditional change of measure
\begin{equation}\label{eq:Qt}
\begin{aligned}
\frac{d\mathbb{Q}^{\bt}_{\bk,\bn}}{d\Pr_\bn}
&=
\frac{d\mathbb{Q}^{\bt}_{\bk,\bn}(\theta)}{d\Pr_\bn}
\coloneqq
\exp\left( \sqrt{\tfrac{\beta}{2}}\zeta_{2^{\bt}}(\bn;r_{\bk},\theta)
-   \tfrac{\beta}{4}\langle \zeta_{2^{\bt}}(\bn;r_{\bk},\theta) \rangle
- \sqrt{\tfrac{\beta}{2}}\varphi_{2^{\bn}}(r_{\bk},\theta)
\right)\\
&= \frac{\mathscr{M}_{\bk,\bt}}{\mathscr{M}_{\bk,\bn}} \exp\left( - S_{\bk,\bn}^{\bt} - L_{\bk,\bn}^{\bt} \right).
\end{aligned}
\end{equation}
Note that $\mathbb{Q}^{\bt}_{\bk,\bn}(\theta)$ is measurable on $\filt_{2^{\bn}}$.

Using \eqref{eq:StLt}, for any $\bt \geq \bn$, we can decompose
\[
\sqrt{\tfrac{\beta}{2}}\varphi_{2^{\bt}}(r_{\bk},\theta) -  \lg_{r_{\bk}}(2^{\bt})
= L_{\bk,\bn}^{\bt} + \sqrt{\tfrac{\beta}{2}}\zeta_{2^{\bt}}(\bn;r_{\bk},\theta) -   \tfrac{\beta}{4}\langle \zeta_{2^{\bt}}(\bn;r_{\bk},\theta) \rangle + S_{\bk,\bn}^{\bt} - \lg_{r_{\bk}}(2^{\bn}).
\]
Thus, by \eqref{eq:Qt}, for bounded random variables $X$, any $\bk, \bt\in \overline{\N}$ with $\bt\ge \bn$,
\begin{equation}\label{eq:com}
\begin{aligned}
\Exp_\bn
\left(
\mathscr{M}_{\bk,\bt}(\theta)
X
\right)
&=
\mathscr{M}_{\bk,\bn}
\mathbb{Q}_{\bk,\bn}^{\bt}
\left( X
\exp( S_{\bk,\bn}^{\bt} + L_{\bk,\bn}^{\bt})
\right).
\end{aligned}
\end{equation}

\subsection{Truncation, linearization, coarse smoothness bounds, and the good event}
\begin{lemma}\label{lem:1st2ndQmoment}
Let  $\bk, \bt\in\overline{\N}$ and $\bn \in \N$.
There are constants \(c_{\rm tr}>0\) and \(C<\infty\)
such that, for every finite \(\bt\geq\bn\),
\begin{equation}\label{eq:tilted-truncation-cost}
 \mathbb Q_{\bk,\bn}^{\bt}\left(
  \operatorname{Tru}_{2^{\bn},2^{\bt}}^{\complement}
\right)
 \leq C2^{-c_{\rm tr}\bn},
\end{equation}
where  $\mathbb{Q}_{\bk,\bn}^{\bt}$ is as in  \eqref{eq:Qt}. Further, under that measure,  the increments $\{Z_j ; 2^{\bn} \leq j <2^{\bt}\}$ are conditionally Gaussian with mean 
\[
\mathbb{Q}_{\bk,\bn}^{\bt}\left( Z_j \middle\vert \filt_j \right) = 
-\frac{\overline{\sigma Q_{j}(r_{\bk},\theta)}}{\sqrt{j+1}} . 
\]
Moreover, for any $\alpha \in (1,3/2)$,
\begin{align}\label{eq:Qgamma}
  \mathbb{Q}_{\bk,\bn}^{\bt}\left( \gamma_j \middle\vert \filt_j \right) &= 
  - \sqrt{\frac{2}{\beta}}\frac{\overline{\sigma Q_{j}(r_{\bk},\theta)}}{{j+1}}
  +\mathcal{O}(j^{-\alpha}) = \mathbb{Q}_{\bk,\bn}^{\bt}\left( \gamma_j 1\{\operatorname{Tru}_{j,j+1}\}\middle\vert \filt_j \right), \\
  \mathbb{Q}_{\bk,\bn}^{\bt}\left( \gamma_j^2 \middle\vert \filt_j \right) &= \mathcal{O}(j^{-\alpha})=  \mathbb{Q}_{\bk,\bn}^{\bt}\left( \gamma_j^2 1\{\operatorname{Tru}_{j,j+1}\}\middle\vert \filt_j \right), \quad
  \mathbb{Q}_{\bk,\bn}^{\bt}\left( |\gamma_j|^2 \middle\vert \filt_j \right) = \frac{2}{\beta j} 
  +\mathcal{O}(j^{-\alpha}).\nonumber
\end{align}
\end{lemma}

\begin{proof}
By \eqref{eq:Z}
\[
\sqrt{\tfrac{\beta}{2}} \zeta_{2^{\bt}}(\bn;r,\theta) =  -   2 \Re\bigg( \sigma \sum_{2^{\bn}\le j <2^{\bt}}
\frac{Q_j(r,\theta)}{\sqrt{j+1}} Z_j \bigg)
+ \sqrt{\tfrac{\beta}{2}}  \varphi_{2^{\bn}}(r,\theta)
\]
where $Z_j$ are Gaussians.
Then, by Girsanov's theorem, under the change of measure $\mathbb{Q}_{\bk,\bn}^{\bt}$, the random 
variables $\{Z_j;2^{\bn}\le j <2^{\bt} \}$ acquire a mean;
\[
\mathbb{Q}_{\bk,\bn}^{\bt}\left( Z_j \middle\vert \filt_j \right) = -\frac{\overline{\sigma Q_{j}(r_{\bk},\theta)}}{\sqrt{j+1}}. 
\]
In particular, we have 
\[
\mathbb{Q}_{\bk,\bn}^{\bt}\left( Z_j^2 \middle\vert \filt_j \right) = \sigma^2\frac{\overline{Q_{j}(r_{\bk},\theta)^2}}{j+1} , \qquad 
\mathbb{Q}_{\bk,\bn}^{\bt}\left( |Z_j|^2 \middle\vert \filt_j \right) = 1+ \frac{|Q_{j}(r_{\bk},\theta)|^2}{j+1}. 
\]

To prove \eqref{eq:tilted-truncation-cost}, write the conditional
one-step density explicitly:
\begin{equation}\label{eq:tilted-one-step-density}
 \frac{d\mathbb Q_{\bk,\bn}^{\bt}
       ((Z_j,\Gamma_j)\in d(z,\gamma)\mid\filt_j)}
      {d\mathbb P((Z_j,\Gamma_j)\in d(z,\gamma))}
 =\exp\left\{-2\Re\frac{\sigma Q_j(r_{\bk},\theta)z}{\sqrt{j+1}}
              -\frac{|Q_j(r_{\bk},\theta)|^2}{j+1}\right\}.
\end{equation}
Thus, conditionally on \(\filt_j\),
\[
 \lawof_{\mathbb Q}(Z_j\mid\filt_j)
 =\lawof\left(G_j-
 \frac{\overline{\sigma Q_j(r_{\bk},\theta)}}{\sqrt{j+1}}\right),
 \qquad
 \lawof_{\mathbb Q}(\Gamma_j\mid\filt_j)
 =\lawof_{\mathbb P}(\Gamma_j),
\]
where \(G_j\) is a standard complex Gaussian independent of
\(\Gamma_j\).  Equivalently, on \([j,j+1]\) the Brownian segment is a
standard segment plus a linear drift of magnitude at most \(Cj^{-1/2}\).
The reflection principle and the  Chernoff bound (applied to $\Gamma_j$) at the two thresholds
in \eqref{eq:goodtruncated} therefore give conditional probability of $ \operatorname{Tru}_{j,j+1}^{\complement}$ which is
at most \(Cj^{-c}\) with universal constants $C>0$ and $c>1$.  Summing over
\(2^{\bn}\leq j<2^{\bt}\) proves
\eqref{eq:tilted-truncation-cost}.

Now recall that
$\gamma_j=Z_j/\sqrt{|Z_j|^2+\Gamma_j}$, where
$\mathbb E\Gamma_j=\beta_j^2$.  Fix $\epsilon>0$. On $\operatorname{Tru}_{j,j+1}$  we have, for $j$ large,
\begin{equation}
\label{eq-250826c}
 |Z_j|^2\leq j^\epsilon,
 \qquad |\Gamma_j-\beta_j^2|\leq j^{1/2+\epsilon}.
\end{equation}
Taylor expansion at $\beta_j^2$ gives that on $\operatorname{Tru}_{j,j+1}$,
\begin{equation}\label{eq:correct-beta-expansion}
 \gamma_j=\frac{Z_j}{\beta_j}\left[
  1-\frac{\Gamma_j-\beta_j^2+|Z_j|^2}{2\beta_j^2}
  +O\!\left(j^{-1}\log j\right)\right].
\end{equation}
On the other hand, the complement of the event in \eqref{eq-250826c}  has probability  smaller than every prescribed negative power of
$j$ (uniformly under the shifted Gaussian law), while $|\gamma_j|\leq1$.
Since
\[
 \mathbb Q(Z_j\mid\filt_j)
 =-\frac{\overline{\sigma Q_j(r_{\bk},\theta)}}{\sqrt{j+1}},
 \quad
 \mathbb Q(Z_j^2\mid\filt_j)
 =\frac{\sigma^2\overline{Q_j(r_{\bk},\theta)^2}}{j+1},
 \quad
 \mathbb Q(|Z_j|^2\mid\filt_j)
 =1+\frac{|Q_j(r_{\bk},\theta)|^2}{j+1},
\]
and $\beta_j^2=\frac\beta2(j+1)$, \eqref{eq:correct-beta-expansion}
yields, uniformly in the conditioning,
\begin{align*}
 \mathbb Q(\gamma_j\mid\filt_j)
 &=-\sqrt{\frac2\beta}
   \frac{\overline{\sigma Q_j(r_{\bk},\theta)}}{j+1}
   +O(j^{-\alpha}),\\
 \mathbb Q(\gamma_j^2\mid\filt_j)&=O(j^{-2}),\\
 \mathbb Q(|\gamma_j|^2\mid\filt_j)
 &=\frac{2}{\beta(j+1)}+O(j^{-2})
  =\frac{2}{\beta j}+O(j^{-2}),
\end{align*}
for every fixed $1<\alpha<3/2$.  These are the unrestricted parts of  \eqref{eq:Qgamma}. On the other hand, using that the probability of $\operatorname{Tru}_{j,j+1}^\complement $ is $O(j^{-8})$,  Cauchy--Schwarz and
$|\gamma_j|\leq1$ give also the parts that include  $\operatorname{Tru}_{j,j+1}$.
\end{proof}

\begin{lemma}[Tilted stopped linearization]
\label{lem:tilted-stopped-linearization}
There are constants \(c,C>0\) such that, with
\(\epsilon_\bu^{\rm lin}=C\bu^2 2^{-\bu/4}\),
\begin{align}
 \sup_{\theta}\sup_{\bt\geq\bu}
 \mathbb Q_{\bk,\bu}^{\bt}\left(
  \one\{\operatorname{Tru}_{2^\bu,2^\bt}\}
  \sup_{\bu\leq\bm\leq\bt}|L_{\bk,\bu}^{\bm}|
  >\epsilon_\bu^{\rm lin}
\right)
 &\leq C\exp\{-c2^{\bu/4}\}, \As.
 \label{eq:tilted-stopped-linearization}
\end{align}
The same estimate holds (with the same constants) if $L^{\bm}_{\bk,\bu}(\theta)$ is replaced by $L^{\bm}_{\bk',\bu}(\theta')$ for any $\theta'$ and any $\bk'\in\overline{\N}$.
\end{lemma}

\begin{proof}
Fix \(\bt\geq\bu\), and in this proof write
\(\mathbb Q=\mathbb Q_{\bk,\bu}^{\bt}(\theta)\).  Let
\begin{equation}\label{eq:tilted-linearization-stopping-index}
 \tau_{\bu,\bt}^{\rm tr}
 :=\inf\left\{2^\bu\leq j<2^\bt:
       \operatorname{Tru}_{j,j+1}^{\complement}\right\}\wedge2^\bt
\end{equation}
and define
\[
 A_j=\one\{j<\tau_{\bu,\bt}^{\rm tr}\}
       \bigl(\log(1-\gamma_jQ_j)+\gamma_jQ_j\bigr),\qquad
 C_j=\one\{j<\tau_{\bu,\bt}^{\rm tr}\}Q_jZ_j
       \left((|Z_j|^2+\Gamma_j)^{-1/2}-\beta_j^{-1}\right).
\]
On \(\operatorname{Tru}_{2^\bu,2^\bt}\), the stopped and original
increments agree, and the identities used in the proof of
Lemma~\ref{lem:zetaphitail} give
\begin{equation}\label{eq:tilted-linearization-decomposition}
 L_{\bk,\bu}^{\bm}
 =\sqrt{2\beta}\,\Re\left(
   \sigma\sum_{j=2^\bu}^{2^\bm-1}(A_j-C_j)\right),
 \qquad \bu\leq\bm\leq\bt.
\end{equation}
For any fixed \(\alpha\in(5/4,3/2)\), the shifted one-step moments in
Lemma~\ref{lem:1st2ndQmoment} give
\begin{align}
 |\mathbb Q(A_j\mid\filt_j)|+|\mathbb Q(C_j\mid\filt_j)|
 &\leq C(j^{-\alpha}+j^{-3/2}\log j),
 \label{eq:tilted-linearization-compensators}\\
 \sum_{j\geq2^\bu}
 \mathbb Q\left(
 |A_j-\mathbb Q(A_j\mid\filt_j)|^2
 +|C_j-\mathbb Q(C_j\mid\filt_j)|^2
 \mathrel{\Big|}\filt_j\right)
 &\leq C\bu^2 2^{-\bu}.
 \label{eq:tilted-linearization-bracket}
\end{align}
The bound on \(\mathbb Q(A_j\mid\filt_j)\) holds because
\(A_j=\one\{j<\tau_{\bu,\bt}^{\rm tr}\}\bigl(-\tfrac12(\gamma_jQ_j)^2+O(|\gamma_j|^3)\bigr)\), and
\eqref{eq:Qgamma} bounds the conditional mean of
\(\gamma_j^2\one\{\operatorname{Tru}_{j,j+1}\}\) by \(Cj^{-\alpha}\),
while \(|\gamma_j|^3\leq C(j^{-1}\log j)^{3/2}\) on
\(\operatorname{Tru}_{j,j+1}\).
For \(C_j\), expand \((|Z_j|^2+\Gamma_j)^{-1/2}-\beta_j^{-1}\) as in
\eqref{eq:correct-beta-expansion}.  Under \(\mathbb Q\), \(\Gamma_j\)
is independent of \(Z_j\) and has mean \(\beta_j^2\), so the term
linear in \(\Gamma_j-\beta_j^2\) has conditional mean zero; the
truncation indicator changes this mean only by \(O(j^{-4})\).  The
remaining terms have conditional mean \(O(j^{-\alpha})\).  Summing \eqref{eq:tilted-linearization-compensators}
over \(j\geq2^\bu\) bounds the total drift by \(C2^{-\bu/4}\).
Freedman's inequality, applied to the centered parts with
\eqref{eq:tilted-linearization-bracket}, and
\eqref{eq:tilted-linearization-decomposition} prove
\eqref{eq:tilted-stopped-linearization}.
\end{proof}

The next lemma bounds the derivative of $Q_k$ on the closed disk by the oscillation \eqref{eq:osc-def} of the polynomial $\Phi_k^*$, and bounds that oscillation.  It is used in the good event \eqref{eq:globalgood} and to control the regularity of $Q$ in $r$ (Proposition \ref{prop:finebracket}).
\begin{lemma}\label{lem:derivativecontroltwo}
For every \(k\geq1\),
\begin{equation}\label{eq:Qprime-osc}
 \max_{|z|\leq1}|Q_k'(z)|\leq1+2k\,e^{\operatorname{Osc}_k}.
\end{equation}
There are constants \(C\) and \(\bn_*\) such that, for all
\(\bn\geq\bn_*\), all integers \(k\geq2^\bn\) and all \(y\geq1\),
\begin{equation}\label{eq:qkprimetail}
 \Prob\Big[\bigl\{\operatorname{Osc}_k>3\operatorname{Osc}_{2^\bn}+C(\log k+y^2)\bigr\}
 \cap\operatorname{Tru}_{2^\bn,k}\,\Big|\,\filt_{2^\bn}\Big]\leq2e^{-y^2}.
\end{equation}
As a consequence, for any \(p\geq1\) there is \(C_p\) such that, for
all \(\bn\geq\bn_*\) and all integers \(k\geq2^\bn\),
\begin{equation}\label{eq:qkprimebound}
 \Prob\Big[\bigl\{\max_{|z|\leq1}|Q_k'(z)|>k^{C_p}e^{3\operatorname{Osc}_{2^\bn}}\bigr\}
 \cap\operatorname{Tru}_{2^\bn,k}\,\Big|\,\filt_{2^\bn}\Big]\leq2k^{-p}.
\end{equation}
The same bounds hold with \(\Prob(\cdot\mid\filt_{2^\bn})\) replaced by
\(\Q^{\bt}_{\bk,\bn}\), for any \(\bk\) and \(\bt\).
\end{lemma}
\begin{proof}
\emph{Proof of \eqref{eq:Qprime-osc}.}  Since \(|\gamma_j|<1\)
for all \(j\), the zeros of \(\Phi_k\) lie in \(\D\), so \(\Phi_k^*\)
has no zeros in the closed disk (see \cite{SimonBook}).  Hence
\(Q_k=z\Phi_k/\Phi_k^*\) is analytic on a neighborhood of the closed
disk, and by the maximum modulus principle it suffices to bound
\(|Q_k'|\) on the unit circle.  There \(|\Phi_k(z)|=|\Phi_k^*(z)|\) by
\eqref{eq:szego}, and
\[
 Q_k'(z)=\frac{\Phi_k(z)}{\Phi_k^*(z)}+\frac{z\Phi_k'(z)}{\Phi_k^*(z)}
 -\frac{z\Phi_k(z)\Phi_k^{*\prime}(z)}{\Phi_k^*(z)^2}.
\]
The first term has modulus one.  The polynomials \(\Phi_k\) and
\(\Phi_k^*\) have degree at most \(k\) and the same modulus on the
circle, so by Bernstein's polynomial inequality \(|\Phi_k'|\) and
\(|\Phi_k^{*\prime}|\) are at most \(k\max_{|w|=1}|\Phi_k^*(w)|\) on the
circle.  Dividing by \(|\Phi_k^*(z)|\geq\min_{|w|=1}|\Phi_k^*(w)|\)
gives \eqref{eq:Qprime-osc}.

\emph{A pointwise bound.}  Fix \(\bn\geq\bn_*\), an integer
\(k\geq2^\bn\) and \(y\geq1\), and write \(\mathbb Q\) for either
\(\Prob(\cdot\mid\filt_{2^\bn})\) or \(\Q^{\bt}_{\bk,\bn}\).  For
\(|z|=1\), \eqref{eq:Phi} gives
\[
 I(z):=\log|\Phi_k^*(z)|-\log|\Phi_{2^\bn}^*(z)|
 =\sum_{j=2^\bn}^{k-1}\Re\log\bigl(1-\gamma_jQ_j(z)\bigr).
\]
Let \(\tau\) be the first \(j\geq2^\bn\) at which
\(\operatorname{Tru}_{j,j+1}\) fails, and let \(I^\tau(z)\) be the same
sum with its \(j\)-th term multiplied by \(\one\{j<\tau\}\).  Then
\(I^\tau(z)=I(z)\) on \(\operatorname{Tru}_{2^\bn,k}\).  We take
\(\bn_*\) so large that \(|\gamma_j|\leq1/64\) on
\(\operatorname{Tru}_{j,j+1}\) for \(j\geq2^{\bn_*}\), so that each term
of \(I^\tau(z)\) is at most \(1/32\) in absolute value.  Since
\(|Q_j(z)|=1\), the conditional second moment of the \(j\)-th term
given \(\filt_j\) is at most \(4\mathbb Q(|\gamma_j|^2\mid\filt_j)\leq
C/j\), by \eqref{eq:Qgamma} under the tilted measure.  Its conditional
mean is zero under \(\Prob\): the event \(\operatorname{Tru}_{j,j+1}\)
is invariant under rotations of the Brownian path on \([j,j+1]\) about
\(B_j\), so the phase of \(\gamma_j\) is uniform given \(\filt_j\) and
\(\operatorname{Tru}_{j,j+1}\).  Under \(\Q^{\bt}_{\bk,\bn}\), expanding
\(\log(1-w)=-w-\tfrac12w^2+O(|w|^3)\) and using \eqref{eq:Qgamma}, the
conditional mean is \(O(1/j)\).  Hence \(I^\tau(z)=D(z)+M(z)\), where
\(|D(z)|\leq C_0V\) and \(M(z)\) is a martingale with increments at most
\(1/16\) and bracket at most \(V\), with
\[
 V:=C_0\bigl(1+\log(k/2^\bn)\bigr).
\]
Freedman's inequality (Lemma~\ref{lem:freedman}, applied to \(16M\))
gives, for \(\mu>0\),
\begin{equation}\label{eq:osc-pointwise}
 \mathbb Q\bigl(|M(z)|>\mu\bigr)
 \leq2\exp\Bigl\{-\frac{\mu^2}{2(\mu/16+V)}\Bigr\}
 \leq2\exp\Bigl\{-\min\Bigl(\frac{\mu^2}{4V},\,4\mu\Bigr)\Bigr\}.
\end{equation}

\emph{Choice of the level.}  Put
\[
 R:=\log(40\pi k)+\operatorname{Osc}_{2^\bn}+2C_0V+y^2,\qquad
 \mu:=16V+\sqrt{8VR}+\tfrac12R,\qquad
 \lambda:=\mu+C_0V .
\]
These are \(\filt_{2^\bn}\)-measurable.  Since \(\mu\geq16V\) and
\(\mu\geq\sqrt{8VR}\), we have \(\mu^2/(4V)\geq2\mu+R\); since
\(\mu\geq R/2\), also \(4\mu\geq2\mu+R\).  As \(|I^\tau(z)|\leq
C_0V+|M(z)|\), \eqref{eq:osc-pointwise} gives, for every \(|z|=1\),
\begin{equation}\label{eq:osc-pointwise-2}
 \mathbb Q\bigl(|I^\tau(z)|>\lambda\bigr)
 \leq2e^{-2\mu-R}
 =2\exp\bigl\{-2\lambda-\log(40\pi k)-\operatorname{Osc}_{2^\bn}-y^2\bigr\}.
\end{equation}

\emph{Meshes.}  Let \(\mathcal Z_1\) consist of \(4k\) equally spaced
points on the unit circle, and \(\mathcal Z_2\) of
\(N:=\lceil10\pi k\,e^{\operatorname{Osc}_{2^\bn}+2\lambda}\rceil\)
equally spaced points.  Since \(N\leq20\pi k\,e^{\operatorname{Osc}_{2^\bn}+2\lambda}\),
a union bound and \eqref{eq:osc-pointwise-2} show that the event
\[
 \mathcal E:=\bigl\{|I^\tau(z)|\leq\lambda\ \text{for all }z\in\mathcal Z_1\cup\mathcal Z_2\bigr\}
\]
satisfies \(\mathbb Q(\mathcal E^\complement)\leq2e^{-y^2}\).  Write
\(M_k:=\max_{|z|=1}|\Phi_k^*(z)|\) and \(m_k:=\min_{|z|=1}|\Phi_k^*(z)|\).
By Bernstein's polynomial inequality, \(|\Phi_k^*(e^{i\theta})-\Phi_k^*(e^{i\theta'})|
\leq kM_k|\theta-\theta'|\).  Every point of the circle is within arc
distance \(\pi/(4k)\) of \(\mathcal Z_1\), so
\(M_k\leq\max_{\mathcal Z_1}|\Phi_k^*|+\tfrac\pi4M_k\), that is,
\(M_k\leq5\max_{\mathcal Z_1}|\Phi_k^*|\).  On
\(\operatorname{Tru}_{2^\bn,k}\cap\mathcal E\), therefore,
\[
 M_k\leq5M_{2^\bn}e^{\lambda}.
\]
Every point of the circle is within arc distance \(\pi/N\) of
\(\mathcal Z_2\), so on the same event
\[
 m_k\geq m_{2^\bn}e^{-\lambda}-\frac{\pi k}{N}M_k
 \geq m_{2^\bn}e^{-\lambda}-\frac{5\pi k}{N}M_{2^\bn}e^{\lambda}
 \geq\tfrac12m_{2^\bn}e^{-\lambda},
\]
by the choice of \(N\) and \(M_{2^\bn}=m_{2^\bn}e^{\operatorname{Osc}_{2^\bn}}\).
Hence, on \(\operatorname{Tru}_{2^\bn,k}\cap\mathcal E\),
\[
 \operatorname{Osc}_k\leq\operatorname{Osc}_{2^\bn}+2\lambda+\log10 .
\]
Finally, \(2\sqrt{8VR}\leq8V+R\) gives \(2\lambda\leq(2C_0+40)V+2R\),
and \(V\leq C_0(1+\log k)\), so
\(\operatorname{Osc}_{2^\bn}+2\lambda+\log10\leq3\operatorname{Osc}_{2^\bn}+C(\log k+y^2)\).
This proves \eqref{eq:qkprimetail}.

\emph{Proof of \eqref{eq:qkprimebound}.}  Take \(y^2=p\log k\) in
\eqref{eq:qkprimetail}.  Outside an event of probability at most
\(2k^{-p}\), \eqref{eq:Qprime-osc} gives
\(\max_{|z|\leq1}|Q_k'(z)|\leq1+2k\,e^{3\operatorname{Osc}_{2^\bn}}k^{C(1+p)}
\leq k^{C_p}e^{3\operatorname{Osc}_{2^\bn}}\), with \(C_p:=C(1+p)+2\).
\end{proof}

We now fix the derivative exponents used throughout the rest of the
proof.  Take \(p=10\) in
\eqref{eq:qkprimebound}, and let
\(C_{10}\) be the corresponding exponent in  \eqref{eq:qkprimebound}.
Next put \(C_{\rm in}:=3C+1\), with \(C\) the constant in
\eqref{eq:qkprimetail}.  Then
\begin{equation}\label{eq:fixed-time-derivative}
 \mathbb P\left(
  \operatorname{Osc}_{2^\bn}>C_{\rm in}\bn\log2\right)
 \longrightarrow_{\bn\to\infty}0.
\end{equation}
Indeed, fix \(\bn_0\geq\bn_*\) and apply \eqref{eq:qkprimetail} with
base \(2^{\bn_0}\), \(k=2^\bn\) and \(y^2=\bn\).  Outside an event of
probability at most
\(\mathbb P(\operatorname{Tru}_{2^{\bn_0}}^\complement)+2e^{-\bn}\),
\[
 \operatorname{Osc}_{2^\bn}
 \leq 3\operatorname{Osc}_{2^{\bn_0}}+C\bn(\log2+1)
 \leq 3\operatorname{Osc}_{2^{\bn_0}}+3C\bn\log2 .
\]
The random variable \(\operatorname{Osc}_{2^{\bn_0}}\) is almost surely
finite and does not depend on \(\bn\), so
\(3\operatorname{Osc}_{2^{\bn_0}}\leq\bn\log2\) with probability tending
to one.  Hence the \(\limsup\) of the left side of
\eqref{eq:fixed-time-derivative} is at most
\(\mathbb P(\operatorname{Tru}_{2^{\bn_0}}^\complement)\), which tends
to zero as \(\bn_0\to\infty\) by \eqref{Truest}.
Finally, choose \(C'>1\) so that
\begin{equation}\label{eq:derivative-constant-hierarchy}
 C'>C_{10}+3C_{\rm in}+2.
\end{equation}
With these choices, \(\mathscr G_\bn\) has probability tending to
one.  Indeed, \(\operatorname{Tru}_{2^\bn}\) does by \eqref{Truest}, the
linearization constraint does by Proposition~\ref{prop:zetaphitail}, and
the constraint on \(\operatorname{Osc}_{2^\bn}\) does by
\eqref{eq:fixed-time-derivative}.  On the latter,
\(k^{C_{10}}e^{3\operatorname{Osc}_{2^\bn}}\leq k^{C_{10}+3C_{\rm in}}\leq k^{C'}\)
for all \(k\geq2^\bn\), so \eqref{eq:qkprimebound} with \(p=10\), summed
over \(k\geq2^\bn\), shows that the remaining derivative constraint
fails on \(\operatorname{Tru}_{2^\bn}\) with probability at most
\(C2^{-9\bn}\).  Moreover, on \(\mathscr G_{\bn,\bn}\),
\eqref{eq:Qprime-osc} gives
\begin{equation}\label{eq:initial-derivative}
 \max_{|z|\leq1}|Q_{2^\bn}'(z)|\leq1+2^{\bn+1}2^{C_{\rm in}\bn}
 \leq2^{(C_{\rm in}+2)\bn}.
\end{equation}

We shall also need to know, under
the tilted measures and scale by scale, how it fails.  For
\(\bm\geq\bn\), \(0\leq r\leq1\) and \(\vartheta\in[0,2\pi]\), put
\begin{equation}\label{eq:complex-linearization-error}
 \mathcal L_{\bn,\bm}^{\mathbb C}(r,\vartheta)
 :=\zeta_{2^\bm}^{\mathbb C}(\bn;r,\vartheta)
   -2\log\Phi_{2^\bm}^*(re^{i\vartheta}),
\end{equation}
so that \(\mathcal L_{\bn,\bn}^{\mathbb C}=0\), and for \(\bu\geq\bn\) let
\begin{equation}\label{eq:scale-failure-event}
 \mathscr N_\bu:=\Big\{\max_{|z|\leq1}|Q_k'(z)|>k^{C'}\ \text{for some }2^\bu\leq k<2^{\bu+1}\Big\}
 \cup\Big\{\sup_{r,\vartheta}\big|\mathcal L_{\bn,\bu+1}^{\mathbb C}(r,\vartheta)
 -\mathcal L_{\bn,\bu}^{\mathbb C}(r,\vartheta)\big|>\epsilon_\bu^{\rm lin}\Big\},
\end{equation}
with \(\epsilon_\bu^{\rm lin}=C\bu^22^{-\bu/4}\) the tolerance of
Lemma~\ref{lem:tilted-stopped-linearization}, be the event that the good
event fails at scale \(2^\bu\), through the derivative bound or through
the linearization increment.

\begin{lemma}[The good event under the tilted measure]
\label{lem:good-event-tilted-tail}
One has
\begin{equation}\label{eq:good-event-cover}
 \mathscr G_{\bn,\bn}\cap\operatorname{Tru}_{2^\bn}
 \cap\mathscr G_\bn^\complement
 \subseteq\bigcup_{\bu\geq\bn}\mathscr N_\bu ,
\end{equation}
and there is \(C<\infty\) such that, on \(\mathscr G_{\bn,\bn}\), for
every \(\bu\geq\bn\) and every \(\bt\in\overline{\mathbb N}\) with
\(\bt\geq\bu+1\),
\begin{equation}\label{eq:good-event-tilted-tail}
 \mathbb Q_{\bk,\bn}^{\bt}\bigl(
  \mathscr N_\bu\cap\operatorname{Tru}_{2^\bn,2^{\bu+1}}\bigr)
 \leq C2^{-9\bu}.
\end{equation}
\end{lemma}

\begin{proof}
The inclusion \eqref{eq:good-event-cover} is read off from
\eqref{eq:globalgood}: on \(\mathscr G_{\bn,\bn}\cap\operatorname{Tru}_{2^\bn}\),
the event \(\mathscr G_\bn\) fails only if the derivative bound fails at
some \(k\geq2^\bn\), or if
\(\sup_{r,\vartheta,\bm\geq\bn}|\mathcal L_{\bn,\bm}^{\mathbb C}|>\epsilon_\bn\);
since \(\mathcal L_{\bn,\bn}^{\mathbb C}=0\) and
\(\sum_{\bu\geq\bn}\epsilon_\bu^{\rm lin}\leq\epsilon_\bn\) by
\eqref{eq:global-linearization-tolerance}, the latter forces
\(|\mathcal L_{\bn,\bu+1}^{\mathbb C}-\mathcal L_{\bn,\bu}^{\mathbb C}|>\epsilon_\bu^{\rm lin}\)
at some point for some \(\bu\geq\bn\).

For the derivative part of \(\mathscr N_\bu\), on \(\mathscr G_{\bn,\bn}\) we have
\(\operatorname{Osc}_{2^\bn}\leq C_{\rm in}\bn\log2\), so by the
choice \eqref{eq:derivative-constant-hierarchy} of \(C'\),
\(k^{C_{10}}e^{3\operatorname{Osc}_{2^\bn}}\leq k^{C'}\) for all
\(k\geq2^\bn\).  Since \(\operatorname{Tru}_{2^\bn,2^{\bu+1}}\subseteq\operatorname{Tru}_{2^\bn,k}\) for \(k<2^{\bu+1}\), \eqref{eq:qkprimebound} with \(p=10\), which holds
under \(\mathbb Q_{\bk,\bn}^{\bt}\), summed over \(2^\bu\leq k<2^{\bu+1}\),
gives at most \(C2^{-9\bu}\).

For the linearization part, write, for \(z=re^{i\vartheta}\),
\begin{equation}\label{eq:linearization-block-increment}
 G_\bu(z):=\mathcal L_{\bn,\bu+1}^{\mathbb C}(r,\vartheta)
 -\mathcal L_{\bn,\bu}^{\mathbb C}(r,\vartheta)
 =-2\sum_{j=2^\bu}^{2^{\bu+1}-1}
 \left(\frac{Z_j}{\beta_j}Q_j(z)+\log(1-\gamma_jQ_j(z))\right),
\end{equation}
by \eqref{eq:Z} and \eqref{eq:Phi}.  Since the derivative part of
\(\mathscr N_\bu\) has already been accounted for, we may assume that
\(\max_{|z|\leq1}|Q_j'(z)|\leq j^{C'}\) for all \(2^\bu\leq j<2^{\bu+1}\).
Since \(|\gamma_jQ_j(z)|<1\) on the closed disk, \(G_\bu\) is analytic
there.  On \(\operatorname{Tru}_{2^\bu,2^{\bu+1}}\), differentiating
term by term and using this derivative bound gives
\[
 \sup_{|z|\leq1}|G_\bu'(z)|\leq C\bu^{1/2}2^{(C'+1/2)\bu}=:D_\bu .
\]
Let \(\Theta_\bu\subset[0,2\pi]\) be a grid of mesh at most
\(\epsilon_\bu^{\rm lin}/(4D_\bu)\), so that \(|\Theta_\bu|\leq e^{C\bu}\).
By the maximum-modulus principle,
\[
 \sup_{|z|\leq1}|G_\bu(z)|\leq\max_{\vartheta\in\Theta_\bu}|G_\bu(e^{i\vartheta})|+\tfrac12\epsilon_\bu^{\rm lin}.
\]
Up to sign and the factor \(\sqrt{\beta/2}\), the real and imaginary
parts of \(G_\bu(e^{i\vartheta})\) are the linearization errors
\(L_{\infty,\bu}^{\bu+1}(\vartheta)\) at radius \(1\), for \(\sigma=1\)
and \(\sigma=i\).  Lemma~\ref{lem:tilted-stopped-linearization}, in the
form with \((\bk',\theta')=(\infty,\vartheta)\), bounds the
\(\mathbb Q_{\bk,\bu}^{\bt}\)-probability that either is large by
\(Ce^{-c2^{\bu/4}}\) at each \(\vartheta\in\Theta_\bu\).  A union bound
over \(\Theta_\bu\) and the tower property at \(\filt_{2^\bu}\) bound
the linearization part of \eqref{eq:good-event-tilted-tail} by
\(e^{C\bu}\cdot Ce^{-c2^{\bu/4}}\leq C2^{-9\bu}\).
\end{proof}

\subsection{The time dilation and the clock deviation}

Fix \(r_{\bk}\) and \(\theta\), and define
\begin{equation}
\label{eq-Rjnew}
 R_j=R_j(r_{\bk},\theta)
 :=\frac{1-|Q_j(r_{\bk},\theta)|^2}{1-r_{\bk}^2}
 =2^{\bk}\bigl(1-|Q_j(r_{\bk},\theta)|^2\bigr).
\end{equation}
We measure the Verblunsky evolution in \textbf{physical time}, the
logarithm of the index.  The step from \(j\) to \(j+1\) advances
physical time by \(1/(j+1)\).  Thus the index \(2^\bm\) corresponds to
physical time \(\bm\log2\) up to a bounded error, and
\(\lg_{r_\bk}(2^\bm)=(\bm\wedge\bk)\log2\) is physical time stopped at
the saturation scale \(\bk\).  Brownian time, the bracket
\eqref{eq-301024}, runs slower: if \(2^{\bn}\leq j\) and
\(s\in[j,j+1)\), then
\begin{equation}\label{eq:R-instantaneous-bracket}
 (j+1)\frac{\beta}{4}\frac{\dif}{\dif s}
 \left\langle\zeta_s(\bn;r_{\bk},\theta)\right\rangle
 =|Q_j(r_{\bk},\theta)|^2=1-2^{-\bk}R_j.
\end{equation}
So \(2^{-\bk}R_j\) is the rate, per unit of physical time, at which
Brownian time falls behind physical time.  We call \(R_j\) the
\textbf{time dilation}.
The time dilation is controlled by the derivative of \(Q_j\) on the
closed disk.  Since \(|Q_j|=1\) on the unit circle and \(Q_j\) is
analytic,
\[
 1-|Q_j(re^{i\theta})|^2
 \leq2\,|Q_j(e^{i\theta})-Q_j(re^{i\theta})|
 \leq2(1-r)\max_{|z|\leq1}|Q_j'(z)|,
\]
and since \(1-r\leq1-r^2\),
\begin{equation}\label{eq:julia-derivative}
 R_j(r_\bk,\theta)\leq2\max_{|z|\leq1}|Q_j'(z)|,
\end{equation}
uniformly in \(\bk\) and \(\theta\).  In particular, on the good event
\(\mathscr G_\bn\) of \eqref{eq:globalgood} one has \(R_j\leq2j^{C'}\)
for every \(j\geq2^\bn\), and \(R_{2^\bn}\leq2\cdot2^{(C_{\rm in}+2)\bn}\)
on \(\mathscr G_{\bn,\bn}\) by \eqref{eq:initial-derivative}; this is the weak bound on \(R\) that
appears as a hypothesis in Lemma~\ref{lem:buffered-one-block} and
Proposition~\ref{prop:brownian-clock-full-tower}.  
\begin{remark}
Since \(Q_j\) is a finite Blaschke product,
\(R_j(r_\bk,\theta)\to|Q_j'(e^{i\theta})|\) as \(\bk\to\infty\), for
fixed \(j\) and \(\theta\).  This is the finite Blaschke case of the
Julia--Carath\'eodory theorem; see \cite[Chapter~5]{AglerMcCarthyYoung}.
\end{remark}

Summing \eqref{eq:R-instantaneous-bracket} over \(j\) shows that,
before the saturation scale \(\bk\), the clock deviation is minus a
weighted sum of the time dilations.  More precisely,
by \eqref{eq-301024}, \eqref{eq:StLt} and \eqref{eq-Rjnew}, for
\(\bu\geq\bn\),
\begin{equation} \label{Qdecompsum}
S_{\bk,\bn}^{\bu}
=-2^{-\bk}
 \sum_{2^{\bn}\leq j<2^{\bu\wedge\bk}}\frac{R_j}{j+1}
 +\sum_{2^{\bk}\leq j<2^{\bu}}
       \frac{|Q_j(r_{\bk},\theta)|^2}{j+1}
 +e_{\bn,\bu\wedge\bk},
\qquad
 e_{\bn,\bm}:=\sum_{i=2^\bn+1}^{2^\bm}\frac1i-(\bm-\bn)\log2,
\end{equation}
where the second sum is empty when \(\bu\leq\bk\).  The error
\(e_{\bn,\bm}\) is deterministic, and comparing the harmonic sum with
\(\int_{2^\bn}^{2^\bm}\dif x/x\) and with
\(\int_{2^\bn+1}^{2^\bm+1}\dif x/x\) gives
\(-2^{-\bn}\leq e_{\bn,\bm}\leq0\).  Since \(R_j\geq0\), it follows
that
\begin{equation}\label{eq:S-sign}
 S_{\bk,\bl}^{\bm}\leq0
 \quad\text{and}\quad
 S_{\bk,\bl}^{\bm'}=S_{\bk,\bl}^{\bm}+S_{\bk,\bm}^{\bm'}\leq S_{\bk,\bl}^{\bm}
 \qquad\text{for }\bl\leq\bm\leq\bm'\leq\bk,
\end{equation}
that is, the deviation is nonpositive and nonincreasing in the upper
index before the saturation scale \(\bk\); after it, the function
\(\lg_{r_\bk}\) is constant and the bracket is nondecreasing, so
\(\bm\mapsto S_{\bk,\bl}^{\bm}\) is nondecreasing for \(\bm\geq\bk\).
In particular, \(\max_{\bl\leq\bm\leq\bu}|S_{\bk,\bl}^{\bm}|=-S_{\bk,\bl}^{\bu}\)
when \(\bu\leq\bk\).

Put
\begin{equation}\label{eq:step2-restoration-A}
A_{\bk,j}= A_{\bk,j}(\theta) :=r_{\bk}^2
 \frac{1-|\gamma_j|^2}
      {|1-\gamma_jQ_j(r_{\bk},\theta)|^2}.
\end{equation}
The Szeg\H{o} recurrence gives
\begin{equation} \label{Qrec}
 |Q_{j+1}(r_{\bk},\theta)|^2
 =r_{\bk}^2\left(
 1-\frac{1-|\gamma_j|^2}
 {|1-\gamma_jQ_j(r_{\bk},\theta)|^2}
 (1-|Q_j(r_{\bk},\theta)|^2)\right)
\end{equation}
and hence
\begin{equation}\label{eq:step2-restoration-affine}
 R_{j+1}=1+A_{\bk,j}R_j,
 \qquad
 R_J=R_N\prod_{s=N}^{J-1}A_{\bk,s}
 +\sum_{i=N}^{J-1}\prod_{s=i+1}^{J-1}A_{\bk,s}.
\end{equation}

The next lemma bounds the positive part of the clock deviation
\(S_{\bk,\bn}\) under a tilted measure \(\mathbb Q_{\bl,\bm}^{\bt}\)
whose radius \(r_\bl\) may differ from \(r_\bk\).

\begin{lemma} \label{lem:QVsmallkold}
Let $p\ge 0$. There are constants $C_p$ and $\bn_0(p)$ so that if
$\bt, \bl\in \overline{\N}$ and $\bm,\bn,\bk\in\N$ with $\bt \geq \bm$, $\bn\geq\bm$ and $\bn\geq\bn_0(p)$, then
\begin{equation}
\label{eq-11012026b}
\Q_{\bl,\bm}^{\bt} \bigl(\exp( p \bigl(2^{(\bn-\bk)_+}\bigr)\sup_{\bu\ge \bn} S_{\bk,\bn}^\bu) \one\{\operatorname{Tru}_{2^{\bn}}\}\,|\, \filt_{2^{\bn}} \bigr)\leq  C_p . 
\end{equation}
\end{lemma}

\begin{proof} Throughout, we write $Q_j$ for $Q_j(r_{\bk},\theta)$, even when working under the measure $\Q^{\bt}_{\bl,\bm}$. We begin with the case $\bn\geq \bk$. Put, for \(j\geq2^{\bn}\),
\[\Delta_j:=|\gamma_j|^2(1-|Q_j|^2)-2\Re(\gamma_j Q_j)-2\Re(\gamma_j^2 Q_j^2),\]
Then,
using  \eqref{Qrec}, on
\(\operatorname{Tru}_{2^{\bn}}\) we have that
 for \(1<\alpha<3/2\),
\begin{equation} \label{Qratiolin}
 \tfrac{1-|\gamma_j|^2}{|1-\gamma_j Q_j|^2} = 1-\Delta_j
+\mathcal{O}(j^{-\alpha}),
\end{equation}
By Lemma~\ref{lem:1st2ndQmoment},
\[
\Q^{\bt}_{\bl,\bm}\left( \Delta_j \vert \filt_j \right) = \mathcal{O}(j^{-1})
\]
and 
\[
|Q_{j+1}|^2 = r_{\bk}^2 |Q_j|^2 + r_{\bk}^2 \Delta_j (1-|Q_j|^2) +\mathcal{O}(j^{-1}). 
\]

This implies that on $\operatorname{Tru}_{2^{\bn}}$,
\[
\sum_{j\ge 2^{\bn}}\frac{|Q_{j+1}|^2}{j+1}
= r_{\bk}^2 \sum_{j\ge 2^{\bn}}\frac{ |Q_{j}|^2}{j+1} +(1- r_{\bk}^2)  M_{\bn,\infty}+ \mathcal{O}(2^{-\bn}) .
\]
where for $t\ge n$, 
\[
M_{\bn,\bt} \coloneqq 
\frac{1}{r_{\bk}^{-2}-1}
\sum_{j = 2^{\bn}}^t 
\frac{\Delta_j- \mathbb{Q}_{\bl,\bm}^{\bt}\left(\Delta_j\vert \filt_j \right)}{j+1} (1-|Q_j|^2) 
\]
is a discrete time martingale.
Using that $\sum_{j\ge 2^{\bn}}\big(\frac{|Q_{j+1}|^2}{j+1} - \frac{ |Q_{j}|^2}{j+1}\big) =  -\frac{ |Q_{2^\bn}|^2}{2^\bn}+ \sum_{j\ge 2^{\bn}}\frac{ |Q_{j}|^2}{j(j+1)} \le 2^{-\bn}$ and  $r_{\bk}^2= 1-2^{-\bk}$, after rearranging,  we obtain
\begin{equation}
\label{eq:sumQjplus1overjplus1}
\sum_{j\ge 2^{\bn}}\frac{|Q_{j+1}|^2}{j+1} 
= M_{\bn,\infty}+ \mathcal{O}(2^{\bk-\bn}) .
\end{equation}

Since $|\Delta_j| \ll |\gamma_j|$, by Lemma~\ref{lem:1st2ndQmoment} again, the quadratic variation of the martingale satisfies  for all $t\ge 2^{\bn}$, 
\[
[M_{\bn,\bt}]
\ll
2^{2\bk}\sum_{j\ge 2^{\bn}} 
\frac{\mathbb{Q}_{\bl,\bm}^{\bt}(\Delta_j^2\vert \filt_j)}{(j+1)^2} \ll 2^{2(\bk-\bn)}.
\]
Moreover, the increments of $M$ are bounded by $2^{\bk}\tfrac{|\gamma_j|}{j} \le 2^{\bk}\tfrac{\sqrt{\log j}}{j^{3/2}}$
on $\operatorname{Tru}_{2^{\bn}}$. 
By Freedman's martingale inequality, there is a constant $c>0$ so that we have the following tail bound; for all $x\geq 0$,
\begin{equation*} 
\Q^{\bt}_{\bl,\bm}\left( \{|M_{\bn,\infty}| > x\} \cap \operatorname{Tru}_{2^{\bn}}\,|\,\filt_{2^{\bn}}
\right)
\ll \exp\big(   -c\min\big\{  x^2 2^{-2(\bk-\bn)},  x  2^{-(\bk-3\bn/2)}/\bn \big\}\big) . 
\end{equation*}
Hence, for all $x\ge 0$,
\[
\Q^{\bt}_{\bl,\bm}\left( \{|M_{\bn,\infty}| > 2^{\bk-\bn} x\} \cap \operatorname{Tru}_{2^{\bn}}\,|\,\filt_{2^{\bn}}
\right)
\ll \exp\big(   -c\min\big\{  x^2  ,  x  2^{\bn/2}/\bn \big\}\big) . 
\]

Going back to \eqref{eq:StLt}, in the regime where $\bn\ge \bk$, using \eqref{eq:sumQjplus1overjplus1}, for any $\bu \geq \bn$,
\[
0\le S_{\bk,\bn}^{\bu} = \sum_{j=2^{\bn}}^{2^{\bu}-1} \frac{|Q_j|^2}{j+1} 
\leq \frac{|Q_{2^{\bn}}|^2}{2^{\bn}+1} + \sum_{j\ge 2^{\bn}}\frac{|Q_{j+1}|^2}{j+1}
\leq
2^{-\bn}+
\sum_{j\ge 2^{\bn}}\frac{|Q_{j+1}|^2}{j+1}
=
M_{\bn,\infty}+ \mathcal{O}(2^{\bk-\bn}).
\]
Note that the bound is independent of $\bu$, and therefore bounds $\sup_{\bu\ge \bn} S_{\bk,\bn}^{\bu}$ by $M_{\bn,\infty}+ \mathcal{O}(2^{\bk-\bn})$.
Then, we conclude that for any $p\ge 0$, if $\bn$ is sufficiently large depending on $p$, and $\bn\geq \bk$, \eqref{eq-11012026b} holds.

Next, assume that \(\bk\geq\bn\).  By \eqref{eq:S-sign},
\(S_{\bk,\bn}^{\bu}\leq 0\)
for $\bu\leq \bk$, and therefore we may and will assume that $\bu\geq \bk$.  By Lemma~\ref{lem:QVsmallkold} with $\bn=\bk$, for any $p\ge 0$,
\[
\Q_{\bl,\bm}^{\bt} \bigl(\exp( p \sup_{\bu\ge \bk} S_{\bk,\bk}^{\bu}) \one\{\operatorname{Tru}_{2^{\bn}}\} \,|\,\filt_{2^{\bn}}\bigr) \le C_p .
\]
On the other hand, by \eqref{Qdecompsum} and using that $R_j\geq 0$ and $e_{\bn,\bk}\leq 0$,
\[ S_{\bk,\bn}^{\bu}\leq \sum_{2^{\bk}\leq  j< 2^{\bu} }\frac{|Q_j|^2}{j+1}=S_{\bk,\bk}^\bu.\]
Therefore, we have proved \eqref{eq-11012026b} in all cases.
\end{proof}

We now take the tilt at the same radius \(r_{\bk}\) and angle
\(\theta\) as the time dilation \(R_j\).  The next lemma bounds the
conditional moments of the factors \(A_{\bk,j}\) in the recursion
\eqref{eq:step2-restoration-affine}.  The bound on
\(\mathbb E_{\mathbb Q}(A_{\bk,j}^p\mid\filt_j)\) for \(p<1\) shows that
the recursion forgets \(R_N\) and brings \(R_J\) back to the scale
\(J\), and the bound on \(\mathbb E_{\mathbb Q}(A_{\bk,j}\mid\filt_j)\)
gives the first moment of \(R_J\); see
Corollary~\ref{cor:step2-local-restoration}.

\begin{lemma}[Conditional moments of \(A_{\bk,j}\)]
\label{lem:step2-local-restoration}
Fix \(p\in(0,1)\).  Throughout this lemma, \(A_{\bk,j}=A_{\bk,j}(\theta)\) is as in \eqref{eq:step2-restoration-A} and
\[
 \mathbb Q=\mathbb Q_{\bk,\bm}^{\bt}(\theta),
 \qquad \bm\in\mathbb N,\quad \bt\in\overline{\mathbb N},\quad
 \bm\leq\bt,
\]
denotes the tilted conditional law defined by \eqref{eq:Qt}.
Then, with  $c_p=2p(1-p)/\beta$, there are $\delta>0$, a
constant $C$ and an integer $j_0$, all depending only on \(p,\beta\), such that
 for every \(j\geq\max(2^{\bm},j_0)\),
\begin{align}
 \mathbb E_{\mathbb Q}(A_{\bk,j}^p\mid\filt_j)
 &\leq1-\frac{c_p}{j}+Cj^{-1-\delta}, \As,
 \label{eq:step2-fractional-contraction}\\
 \mathbb E_{\mathbb Q}(A_{\bk,j}\mid\filt_j)
 &\leq1+Cj^{-1-\delta}, \As,
 \label{eq:step2-first-moment}\\
 \mathbb E_{\mathbb Q}(A_{\bk,j}^{-p}\mid\filt_j)
 &\leq1+\frac{C}{j}+C2^{-\bk}, \As.
 \label{eq:step2-negative-moment}
\end{align}
\end{lemma}
\noindent
Note that when \(j\geq2^{\bt}\), the conditional one-step law under $\mathbb Q_{\bk,\bm}^{\bt}$ is the
original untilted law.
\begin{proof}
Conditionally on \(\filt_j\), write
\(q=Q_j(r_{\bk},\theta)\).  When $j<2^{\bt}$,
\(Z_j\) is
a complex Gaussian with its original covariance and mean
\(-\overline{\sigma q}/\sqrt{j+1}\), while the \(\Gamma_j\)
variable retains its original gamma law and is independent of the shifted
Gaussian.  
On \(\{|\gamma_j|\leq1/2\}\), Taylor expansion of
\eqref{eq:step2-restoration-A} gives
\begin{align*}
 A_{\bk,j}^p=r_{\bk}^{2p}\bigl(&1-p|\gamma_j|^2
 +2p\Re(\gamma_jq)+p\Re(\gamma_j^2q^2)
+2p^2(\Re(\gamma_jq))^2+O_p(|\gamma_j|^3)\bigr).
\end{align*}
The conditional moment formulas in
Lemma~\ref{lem:1st2ndQmoment} therefore give, uniformly in \(|q|\leq1\),
\[
 \mathbb E_{\mathbb Q}(A_{\bk,j}^p\mid\filt_j)
 \leq1-\frac{2p}{\beta j}
       +\frac{2p^2|q|^2}{\beta j}
       -2p\sqrt{\frac2\beta}\frac{\Re(\sigma)|q|^2}{j}
       +O(j^{-1-\delta}).
\]
The last term of order \(j^{-1}\) is nonpositive for
\(\sigma\in\{1,i\}\), so it may be dropped, and the second term of order $j^{-1}$ is largest when $|q|=1$.  Hence there is a constant $C$ depending only on $p,\beta$ with
\[
 \mathbb E_{\mathbb Q}(A_{\bk,j}^p\mid\filt_j)
 \leq1-\frac{2p(1-p)}{\beta j}+Cj^{-1-\delta},
\]
which is \eqref{eq:step2-fractional-contraction}.  The same expansion at exponent
one proves \eqref{eq:step2-first-moment}. For $j\geq2^{\bt}$,
the Gaussian mean is zero; deleting the
displayed nonpositive drift proves the same estimates for the untilted
one-step law.

For the negative moment, the same Taylor expansion on
\(\{|\gamma_j|\leq1/2\}\), now at exponent \(-p\), gives
\begin{align*}
 A_{\bk,j}^{-p}=r_{\bk}^{-2p}\bigl(&1+p|\gamma_j|^2
 -2p\Re(\gamma_jq)-p\Re(\gamma_j^2q^2)
+2p^2(\Re(\gamma_jq))^2+O_p(|\gamma_j|^3)\bigr).
\end{align*}
Each term of order \(|\gamma_j|^2\) or \(|\gamma_j|\) has conditional
expectation of order \(j^{-1}\) by Lemma~\ref{lem:1st2ndQmoment}, uniformly in \(|q|\leq1\), whether or not
the tilt is active, and the remainder has conditional expectation \(O(j^{-3/2})\).  Hence
\[
 \mathbb E_{\mathbb Q}(A_{\bk,j}^{-p}\one\{|\gamma_j|\leq1/2\}\mid\filt_j)
 \leq r_{\bk}^{-2p}\left(1+\frac{C}{j}\right)
 \leq 1+\frac{C}{j}+C2^{-\bk},
\]
where in the last step we used \(r_{\bk}^{-2p}=(1-2^{-\bk})^{-2p}\leq1+C2^{-\bk}\).
Together with the bound on \(\{|\gamma_j|>1/2\}\) established next, this proves
\eqref{eq:step2-negative-moment}.

For the estimate on the event \(\{|\gamma_j|>1/2\}\), recall from
\eqref{eq:tilted-one-step-density} that the conditional law of
\((Z_j,\Gamma_j)\) given \(\filt_j\) under \(\mathbb Q\) has density
\[
 L_j=\exp\left\{-2\Re\frac{\sigma Q_jZ_j}{\sqrt{j+1}}-\frac{|Q_j|^2}{j+1}\right\}
\]
with respect to its law under \(\mathbb P\), where \(Q_j=Q_j(r_{\bk},\theta)\), and that \(L_j=1\) when \(j\geq2^{\bt}\).  Since
\(\Re(\sigma Q_jZ_j)\) is a centered Gaussian with variance \(|Q_j|^2/2\) under \(\mathbb P\),
\[
 \mathbb E_{\mathbb P}(L_j^2\mid\filt_j)
 =\exp\left(\frac{2|Q_j|^2}{j+1}\right)\leq e^2.
\]
On the other hand, \(|1-\gamma_jq|\leq2\) and
\(|1-\gamma_jq|\geq1-|\gamma_j|\geq\tfrac12(1-|\gamma_j|^2)\), so for every
\(h\in(0,1]\),
\[
 A_{\bk,j}^{h}\leq\frac{4^h}{(1-|\gamma_j|^2)^{h}},
 \qquad
 A_{\bk,j}^{-h}\leq\frac{4^h}{r_{\bk}^{2h}(1-|\gamma_j|^2)^{h}}.
\]
Under \(\mathbb P\), \(|\gamma_j|^2\) has density \(a_j(1-x)^{a_j-1}\) on
\([0,1]\) with \(a_j=\beta(j+1)/2\), so if \(a_j\geq4\) then \(a_j-2h\geq a_j/2\) and
\[
 \mathbb E_{\mathbb P}\left((1-|\gamma_j|^2)^{-2h}\one\{|\gamma_j|>1/2\}\right)
 =\frac{a_j}{a_j-2h}\left(\frac34\right)^{a_j-2h}
 \leq2\left(\frac34\right)^{a_j/2}.
\]
Let \(j_0\) be such that \(a_j\geq4\) for all \(j\geq j_0\).  Then, by Cauchy--Schwarz and \(r_{\bk}\geq1/2\), for every \(j\geq j_0\) and \(h\in\{p,1\}\),
\[
 \mathbb E_{\mathbb Q}\left(A_{\bk,j}^{\pm h}\one\{|\gamma_j|>1/2\}\mid\filt_j\right)
 \leq e\cdot4^h r_{\bk}^{-2h}\sqrt2\left(\frac34\right)^{a_j/4}
 \leq Ce^{-c\beta j},
\]
with \(C,c\) absolute constants, uniformly in \(\filt_j\) and in the tilt.  This is
smaller than every error term of order \(j^{-1-\delta}\) above, which completes the proof.
\end{proof}

\begin{corollary}[Restoration of \(R_J\)]
\label{cor:step2-local-restoration}
Fix \(p\in(0,1)\) and \(N=2^{\bl}\), and let \(J>N\) be an integer.  Let \(A_{\bk,j}\), \(\mathbb Q=\mathbb Q_{\bk,\bm}^{\bt}(\theta)\) with \(\bm\leq\bl\), \(c_p\), $\delta$, $C$ and $j_0$ be as in Lemma~\ref{lem:step2-local-restoration}, and let \(R_j=R_j(r_{\bk},\theta)\) be as in
\eqref{eq-Rjnew}.  Then, provided \(N\geq j_0\),
\begin{equation}\label{eq:step2-first-moment-evolution}
 \mathbb E_{\mathbb Q}(R_J\mid\filt_N)
 \leq (R_N+J-N)(1+CN^{-\delta}),
\end{equation}
and, with \(P_N=1\) and
\(P_{j+1}/P_j=1+Cj^{-1-\delta}\) for \(j\geq N\), one has
\begin{equation}\label{eq-11012026f}
 P_{j+1}^{-1}
 \mathbb E_{\mathbb Q}\left(
  R_{j+1}\one\{\operatorname{Tru}_{N,j+1}\}\mid\filt_j\right)
 \leq 1+P_j^{-1}R_j
       \one\{\operatorname{Tru}_{N,j}\}.
\end{equation}
Finally, if \(J=2^{\bu}\), then for every fixed \(B>0\), there exists $C_B=C_B(p,\beta)$ so that
\begin{equation}\label{eq:step2-restoration-tail}
 \mathbb Q\left(R_J>C_BJ\bu^B\mid\filt_N\right)
 \leq C_B\bu^{-B}
 +C_B\bu^{-Bp}\left(\frac{R_N}{J}\right)^p
                    \left(\frac NJ\right)^{c_p}.
\end{equation}
\end{corollary}
\begin{proof}
Since \(A_{\bk,s}\) is \(\filt_{s+1}\)-measurable and \(R_N\) is \(\filt_N\)-measurable, conditioning successively on \(\filt_{J-1},\filt_{J-2},\dots,\filt_N\) and applying \eqref{eq:step2-first-moment} at each step gives, for every \(N\leq i<J\),
\begin{equation}\label{eq:step2-product-first-moment}
 \mathbb E_{\mathbb Q}\Bigl(\prod_{s=i}^{J-1}A_{\bk,s}\mathrel{\Big|}\filt_N\Bigr)
 \leq\prod_{s=i}^{J-1}(1+Cs^{-1-\delta})
 \leq\prod_{s=N}^{\infty}(1+Cs^{-1-\delta})
 \leq1+CN^{-\delta}.
\end{equation}
Taking conditional expectations in \eqref{eq:step2-restoration-affine} and applying \eqref{eq:step2-product-first-moment} to each of the \(J-N+1\) products appearing there proves \eqref{eq:step2-first-moment-evolution}.

For \eqref{eq-11012026f}, since \(\operatorname{Tru}_{N,j+1}\subseteq\operatorname{Tru}_{N,j}\), the latter being \(\filt_j\)-measurable, the one-step recursion \(R_{j+1}=1+A_{\bk,j}R_j\) and \eqref{eq:step2-first-moment} give
\[
 \mathbb E_{\mathbb Q}\left(
  R_{j+1}\one\{\operatorname{Tru}_{N,j+1}\}\mid\filt_j\right)
 \leq\one\{\operatorname{Tru}_{N,j}\}\bigl(1+(1+Cj^{-1-\delta})R_j\bigr)
 =\one\{\operatorname{Tru}_{N,j}\}\Bigl(1+\frac{P_{j+1}}{P_j}R_j\Bigr).
\]
Dividing by \(P_{j+1}\geq1\) proves \eqref{eq-11012026f}.

For \eqref{eq:step2-restoration-tail}, the same successive conditioning with \eqref{eq:step2-fractional-contraction} in place of \eqref{eq:step2-first-moment} gives
\[
 \mathbb E_{\mathbb Q}\Bigl(\Bigl(\prod_{s=N}^{J-1}A_{\bk,s}\Bigr)^p\mathrel{\Big|}\filt_N\Bigr)
 \leq\prod_{s=N}^{J-1}\Bigl(1-\frac{c_p}{s}+Cs^{-1-\delta}\Bigr)
 \leq C\exp\Bigl(-c_p\sum_{s=N}^{J-1}\frac1s\Bigr)
 \leq C\left(\frac NJ\right)^{c_p}.
\]
By \eqref{eq:step2-restoration-affine}, on the event \(\{R_J>C_BJ\bu^B\}\) either the term \(R_N\prod_{s=N}^{J-1}A_{\bk,s}\) or the sum \(\sum_{i=N}^{J-1}\prod_{s=i+1}^{J-1}A_{\bk,s}\) exceeds \(C_BJ\bu^B/2\).  By \eqref{eq:step2-product-first-moment}, the sum has conditional expectation at most \((J-N)(1+CN^{-\delta})\leq CJ\), so Markov's inequality bounds the conditional probability that the sum exceeds \(C_BJ\bu^B/2\) by \(2C\bu^{-B}/C_B\).  Markov's inequality applied to the \(p\)-th power of the term, together with the previous display, bounds the conditional probability that the term exceeds \(C_BJ\bu^B/2\) by
\[
 \left(\frac{2}{C_BJ\bu^B}\right)^pR_N^p\,C\left(\frac NJ\right)^{c_p}
 =\frac{2^pC}{C_B^p}\,\bu^{-Bp}\left(\frac{R_N}J\right)^p\left(\frac NJ\right)^{c_p}.
\]
Combining the two bounds proves \eqref{eq:step2-restoration-tail}.
\end{proof}

For \(B>0\) and \(C_B\) as in \eqref{eq:step2-restoration-tail}, the event of typical time dilation at scale \(2^\bl\) is
\begin{equation}\label{eq:step2-radial-window}
\mathscr W_\bl^B
:=\left\{R_{2^\bl}(r_\bk,\theta)
\leq C_B2^\bl\bl^B\right\}.
\end{equation}

\begin{lemma}\label{lem:QVsmallk}
There exists a constant $C>0$ so that  for $\bt\in \overline{\N}$,  $\bn,\bk\in\N$, 
\begin{equation}
\label{eq-11012026dnew}
\Q_{\bk,\bn}^{\bt} \Bigl(\sup_{\bu\ge \bn}
\left(|S_{\bk,\bn}^\bu|
\one\{\operatorname{Tru}_{2^\bn,2^\bu}\}\right)
\Bigr) \le C\left( 1 + R_{2^{\bn}}2^{-\bk}\bk\right)
\end{equation}
and if $\bk\ge \bt\ge \bn$ then
\begin{equation}
\label{eq-11012026enew}
S_{\bk,\bn}^{\bt} \leq 0
\quad \text{and} \quad
\Q_{\bk,\bn}^{\bt} \bigl((|S_{\bk,\bn}^{\bt}|-2^{-\bn})_+
\one\{\operatorname{Tru}_{2^\bn,2^\bt}\}
\bigr) \le C\left(2^{\bt} + R_{2^{\bn}}\bt\right) 2^{-\bk}.
\end{equation}
\end{lemma}
The slack \(2^{-\bn}\) in \eqref{eq-11012026enew} is the
deterministic error \(e_{\bn,\bt}\) of \eqref{Qdecompsum}.  In every
application the deviation is tested at a level much larger than
\(2^{-\bn}\).

\begin{proof}
Since $R_j\ge 0$ and $-2^{-\bn}\leq e_{\bn,\bm}\leq 0$, and since
\(\operatorname{Tru}_{2^\bn,2^\bu}\subseteq\operatorname{Tru}_{2^\bn,j}\)
for \(j\leq2^\bu\), \eqref{Qdecompsum} gives
\begin{equation}
\label{eq-11012026c}
\Q_{\bk,\bn}^{\bt} \Bigl(\sup_{\bu\ge \bn}
\left((-S_{\bk,\bn}^{\bu})_+
\one\{\operatorname{Tru}_{2^\bn,2^\bu}\}\right)
\Bigr)
\le 2^{-\bn}+2^{-\bk} \sum_{2^{\bn}\le  j < 2^{\bk}}
\frac{\Q_{\bk,\bn}^{\bt}(R_{j}
\one\{\operatorname{Tru}_{2^{\bn},j}\}
)}{j+1}.
\end{equation}

Put \(N=2^{\bn}\) and apply
Corollary~\ref{cor:step2-local-restoration} with \(p=1/2\):
\eqref{eq:step2-first-moment-evolution}, with the truncation indicator
dropped, gives for \(N\leq j<2^{\bk}\)
\[
 \mathbb Q_{\bk,\bn}^{\bt}\left(
  R_{j+1}\one\{\operatorname{Tru}_{N,j+1}\}
  \mathrel{\Big|}\filt_N\right)
 \leq (R_N+j+1-N)(1+CN^{-\delta}).
\]
Consequently, separating the term \(j=N\), which equals \(R_N/(N+1)\),
\[
\sum_{N\leq j<2^{\bk}}
 \frac{\mathbb Q_{\bk,\bn}^{\bt}\left(
 R_{j}\one\{\operatorname{Tru}_{N,j}\}
 \mathrel{\big|}\filt_N\right)}{j+1}
\leq \frac{R_N}{N+1}+\sum_{N< j<2^{\bk}}\frac{(R_N+j-N)(1+CN^{-\delta})}{j+1}
\leq C\left(2^{\bk}+R_N\bk\right).
\]
Further, on the event $\widehat{\mathscr{G}}_{\bn}$, we have $R_{2^{\bn}} \ll 2^{C'\bn}$ with $C'$ the constant from \eqref{eq:globalgood}.
The proof of Lemma~\ref{lem:QVsmallkold}, applied only while
\(\operatorname{Tru}_{2^\bn,j}\) holds, controls the positive part of
the same supremum.  Together with \eqref{eq-11012026c}, this proves
\eqref{eq-11012026dnew}.  
Since
\(\operatorname{Tru}_{2^\bn}\subseteq
\operatorname{Tru}_{2^\bn,2^\bu}\) for every \(\bu\), the corresponding
bound with \(\operatorname{Tru}_{2^\bn}\) follows as a corollary.
Moreover, if $\bk\ge\bt\geq\bn$, the first part of
\eqref{eq-11012026enew} is \eqref{eq:S-sign}, and by
\eqref{Qdecompsum},
\((|S_{\bk,\bn}^{\bt}|-2^{-\bn})_+\leq2^{-\bk}\sum_{2^\bn\leq j<2^\bt}R_j/(j+1)\),
so the same argument as used for \eqref{eq-11012026dnew}, restricted
to the indices \(2^\bn\leq j<2^\bt\), gives the second part of \eqref{eq-11012026enew}.
\end{proof}

\subsection{Time regularity of the field}

From deviation results for the quadratic variation (and estimates on the linearization error) we can also give tail estimates for increments of the process $\varphi_{2^{\bj}}$ over adjacent $\bj$ under the tilted measure.
\begin{proposition}\label{prop:fieldincrement}
For any $p > 0$, there is a $C_p > 0$ and an $\bn_*$ sufficiently large so that
for any $\bt,\bj \geq \bn \geq\bn_*$
and any $\bk,\bl$
\begin{equation}\label{eq:fridaythe13}
\mathbb{Q}_{\bk,\bn}^{\bt}
\left(
\one\{\operatorname{Tru}_{2^{\bn}}\}
\exp\left(
  p| \varphi_{2^{\bj+1}}(r_{\bl},\theta)- \varphi_{2^{\bj}}(r_{\bl},\theta)|
\right)
\right)
\leq C_p.
\end{equation}
Further, if $\bj \geq \bl$, there is a $C>0$ so that
\begin{equation}\label{eq:fridaythe14}
\mathbb{Q}_{\bk,\bn}^{\bt}
\left(
\one\{\operatorname{Tru}_{2^{\bn}}\}
| \varphi_{2^{\bj+1}}(r_{\bl},\theta)- \varphi_{2^{\bj}}(r_{\bl},\theta)|\,
\right)
\leq C 2^{-(\bj-\bl)/6}.
\end{equation}
\end{proposition}

\begin{proof}
We work with a fixed $\theta$ and suppress the dependence on $\theta$ throughout. By the recurrence \eqref{eq:upsilon}, we have
\[
\varphi_{2^{\bj+1}}(r_{\bl}) - \varphi_{2^{\bj}}(r_{\bl})
= \sum_{s=2^{\bj}}^{2^{\bj+1}-1} 2\Re\{\sigma \log(1-\gamma_s Q_s(r_{\bl}))\}.
\]
Using $\log(1-x) = -x - \tfrac{1}{2}x^2 + \mathcal{O}(|x|^3)$ and $\gamma_s = Z_s/\sqrt{|Z_s|^2 + \Gamma_s}$, on the event $\operatorname{Tru}_{2^{\bn}}$ we have $|\gamma_j|^2\ll \tfrac{\log j}{j}$, so
\[
\varphi_{2^{\bj+1}}(r_{\bl}) - \varphi_{2^{\bj}}(r_{\bl})
= -2\Re\bigg\{\sigma \sum_{s=2^{\bj}}^{2^{\bj+1}-1} \gamma_s Q_s(r_{\bl})\bigg\}
- \Re\bigg\{\sigma \sum_{s=2^{\bj}}^{2^{\bj+1}-1} (\gamma_s Q_s(r_{\bl}))^2\bigg\}
+ \mathcal{O}(\bj^{3/2}2^{-\bj/2}).
\]
Write $\Q$ for $\mathbb Q_{\bk,\bn}^{\bt}$; all estimates below
are conditional on $\filt_{2^\bn}$.  As in the proof of
Lemma~\ref{lem:tilted-stopped-linearization}, each martingale is stopped
at the first index at which the truncation fails, as in
\eqref{eq:tilted-linearization-stopping-index}.  On
$\operatorname{Tru}_{2^\bn}$ the stopped and original sums agree.
We decompose each sum into its predictable and martingale parts under $\Q$. Writing $Q_s = Q_s(r_{\bl})$, we have $\sum_{s} \gamma_s Q_s = A^{(1)} + M^{(1)}$ where
\[
A^{(1)} \coloneqq \sum_{s=2^{\bj}}^{2^{\bj+1}-1} \Q(\gamma_s | \filt_s) Q_s(r_{\bl}),
\qquad
M^{(1)} \coloneqq \sum_{s=2^{\bj}}^{2^{\bj+1}-1} \big(\gamma_s - \Q(\gamma_s | \filt_s)\big) Q_s(r_{\bl}).
\]
By Lemma~\ref{lem:1st2ndQmoment}, $\Q(\gamma_s | \filt_s) = -\sqrt{\tfrac{2}{\beta}}\frac{\overline{\sigma Q_s(r_{\bk})}}{s+1} + \mathcal{O}(s^{-\alpha})$ for $\alpha > 1$, so
\begin{equation}\label{eq:A1}
A^{(1)} = -\sqrt{\tfrac{2}{\beta}}\overline{\sigma} \sum_{s=2^{\bj}}^{2^{\bj+1}-1} \frac{\overline{Q_s(r_{\bk})} Q_s(r_{\bl})}{s+1} + \mathcal{O}(2^{-\bj(\alpha-1)}).
\end{equation}
The martingale $M^{(1)}$ has increments $(\gamma_s - \Q(\gamma_s|\filt_s))Q_s$ with conditional variance $\asymp |Q_s|^2/(s+1)$.

For the quadratic term, we write $\sum_{s} (\gamma_s Q_s(r_{\bl}))^2 = A^{(2)} + M^{(2)}$ where
\[
A^{(2)} \coloneqq \sum_{s=2^{\bj}}^{2^{\bj+1}-1} \Q(\gamma_s^2 | \filt_s) Q_s(r_{\bl})^2,
\qquad
M^{(2)} \coloneqq \sum_{s=2^{\bj}}^{2^{\bj+1}-1} \big(\gamma_s^2 - \Q(\gamma_s^2 | \filt_s)\big) Q_s(r_{\bl})^2.
\]
By Lemma~\ref{lem:1st2ndQmoment}, $\Q(\gamma_s^2 | \filt_s) = \mathcal{O}(s^{-\alpha})$, so $A^{(2)} = \mathcal{O}(2^{-\bj(\alpha-1)})$.

For \eqref{eq:fridaythe13}, we use that $|Q| \leq 1$ from which it follows that both $A^{(1)}$ and $A^{(2)}$ are bounded by $C_p$ (deterministically).  For the martingales, we use Freedman's inequality Lemma \ref{lem:freedman},
noting that the quadratic variations of $M^{(1)}$ and $M^{(2)}$ are bounded  and their increments are uniformly small (and in fact vanishing with $\bn$) on $\operatorname{Tru}_{2^{\bn}}$.

For \eqref{eq:fridaythe14}, we use
Lemma~\ref{lem:QVsmallkold} conditionally on $\filt_{2^\bn}$ to bound the
quadratic variation $[M^{(i)}]$ of the stopped martingales, $i=1,2$.  For
$i=1$,
\[
  \Q([M^{(1)}]) \ll
  \Q
  \left(
   \sum_{s=2^{\bj}}^{2^{\bj+1}-1} \frac{|Q_s(r_{\bl})|^2}{s+1}
  \right)\ll 2^{-(\bj-\bl)}.
\]
For $A^{(1)}$, we use Cauchy-Schwarz and $|Q|\leq 1$ to give (taking $\alpha=4/3$ in \eqref{eq:A1})
\[
|A^{(1)}|^2 \ll \sum_{s=2^{\bj}}^{2^{\bj+1}-1} \frac{|Q_s(r_{\bl})|^2}{s+1} + \mathcal{O}(2^{-\bj/3}) \ll 2^{-(\bj-\bl)} + 2^{-\bj/3}.
\]
Thus from Doob's $L^2$ martingale inequality, we have
\[
\Q
\left(
(|A^{(1)}| + |M^{(1)}|)
\right)
\leq C 2^{-(\bj-\bl)/4} + C 2^{-\bj/6}.
\]
Stronger bounds can be derived for $A^{(2)}$ and $M^{(2)}$ by similar arguments.
\end{proof}
\begin{corollary}
  \label{cor:ezfieldincrement}
  There are $C,c > 0$ and an $\bn_*$ sufficiently large such that 
  for any $x \geq 1$, any $\bl \geq \bn \geq \bn_*$ and
  any $\bt,\bk$,
  \[
    \mathbb{Q}_{\bk,\bn}^{\bt}
    \left(
    \operatorname{Tru}_{2^{\bn}}
    \cap
    \{
    \sup_{\bu \geq \bl}| \varphi_{\infty}(r_{\bl},\theta)- \varphi_{2^{\bu}}(r_{\bl},\theta)|
    \geq x
    \}\,
    \right)
    \leq Ce^{-cx^{1/3}}.
\]
\end{corollary}
\begin{proof}
We split the supremum into dyadic increments.  The first $x^{1/3}$ of them are bounded using the exponential moment \eqref{eq:fridaythe13}, and the rest using Markov's inequality and \eqref{eq:fridaythe14}.
Fixing a sequence $(10 + \bu-\bl)^{-2}$ whose sum over $\bu \geq \bl$ is less than $1$ we can write 
\[
  \{
    \sup_{\bu \geq \bl}| \varphi_{\infty}(r_{\bl},\theta)- \varphi_{2^{\bu}}(r_{\bl},\theta)|
    \geq x
    \}
  \subset
  \bigcup_{\bu \geq \bl}
  \left\{
    | \varphi_{2^{\bu+1}}(r_{\bl},\theta)- \varphi_{2^{\bu}}(r_{\bl},\theta)|
    \geq x (10 + \bu-\bl)^{-2}
  \right\}.
\]
Hence applying the exponential bound to the $\bu$ with  $\bu - \bl \leq x^{1/3}$ and the Markov bound to the remainder, we get
\[
  \mathbb{Q}_{\bk,\bn}^{\bt}
  \left(
  \operatorname{Tru}_{2^{\bn}}
  \cap
  \{
  \sup_{\bu \geq \bl}| \varphi_{\infty}(r_{\bl},\theta)- \varphi_{2^{\bu}}(r_{\bl},\theta)|
  \geq x
  \}\,
  \right)
  \ll e^{-cx^{1/3}}.
\] 
\end{proof}
\begin{corollary}\label{cor:fieldincrement}
There is a $C>0$ and an $\bn_*$ so that
for any $\bt,\bj \geq \bn \geq \bn_*$, any $\bk,\bl$, and all $x > 0$
\[
\mathbb{Q}_{\bk,\bn}^{\bt}
\left(
\operatorname{Tru}_{2^{\bn}}
\cap \{
| \varphi_{2^{\bj+1}}(r_{\bl},\theta)- \varphi_{2^{\bj}}(r_{\bl},\theta)|
> x
\}\,
\right)
\leq
C
\min \left\{
e^{-x}, \frac{2^{(\bl-\bj)/6}}{x}
\right\}.
\]
\end{corollary}

\subsection{Radial regularity of the field}
\label{subsec-2.3}
This section quantifies the second part of Remark~\ref{rem-almostgaussian}: for indices $\ell\leq 1/(1-r^2)$, the recursion $Q_\ell(r)$ stays close to $Q_\ell(1)$.  The key estimate, Proposition~\ref{prop:finebracket}, is a tail bound for $\sum_\ell|Q_{\ell+1}(r_\bv)-Q_{\ell+1}(1)|/(\ell+1)$ over $\ell\leq2^\bv=1/(1-r_\bv^2)$. Throughout this section and the rest of the paper, we set \(\mathfrak c:=4C'\) with $C'$ as in \eqref{eq:derivative-constant-hierarchy}.
\begin{proposition}\label{prop:finebracket}
There exist \(\delta>0\) and \(\bn_*>0\)
such that for all $\bu \geq \bv$, all $x > 1$, all $\bv \geq \mathfrak{c}\bn \geq \mathfrak{c}\bn_*$
\[
\Q_{\bu,\bn}^{\bv} \left( 
  \widehat{\mathscr{G}}_{\bn}
  \cap
  \left\{
    \sum_{\ell =2^{\bn}}^{2^{\bv}}
    \frac{\left|Q_{\ell+1}(r_{\bv},\theta) -  Q_{\ell+1}(1,\theta)\right|}
    {\ell+1}
    > x
   \right\}
    ~\middle|~ \filt_{2^{\bn}}
   \right)
  \ll 2^{-\delta x}.
\]
\end{proposition}
\begin{proof}
We work $\theta$ by $\theta$ and suppress $\theta$ from the notation throughout, writing $Q_\ell(r)$ for $Q_\ell(r,\theta)$.
Observe that on the unit circle, $|Q_\ell(1)| = 1$ for all $\ell$.
We introduce the quantity
\begin{equation}
\label{eq-theotherDelta}
\Delta_{\bk,\ell} \coloneqq Q_\ell(r_{\bk}) \overline{Q_\ell(1)} - 1,
\end{equation}
where $\bk$ will be either $\bu$ or $\bv$; we write $\bk$ to illustrate the parts of the argument that are common to both cases.
Define the sum
\[
S \coloneqq \sum_{\ell = 2^{\bn}}^{2^{\bv}} \frac{|Q_{\ell+1}(r_{\bv}) - Q_{\ell+1}(1)|}{\ell+1}
= \sum_{\ell = 2^{\bn}}^{2^{\bv}} \frac{|\Delta_{\bv,\ell+1}|}{\ell+1}.
\]
Since $|Q_\ell(z)| \leq 1$ for all $\ell$ and $|z| \leq 1$, we have $|Q_{\ell+1}(r_{\bv}) - Q_{\ell+1}(1)| \leq 2$. In particular, $S \leq 4 \sum_{\ell=2^{\bn}}^{2^{\bv}} \frac{1}{\ell+1} \leq 4 (\bv - \bn)\log 2 $, so we only need to prove the tail bound for $x \leq 4  (\bv - \bn) \log 2 \ll \bv$. This means we may ignore events of probability $e^{-c\bv}$ for any fixed $c > 0$, and further by monotonicity it in fact suffices to prove the bound for $x < c_1 \bv$ for any fixed $c_1 > 0$ sufficiently small.

Fix a small constant $\epsilon<1/(32C)$, with $C$ as in \eqref{eq:qkprimetail}, to be determined below; we may suppose $\epsilon \bv$ is an integer.  We split $S$ into an early part and a late part:
\[
S = S_{\mathrm{early}} + S_{\mathrm{late}}, \quad
S_{\mathrm{early}} \coloneqq \sum_{\ell = 2^{\bn}}^{2^{\epsilon\bv}} \frac{|\Delta_{\bv,\ell+1}|}{\ell+1}, \quad
S_{\mathrm{late}} \coloneqq \sum_{\ell = 2^{\epsilon\bv}+1}^{2^{\bv}} \frac{|\Delta_{\bv,\ell+1}|}{\ell+1}.
\]
For the early part, we use the mean value theorem: $|Q_{\ell+1}(r_{\bv}) - Q_{\ell+1}(1)| \leq \sup_{|z| \leq 1} |Q_{\ell+1}'(z)| \cdot (1 - r_{\bv})$. 

Let $\mathcal{B}$ be the event that
\[
 \operatorname{Osc}_\ell\leq\tfrac12\bv\log2
 \qquad\text{for all }2^\bn<\ell\leq2^{\epsilon\bv}+1 .
\]
On $\widehat{\mathscr G}_\bn$ we have
$\operatorname{Osc}_{2^\bn}\leq C_{\rm in}\bn\log2\leq\tfrac1{12}\bv\log2$,
since $\bv\geq\mathfrak c\bn=4C'\bn\geq12C_{\rm in}\bn$ by
\eqref{eq:derivative-constant-hierarchy}.  Apply \eqref{eq:qkprimetail},
which holds under $\Q_{\bu,\bn}^{\bv}$, with $y^2:=\bv\log2/(16C)$.
Outside the exceptional event there,
\[
 \operatorname{Osc}_\ell
 \leq\tfrac14\bv\log2+C(\epsilon\bv+1)\log2+\tfrac1{16}\bv\log2
 \leq\tfrac12\bv\log2,
\]
since $C\epsilon<1/32$ and $\bv$ is large.  A union bound over the at
most $2^{\epsilon\bv+1}$ values of $\ell$ therefore gives
\begin{equation}\label{eq:Bbound}
\Q_{\bu,\bn}^{\bv}(\mathcal{B}^{\complement} \cap \widehat{\mathscr{G}}_{\bn})
\leq2^{\epsilon\bv+2}\,2^{-\bv/(16C)}
\leq4\cdot2^{-\bv/(32C)} .
\end{equation}
On $\mathcal B$, \eqref{eq:Qprime-osc} gives, for
$2^\bn\leq\ell\leq2^{\epsilon\bv}$,
\[
 \max_{|z|\leq1}|Q_{\ell+1}'(z)|\leq1+2(\ell+1)2^{\bv/2}\leq2^{(1/2+\epsilon)\bv+3}.
\]
Since $1-r_\bv\leq1-r_{\bv}^2= 2^{-\bv}$, on $\mathcal B$,
\[
|Q_{\ell+1}(r_{\bv}) - Q_{\ell+1}(1)| \leq 2^{-(1/2-\epsilon)\bv+3}.
\]
Therefore, as $\epsilon<1/4$,
\[
S_{\mathrm{early}} \leq \sum_{\ell = 2^{\bn}}^{2^{\epsilon\bv}} \frac{2^{-(1/2-\epsilon)\bv+3}}{\ell+1} \leq 8(\epsilon\bv+1)\,
2^{-(1/2-\epsilon)\bv}
=\mathcal{O}(\bv2^{-\bv/4}),
\]
as needed.

It remains to bound $S_{\mathrm{late}}$. Splitting the sum at
\begin{equation}
\label{eq-L}
 L \coloneqq 2^{\bv} \cdot 2^{-cx} \quad \mbox{\rm  with $c=1/(12 \log2)$,}
\end{equation}
the contribution from $\ell \geq L$ is bounded deterministically by
\[
\sum_{\ell = L}^{2^{\bv}} \frac{4}{\ell+1} \leq 4\log(2^{\bv}/L) = 4\log 2 \cdot cx\leq x/3.
\]
Thus,
it suffices to show that
\begin{equation}\label{eq:Sprimedef}
  \Q_{\bu,\bn}^{\bv}(\widehat{\mathscr{G}}_{\bn} \cap \{S '> x/3\})
  \leq 2^{-\delta x}
  \quad \text{where} \quad
  S' \coloneqq \sum_{\ell = 2^{\epsilon\bv}}^{L} \frac{|Q_{\ell+1}(r_{\bv}) - Q_{\ell+1}(1)|}{\ell+1}.
\end{equation}
Then on the event $S' \leq \frac{x}{3}$ with $x\geq 1$,
\[
S = S_{\mathrm{early}} + S_{\mathrm{late}} \leq \mathcal{O}(\bv2^{-\bv/4}) + \frac{x}{3}+ S' \leq x.
\]

It thus remains to bound $S'$ as in ~\eqref{eq:Sprimedef}.
By the recurrence~\eqref{eq:Q}, we have
\begin{align*}
Q_{\ell+1}(r_{\bk}) = r_{\bk}e^{i\theta} \frac{Q_\ell(r_{\bk}) - \overline{\gamma}_\ell}{1 - \gamma_\ell Q_\ell(r_{\bk})},
\quad \text{and} \quad
Q_{\ell+1}(1) = e^{i\theta} \frac{Q_\ell(1) - \overline{\gamma}_\ell}{1 - \gamma_\ell Q_\ell(1)}.
\end{align*}
Multiplying, we obtain
\begin{align*}
\Delta_{\bk,\ell+1} + 1 &= Q_{\ell+1}(r_{\bk}) \overline{Q_{\ell+1}(1)}
= r_{\bk} \frac{Q_\ell(r_{\bk}) - \overline{\gamma}_\ell}{1 - \gamma_\ell Q_\ell(r_{\bk})}
\cdot \frac{\overline{Q_\ell(1)} - \gamma_\ell}{1 - \overline{\gamma}_\ell \overline{Q_\ell(1)}}\\
&= r_{\bk} \frac{(Q_\ell(r_{\bk}) - \overline{\gamma}_\ell)(\overline{Q_\ell(1)} - \gamma_\ell)}{(1 - \gamma_\ell Q_\ell(r_{\bk}))(1 - \overline{\gamma}_\ell \overline{Q_\ell(1)})}
= r_{\bk} \frac{(\Delta_{\bk,\ell} + 1) + \omega_{\bk,\ell} + |\gamma_\ell|^2}{1 + \omega_{\bk,\ell} + |\gamma_\ell|^2 (\Delta_{\bk,\ell} + 1)},
\end{align*}
where we define
\[
\omega_{\bk,\ell} \coloneqq -\gamma_\ell Q_\ell(r_{\bk}) - \overline{\gamma}_\ell \overline{Q_\ell(1)} = \mathcal{O}(\gamma_\ell).
\]
Expanding the denominator in powers of $\gamma_\ell$, we have
\begin{align}
\Delta_{\bk,\ell+1} + 1
&= r_{\bk} \bigl((\Delta_{\bk,\ell} + 1) + \omega_{\bk,\ell} + |\gamma_\ell|^2\bigr)
\bigl(1 - \omega_{\bk,\ell} - |\gamma_\ell|^2 (\Delta_{\bk,\ell} + 1) + \omega_{\bk,\ell}^2 + \mathcal{O}(\gamma_\ell^3)\bigr) \notag\\
&= r_{\bk} \Bigl[ (\Delta_{\bk,\ell} + 1) - (\Delta_{\bk,\ell} + 1)\omega_{\bk,\ell} + \omega_{\bk,\ell} \notag\\
&\qquad - (\Delta_{\bk,\ell} + 1)^2 |\gamma_\ell|^2 + (\Delta_{\bk,\ell} + 1)\omega_{\bk,\ell}^2 - \omega_{\bk,\ell}^2 + |\gamma_\ell|^2 + \mathcal{O}(\gamma_\ell^3) \Bigr] \notag\\
&= r_{\bk} \Bigl[ (\Delta_{\bk,\ell} + 1) - \Delta_{\bk,\ell} \omega_{\bk,\ell}
+ \Delta_{\bk,\ell} \omega_{\bk,\ell}^2 - (\Delta_{\bk,\ell} + 1)^2 |\gamma_\ell|^2 + |\gamma_\ell|^2 + \mathcal{O}(\gamma_\ell^3) \Bigr] \notag\\
&= r_{\bk} + \Delta_{\bk,\ell} \cdot r_{\bk}\bigl(1 - \omega_{\bk,\ell} + \omega_{\bk,\ell}^2 - (2 + \Delta_{\bk,\ell})|\gamma_\ell|^2\bigr) + \mathcal{O}(\gamma_\ell^3). \label{eq:Deltaexpansion}
\end{align}

We define $P_{\bk,\ell}\neq 0$ by the ratio
\begin{equation}\label{eq:Pdef}
\frac{P_{\bk,\ell+1}}{P_{\bk,\ell}} \coloneqq r_{\bk}\bigl(1 - \omega_{\bk,\ell} + \omega_{\bk,\ell}^2 - (2 + \Delta_{\bk,\ell})|\gamma_\ell|^2\bigr),
\end{equation}
so that~\eqref{eq:Deltaexpansion} becomes
\[
\Delta_{\bk,\ell+1} = (r_{\bk}-1) + \Delta_{\bk,\ell} \cdot \frac{P_{\bk,\ell+1}}{P_{\bk,\ell}} + \mathcal{O}(\gamma_\ell^3).
\]
Dividing by $P_{\bk,\ell+1}$,
\[
\frac{\Delta_{\bk,\ell+1}}{P_{\bk,\ell+1}} = \frac{r_{\bk}-1}{P_{\bk,\ell+1}} + \frac{\Delta_{\bk,\ell}}{P_{\bk,\ell}} + \mathcal{O}\Bigl(\frac{\gamma_\ell^3}{P_{\bk,\ell+1}}\Bigr).
\]
Note that $P_{\bk,\ell+1}$ depends on $\Delta_{\bk,\ell}$.
Iterating from $j = 2^{\epsilon\bv}$ to $\ell-1$, we obtain
\begin{equation}\label{eq:Deltasolution}
\frac{\Delta_{\bk,\ell}}{P_{\bk,\ell}} = \frac{\Delta_{\bk,2^{\epsilon \bv}}}{P_{\bk,2^{\epsilon\bv}}} + \sum_{j=2^{\epsilon \bv}}^{\ell-1} \frac{r_{\bk}-1}{P_{\bk,j+1}} + \mathcal{O}\Bigl(\sum_{j=2^{\epsilon \bv}}^{\ell-1} \frac{\gamma_j^3}{P_{\bk,j+1}}\Bigr).
\end{equation}

\begin{claim}\label{claim:Pmoment}  
 For any $\eta > 0$ there exists $\delta=\delta(c_1,\eta,\beta) > 0$ such that with $L$ as in \eqref{eq-L} and $x < c_1 \bv$, for $j \leq L$ and $h\ll \bv$, $\bk \in \{\bu, \bv\}$,  $\bu \geq \bv \geq \mathfrak{c}\bn$,
\begin{equation}\label{eq:Pmoment}
\Q_{\bu,\bn}^{\bv}\left[\widehat{\mathscr{G}}_\bn 
\bigcap 
\Bigg\{
\exists~(j,\ell) : 2^{\epsilon \bv}\leq j <\ell\leq L, \frac{|P_{\bk,\ell}|}{|P_{\bk,j}|}\geq \left(\frac{j}{\ell}\right)^{1/(2\beta)} \left(\frac{L}{j}\right)^\eta 2^{h} \Bigg\}
~\middle|~ \filt_{2^{\bn}}\right] \ll e^{-\delta h} \vee e^{-\delta x}.
\end{equation}
\end{claim}
\begin{proof}[Proof of Claim \ref{claim:Pmoment}]
Throughout this proof, $\Q$ denotes $\Q_{\bu,\bn}^{\bv}$, conditional
on $\filt_{2^\bn}$.  Besides the stopping time $\tau$ below, all processes
are stopped at the first index at which the truncation fails, as in the
proof of Lemma~\ref{lem:tilted-stopped-linearization}.  Since
$\widehat{\mathscr G}_\bn\subseteq\operatorname{Tru}_{2^\bn}$, the stopped
and original variables in \eqref{eq:Pmoment} agree on
$\widehat{\mathscr G}_\bn$.
We compute the conditional mean and conditional second moment of $\log|P_{\bk,\ell+1}/P_{\bk,\ell}|$ under the $\Q$ measure. Note that $P_{\bk,\ell}$ is complex-valued (since $\Delta_{\bk,\ell} \in \C$), so we work with the modulus.

We introduce a stopping time $\tau$ to control the dependence of $P_{\bk,\ell}$ on $\Delta_{\bk,\ell}$: let $\epsilon' > 0$ be a small constant, and define
\[
\tau \coloneqq \inf\{\ell \geq 2^{\epsilon\bv} : |\Delta_{\bu,\ell}| > \epsilon' \text{ or } |\Delta_{\bv,\ell}| > \epsilon'\}.
\]
Throughout the proof we may assume that $x$ is a sufficiently large constant depending only on $\epsilon'$.

Recalling $\mathcal{B}$ from \eqref{eq:Bbound}, note that $\Delta_{\bk,2^{\epsilon\bv}}$ with $\bk \in \{\bu, \bv\}$ is bounded (on $\mathcal{B}\cap \mathscr{G}_\bn$) by
\begin{equation}\label{eq:InitialDeltaBound}
|\Delta_{\bk,2^{\epsilon\bv}}| 
= |Q_{2^{\epsilon\bv}}(r_{\bk})\overline{Q_{2^{\epsilon\bv}}(1)} -1| 
\leq (1-r_{\bk})\sup_{|z| \leq 1} |Q'_{2^{\epsilon\bv}}(z)|  \ll 2^{-\bv/4}
\end{equation}
And moreover, this holds with probability better than $1-e^{-\delta \bv}.$

From~\eqref{eq:Pdef}, we have
\[
\frac{P_{\bk,\ell+1}}{P_{\bk,\ell}} = r_{\bk}(1 - \omega_{\bk,\ell} + \omega_{\bk,\ell}^2 - (2+\Delta_{\bk,\ell})|\gamma_\ell|^2).
\]
Writing $X_{\bk,\ell} \coloneqq -\omega_{\bk,\ell} + \omega_{\bk,\ell}^2 - (2+\Delta_{\bk,\ell})|\gamma_\ell|^2$, we have $|X_{\bk,\ell}| = \mathcal{O}(|\gamma_\ell|)$. For the modulus,
\[
\log\left|\frac{P_{\bk,\ell+1}}{P_{\bk,\ell}}\right| = \log r_{\bk} + \frac{1}{2}\log|1 + X_{\bk,\ell}|^2 = \log r_{\bk} + \frac{1}{2}\log(1 + 2\Re(X_{\bk,\ell}) + |X_{\bk,\ell}|^2).
\]
(Since $|\Delta_{\bk,\ell}| \leq 2$ and is multiplied by $|\gamma_\ell|^2$,
on $\operatorname{Tru}_{2^\bn}$ we are ensured that the argument of the logarithm in the last display
 is bounded away from zero.)
Expanding for small $|X_{\bk,\ell}|$,
\begin{equation}\label{eq:logPexpand}
\log\left|\frac{P_{\bk,\ell+1}}{P_{\bk,\ell}}\right| = \log r_{\bk} + \Re(X_{\bk,\ell}) + \frac{|X_{\bk,\ell}|^2 - 2(\Re X_{\bk,\ell})^2}{2} + \mathcal{O}(|\gamma_\ell|^3).
\end{equation}
Since $X_{\bk,\ell}=-\omega_{\bk,\ell}+O(|\gamma_\ell|^2)$,
\[
 \frac{|X_{\bk,\ell}|^2-2(\Re X_{\bk,\ell})^2}{2}
 =-\frac12\Re(\omega_{\bk,\ell}^2)+O(|\gamma_\ell|^3).
\]
Substituting in \eqref{eq:logPexpand} and using the expression for $X_{\bk,\ell}$, we obtain that
\begin{equation}
\label{eq:logPexpandcode}
\log\left|\frac{P_{\bk,\ell+1}}{P_{\bk,\ell}}\right|
= \log r_{\bk} - \Re(\omega_{\bk,\ell})
+\frac12\Re(\omega_{\bk,\ell}^2)
-(2+\Re(\Delta_{\bk,\ell}))|\gamma_\ell|^2
+\mathcal{O}(|\gamma_\ell|^3).
\end{equation}
Note that
$\frac12\mathbb Q\Re(\omega_{\bk,\ell}^2)
=2\Re(\Delta_{\bk,\ell}+1)/(\beta\ell)+O(\ell^{-\alpha})$.

\textbf{Conditional moment computations for $\omega_{\bk,\ell}$:} We now compute the $\Q$-conditional expectations using Lemma~\ref{lem:1st2ndQmoment}. The key estimates are:
\begin{align*}
\Q[\gamma_\ell \mid \filt_\ell] &= -\sqrt{\frac{2}{\beta}} \frac{\overline{\sigma Q_\ell(r_{\bu})}}{\ell+1} + O(\ell^{-\alpha}), \\
\Q[\gamma_\ell^2 \mid \filt_\ell] &= O(\ell^{-\alpha}), \\
\Q[|\gamma_\ell|^2 \mid \filt_\ell] &= \frac{2}{\beta \ell} + O(\ell^{-\alpha}),
\end{align*}
where $\alpha>1$; note that the bias in the first moment is determined by $r_{\bu}$, the point at which the measure $\Q_{\bu,\bn}^{\bv}$ is tilted.

Recall that $\omega_{\bk,\ell} = -\gamma_\ell Q_\ell(r_{\bk}) - \overline{\gamma}_\ell \overline{Q_\ell(1)}$. Taking the complex conjugate of the first moment formula, we have $\Q[\overline{\gamma}_\ell \mid \filt_\ell] = -\sqrt{\frac{2}{\beta}} \frac{\sigma Q_\ell(r_{\bu})}{\ell+1} + O(\ell^{-\alpha})$. Thus
\begin{align*}
\Q[\omega_{\bk,\ell} \mid \filt_\ell] &= -Q_\ell(r_{\bk}) \Q[\gamma_\ell \mid \filt_\ell] - \overline{Q_\ell(1)} \Q[\overline{\gamma}_\ell \mid \filt_\ell] \\
&= \sqrt{\frac{2}{\beta}} \frac{1}{\ell+1} \bigl( \overline{\sigma} Q_\ell(r_{\bk}) \overline{Q_\ell(r_{\bu})} + \sigma \overline{Q_\ell(1)} Q_\ell(r_{\bu}) \bigr) + \mathcal{O}(\ell^{-\alpha})\\
&= \sqrt{\frac{2}{\beta}} \frac{1}{\ell+1} \bigl( \overline{\sigma} (\Delta_{\bk,\ell} + 1)\overline{(\Delta_{\bu,\ell} + 1)} + \sigma (\Delta_{\bu,\ell} + 1) \bigr) + \mathcal{O}(\ell^{-\alpha}).
\end{align*}
Taking the real part, $\Re(\Q[\omega_{\bk,\ell} \mid \filt_\ell]) = \sqrt{\frac{2}{\beta}} \frac{1}{\ell+1} \Re\bigl( \overline{\sigma} (\Delta_{\bk,\ell} + 1)\overline{(\Delta_{\bu,\ell} + 1)} + \sigma (\Delta_{\bu,\ell} + 1) \bigr) + \mathcal{O}(\ell^{-\alpha})$.

Since $\omega_{\bk,\ell}$ is complex, we compute $\Q[|\omega_{\bk,\ell}|^2 \mid \filt_\ell]$. Expanding,
\[
|\omega_{\bk,\ell}|^2 = |\gamma_\ell|^2 |Q_\ell(r_{\bk})|^2 + \gamma_\ell^2 Q_\ell(r_{\bk}) Q_\ell(1) + \overline{\gamma_\ell^2 Q_\ell(r_{\bk}) Q_\ell(1)} + |\gamma_\ell|^2 |Q_\ell(1)|^2.
\]
Using $|Q_\ell(1)|^2 = 1$, $|Q_\ell(r_{\bk})|^2 = |\Delta_{\bk,\ell} + 1|^2$, and $Q_\ell(r_{\bk}) Q_\ell(1) = (\Delta_{\bk,\ell} + 1) Q_\ell(1)^2$,
\[
|\omega_{\bk,\ell}|^2 = |\gamma_\ell|^2 (|\Delta_{\bk,\ell} + 1|^2 + 1) + 2\Re(\gamma_\ell^2 (\Delta_{\bk,\ell} + 1) Q_\ell(1)^2).
\]
Since $\Q[\gamma_\ell^2 \mid \filt_\ell] = O(\ell^{-\alpha})$, we obtain
\[
\Q[|\omega_{\bk,\ell}|^2 \mid \filt_\ell] = (|\Delta_{\bk,\ell} + 1|^2 + 1) \cdot \frac{2}{\beta\ell} + O(\ell^{-\alpha}).
\]
Expanding $\omega_{\bk,\ell}^2 = \gamma_\ell^2 Q_\ell(r_{\bk})^2 + 2|\gamma_\ell|^2 (\Delta_{\bk,\ell} + 1) + \overline{\gamma_\ell^2} \overline{Q_\ell(1)}^2$ and using that $\Q[\gamma_\ell^2 \mid \filt_\ell] = O(\ell^{-\alpha})$,
\[
\Q[\omega_{\bk,\ell}^2 \mid \filt_\ell] = \frac{4(\Delta_{\bk,\ell} + 1)}{\beta \ell} + O(\ell^{-\alpha}).
\]
Combining the above,
\begin{align}
\Q\left[\log \left|\frac{P_{\bk,\ell+1}}{P_{\bk,\ell}}\right| \,\Big|\, \filt_\ell\right] &= \log r_{\bk} - \sqrt{\frac{2}{\beta}} \frac{\Re\bigl( \overline{\sigma} (\Delta_{\bk,\ell} + 1)\overline{(\Delta_{\bu,\ell} + 1)} + \sigma (\Delta_{\bu,\ell} + 1) \bigr)}{\ell} \notag\\
&\quad + \frac{2\Re(\Delta_{\bk,\ell} + 1)}{\beta\ell} - (2+\Re(\Delta_{\bk,\ell})) \cdot \frac{2}{\beta\ell} + O(\ell^{-\alpha}) \notag\\
&= \log r_{\bk} - \sqrt{\frac{2}{\beta}} \frac{\Re\bigl( \overline{\sigma} (\Delta_{\bk,\ell} + 1)\overline{(\Delta_{\bu,\ell} + 1)} + \sigma (\Delta_{\bu,\ell} + 1) \bigr)}{\ell} - \frac{2}{\beta\ell} + O(\ell^{-\alpha}). \label{eq:logPmean}
\end{align}

The dominant contribution to the fluctuations comes from $\omega_{\bk,\ell}$. We have
\[
\Q[|\omega_{\bk,\ell} - \Q[\omega_{\bk,\ell} \mid \filt_\ell]|^2 \mid \filt_\ell] = \Q[|\omega_{\bk,\ell}|^2 \mid \filt_\ell] - |\Q[\omega_{\bk,\ell} \mid \filt_\ell]|^2 = \frac{2(|\Delta_{\bk,\ell} + 1|^2 + 1)}{\beta\ell} + O(\ell^{-\alpha}),
\]
since $|\Q[\omega_{\bk,\ell} \mid \filt_\ell]|^2 = O(\ell^{-2})$. The contributions from $|\gamma_\ell|^2$ and higher-order terms are $O(\ell^{-2})$, so
\begin{equation}\label{eq:logPvar}
\Q\left[\left|\log \left|\frac{P_{\bk,\ell+1}}{P_{\bk,\ell}}\right| - \Q\left[\log \left|\frac{P_{\bk,\ell+1}}{P_{\bk,\ell}}\right| \,\Big|\, \filt_\ell\right]\right|^2 \,\Big|\, \filt_\ell\right] \leq \frac{2(|\Delta_{\bk,\ell} + 1|^2 + 1)}{\beta\ell} + O(\ell^{-\alpha}).
\end{equation}

\textbf{Martingale bounds:} Define the martingale  $M_{\bk,\ell} \coloneqq \log(|P_{\bk,\ell}|/|P_{\bk,2^{\epsilon\bv}}|)- A_{\bk,\ell}$ where
\[
A_{\bk,\ell} \coloneqq \sum_{j=2^{\epsilon\bv}}^{\ell-1} \Q\left[\log \left|\frac{P_{\bk,j+1}}{P_{\bk,j}}\right| \,\Big|\, \filt_j\right]
\]
is the predictable (drift) part.
Since on $\{\tau > \ell\}$ we have $|\Delta_{\bk,\ell}|, |\Delta_{\bu,\ell}| \leq \epsilon'$, the quantity
\[-\Re\bigl( \overline{\sigma} (\Delta_{\bk,\ell} + 1)\overline{(\Delta_{\bu,\ell} + 1)} + \sigma (\Delta_{\bu,\ell} + 1) \bigr)\leq -\Re(\sigma+\overline{\sigma})+O(\epsilon')\leq
O(\epsilon'),\] and we have from
\eqref{eq:logPmean}, by selecting $\epsilon'$ small enough (depending only on $\beta$) that for $\tau \geq \ell\geq \ell'$,
\[
A_{\bk,\ell} -A_{\bk,\ell'}< (\ell-\ell')2^{-\bk}+\sum_{j=\ell'}^{\ell-1} -\frac{1}{\beta j} \leq -\frac{1}{\beta} \log(\ell / \ell')+(\ell-\ell')2^{-\bk}
\leq -\frac{1}{\beta} \log(\ell /\ell')+1.
\]

The stopped martingale $M^{\tau}_{\bk,\ell} \coloneqq M_{\bk,\ell \wedge \tau}$ therefore satisfies the hypotheses of Lemma \ref{lem:appendixmartingale}.  Hence we have that for any $\eta> 0$ and any integer $\bm \geq 0$, there is a $\delta>0$ so that for all $y \geq 0$,
\begin{equation}\label{eq:Mtailbound}
\Q\left(
\widehat{\mathscr{G}}_{\bn}
\cap \Big\{
\exists 2^{\epsilon \bv}\leq j <\ell\leq L \cdot 2^{-\bm} :  \exp(M^{\tau}_{\bk,\ell} - M^{\tau}_{\bk,j}) \geq \left(\frac{\ell}{j}\right)^{1/(2\beta)} \left(\frac{L \cdot 2^{-\bm}}{\ell}\right)^{\eta} 2^{y}
\Big\}
\middle\vert
\filt_{2^{\bn}}
\right) \leq 2e^{-\delta y}.
\end{equation}

We now argue that the stopping time $\tau$ is unlikely to occur before $L$. The argument applies to both $\Delta_{\bu,\ell}$ and $\Delta_{\bv,\ell}$ with the same bounds, so we write it for a general $\bk \in \{\bu, \bv\}$, and we apply \eqref{eq:Mtailbound} twice, once for each $\bk$ as well. Recall from~\eqref{eq:Deltasolution} that for any $\ell \leq \tau$,
\[
\Delta_{\bk,\ell} = \frac{\Delta_{\bk,2^{\epsilon \bv}} P_{\bk,\ell}}{P_{\bk,2^{\epsilon\bv}}} + \sum_{j=2^{\epsilon \bv}}^{\ell-1} \frac{(r_{\bk}-1) P_{\bk,\ell}}{P_{\bk,j+1}} + O\Bigl(\sum_{j=2^{\epsilon \bv}}^{\ell-1} \frac{\gamma_j^3 P_{\bk,\ell}}{P_{\bk,j+1}}\Bigr).
\]
On the event $\{L \cdot 2^{-\bm-1} < \tau \leq L \cdot 2^{-\bm}\}$, we have $\max(|\Delta_{\bu,\tau}|, |\Delta_{\bv,\tau}|) \geq \epsilon'$ by definition of $\tau$. We apply~\eqref{eq:Mtailbound} with $y$ to be determined below. On the complement of the exceptional event in~\eqref{eq:Mtailbound}, combining with the drift bound $A_{\bk,\ell} - A_{\bk,j} \leq -\frac{1}{\beta}\log(\ell/j)+1$ from above, we obtain for $2^{\epsilon\bv} \leq j< \tau$,
\begin{equation}\label{eq:Pratio}
\frac{|P_{\bk,\tau}|}{|P_{\bk,j}|} \ll \left(\frac{j}{\tau}\right)^{1/(2\beta)} \left(\frac{L \cdot 2^{-\bm}}{\tau}\right)^{\eta} 2^{y} \ll \left(\frac{j}{\tau}\right)^{1/(2\beta)} 2^{y}.
\end{equation}
We now bound each term in~\eqref{eq:Deltasolution}.

For the initial term, we recall \eqref{eq:InitialDeltaBound}, which gives $|\Delta_{\bk,2^{\epsilon\bv}}| \ll 2^{-\bv/4}.$  By~\eqref{eq:Pratio} with $j = 2^{\epsilon\bv}$ and $\ell = \tau$, we have $|P_{\bk,\tau}/P_{\bk,2^{\epsilon\bv}}| \ll 2^{y}$. Thus the initial term satisfies
\begin{equation}\label{eq:initialterm}
\left|\frac{\Delta_{\bk,2^{\epsilon\bv}} P_{\bk,\tau}}{P_{\bk,2^{\epsilon\bv}}}\right| \ll 2^{-\bv/4+y}. 
\end{equation}

The middle sum satisfies
\begin{equation}\label{eq:middleterm}
\left| \sum_{j=2^{\epsilon \bv}}^{\tau-1} \frac{(r_{\bk}-1) P_{\bk,\tau}}{P_{\bk,j+1}} \right| \ll 2^{-\bk + y} \sum_{j=2^{\epsilon \bv}}^{\tau-1} \left(\frac{j+1}{\tau}\right)^{1/(2\beta)} \ll 2^{-\bk +y} \tau \ll L 2^{-\bk -\bm + y} \leq 2^{-\bk + \bv - \bm - cx + y}.
\end{equation}

For the error term, on $\mathscr{G}_{\bn}$ we have $|\gamma_j|^3 \ll (\log j)^{3/2} / j^{3/2}$. By~\eqref{eq:Pratio},
\begin{equation}\label{eq:errorterm}
\sum_{j=2^{\epsilon \bv}}^{\tau-1} \frac{|\gamma_j|^3 |P_{\bk,\tau}|}{|P_{\bk,j+1}|} \ll 2^{y} \sum_{j=2^{\epsilon \bv}}^{\tau-1} \frac{(\log j)^{3/2}}{j^{3/2}} \left(\frac{j+1}{\tau}\right)^{1/(2\beta)}
\ll 2^{-(\epsilon/4)\bv +y}.
\end{equation}

We now choose $y = \nu(cx + \bm)$ and $\nu < \min(\epsilon/4, 1)$. With this choice, and using that $\bm + cx \leq \bv$ and $\bk\geq\bv$, the bounds~\eqref{eq:initialterm} and~\eqref{eq:errorterm} are $o(1)$ as $\bv \to \infty$, and~\eqref{eq:middleterm} is at most $2^{-(1-\nu)cx}$, which is small since $x$ is large.  Hence all three are smaller than $\epsilon'/3$.

Combining these bounds, on $\mathscr{G}_{\bn} \cap \{L \cdot 2^{-\bm-1} < \tau \leq L \cdot 2^{-\bm}\}$ and outside the exceptional event in~\eqref{eq:Mtailbound}, we obtain $|\Delta_{\bk,\tau}| < \epsilon'$ for both $\bk = \bu$ and $\bk = \bv$, but $\max(|\Delta_{\bu,\tau}|, |\Delta_{\bv,\tau}|) > \epsilon'$ by definition of $\tau$, a contradiction. With $y$ as chosen, using \eqref{eq:Mtailbound},
\[
\Q\bigl(\{L \cdot 2^{-\bm-1} < \tau \leq L \cdot 2^{-\bm}\} \,\big|\, \filt_{2^{\bn}}\bigr) \leq 2e^{-\delta \nu (cx + \bm)/2}.
\]
Summing over $\bm \geq 0$,
\begin{equation}\label{eq:tauL}
\Q\bigl(\{\tau \leq L\} \,\big|\, \filt_{2^{\bn}}\bigr) \leq \sum_{\bm=0}^{\infty} 2e^{-\delta \nu (cx + \bm)/2} = e^{-\delta \nu cx/2} \sum_{\bm=0}^{\infty} 2e^{-\delta \nu \bm/2} \ll e^{-\delta' cx}
\end{equation}
for some $\delta' > 0$. On the event $\{\tau > L\}$, the stopped and unstopped processes coincide, and using \eqref{eq:Mtailbound} with $y=h$ shows
\[
  \Q_{\bu,\bn}^{\bv}\left[\widehat{\mathscr{G}}_\bn 
  \bigcap 
  \Bigg\{
  \exists~(j,\ell) : 2^{\epsilon \bv}\leq j <\ell\leq L, \frac{|P_{\bk,\ell}|}{|P_{\bk,j}|}\geq \left(\frac{j}{\ell}\right)^{1/2\beta} \left(\frac{L}{j}\right)^\eta 2^{{h}} \Bigg\}
  \bigcap \{ \tau > L\}
  ~\middle|~ \filt_{2^{\bn}}\right] \ll e^{-\delta {h}}.
\]
Combining this with \eqref{eq:tauL}, this shows Claim \ref{claim:Pmoment}.  
\end{proof}

We now bound $S'$ from~\eqref{eq:Sprimedef}. Substituting~\eqref{eq:Deltasolution} (with $\bk = \bv$) we obtain
\[
S' \leq \sum_{\ell = 2^{\epsilon\bv}}^{L} \frac{1}{\ell+1} \left| P_{\bv,\ell+1} \left( \frac{\Delta_{\bv,2^{\epsilon\bv}}}{P_{\bv,2^{\epsilon\bv}}} + \sum_{j=2^{\epsilon\bv}}^{\ell} \frac{r_{\bv}-1}{P_{\bv,j+1}} + \mathcal{O}\Bigl(\sum_{j=2^{\epsilon\bv}}^{\ell} \frac{\gamma_j^3}{P_{\bv,j+1}}\Bigr) \right) \right|.
\]
By the triangle inequality, $S' \ll (S_1 + S_2 + S_3)$ where
\[
S_1 \coloneqq \sum_{\ell = 2^{\epsilon\bv}}^{L} \frac{|P_{\bv,\ell+1}|}{\ell+1} \cdot \frac{|\Delta_{\bv,2^{\epsilon\bv}}|}{|P_{\bv,2^{\epsilon\bv}}|}  , \quad
S_2 \coloneqq \sum_{\ell = 2^{\epsilon\bv}}^{L} \frac{|P_{\bv,\ell+1}|}{\ell+1} \left| \sum_{j=2^{\epsilon\bv}}^{\ell} \frac{r_{\bv}-1}{P_{\bv,j+1}} \right|, \quad
S_3 \coloneqq \sum_{\ell = 2^{\epsilon\bv}}^{L} \frac{|P_{\bv,\ell+1}|}{\ell+1} \left| \sum_{j=2^{\epsilon\bv}}^{\ell} \frac{\gamma_j^3}{P_{\bv,j+1}} \right|.
\]
By Claim~\ref{claim:Pmoment} with $h=cx/2$ and $\eta > 0$ to be determined, on $\widehat{\mathscr{G}}_{\bn}$ and outside an exceptional event of probability at most $2e^{-c\delta x/2}$ (where we recall $\delta > 0$ depends on $\eta$), we have $|P_{\bv,\ell+1}|/|P_{\bv,j}| \ll (j/(\ell+1))^{1/(2\beta)} (L/j)^{\eta} 2^{cx/2}$ for all $2^{\epsilon\bv} \leq j < \ell \leq L$.

For $S_1$, we have using the bound \eqref{eq:InitialDeltaBound} on $\Delta_{\bv,2^{\epsilon\bv}}$ which holds with probability $1-e^{-\delta \bv}$,
\[
S_1 \leq 2|\Delta_{\bv,2^{\epsilon\bv}}| \sum_{\ell = 2^{\epsilon\bv}}^{L} \frac{1}{\ell+1} \cdot \frac{|P_{\bv,\ell+1}|}{|P_{\bv,2^{\epsilon\bv}}|}  \ll 2^{-\bv/4} \sum_{\ell = 2^{\epsilon\bv}}^{L} \frac{1}{\ell+1} \cdot \frac{|P_{\bv,\ell+1}|}{|P_{\bv,2^{\epsilon\bv}}|}.
\]
We have $|P_{\bv,\ell+1}|/|P_{\bv,2^{\epsilon\bv}}| \ll (2^{\epsilon\bv}/(\ell+1))^{1/(2\beta)} (L/2^{\epsilon\bv})^{\eta} 2^{cx/2}$, so
\[
S_1 \ll 2^{-\bv/4+cx/2} L^{\eta} \sum_{\ell = 2^{\epsilon\bv}}^{L} \frac{1}{\ell+1} \left(\frac{2^{\epsilon\bv}}{\ell+1}\right)^{1/(2\beta)} (2^{\epsilon\bv})^{-\eta}  \ll 2^{-\bv/4+cx/2} L^{\eta} (2^{\epsilon\bv})^{-\eta}.
\]
This is $\ll 2^{cx/2 + \eta(\bv - cx - \epsilon\bv) - \bv/4} \ll 2^{cx/2 - \delta_1 \bv}$ for some $\delta_1 > 0$, provided $\epsilon$ and $\eta$ are chosen small enough; this suffices as we may assume $cx/2 < \delta_1 \bv/2$ by choosing $c_1$ small.

For $S_2$, since $|r_{\bv} - 1| \ll 2^{-\bv}$, we have
\[
S_2 \ll 2^{-\bv} \sum_{\ell = 2^{\epsilon\bv}}^{L} \frac{1}{\ell+1} \sum_{j=2^{\epsilon\bv}}^{\ell} \frac{|P_{\bv,\ell+1}|}{|P_{\bv,j+1}|}.
\]
Using $|P_{\bv,\ell+1}|/|P_{\bv,j+1}| \ll (j/(\ell+1))^{1/(2\beta)} (L/j)^{\eta} 2^{cx/2}$, the inner sum is bounded by $2^{cx/2} L^{\eta} \sum_{j=2^{\epsilon\bv}}^{\ell} (j/(\ell+1))^{1/(2\beta)} j^{-\eta} \ll 2^{cx/2} L^{\eta} \ell^{1-\eta}$. Hence
\[
S_2 \ll 2^{-\bv +cx/2} L^{\eta} \sum_{\ell = 2^{\epsilon\bv}}^{L} \frac{\ell^{1-\eta}}{\ell+1} \ll 2^{-\bv + cx/2} L^{\eta} L^{1-\eta} = 2^{- cx/2}.
\]

For $S_3$, on the event $\mathscr{G}_{\bn}$ we have $|\gamma_j|^3 \ll (\log j)^{3/2} / j^{3/2}$. Interchanging the order of summation,
\[
S_3 \ll \sum_{j=2^{\epsilon\bv}}^{L} \frac{(\log j)^{3/2}}{j^{3/2}} \sum_{\ell=j}^{L} \frac{1}{\ell+1} \cdot \frac{|P_{\bv,\ell+1}|}{|P_{\bv,j+1}|}.
\]
Using $|P_{\bv,\ell+1}|/|P_{\bv,j+1}| \ll (j/(\ell+1))^{1/(2\beta)} (L/j)^{\eta} 2^{cx/2}$, the inner sum is bounded by
\[
2^{cx/2} (L/j)^{\eta} \sum_{\ell=j}^{L} \frac{1}{\ell+1} \left(\frac{j}{\ell+1}\right)^{1/(2\beta)} = 2^{cx/2} (L/j)^{\eta} j^{1/(2\beta)} \sum_{\ell=j}^{L} (\ell+1)^{-1-1/(2\beta)} \ll 2^{cx/2} (L/j)^{\eta}.
\]
Summing over $j \geq 2^{\epsilon\bv}$ with the factor $(\log j)^{3/2}/j^{3/2}$ gives $S_3 \ll 2^{cx/2} L^{\eta} \cdot (\epsilon\bv)^{3/2} 2^{-\epsilon\bv(1/2 + \eta)} \ll 2^{cx/2 + \eta(\bv - cx) - \epsilon \bv/2}$.  Choosing $\eta < \epsilon/2$ and choosing $c_1$ small gives $S_3 \ll 2^{-\delta_3 \bv}$ for some $\delta_3 >0$.

Combining the bounds for $S_1$, $S_2$, and $S_3$, on the event from Claim~\ref{claim:Pmoment} we have $S' \ll 2^{-c x/2}$ provided $c_1$ is sufficiently small. Recall the exceptional event has probability at most $2e^{-c \delta x/2}$. 
\end{proof}
We finally compare the two $\zeta$ processes evaluated at $r_\bl$ and at $1$.
\begin{corollary}\label{cor:radialincrement}
For any $\bk\geq \bl\geq \mathfrak{c}\bn$, $\bl\geq \bu\geq \bn$,
%
\[
\mathbb{Q}_{\bk,\bn}^{\infty}
\left(
  \widehat{\mathscr{G}}_{\bn}
  \cap
  \{
  |
  \zeta_{2^{\bu}}(\bn;r_{\bl},\theta)
  -
  \zeta_{2^{\bu}}(\bn;1,\theta)
  |
  > x 
  \}
\right)
\leq 
2e^{-c x}.
\]
\end{corollary}
\begin{proof}
Under the measure $\mathbb{Q}_{\bk,\bn}^{\infty}$ conditioned on $\filt_{2^\bn}$, we decompose the difference $\zeta_{2^{\bu}}(\bn;r_{\bl},\theta) - \zeta_{2^{\bu}}(\bn;1,\theta)$ into a drift part and a martingale part. By \eqref{eq:Z} and the Girsanov theorem, the drift contribution is bounded by
\[
\sum_{\ell=2^{\bn}}^{2^{\bu}-1} \frac{|Q_\ell(r_{\bl},\theta) - Q_\ell(1,\theta)|}{\ell+1} = \sum_{\ell=2^{\bn}}^{2^{\bu}-1} \frac{|\Delta_{\bl,\ell}|}{\ell+1},
\]
while the martingale part has quadratic variation
\[
\langle M \rangle = \sum_{\ell=2^{\bn}}^{2^{\bu}-1} \frac{|Q_\ell(r_{\bl},\theta) - Q_\ell(1,\theta)|^2}{\ell+1} = \sum_{\ell=2^{\bn}}^{2^{\bu}-1} \frac{|\Delta_{\bl,\ell}|^2}{\ell+1}.
\]
By Proposition~\ref{prop:finebracket}, and since $|\Delta_{\bl,\ell}| \leq 2$ the quadratic variation satisfies $\mathbb{Q}_{\bk,\bn}^{\infty}(\widehat{\mathscr{G}}_{\bn} \cap \{\langle M \rangle > x\}) \ll 2^{-\delta x}$. Proposition~\ref{prop:finebracket} implies that the drift also has exponential tails. The corollary then follows from standard Gaussian tail bounds for the martingale part.
\end{proof}

\section{The limiting clock deviation}
\label{sec:identification-constant}

Fix the radius \(r_\bk\).  Since \(\lg_{r_\bk}(2^\bt)=\bk\log2\)
for \(\bt\geq\bk\), the typical Brownian time
\(\lg_{r_\bk}(2^\bt)-\lg_{r_\bk}(2^\bn)\) stops growing once \(\bt\)
reaches the scale \(\bk\).  Brownian time itself keeps growing, so the
clock deviation \(S_{\bk,\bn}^{\bt}\) keeps changing after scale
\(\bk\).  Its limit \(S_{\bk,\bn}^{\infty}\) as \(\bt\to\infty\) is the
\textbf{limiting clock deviation}.  By Lemma~\ref{lem:QVsmallkold}, it is
finite and has exponential moments on the truncation event.

The weight \(e^{S^{\infty}}\) enters through the change of measure
\eqref{eq:com} with \(\bt=\infty\).  In
Section~\ref{sec:brownian-clock-approximate-martingale} we need its
conditional expectation given \(\filt_{2^\bt}\), for a scale \(\bt\)
slightly below \(\bk\), in the limit \(\bk\to\infty\).  This section
computes it.  Lemma~\ref{lem:step3-late-forgetting} shows that this
conditional expectation is asymptotically deterministic, because the
Markov chain of time dilations after scale \(\bt\) forgets its initial
condition.  Proposition~\ref{lem:averaged-FK-transfer} identifies the
resulting constant \(c_\bk\) by averaging over \(\theta\) against the
martingale \(\mathscr M\).

Choose
\begin{equation}\label{eq:step3-near-saturation-scales}
 1<\alpha'<2,\qquad
 \alpha>4\alpha',\qquad
 \bt_\bk=\left\lfloor\bk-(\log(\bk+2))^\alpha\right\rfloor,\quad
 \bt_\bk'=\left\lfloor\bk-(\log(\bk+2))^{\alpha'}\right\rfloor .
\end{equation}
Let 
\[ \mathcal Z_{\bk,\bt}^A :=\{z\in\mathbb D: (1-|z|^2)/(1-r_\bk^2)\leq C_A2^\bt\bt^A\}.\]
With this notation, the set \(\mathscr W_\bt^A=\{R_{2^\bt}(r_\bk,\theta)\leq C_A2^\bt\bt^A\}\) as in 
\eqref{eq:step2-radial-window} is such that 
 \(\mathscr W_\bt^A=\{Q_{2^\bt}\in\mathcal Z_{\bk,\bt}^A\}\).
\begin{lemma}[Mixing of the limiting clock deviation]
\label{lem:step3-late-forgetting}
Let \(A<\infty\) be fixed, let \(\bt=\bt_\bk\) and let
\(\mathbb Q_{\bk,\bt}^{\infty}\) be as in \eqref{eq:Qt}.
Then there are  deterministic $c,C$, \(c_\bk\), with
\begin{equation}\label{eq:step3-ck-bounds}
 0<c\leq c_\bk\leq C<\infty,
\end{equation}
such that, almost surely on \(\mathscr W_\bt^A\),
\begin{equation}\label{eq:step3-FK-forgetting}
 \left|
 \mathbb Q_{\bk,\bt}^{\infty}\left(
  \operatorname{Tru}_{2^\bt}e^{S_{\bk,\bt}^{\infty}}
\right)
 -c_\bk\right|=o_\bk(1),
\end{equation}
where the \(o_\bk(1)\) term is deterministic and depends only on \(A\) and \(\beta\).
\end{lemma}

\begin{proof}
Throughout the proof write \(\mathbb Q=\mathbb Q_{\bk,\bt}^{\infty}\) and 
write \(\mathbb E\) for expectation with respect to \(\mathbb Q\) (both are conditional on $\filt_{2^{\bt}}$).
Put \(Q_j=Q_j(r_\bk,\theta)\) and \(R_j=R_j(r_\bk,\theta)\), the latter
as in \eqref{eq-Rjnew}. 
By \eqref{eq:tilted-one-step-density} and \eqref{eq:Q}, under
\(\mathbb Q\) the conditional law of \((\gamma_j)_{j\geq2^\bt}\) given
\(\filt_{2^\bt}\) depends on \(\filt_{2^\bt}\) only through
\(Q_{2^\bt}\).  For \(z\in\mathbb D\), let
\(\mathbb Q_{\bk,\bt}^{\infty,z}\) denote this conditional law on the
event \(\{Q_{2^\bt}=z\}\).  Since the conditional expectation in
\eqref{eq:step3-FK-forgetting} equals
\(\mathbb Q_{\bk,\bt}^{\infty,z}(\operatorname{Tru}_{2^\bt}
e^{S_{\bk,\bt}^{\infty}})\) at \(z=Q_{2^\bt}\), it suffices to prove
\begin{equation}\label{eq:step3-FK-coupling-target}
 \sup_{z,z'\in\mathcal Z_{\bk,\bt}^A}
 \left|
 \mathbb Q_{\bk,\bt}^{\infty,z}\left(
  \operatorname{Tru}_{2^\bt}e^{S_{\bk,\bt}^{\infty}}\right)
 -\mathbb Q_{\bk,\bt}^{\infty,z'}\left(
  \operatorname{Tru}_{2^\bt}e^{S_{\bk,\bt}^{\infty}}\right)
 \right|=o_\bk(1).
\end{equation}

By \eqref{eq:tilted-one-step-density}, under \(\mathbb Q\) and
conditionally on \(\filt_j\), \(Z_j\) is a complex Gaussian with mean
\(-\overline{\sigma Q_j}/\sqrt{j+1}\), while \(\Gamma_j\) is
independent of \(Z_j\) and retains its law.  Hence, for
\(\phi\in\mathbb R\), replacing \(Q_j\) by \(e^{i\phi}Q_j\) replaces
the conditional law of \(\gamma_j\) by that of \(e^{-i\phi}\gamma_j\),
so the conditional law of \(\gamma_jQ_j\) given \(\filt_j\) depends on
\(Q_j\) only through \(|Q_j|\), that is, only through \(R_j\).  By
\eqref{eq:step2-restoration-affine}, the conditional law of
\(R_{j+1}\) given \(\filt_j\) therefore depends on the past only
through \(R_j\), and \((R_j)_{j\geq2^\bt}\) is a time-inhomogeneous
Markov chain under \(\mathbb Q\).

We establish three estimates. 

\smallskip
\noindent\emph{Step 1. Two-sided bounds on \(R_j\) for
\(2^{\bt}\leq j\leq2^{\bt_\bk'}\).}
Fix \(p_0\in(0,1)\), and let \(A_{\bk,j}\) be as in
\eqref{eq:step2-restoration-A}.  Since \(\bt=\bt_\bk\to\infty\), for all \(\bk\) large we have \(2^{\bt}\geq j_0\), with \(j_0\) as in Lemma~\ref{lem:step2-local-restoration}, and  \(2^{-\bk}\leq j^{-1}\) for \(j\leq2^{\bt_\bk'}\leq2^{\bk}\).  Hence
\eqref{eq:step2-negative-moment} with \(p=p_0\) gives a
constant \(C<\infty\) depending only on \(p_0,\beta\) such that, almost surely, for \(2^{\bt}\leq j\leq2^{\bt_\bk'}\),
\begin{equation}
 \mathbb E(A_{\bk,j}^{-p_0}\mid \filt_j)
 \leq1+C/j .
 \label{eq:step3-lower-restoration}
\end{equation}

Let \(N=2^\bt\) and \(J\geq2N\).  Keeping only the terms
\(J/2\leq i<J\) of the second sum in
\eqref{eq:step2-restoration-affine}, and applying Jensen's inequality to
the convex map \(x\mapsto x^{-p_0}\), gives
\begin{align}
 R_J^{-p_0}
 &\leq(J/2)^{-p_0-1}\sum_{J/2\leq i<J}
 \left(\prod_{s=i+1}^{J-1}A_{\bk,s}\right)^{-p_0},\notag\\
 \mathbb E(R_J^{-p_0}\mid\filt_N)&\leq CJ^{-p_0},
 \label{eq:step3-endpoint-negative-moment}
\end{align}
where the second line follows from the Markov property and \eqref{eq:step3-lower-restoration}.
Now let \(J=2^\bm\) with \(\bm\geq\bt+1\), and let \(\tau\) be the
first \(j\geq J\) with \(R_j<j(\bk+2)^{-C_0}\).  By \eqref{eq:step2-restoration-affine} and \eqref{eq:step3-lower-restoration}, the process
\begin{equation}
\label{eq-011026a}
 \frac{R_{j\wedge\tau}^{-p_0}}
 {\prod_{s=J}^{(j\wedge\tau)-1}(1+C/s)}
\end{equation}
is a nonnegative supermartingale.  Note that for $j\in [J,2J]$, the denominator in \eqref{eq-011026a}
is uniformly bounded above. Hence, Doob's maximal inequality followed by
\eqref{eq:step3-endpoint-negative-moment} yields
\begin{equation}\label{eq:step3-maximal-lower-exit}
 \mathbb Q\left(\min_{J\leq j\leq2J}\frac{R_j}{j}
 <(\bk+2)^{-C_0}\ \middle|\ \filt_N\right)
 \leq C(\bk+2)^{-p_0C_0}.
\end{equation}
For the upper bound on the maximum of $R_j/j$, let \(J=2^\bm\) with \(\bm\geq\bt\), and let
\(P_j=\prod_{s=J}^{j-1}(1+Cs^{-1-\delta})\leq C\).  By
\eqref{eq:step2-first-moment} and \eqref{eq:step2-restoration-affine},
as in \eqref{eq-11012026f}, the process \(P_j^{-1}R_j-j\) is a
supermartingale for \(j\geq J\).  On \(\mathscr W_\bt^A\), by
\eqref{eq:step2-first-moment-evolution},
\(\mathbb E(R_J\mid\filt_N)\leq CJ\bt^A\), so Doob's maximal inequality
gives
\begin{equation}\label{eq:step3-maximal-upper-exit}
 \mathbb Q\left(\max_{J\leq j\leq2J}\frac{R_j}{j}
 >(\bk+2)^{C_0}\ \middle|\ \filt_N\right)
 \leq C(\bk+2)^{A-C_0}.
\end{equation}
There are only \(O((\log\bk)^\alpha)\) blocks \([2^\bm,2^{\bm+1}]\)
with \(\bt\leq\bm<\bt_\bk'\).  Fix \(M>0\).  A union bound over these blocks, with
\(C_0\) sufficiently large depending on \(A\) and \(M\), shows that on
\(\mathscr W_\bt^A\), with $\mathbb Q$-conditional probability \(1-o(\bk^{-M})\),
\begin{equation}\label{eq:step3-two-sided-window}
 j(\bk+2)^{-C_0}\leq R_j\leq j(\bk+2)^{C_0}
 \quad\text{for }2^{\bt+1}\leq j\leq2^{\bt_\bk'},
 \qquad
 R_j\leq j(\bk+2)^{C_0}
 \quad\text{for }2^{\bt}\leq j<2^{\bt+1}.
\end{equation}
Note that on the event of \eqref{eq:step3-two-sided-window}, for every fixed $\widetilde M$,
\begin{equation}\label{eq:step3-late-small-parameter}
 \sup_{2^{\bt}\leq j\leq2^{\bt_\bk'}}2^{-\bk}R_j
 \leq2^{\bt_\bk'-\bk}(\bk+2)^{C_0}=o(\bk^{-\widetilde M}).
\end{equation}

By H\"older's inequality and Lemma~\ref{lem:QVsmallkold}, for any fixed
\(p>0\), the complement of \eqref{eq:step3-two-sided-window} also
satisfies, on \(\mathscr W_\bt^A\),
\[
 \mathbb Q\left(
 e^{pS_{\bk,\bt}^{\infty}}\one\{\operatorname{Tru}_{2^\bt}\}
 \one\{\text{\eqref{eq:step3-two-sided-window} fails}\}
\right)=o(\bk^{-M/2}).
\]

\smallskip
\noindent\emph{Step 2. Crossing of two independent copies of \(U_j=\log(R_j/j)\).}
Put \(U_j=\log(R_j/j)\) and write the Doob decomposition of its
increments as
\begin{equation}\label{eq:step3-log-transition}
 U_{j+1}-U_j=m_{\bk,j}(U_j)+\eta_{\bk,j+1},\qquad
 m_{\bk,j}(u)=\mathbb E(U_{j+1}-U_j\mid U_j=u),\quad
 \mathbb E(\eta_{\bk,j+1}\mid\filt_j)=0.
\end{equation}
Since \((R_j)\) is a Markov chain under \(\mathbb Q\) and \(|Q_j|^2=1-2^{-\bk}je^{U_j}\), the conditional law of \(U_{j+1}\) given \(\filt_j\) depends only on \(U_j\), so \(m_{\bk,j}\) is a function of \(u\) alone.
By \eqref{eq:step2-restoration-affine},
\begin{equation}\label{eq:step3-U-increment}
 U_{j+1}-U_j
 =\log A_{\bk,j}+\log\left(1+\frac{e^{-U_j}}{jA_{\bk,j}}\right)
  -\log(1+1/j).
\end{equation}
We need two routine consequences of this identity: the drift \(m_{\bk,j}\) is nearly decreasing, and the martingale increments have variance of exact order \(1/j\).  Precisely, for \(u\geq v\) both in the range permitted by
\eqref{eq:step3-two-sided-window} and \(2^{\bt+1}\leq j<2^{\bt_\bk'}\),
\begin{equation}\label{eq:step3-monotone-drift}
 m_{\bk,j}(u)-m_{\bk,j}(v)\leq\varepsilon_{\bk,j}
 :=C\left(j^{-2}(\bk+2)^{2C_0}
 +j^{-1}2^{\bt_\bk'-\bk}(\bk+2)^{C_0}
 +j^{-1-\delta}\right),
 \end{equation}
and, when the state at time \(j\) satisfies \eqref{eq:step3-two-sided-window},
\begin{equation}\label{eq:step3-full-increment-moments}
 \frac cj\leq
 \mathbb E(\eta_{\bk,j+1}^2\mid\filt_j)\leq\frac Cj,\qquad
 \mathbb E(|\eta_{\bk,j+1}|^{2+\rho}\mid\filt_j)
 \leq Cj^{-1-\rho/2}.
\end{equation}
Note that 
summing the three terms defining \(\varepsilon_{\bk,j}\) and using \eqref{eq:step3-near-saturation-scales}, one obtains that
\begin{equation}
\label{eq-071026a}
 \sum_{2^{\bt+1}\leq j<2^{\bt_\bk'}}\varepsilon_{\bk,j}=o(1).\end{equation}

The heuristic behind \eqref{eq:step3-monotone-drift} and \eqref{eq:step3-full-increment-moments} is that \(\log A_{\bk,j}\approx2\Re(\gamma_jQ_j)\) is, under $\mathbb Q$,  a centered Gaussian step of variance \(4|Q_j|^2/(\beta j)\), plus a drift of order \(1/j\) that depends on the state only through \(|Q_j|^2=1-2^{-\bk}je^{u}\), which is within \(2^{\bt_\bk'-\bk}(\bk+2)^{C_0}\) of one on the window, while the middle term of \eqref{eq:step3-U-increment} contributes the drift \(e^{-u}/j\), which is decreasing in \(u\).  To make this precise, one expands \(\log A_{\bk,j}\) to second order in \(\gamma_j\) on \(\{|\gamma_j|\leq1/2\}\), takes conditional expectations using \eqref{eq:Qgamma}, and controls \(\{|\gamma_j|>1/2\}\) by Cauchy--Schwarz, exactly as in the proof of Lemma~\ref{lem:step2-local-restoration}; we omit the details.

Let \(R^{(1)}\) and \(R^{(2)}\) be independent copies of the chain,
with laws \(\mathbb Q_{\bk,\bt}^{\infty,z}\) and
\(\mathbb Q_{\bk,\bt}^{\infty,z'}\), and put
\(U_j^{(i)}=\log(R_j^{(i)}/j)\),
\(D_j=U_j^{(1)}-U_j^{(2)}\), and
\(\Delta M_{j+1}:=\eta_{\bk,j+1}^{(1)}-\eta_{\bk,j+1}^{(2)}\), where
\(\eta^{(i)}\) is the martingale increment of \(U^{(i)}\) in
\eqref{eq:step3-log-transition}.  The two copies are driven by independent noise, so with
\(\filt_j\) now denoting the filtration generated by both copies,
\eqref{eq:step3-full-increment-moments} gives, whenever both states at time \(j\) satisfy \eqref{eq:step3-two-sided-window},
\begin{equation}\label{eq:step3-separation-moments}
 \frac cj\leq
 \mathbb E((\Delta M_{j+1})^2\mid\filt_j)\leq\frac Cj,\qquad
 \mathbb E(|\Delta M_{j+1}|^{2+\rho}\mid\filt_j)
 \leq Cj^{-1-\rho/2}.
\end{equation}
Choose
\begin{equation}\label{eq:step3-short-log-length}
 L_\bk=(\log(\bk+2))^{2+\delta_0},\qquad
 0<\delta_0<\alpha-2.
\end{equation}

\smallskip
\noindent\emph{Crossing within one block.}
Fix \(\bs\) with \(\bt+1\leq \bs\) and \(\bs+L_\bk\leq\bt_\bk'\), and put
\(j_1=2^\bs\) and \(j_2=2^{\bs+L_\bk}\).  Let \(\mathscr Y_\bs\) be the event
that both copies satisfy \eqref{eq:step3-two-sided-window} for all
\(j_1\leq j\leq j_2\), and let \(\tau\) be the first \(j\geq j_1\) with
\(D_j\leq0\).  By symmetry we may assume \(D_{j_1}\geq0\), and then
\(\tau=j_1\) if \(D_{j_1}=0\).  We claim that
\begin{equation}\label{eq:step3-brownian-crossing}
 \mathbb Q\left(\tau>j_2,\ \mathscr Y_\bs\mathrel{\big|}\filt_{j_1}\right)=o(1),
\end{equation}
uniformly in \(\bs\) and in the states of the two copies at time \(j_1\).
Let \(\sigma_W\) be the first \(j\geq j_1\) at which one of the copies
violates \eqref{eq:step3-two-sided-window}, and put
\[
 M_j=\sum_{j_1\leq i<j}\Delta M_{i+1}\one\{\sigma_W>i\},
 \qquad j_1\leq j\leq j_2,
\]
which is a martingale under \(\mathbb Q(\cdot\mid\filt_{j_1})\) with
\(M_{j_1}=0\), since \(\{\sigma_W>i\}\in\filt_i\).  Every nonzero
increment of \(M\) is taken from states satisfying
\eqref{eq:step3-two-sided-window}, so \eqref{eq:step3-separation-moments}
applies to each of them.  Hence
\[
 \kappa:=\sum_{j_1\leq i<j_2}\mathbb E\bigl(|M_{i+1}-M_i|^{2+\rho}\mathrel{\big|}\filt_{j_1}\bigr)
 \leq C\sum_{i\geq2^\bs}i^{-1-\rho/2}\leq C2^{-\bs\rho/2},
\]
and on \(\mathscr Y_\bs\), where \(\sigma_W>j_2\),
\[
 \langle M\rangle_{j_2}
 \geq c\sum_{j_1\leq i<j_2}\frac1i
 \geq cL_\bk
\]
after decreasing \(c\).  On \(\{\tau>j\}\) we have
\(U^{(1)}_i>U^{(2)}_i\) for all \(j_1\leq i<j\), so summing
\eqref{eq:step3-log-transition} for both copies and applying
\eqref{eq:step3-monotone-drift} and \eqref{eq-071026a} gives, on \(\{\tau>j\}\cap\mathscr Y_\bs\),
\begin{equation}\label{eq:step3-D-below-M}
 D_j=D_{j_1}+M_j+\sum_{j_1\leq i<j}\bigl(m_{\bk,i}(U^{(1)}_i)-m_{\bk,i}(U^{(2)}_i)\bigr)
 \leq D_{j_1}+M_j+1
\end{equation}
for all \(\bk\) large.  Hence on \(\{\tau>j_2\}\cap\mathscr Y_\bs\) the martingale \(M\) stays above \(-(2C_0\log(\bk+2)+1)\) throughout \([j_1,j_2]\) while its bracket reaches \(cL_\bk\gg(\log(\bk+2))^2\), and optional stopping, in the form of Lemma~\ref{lem:appendix-crossing} with \(x=2C_0\log(\bk+2)+1\) and \(V=cL_\bk\),
gives \eqref{eq:step3-brownian-crossing}.

\smallskip
\noindent\emph{The overshoot at the crossing.}
We next show that at the crossing time the two copies are within
\(O(\tau^{-1/2})\) of each other:  for every fixed \(K\geq6\),
\begin{equation}\label{eq:step3-crossing-overshoot}
 \mathbb Q\left(\tau\leq j_2,\ |D_\tau|>K\tau^{-1/2},\ \mathscr Y_\bs\right)
 \leq e^{-\beta(K-6)^2/16}+o(1),
\end{equation}
uniformly in \(\bs\) and in the states at time \(j_1\).  This rests on the
Gaussian structure of the one-step increment.  Let
\(\mathcal X_j=\operatorname{Tru}^{(1)}_{j,j+1}\cap\operatorname{Tru}^{(2)}_{j,j+1}\)
be the event that both copies satisfy the one-step truncation
\eqref{eq:goodtruncated} at step \(j\), and let \(\mathcal Y_j\) be the
\(\filt_j\)-measurable event that both states at time \(j\) satisfy
\eqref{eq:step3-two-sided-window}.  On \(\mathcal X_j\), by \eqref{eq-250826c}
with \(\epsilon=1/6\) and \eqref{eq:correct-beta-expansion}, we have
\(|\gamma^{(i)}_j|\leq j^{-1/3}\) and
\(\gamma^{(i)}_j=Z^{(i)}_j/\beta_j+O(j^{-3/4})\) for \(i=1,2\).
Inserting this in the second-order expansion of \(\log A_{\bk,j}\) in \(\gamma_j\) from the proof of Lemma~\ref{lem:step2-local-restoration},
and using \(A_{\bk,j}\geq1/2\) together with
\(0\leq\log(A_{\bk,j}+e^{-U_j}/j)-\log A_{\bk,j}\leq2e^{-U_j}/j\) and
\(e^{-U_j}\leq(\bk+2)^{C_0}\) on \(\mathcal Y_j\), we obtain on \(\mathcal X_j\cap \mathcal Y_j\)
\begin{equation}\label{eq:step3-gaussian-increment}
 D_{j+1}-D_j=G_{j+1}+\xi_{j+1},
 \qquad
 G_{j+1}:=\frac{2}{\beta_j}\Bigl(\Re\bigl(Z^{(1)}_jQ^{(1)}_j\bigr)-\Re\bigl(Z^{(2)}_jQ^{(2)}_j\bigr)\Bigr),
 \qquad
 |\xi_{j+1}|\leq j^{-1/2},
\end{equation}
for all \(\bk\) large, because \(j^{-3/4}\), \(j^{-2/3}\),
\((\bk+2)^{C_0}j^{-1}\) and \(2^{-\bk}\) are all \(o(j^{-1/2})\)
uniformly in \(2^{\bt}\leq j\leq2^{\bt_\bk'}\).  Conditionally on
\(\filt_j\), by Lemma~\ref{lem:1st2ndQmoment} and the independence of
the two copies, \(G_{j+1}\) is Gaussian with mean \(m_j\) and variance
\(\sigma_j^2\) satisfying
\[
 |m_j|\leq\frac{4}{\beta_j\sqrt{j+1}}\leq j^{-1/2},
 \qquad
 \sigma_j^2=\frac{2}{\beta_j^2}\bigl(|Q^{(1)}_j|^2+|Q^{(2)}_j|^2\bigr)\leq\frac{8}{\beta(j+1)}.
\]

Fix \(j_1\leq j<j_2\) and work on \(\{\tau>j\}\cap \mathcal Y_j\), on which
\(D_j>0\).  By \eqref{eq:step3-gaussian-increment}, on \(\mathcal X_j\) a crossing at step \(j+1\) is guaranteed by \(G_{j+1}\leq-D_j-j^{-1/2}\), while an overshoot of more than \(K(j+1)^{-1/2}\) requires \(G_{j+1}<-D_j-K(j+1)^{-1/2}+j^{-1/2}\).  These are lower-tail events for the Gaussian \(G_{j+1}\) at two levels below its mean whose separation, after accounting for \(|m_j|\leq j^{-1/2}\), is at least \((K-6)(j+1)^{-1/2}\geq(K-6)\sqrt{\beta/8}\,\sigma_j\).  Hence, by a Gaussian computation, namely the tail bound \(\Phi(-t-b)\leq e^{-b^2/2}\Phi(-t)\) for \(t,b\geq0\), we have on \(\{\tau>j\}\cap \mathcal Y_j\)
\begin{equation}\label{eq:step3-tail-ratio}
 \mathbb Q\bigl(\{D_{j+1}<-K(j+1)^{-1/2}\}\cap \mathcal X_j\mid\filt_j\bigr)
 \leq e^{-\beta(K-6)^2/16}\Bigl(\mathbb Q(\tau=j+1\mid\filt_j)+\mathbb Q({\mathcal X_j}^{\complement}\mid\filt_j)\Bigr).
\end{equation}
The event in \eqref{eq:step3-crossing-overshoot} is the union over
\(j_1\leq j<j_2\) of \(\{\tau>j\}\cap\{D_{j+1}<-K(j+1)^{-1/2}\}\cap\mathscr Y_\bs\),
and \(\mathscr Y_\bs\subseteq \mathcal Y_j\).  Splitting each of these events
according to \(\mathcal X_j\), multiplying \eqref{eq:step3-tail-ratio} by
\(\one\{\tau>j\}\one\{\mathcal Y_j\}\), taking expectations and summing over
\(j\) gives
\[
 \mathbb Q\left(\tau\leq j_2,\ |D_\tau|>K\tau^{-1/2},\ \mathscr Y_\bs\right)
 \leq e^{-\beta(K-6)^2/16}\,\mathbb Q(\tau\leq j_2)
 +2\sum_{j\geq2^\bs}\mathbb Q({\mathcal X_j}^{\complement}).
\]
As shown in the proof of Lemma~\ref{lem:1st2ndQmoment},
\(\mathbb Q(\operatorname{Tru}_{j,j+1}^{\complement}\mid\filt_j)\leq Cj^{-c}\)
with \(c>1\), so the last sum is \(O(2^{-\bs(c-1)})=o(1)\), which proves
\eqref{eq:step3-crossing-overshoot}.

\smallskip
\noindent\emph{Step 3. Coupling after a crossing.}
Let \(P_{\bk,j}(u,dy)\) denote the conditional law of \(U_{j+1}\) given \(\filt_j\) on \(\{U_j=u\}\).  We claim that for each fixed \(K\) there is \(\varepsilon_K>0\) such that, for all \(\bk\) large, all \(2^{\bt+1}\leq j\leq2^{\bt_\bk'}\) and all \(u,v\) in the range permitted by \eqref{eq:step3-two-sided-window},
\begin{equation}\label{eq:step3-local-minorization}
 |u-v|\leq Kj^{-1/2}
 \quad\Longrightarrow\quad
 \int P_{\bk,j}(u,dy)\wedge P_{\bk,j}(v,dy)
 \geq\varepsilon_K .
\end{equation}
This is a density computation, which we outline.  Since the conditional law of \(\gamma_jQ_j\) depends on \(Q_j\) only through \(|Q_j|\), we may replace \(Z_j\) by \(\widetilde Z_j=e^{i\arg Q_j}Z_j=x+iy\) and \(Q_j\) by \(|Q_j|\); freezing \(y\) and \(\Gamma_j\), which have the same law under both kernels and are coupled identically, leaves \(U_{j+1}=f_u(x)\) for an explicit map \(f_u\) given by \eqref{eq:step2-restoration-affine}, with \(x\) a Gaussian of variance \(1/2\) and mean \(O(j^{-1/2})\) by \eqref{eq:tilted-one-step-density}.  On the event \(\{|x|,|y|\leq K',\,|\Gamma_j-\beta_j^2|\leq K'\sqrt j\}\), which has conditional probability at least \(1/2\) for \(K'\) large, direct differentiation gives \(f_u'(x)=\frac{2|Q_j|}{\beta_j}(1+o(1))\) and \(\partial_uf_u(x)=1+o(1)\), so for \(|u-v|\leq Kj^{-1/2}\) the map \(f_v^{-1}\circ f_u\) is a shift by \(O_K(1)\) with derivative \(1+o(1)\).  The two kernels are therefore the pushforwards under \(f_u\) of two laws on \(\mathbb R\) whose densities on \([-K',K']\) are comparable to Gaussians of variance \(1/2\) shifted by \(O_K(1)\), and these overlap by a positive amount depending only on \(K\), which proves \eqref{eq:step3-local-minorization}.

\smallskip
\noindent\emph{The coupling.}
We now define the coupling \(\widehat{\mathbb Q}^{z,z'}\) of \(\mathbb Q_{\bk,\bt}^{\infty,z}\) and \(\mathbb Q_{\bk,\bt}^{\infty,z'}\) used in the rest of the proof.  Choose \(K\geq6\) with \(e^{-\beta(K-6)^2/16}\leq1/8\), and consider the consecutive blocks \([2^{\bs_i},2^{\bs_i+L_\bk}]\) with \(\bs_0=\bt+1\) and \(\bs_{i+1}=\bs_i+L_\bk\).  Run the two copies with independent noise, except at the first crossing time \(\tau\) of each block, where \((U^{(1)}_{\tau+1},U^{(2)}_{\tau+1})\) is drawn from a maximal coupling of \(P_{\bk,\tau}(U^{(1)}_\tau,\cdot)\) and \(P_{\bk,\tau}(U^{(2)}_\tau,\cdot)\); after a failed attempt the noise is independent again, so that the hypotheses of Step~2 hold at the start of the next block.  
From the first time \(j\) with \(R^{(1)}_j=R^{(2)}_j\) onward, the two copies have \(Q^{(2)}_j=e^{i\phi}Q^{(1)}_j\), and by the rotation covariance of \eqref{eq:Q} and of the tilted law \eqref{eq:tilted-one-step-density} we may couple the remaining noise of the second copy to that of the first by the rotation \(e^{-i\phi}\), so that \(|Q^{(2)}_i|=|Q^{(1)}_i|\) for all \(i\geq j\).

Within one block starting at \(2^\bs\), by \eqref{eq:step3-brownian-crossing}, \eqref{eq:step3-crossing-overshoot}, the choice of \(K\), and \eqref{eq:step3-local-minorization} applied at the stopping time \(\tau\), the conditional probability given \(\filt_{2^\bs}\) that the two copies of \(R\) do not coincide by the end of the block while \(\mathscr Y_\bs\) holds is at most \(1-c_*\), where \(c_*:=\varepsilon_K/2\), for all \(\bk\) large.  Any fixed number \(N_0\) of blocks fits before \(2^{\bt_\bk'}\), because \(N_0L_\bk=o(\bt_\bk'-\bt_\bk)\) by \eqref{eq:step3-near-saturation-scales} and \eqref{eq:step3-short-log-length}.  Iterating over the first \(N_0\) blocks and using Step~1 for the window, we get
\begin{equation}\label{eq:step3-repeated-coupling}
 \sup_{z,z'\in\mathcal Z_{\bk,\bt}^A}\widehat{\mathbb Q}^{z,z'}\bigl(R^{(1)}_j\neq R^{(2)}_j\text{ for all }j\leq2^{\bt_\bk'}\bigr)
 \leq(1-c_*)^{N_0}+o_\bk(1).
\end{equation}
Since \(c_*\) depends neither on \(\bk\), \(N_0\), \(z\) nor \(z'\), letting \(\bk\to\infty\) and then \(N_0\to\infty\) shows that the left side of \eqref{eq:step3-repeated-coupling} is \(o(1)\).

Fix \(z,z'\in\mathcal Z_{\bk,\bt}^A\) and put
\[
 Y^z:=\one\{\operatorname{Tru}_{2^\bt}^z\}
       e^{S_{\bk,\bt}^{\infty,z}},
 \qquad
 Y^{z'}:=\one\{\operatorname{Tru}_{2^\bt}^{z'}\}
       e^{S_{\bk,\bt}^{\infty,z'}}.
\]
Let \(\operatorname{Couple}_\bk\) be the event that
\(R_j^{(1)}=R_j^{(2)}\) for some \(j\leq2^{\bt_\bk'}\) and that both
copies satisfy \eqref{eq:step3-two-sided-window}.  By Steps~1--3,
\begin{equation}
\label{eq-071026b}
\widehat{\mathbb Q}^{z,z'}(\operatorname{Couple}_\bk^\complement)=o(1),\end{equation}
uniformly in \(z,z'\in\mathcal Z_{\bk,\bt}^A\), and
\eqref{eq:step3-FK-coupling-target} will follow from
\begin{align}
 \sup_{z,z'\in\mathcal Z_{\bk,\bt}^A}
 \widehat{\mathbb Q}^{z,z'}\left(
  (Y^z+Y^{z'})\one\{\operatorname{Couple}_\bk^\complement\}\right)
 &=o_\bk(1),
 \label{eq:step3-FK-coupling-bad}\\
 \sup_{z,z'\in\mathcal Z_{\bk,\bt}^A}
 \widehat{\mathbb Q}^{z,z'}\left(
  |Y^z-Y^{z'}|\one\{\operatorname{Couple}_\bk\}\right)
 &=o_\bk(1).
 \label{eq:step3-FK-coupling-good}
\end{align}
Split the truncation event and the clock at time \(2^{\bt_\bk'}\):
\[
 \operatorname{Tru}_{2^\bt}
 =\operatorname{Tru}_{2^\bt,2^{\bt_\bk'}}
  \cap\operatorname{Tru}_{2^{\bt_\bk'}},
 \qquad
 S_{\bk,\bt}^{\infty}=S_{\bk,\bt}^{\bt_\bk'}+S_{\bk,\bt_\bk'}^{\infty}.
\]
Since \(\mathbb Q_{\bk,\bt}^{\infty}\) and
\(\mathbb Q_{\bk,\bt}^{\bt_\bk'}\) agree on \(\filt_{2^{\bt_\bk'}}\),
\eqref{eq:tilted-truncation-cost} bounds the probability of
\(\operatorname{Tru}_{2^\bt,2^{\bt_\bk'}}^\complement\) by
\(C2^{-c_{\rm tr}\bt}\).  Using
\(S_{\bk,\bt}^{\bt_\bk'}\leq 0\) from
\eqref{eq:S-sign}, and H\"older's inequality together with
Lemma~\ref{lem:QVsmallkold} for
\(e^{S_{\bk,\bt_\bk'}^{\infty}}\one\{\operatorname{Tru}_{2^{\bt_\bk'}}\}\),
we may therefore drop the indicator
\(\one\{\operatorname{Tru}_{2^\bt,2^{\bt_\bk'}}\}\) from both \(Y^z\)
and \(Y^{z'}\) at cost \(o(1)\) in
\eqref{eq:step3-FK-coupling-bad} and \eqref{eq:step3-FK-coupling-good}.
On \(\operatorname{Couple}_\bk\), the two copies are coupled by the rotation \(B\mapsto e^{-i\phi}B\) from the meeting time \(j\leq2^{\bt_\bk'}\) onward, so that \(|Q_i^{(2)}|=|Q_i^{(1)}|\) for all \(i\geq j\).  The event \(\operatorname{Tru}_{2^{\bt_\bk'}}\) depends on the Brownian path only through \(|B_t-B_i|\) and on the \(\Gamma_i\), which the rotation leaves unchanged, and \(S_{\bk,\bt_\bk'}^{\infty}\) depends on the coefficients only through \(|Q_i|\) by \eqref{Qdecompsum}, so both are invariant under this rotation and the variable
\[
 W:=\one\{\operatorname{Tru}_{2^{\bt_\bk'}}\}
     e^{S_{\bk,\bt_\bk'}^\infty}
\]
is common to the two copies.  After the prefix truncation indicators have
been removed, the two variables in the coupling comparison factor as
\[
 Y^z=e^{S_{\bk,\bt}^{\bt_\bk',z}}W,
 \qquad
 Y^{z'}=e^{S_{\bk,\bt}^{\bt_\bk',z'}}W.
\]

It remains to show that the two factors
\(e^{S_{\bk,\bt}^{\bt_\bk',z}}\) and \(e^{S_{\bk,\bt}^{\bt_\bk',z'}}\)
are \(1+o(1)\).  On \eqref{eq:step3-two-sided-window}, by
\eqref{Qdecompsum} and \eqref{eq:step3-late-small-parameter},
\[
 |S_{\bk,\bt}^{\bt_\bk'}|
 \leq C2^{\bt_\bk'-\bk}(\bk+2)^{C_0}+C2^{-\bt}=o(1),
\]
while the complement of \eqref{eq:step3-two-sided-window} has weighted
cost \(o(\bk^{-M})\) by Step~1.  Together with the bound \eqref{eq-071026b} on
\(\widehat{\mathbb Q}^{z,z'}(\operatorname{Couple}_\bk^\complement)\),
H\"older's inequality, and Lemma~\ref{lem:QVsmallkold}, these
estimates prove \eqref{eq:step3-FK-coupling-bad} and
\eqref{eq:step3-FK-coupling-good}, hence
\eqref{eq:step3-FK-coupling-target}.  Define
\[
 c_\bk:=\mathbb Q_{\bk,\bt}^{\infty,z_0}\left(
  \operatorname{Tru}_{2^\bt}e^{S_{\bk,\bt}^{\infty}}\right),
 \qquad z_0:=r_\bk,
\]
which lies in \(\mathcal Z_{\bk,\bt}^A\) since
\((1-r_\bk^2)/(1-r_\bk^2)=1\).  Lemma~\ref{lem:QVsmallkold} gives the
upper bound in \eqref{eq:step3-ck-bounds}.  The first-moment bound
\eqref{eq-11012026dnew} and \eqref{eq:tilted-truncation-cost} give
constants \(C,p>0\) such that, on \(\mathscr W_\bt^A\),
\[
 \mathbb Q_{\bk,\bt}^{\infty}\left(
  \operatorname{Tru}_{2^\bt}\cap\{|S_{\bk,\bt}^{\infty}|\leq C\}
\right)\geq p.
\]
The conditional expectation in \eqref{eq:step3-FK-forgetting} is
consequently at least \(pe^{-C}\) on \(\mathscr W_\bt^A\), and in
particular \(c_\bk\geq pe^{-C}\), proving the lower bound in
\eqref{eq:step3-ck-bounds} and completing the proof.
\end{proof}

In the next proposition, we revert to our standard usage of $\mathbb{P}$ and $\mathbb{E}$ as the original (untilted) law and expectation over the Verblunsky coefficients $\gamma_j$.
\begin{proposition}[The expected clock weight]
\label{lem:averaged-FK-transfer}
Put
\begin{equation}\label{eq:averaged-FK-factor}
 Z_{\bk,\bn}
 :=\mathbb Q_{\bk,\bn}^{\infty}\left(
  e^{S_{\bk,\bn}^{\infty}+L_{\bk,\bn}^{\infty}}
\right).
\end{equation}
There is a deterministic \(\rho_\bn\downarrow0\) such that
\begin{equation}\label{eq:averaged-FK-transfer}
 \mathbb E\left[
  \mathscr M_{\bk,\bn}|Z_{\bk,\bn}-c_\bk|\right]
 \leq(\rho_\bn+o_\bk(1))\,
       \mathbb E\mathscr M_{\bk,\bn}.
\end{equation}
Further, in the iterated limit \(\bk\to\infty\), then
\(\bn\to\infty\),
\begin{equation}\label{eq:identify-step3-ck}
 2\pi c_\bk\,\mathbb E\mathscr M_{\bk,\bn}
 =1+o_{\bk,\bn}(1).
\end{equation}
\end{proposition}

\begin{proof}
Since \(Z_{\bk,\bn}\) involves no barrier event, we may condition
directly at time \(2^\bt\), \(\bt=\bt_\bk\), and apply
Lemma~\ref{lem:step3-late-forgetting}.  Put
\begin{equation}\label{eq:averaged-FK-access-event}
 \mathscr I_{\bk,\bn}
 :=\{R_{2^\bn}(r_\bk,\theta)\leq2^{\bk/2}\}.
\end{equation}
For fixed \(\bn\), \(Q_{2^\bn}\) is a finite Blaschke product, so by
\eqref{eq:julia-derivative},
\(R_{2^\bn}(r_\bk,\theta)\leq2\max_{|z|\leq1}|Q_{2^\bn}'(z)|<\infty\),
uniformly in \(\bk\) and \(\theta\).
For fixed \(\bn\), \(\mathscr M_{\bk,\bn}\leq C_\bn\) almost surely,
and the conditional expectation
\(\mathbb Q_{\bk,\bn}^\infty(e^{S_{\bk,\bn}^{\infty}+|L_{\bk,\bn}^{\infty}|})\)
has a bounded moment of some fixed positive order by
Lemmas~\ref{lem:QVsmallkold} and \ref{lem:zetaphitail}.  Dominated
convergence, followed by H\"older's inequality, with
\eqref{eq:tilted-truncation-cost} used to remove the indicator
\(\one\{\operatorname{Tru}_{2^\bn}\}\) present in those lemmas, gives
\begin{equation}\label{eq:averaged-FK-access-cost}
 \frac{\mathbb E\left[
  \mathscr M_{\bk,\bn}\one\{\mathscr I_{\bk,\bn}^\complement\}
  \mathbb Q_{\bk,\bn}^{\infty}\left(
  e^{S_{\bk,\bn}^{\infty}+|L_{\bk,\bn}^{\infty}|}
\right)\right]}
 {\mathbb E\mathscr M_{\bk,\bn}}
 =o_\bk(1).
\end{equation}

On \(\mathscr I_{\bk,\bn}\), let
\[
 a_{\bk,\bn}=C(2^{\bt-\bk}+\bt 2^{-\bk/2}),\qquad
 \delta_{\bk,\bn}=a_{\bk,\bn}^{1/2}=o_\bk(1).
\]
Equation~\eqref{eq-11012026enew} and Markov's inequality give
\begin{equation}\label{eq:averaged-FK-prefix-clock}
 \one\{\mathscr I_{\bk,\bn}\}
 \mathbb Q_{\bk,\bn}^{\bt}\left(
  \one\{\operatorname{Tru}_{2^\bn,2^\bt}\}
  \one\{|S_{\bk,\bn}^{\bt}|>2^{-\bn}+\delta_{\bk,\bn}\}
\right)
 \leq\delta_{\bk,\bn}.
\end{equation}
Applying \eqref{eq:step2-restoration-tail} with \(N=2^\bn\),
\(J=2^\bt\) and \(R_N\leq2^{\bk/2}\) also yields, for fixed large
\(B\),
\begin{align}
 &\one\{\mathscr I_{\bk,\bn}\}
 \mathbb Q_{\bk,\bn}^{\bt}\left(
  R_{2^\bt}>C_B2^\bt\bt^B\right)\leq C_B\bt^{-B}
 +C_B\bt^{-Bp}2^{p(\bk/2-\bt)+c_p(\bn-\bt)}
 =o_\bk(1).
 \label{eq:averaged-FK-direct-restoration}
\end{align}
The second term is exponentially small in \(\bk\) because
\(\bt-\bk/2=\bk/2-o(\bk)\).  By the conditional H\"older inequality
with Lemmas~\ref{lem:QVsmallkold} and \ref{lem:zetaphitail}, both
\eqref{eq:averaged-FK-prefix-clock} and
\eqref{eq:averaged-FK-direct-restoration} remain \(o_\bk(1)\) when the
weight \(e^{S_{\bk,\bn}^{\bt}+|L_{\bk,\bn}^{\bt}|}\) is inserted.

By the tower property and the additivity
\(S_{\bk,\bn}^{\infty}=S_{\bk,\bn}^{\bt}+S_{\bk,\bt}^{\infty}\),
\(L_{\bk,\bn}^{\infty}=L_{\bk,\bn}^{\bt}+L_{\bk,\bt}^{\infty}\), which
follows from \eqref{eq:StLt},
\begin{align}\label{eq:averaged-FK-terminal-factorization}
 Z_{\bk,\bn}
 =\mathbb Q_{\bk,\bn}^{\bt}\Big(
 e^{S_{\bk,\bn}^{\bt}+L_{\bk,\bn}^{\bt}}
 \mathbb Q_{\bk,\bt}^{\infty}\big(
 e^{S_{\bk,\bt}^{\infty}+L_{\bk,\bt}^{\infty}}
\big)
\Big).
\end{align}
On the event \(\{R_{2^\bt}\leq C_B2^\bt\bt^B\}=\mathscr W_\bt^B\),
insert the indicator \(\one\{\operatorname{Tru}_{2^\bt}\}\) into the
inner expectation at cost \(o_\bk(1)\), using
\eqref{eq:tilted-truncation-cost} and H\"older's inequality.  By
Lemma~\ref{lem:tilted-stopped-linearization} with \(\bu=\bt\), on
\(\operatorname{Tru}_{2^\bt}\) we have
\(|L_{\bk,\bt}^{\infty}|\leq\epsilon_\bt^{\rm lin}=C\bt^2 2^{-\bt/4}=o_\bk(1)\) outside
an event of conditional probability \(C\exp\{-c2^{\bt/4}\}\), so after
a further application of H\"older's inequality the factor
\(e^{L_{\bk,\bt}^{\infty}}\) may be replaced by \(1\).
Lemma~\ref{lem:step3-late-forgetting} with \(A=B\) then shows that the
inner expectation equals \(c_\bk+o_\bk(1)\) on \(\mathscr W_\bt^B\).
Lemma~\ref{lem:tilted-stopped-linearization} with \(\bu=\bn\), together
with \eqref{eq:averaged-FK-prefix-clock}, similarly gives
\[
 e^{S_{\bk,\bn}^{\bt}+L_{\bk,\bn}^{\bt}}
 =1+O(2^{-\bn})+O(\epsilon_\bn^{\rm lin})+o_\bk(1)
\]
outside an event of conditional probability
\(C\exp\{-c2^{\bn/4}\}+\delta_{\bk,\bn}\), whose weighted cost is
controlled by H\"older's inequality; the terms depending only on
\(\bn\) are collected into \(\rho_\bn\).  Finally,
\eqref{eq:averaged-FK-access-cost} removes the indicator
\(\one\{\mathscr I_{\bk,\bn}\}\).  Multiplying by
\(\mathscr M_{\bk,\bn}\) and taking expectation proves
\eqref{eq:averaged-FK-transfer}.  By Theorem~\ref{thm:fieldconverges},
\(\varphi_\infty(r_\bk,\theta)\) is a centered Gaussian variable with
variance \(-\frac4\beta\log(1-r_\bk^2)\), so that, with
\(\mathscr M_{\bk,\infty}\) defined by \eqref{eq:Mn} with
\(\bt=\infty\), \(\mathbb E\mathscr M_{\bk,\infty}=\frac1{2\pi}\).
Together with \eqref{eq:com} for \(\bt=\infty\) and \(X=1\), this
gives
\[
 \frac1{2\pi}
 =\mathbb E\mathscr M_{\bk,\infty}
 =\mathbb E\left[\mathscr M_{\bk,\bn}Z_{\bk,\bn}\right].
\]
Substitution of \eqref{eq:averaged-FK-transfer}, and the two-sided bounds
\eqref{eq:step3-ck-bounds}, give
\eqref{eq:identify-step3-ck}.
\end{proof}

\section{Approximate martingale behavior under the Brownian barrier}
\label{sec:brownian-clock-approximate-martingale}

The idea behind this section is a fact about Brownian motion.  Let
\(B\) be a standard Brownian motion, let \(x>0\), and let \(\tau\) be
the first time at which \(x+B_s\) reaches \(0\).  Then
\begin{equation}\label{eq:killed-bm-martingale}
 \mathbb E\bigl[(x+B_t)\one\{\tau>t\}\bigr]=x
 \qquad\text{for every }t\geq0 ,
\end{equation}
that is, Brownian motion killed at a straight-line barrier is a
martingale, and the expected surviving endpoint is the starting point,
however long one waits.

For our process, fix \(\bn\leq\bl\leq\bu\leq\bk\).  For
\(\bn\leq\bm\leq\bu\), define the \textbf{deficit}
\begin{equation}\label{eq:buffered-gap-definition}
\Delta_{\bk,\bm}^{[\bn]}
:=\sqrt2\lg_{r_\bk}(2^\bm)
-\sqrt{\frac{\beta}{4}}\,
\zeta_{2^\bm}(\bn;r_\bk,\theta),
\end{equation}
not to be confused with \(\Delta_{\bk,\ell}\) of
\eqref{eq-theotherDelta}.  Recall that the Gaussian field is
initialized at the left endpoint of a dyadic block,
\begin{equation}\label{eq:zeta-initial-is-field}
\zeta_{2^\bl}(\bl;r_\bk,\theta)
=\varphi_{2^\bl}(r_\bk,\theta).
\end{equation}
We call an interval \([\bl,\bu]\) of dyadic exponents a \textbf{block}.  On a block beginning at \(\bl\geq\bn\), introduce the
\textbf{Brownian clock}
\begin{equation}\label{eq:brownian-clock}
q_{\bk,\bl}(s)
:=\lg_{r_\bk}(2^\bl)
+\frac{\beta}{4}\left(
\langle\zeta(\bn;r_\bk,\theta)\rangle_s
-\langle\zeta(\bn;r_\bk,\theta)\rangle_{2^\bl}
\right).
\end{equation}
Its increment is independent of the initial exponent \(\bn\leq\bl\),
and at dyadic endpoints
\begin{equation}\label{eq:clock-at-dyadic-point}
q_{\bk,\bl}(2^\bm)
=\lg_{r_\bk}(2^\bm)+S_{\bk,\bl}^{\bm},
\qquad \bl\leq\bm\leq\bk.
\end{equation}
Define the \textbf{Brownian barrier} on this block by
\begin{equation}\label{eq:clock-adapted-barrier}
\widetilde{\mathscr E}_{\bk,\bl}^{\bu,L,a,[\bn]}
:=\left\{
\sqrt{\frac{\beta}{8}}\,
\zeta_{2^\bm}(\bn;r_\bk,\theta)
\leq q_{\bk,\bl}(2^\bm)
+a q_{\bk,\bl}(2^\bm)^{1/8}+L:
\bl\leq\bm\leq\bu
\right\}.
\end{equation}
The estimate we want on this block is the analogue of
\eqref{eq:killed-bm-martingale} for the deficit: under the tilted
measure \(\mathbb Q^{\bu}_{\bk,\bl}\), the expectation of the deficit
at the end of the block, weighted by \(e^{S}\) and killed when the
barrier fails, should be the deficit at the start,
\[
 \mathbb Q_{\bk,\bl}^{\bu}\Big(
 \one\{\operatorname{Tru}_{2^\bl,2^\bu}\}
 \one\{\widetilde{\mathscr E}_{\bk,\bl}^{\bu,L,a,[\bn]}\}
 e^{S_{\bk,\bl}^{\bu}}\Delta_{\bk,\bu}^{[\bn]}\Big)
 =(1+o(1))\,\Delta_{\bk,\bl}^{[\bn]}.
\]
The barrier
\eqref{eq:clock-adapted-barrier} is not a straight line: it is bent by
the eighth root of the clock.  This curvature does not change
the estimate, provided the deficit at the start of the block is of the
order of the square root of the clock elapsed over the block.  We call
this condition the \emph{buffer event}; it is defined in
\eqref{eq:buffer-event}, and the estimate is
Lemma~\ref{lem:buffered-one-block}.

Our goal in this section is to prove a barrier estimate on the whole time axis for which the
same holds for the full object: with a suitable good barrier
\(\widetilde{\mathscr E}\), defined in Section~\ref{sec:full-barrier},
\[
 \widetilde{\mathcal U}
 :=\mathbb Q_{\bk,\bn}^{\infty}\Big(
 \one\{\operatorname{Tru}_{2^\bn}\}
 \one\{\widetilde{\mathscr E}\}
 e^{S_{\bk,\bn}^{\infty}}
 (\Delta_{\bk,\infty}^{[\bn]})_+\Big)
 =(c_\bk+o(1))\,\Delta_{\bk,\bn}^{[\bn]} .
\]
This is Proposition~\ref{prop:brownian-clock-full-tower}, with
\(\widetilde{\mathcal U}\) defined precisely in
\eqref{eq:brownian-clock-full-operator}.  The constant \(c_\bk\) of
Lemma~\ref{lem:step3-late-forgetting} replaces \(1\) because of the
terminal block after \(\bt_\bk\), on which the Brownian time saturates
and the weight \(e^{S}\) averages to \(c_\bk\).  We prove the block
version first, in Section~\ref{sec:one-block}; it is the easier of the
two, since on one block the clock runs only for a time comparable to
the starting scale.  The full version is proved in
Section~\ref{sec:full-barrier}, by iterating the block version along
a partition of the time axis.

\subsection{Preliminaries}
\label{sec:prelim-brownian}

The following lemma collects the estimates for Brownian motion killed
at a curved barrier that we need.  Parts~(a) and~(b) remove the
curvature: for a barrier that deviates from a straight line by at most
a constant times \(\bl^{1/8}\), the expected surviving endpoint is the
starting point up to an error of that order, and the expected absolute
value of the surviving endpoint is no larger.  Part~(c) says that a
Brownian motion started in the window
\([\bl^{1/2-\nu},\bl^{1/2+\nu}]\) and surviving for a time of order
\(\bl\) ends in the corresponding window, apart from a fraction
\(\bl^{-3\nu}\) of its mass; this is what lets us stay inside the
buffer event \eqref{eq:buffer-event} from one block to the next.
Part~(d) is a weighted form of the hitting estimate, run in a random
clock, which is used only on the terminal window after \(\bt_\bk\).

\begin{lemma}[Brownian barrier estimates]\label{lem:buffered-brownian}
Fix \(0<\nu<3/8\), \(0\leq A<\infty\),
\(0<c_0<C_0<\infty\), and \(1<\lambda_0<\lambda_1<\infty\).
There is
\(\bl_0=\bl_0(A,c_0,C_0,\nu,\lambda_0,\lambda_1)\) such that the
following holds for every \(\bl\geq\bl_0\).  Let \(B\) be a standard
Brownian motion, let \(T\) be a stopping time with
\(T\in[c_0\bl,C_0\bl]\), and choose
\[
 v\in[c_0\bl,C_0\bl],\qquad |b|\leq A\bl^{1/8},\qquad
 h(s)=d\big((v+s)^{1/8}-v^{1/8}\big),\quad |d|\leq A.
\]
For \(x\in[\bl^{1/2-\nu},\bl^{1/2+\nu}]\), let
\(\tau=\inf\{s\geq0:x+B_s=b+h(s)\}\).  Then there exists $C=C(A,\lambda_0,\lambda_1,c_0,C_0,\nu)$ so that:
\begin{enumerate}
\item[(a)]
\begin{align}
 \left|\mathbb E\left[(x+B_T)\one\{\tau>T\}\right]-x\right|
 &\leq C\bl^{1/8}.
 \label{eq:buffered-brownian-main}
\end{align}
\item[(b)]
\begin{align}
 \mathbb E\left[\one\{\tau>T\}|x+B_T|\right]
 &\leq x+C\bl^{1/8}.
 \label{eq:buffered-brownian-upper}
\end{align}
\item[(c)] If, in addition, \(T\) is deterministic and
\(\lambda_0\bl\leq\bu\leq\lambda_1\bl\), put
\(y_-=\bu^{1/2-\nu}\), \(y_+=\bu^{1/2+\nu}\).
Then,
\begin{align}
 &\mathbb E\left[(1+|x+B_T|)\one\{\tau>T\}
 \one\{x+B_T<y_-+2\ \text{or}\ x+B_T>y_+-2\}\right]
 \leq C\bl^{-3\nu}x .
 \label{eq:buffered-brownian-buffer-tail}
\end{align} 
\item[(d)]
Let \(M\) be a continuous local
martingale starting from zero, with bracket \(\langle M\rangle_s\), and let \(\ell_s\)
be a deterministic increasing reference clock.  Suppose
\[
 \langle M\rangle_0=\ell_0=0,\qquad \langle M\rangle_\infty=\overline T,\qquad
 \ell_\infty=T_0,\qquad 1\leq T_0,
\]
and, for some event \(\mathscr T\), weight \(W\geq0\), and \(p_0>0\),
\begin{equation}\label{eq:weighted-brownian-clock-inputs}
 \mathbb E\left(
  \one\{\mathscr T\}e^{p_0(\overline T-T_0)_+}\right)
  \leq C_{p_0},\qquad
 \mathbb E\left(
  \one\{\mathscr T\}W^2\right)^{1/2}\leq K.
\end{equation}
Let 
$y>0$ and let $h$ be a continuous function with $h(0)=0$ and $|h(u)|\leq Au^{1/8}$. Then
\begin{align}
 \mathbb E\left(
  \one\{\mathscr T\}
  \one\{\overline T\leq2T_0\}
  \one\{\exists s:\ y+M_s\leq h(\langle M\rangle_s)\}W\right)
 &\leq C_AK
  \exp\left\{
       -c\frac{y^2}{T_0}\right\}.
 \label{eq:weighted-brownian-clock-DDS}
\end{align}
Moreover, for every \(D>2T_0\),
\begin{align}
 \mathbb E\left(
  \one\{\mathscr T\}
  \one\left\{\sup_s|\langle M\rangle_s-\ell_s|>D\right\}W\right)
 &\leq C_{p_0}Ke^{-cT_0}.
 \label{eq:weighted-brownian-clock-deviation}
\end{align}
\end{enumerate}
The constants $C_A$, $c$ are uniform over the displayed ranges and may depend on
\(A,c_0,C_0,\nu,\lambda_0,\lambda_1\).
\end{lemma}

\begin{proof}

We write $X_s=x+B_s$. Note that $X_s$ is a martingale with $X_0=x$. Now,
\begin{equation}
\label{eq-260826a}
\mathbb E[ X_T\one \{\tau>T\}]=\mathbb E[ X_T]-\mathbb E [ X_T\one\{\tau\leq T\}]=
x- \mathbb E [X_{\tau}\one\{\tau\leq T\}].
\end{equation}
Since $|X_{\tau}\one\{\tau\leq T\}|\leq C \bl^{1/8}$, we obtain \eqref{eq:buffered-brownian-main}.

Next, put $m_\bl:=A\bl^{1/8}+A(2C_0\bl)^{1/8}\leq C\bl^{1/8}$, which bounds $|b|+\sup_{s\leq T}|h(s)|$. Note that $|X_T|\one\{\tau>T\}\leq (X_T+m_\bl)\one\{\tau>T\}$ since $X_s\geq -m_\bl$ for $s\leq \tau$. But then,
\begin{equation}
\mathbb E[|X_T|\one\{\tau>T\}]\leq \mathbb E[X_T\one\{\tau>T\}]+m_\bl,
\end{equation}
and \eqref{eq:buffered-brownian-upper} then follows from  \eqref{eq:buffered-brownian-main}.

Turning to the proof of \eqref{eq:buffered-brownian-buffer-tail}, note that
\begin{eqnarray}\nonumber
&&\mathbb E[|X_T| \one\{\tau>T\}\one\{X_T\not \in [y_-+2,y_+-2]\}]\\
&\leq & \mathbb E[|X_T| \one\{X_s\geq -m_\bl, s\in[0,T]\}\one\{X_T\not \in [y_-+2,y_+-2]\}]\nonumber\\
&= &\int_{y\geq -m_\bl, y\not\in [y_-+2,y_+-2]} |y|  [  p_T(x,y)-p_T(x,-y-2m_\bl) ] dy,
\label{eq-260826c}
\end{eqnarray}
where $p_T(x,y)=(2\pi T)^{-1/2} e^{-(x-y)^2/2T}$ is the standard heat kernel, and we have used here the reflection principle. Using that $T\leq C_0\bl$ and that $y_+-x-2\geq (\lambda_0-1)\bl^{1/2+\nu}$, the contribution to the integral of the domain $y>y_+-2$ decays exponentially in a power of $\bl$, while for $y<y_-$ we have that $xy/T$ is bounded by the constant $\lambda_1/c_0$. It follows by linearizing the difference of the heat kernels that the 
integral in the right hand side of \eqref{eq-260826c} is bounded (since $x, y_->\bl^{1/8}$) by
 $$ ce^{-\bl^{c'}}+C\frac{x}{\bl^{3/2}}\int_{-m_\bl}^{y_-+2}  |y| (y+m_\bl) dy\leq  ce^{-\bl^{c'}}+C' x\frac{y_-^3}{\bl^{3/2}}\leq ce^{-\bl^{c'}}+C' x\bl^{-3\nu}.$$
This yields $\eqref{eq:buffered-brownian-buffer-tail}$.

It remains to prove the weighted random-clock conclusion. By Cauchy--Schwarz, we may and will assume that $y=\sqrt{qT_0}$ with $q\geq q_A:=\max\{1,4\cdot2^{1/4}A^2\}$ (otherwise there is nothing to prove, after enlarging $C_A$); since $T_0\geq1$, this gives $A(2T_0)^{1/8}\leq y/2$. Using 
\eqref{eq:weighted-brownian-clock-inputs}, we have with  $\mathcal{A}=\{\bar T\leq q+T_0\}$ that
\[ \mathbb E\left(
  \one\{\mathscr T\}\one\{\mathcal{A}^\complement\}
W\right)\leq KC_{p_0}^{1/2} e^{-p_0 q/2}.\] 
On the other hand, Dambis--Dubins--Schwarz gives \(M_s=B_{\langle M\rangle_s}\) for a Brownian motion $B$.  On $\mathcal A\cap\{\overline T\leq2T_0\}$ we have $\langle M\rangle_s\leq(q+T_0)\wedge2T_0$ for all $s$, and on $\{\exists s:\ y+M_s\leq h(\langle M\rangle_s)\}$, since $|h(\langle M\rangle_s)|\leq A(2T_0)^{1/8}\leq y/2$, some $|B_{\langle M\rangle_s}|$ is at least $y/2$.  Hence, by Cauchy--Schwarz,
\begin{eqnarray}
&&\mathbb E\left(
  \one\{\mathscr T\}\one\{\mathcal{A}\} \one\{\overline T\leq2T_0\}
  \one\{\exists s: y+M_s\leq h(\langle M\rangle_s)\}W\right)\\
\nonumber
&&\leq   \mathbb E\left( \one\{\mathscr T\} \one\{\max_{0\leq s\leq (q+T_0)\wedge 2T_0} |B_s|\geq \sqrt{qT_0}/2\} W\right)
\leq K e^{-cqT_0/((T_0+q)\wedge 2T_0)}\leq K e^{-cq}.
\end{eqnarray}
Combining the last two displays yields \eqref{eq:weighted-brownian-clock-DDS}.

For the final assertion, exponential Markov and
\eqref{eq:weighted-brownian-clock-inputs} give
\[
 \mathbb P\left(
  \mathscr T\cap\{\overline T>2T_0\}\right)
 \leq C_{p_0}e^{-p_0T_0}.
\]
On \(\{\overline T\leq2T_0\}\), monotonicity gives
\(0\leq \langle M\rangle_s\leq2T_0\) and \(0\leq\ell_s\leq T_0\), hence
\(\sup_s|\langle M\rangle_s-\ell_s|\leq2T_0\).  The event in
\eqref{eq:weighted-brownian-clock-deviation} is therefore contained in the
same clock-overrun event.  Another Cauchy--Schwarz estimate proves that
bound.
\end{proof}

\subsection{Brownian evaluation of one block}
\label{sec:one-block}
\label{sec-mainpropproof}

We now prove the one-block estimate described at the start of the
section, Lemma~\ref{lem:buffered-one-block}, and its one-sided
companion, Lemma~\ref{lem:one-sided-one-block}.
Fix \(\bn\leq\bl\leq\bu\leq\bk\) as at the start of the section.
For Brownian comparison we also use the continuous counterpart of
\eqref{eq:clock-adapted-barrier},
\begin{equation}\label{eq:continuous-clock-barrier}
\widetilde{\mathscr R}_{\bk,\bl}^{\bu,L,a,[\bn]}
:=\left\{
\sqrt{\frac{\beta}{8}}\,
\zeta_s(\bn;r_\bk,\theta)
\leq q_{\bk,\bl}(s)+a q_{\bk,\bl}(s)^{1/8}+L:
2^\bl\leq s\leq2^\bu
\right\}.
\end{equation}
The terminal restriction is the increasing-horizon limit
\begin{equation}\label{eq:terminal-clock-barrier}
 \widetilde{\mathscr E}_{\bk,\bl}^{\infty,L,a,[\bn]}
 :=\bigcap_{\bu\geq\bl}
 \widetilde{\mathscr E}_{\bk,\bl}^{\bu,L,a,[\bn]}.
\end{equation}
For \(0<\nu<3/8\), the \textbf{deficit-buffer event} is
\begin{equation}\label{eq:buffer-event}
\mathscr V_{\bk,\bl}^{\nu,[\bn]}
:=\left\{
\bl^{1/2-\nu}
\leq\Delta_{\bk,\bl}^{[\bn]}
\leq\bl^{1/2+\nu}
\right\}.
\end{equation}
This is typical Brownian behavior: a Brownian motion started in this
range at scale \(\bl\) and killed at the barrier is, at the end of a
block of comparable length, again in the corresponding range, up to a
small fraction of its mass; see part~(c) of
Lemma~\ref{lem:buffered-brownian}.

Corollary~\ref{cor:step2-local-restoration} controls \(R\) at the
end of a block.  With \(C'\) and \(\mathfrak c\) fixed as in
Section~\ref{subsec-2.3}, choose \(0<p<1\) and constants
\(1<\lambda_-<\lambda_+<\infty\) so that
\begin{equation}\label{eq:step2-ratio-choice}
(p+c_p)\lambda_->pC'+c_p,
\qquad c_p=\frac{2p(1-p)}{\beta}.
\end{equation}
Then, for \(B\) sufficiently large, \eqref{eq:step2-restoration-tail}
gives, uniformly for \(\lambda_-\bl\leq\bu\leq\lambda_+\bl\),
\begin{equation}\label{eq:step2-radial-failure}
 \mathbb Q_{\bk,\bl}^{\bu}\left(
  (\mathscr W_\bu^B)^{\complement}
\right)
 \leq C_B\bl^{-B}+Ce^{-c\bl}
\end{equation}
whenever \(R_{2^\bl}\leq C2^{C'\bl}\).  If
\(\mathscr W_\bl^B\) holds, the same bound holds whenever
\(\bu/\bl\) lies in a fixed compact subset of \((1,\infty)\).

\begin{lemma}[Short-time approximate martingale property]\label{lem:buffered-one-block}
Fix \(a<0\), \(L<\infty\), \(0<\nu<3/8\), and take \(B\)
sufficiently large.  Fix also \(1<\lambda_0<\lambda_1<\infty\).
There is \(\bn_*<\infty\) such that the following holds for every \(\bn\geq\bn_*\) and
\(\bn\leq\bl\leq\bu\leq\bk\).  Suppose that either
\begin{enumerate}
\item[(i)]
\[
 \lambda_-\bl\leq\bu\leq\lambda_+\bl,\qquad
 \mathfrak c\bu\leq\bk,\qquad
 R_{2^\bl}\leq C2^{C'\bl};
\]
in which case put \(\mathscr C_{\bk,\bl}^{\bu}=\Omega\) and
\(r_\bl^*=0\); or
\item[(ii)] \(\lambda_0\bl\leq\bu\leq\lambda_1\bl\),
\(\mathscr W_\bl^B\) holds, and there are
\(\mathscr C_{\bk,\bl}^{\bu}\in\filt_{2^\bu}\) such that
\begin{align}
 \max_{\bl\leq\bm\leq\bu}|S_{\bk,\bl}^{\bm}|
 &\leq\bl^{-10}
 \quad\hbox{on }\mathscr C_{\bk,\bl}^{\bu},
 &
 \mathbb Q_{\bk,\bl}^{\bu}\left(
  \one\{\operatorname{Tru}_{2^\bl,2^\bu}\}
  \one\{(\mathscr C_{\bk,\bl}^{\bu})^\complement\}
\right)
 &\leq C\bl^{-40}.
 \label{eq:one-block-clock-input}
\end{align}
In this case, let \(r_\bl^*\in[0,\bl^{1/16}]\) be a constant.
\end{enumerate}

Then, on
\(
 \mathscr V_{\bk,\bl}^{\nu,[\bn]}
 \cap\widetilde{\mathscr E}_{\bk,\bl}^{\bl,L,a,[\bn]},
\)
conclusions (a)--(c) hold, with $C=C_{a,L,\nu,\lambda_0,\lambda_1}$:
\begin{enumerate}
\item[(a)]
\begin{align}
 \left|\mathbb Q_{\bk,\bl}^{\bu}\left(
   \one\{\operatorname{Tru}_{2^\bl,2^\bu}\}
   \one\{\mathscr C_{\bk,\bl}^{\bu}\}
   \one\{\widetilde{\mathscr E}_{\bk,\bl}^{\bu,L,a,[\bn]}\}
   e^{S_{\bk,\bl}^{\bu}}\Delta_{\bk,\bu}^{[\bn]}
\right)
       -\Delta_{\bk,\bl}^{[\bn]}\right|
 &\leq
 C(\bl^{1/8}+Ce^{-c\bl^{1/16}}).
 \label{eq:buffered-one-block-main}
\end{align}
\item[(b)]
\begin{align}
 \mathbb Q_{\bk,\bl}^{\bu}\left(
   \one\{\operatorname{Tru}_{2^\bl,2^\bu}\}
   \one\{\mathscr C_{\bk,\bl}^{\bu}\}
   \one\{\widetilde{\mathscr E}_{\bk,\bl}^{\bu,L,a,[\bn]}\}
   e^{S_{\bk,\bl}^{\bu}}(1+|\Delta_{\bk,\bu}^{[\bn]}|)
\right)
 &\leq \Delta_{\bk,\bl}^{[\bn]}+C\bl^{1/8} .
 \label{eq:buffered-one-block-upper}
\end{align}
\item[(c)]
\begin{align}
 &\mathbb Q_{\bk,\bl}^{\bu}\left(
 \one\{\operatorname{Tru}_{2^\bl,2^\bu}\}
 \one\{\widetilde{\mathscr E}_{\bk,\bl}^{\bu,L,a,[\bn]}\}
 e^{S_{\bk,\bl}^{\bu}}
 (1+|\Delta_{\bk,\bu}^{[\bn]}|+r_\bl^*)
 \one\left\{(\mathscr C_{\bk,\bl}^{\bu}
              \cap\mathscr V_{\bk,\bu}^{\nu,[\bn]}
              \cap\mathscr W_\bu^B)^\complement\right\}
\right)
 \leq C\bl^{-3\nu}\Delta_{\bk,\bl}^{[\bn]} .
 \label{eq:step2-retention}
\end{align}
\end{enumerate}
\end{lemma}

For later use, set
\begin{equation}\label{eq:kappa-rho}
 \kappa:=\min\{3\nu,\tfrac38-\nu\}>0,
 \qquad
 \rho_\bl:=C\bl^{-3\nu},
\end{equation}
with \(C\) the constant in \eqref{eq:step2-retention}.  On
\(\mathscr V_{\bk,\bl}^{\nu,[\bn]}\) one has
\(\Delta_{\bk,\bl}^{[\bn]}\geq\bl^{1/2-\nu}\), so the right sides of
\eqref{eq:buffered-one-block-main}, \eqref{eq:buffered-one-block-upper}
and \eqref{eq:step2-retention} are at most
\(C\bl^{-\kappa}\Delta_{\bk,\bl}^{[\bn]}\),
\((1+C\bl^{-\kappa})\Delta_{\bk,\bl}^{[\bn]}\) and
\(\rho_\bl\Delta_{\bk,\bl}^{[\bn]}\), respectively.  Note also that,
for \(\bl\geq\bn_*\),
\begin{equation}\label{eq:buffer-implies-barrier}
 \mathscr V_{\bk,\bl}^{\nu,[\bn]}
 \subseteq
 \widetilde{\mathscr E}_{\bk,\bl}^{\bl,L,a,[\bn]}
 =\left\{\Delta_{\bk,\bl}^{[\bn]}
 \geq-\sqrt2a\,\lg_{r_\bk}(2^\bl)^{1/8}-\sqrt2L\right\},
\end{equation}
since \(\bl^{1/2-\nu}\geq-\sqrt2a(\bl\log2)^{1/8}-\sqrt2L\) for all
large \(\bl\).

\begin{proof}
For \(2^\bl\leq s\leq2^\bu\), define
\[
 \widehat\Delta_{\bk}^{[\bn]}(s)
 :=\sqrt2q_{\bk,\bl}(s)
   -\sqrt{\frac\beta4}\,\zeta_s(\bn;r_\bk,\theta).
\]
Under \(\mathbb Q_{\bk,\bl}^{\bu}\), this process is a continuous
martingale with initial value \(\Delta_{\bk,\bl}^{[\bn]}\) and quadratic
variation \(q_{\bk,\bl}(s)-q_{\bk,\bl}(2^\bl)\).  Thus,
Dambis--Dubins--Schwarz gives a Brownian motion \(B\) adapted to its filtration $\mathcal{G}_t$ such that
\begin{equation}\label{eq:clock-gap-Brownian}
 \widehat\Delta_{\bk}^{[\bn]}(s)
 =\Delta_{\bk,\bl}^{[\bn]}
  -B_{q_{\bk,\bl}(s)-q_{\bk,\bl}(2^\bl)}.
\end{equation}
We use \eqref{eq:clock-gap-Brownian} only for \(s\leq2^\bu\).  Set
\(T=q_{\bk,\bl}(2^\bu)-q_{\bk,\bl}(2^\bl)\), the Brownian time
corresponding to the endpoint \(s=2^\bu\).  
  If
necessary, extend \(B\) and its filtration $\mathcal{G}_t$ beyond \(T\) by a Brownian motion independent of
\(\mathcal G_T\).
At \(s=2^\bm\),
\[
 \widehat\Delta_{\bk}^{[\bn]}(2^\bm)
 =\Delta_{\bk,\bm}^{[\bn]}+\sqrt2S_{\bk,\bl}^{\bm}.
\]
Set
\[
 \tau_\bm=q_{\bk,\bl}(2^\bm)-q_{\bk,\bl}(2^\bl),
 \qquad\bl\leq\bm\leq\bu.
\]
Each \(\tau_\bm\) is a \(\mathcal G_t\)-stopping time.  Moreover, since
\(|Q_j|\leq1\),
\begin{equation}\label{eq:dyadic-bracket-step}
 0\leq\tau_{\bm+1}-\tau_\bm
 =q_{\bk,\bl}(2^{\bm+1})-q_{\bk,\bl}(2^\bm)\leq C .
\end{equation}

Put \(w_\bl=\bl^{1/32}\).  Let 
\begin{equation}
\label{eq-300826a}
\mathscr O_{\bl,\bu}=\{ \max_{t\in[\tau_\bm,\tau_{\bm+1}]} B_t- \min_{t\in[\tau_\bm,\tau_{\bm+1}]} B_t\leq w_\bl, \quad \bl\leq \bm< \bu \} 
\end{equation}
be the event on
which the Brownian oscillation in every
\([\tau_\bm,\tau_{\bm+1}]\) is at most \(w_\bl\), and let
\begin{equation}\label{eq:buffered-good-clock-event}
 \mathscr S_{\bk,\bl}^{\bu}
 :=\left\{\max_{\bl\leq\bm\leq\bu}
 |S_{\bk,\bl}^{\bm}|\leq e^{-c\bl/2}\right\},
\end{equation}
for a fixed \(c\in(0,1]\).
On \(\mathscr S_{\bk,\bl}^{\bu}\), we have
\(T=\lg_{r_\bk}(2^\bu)-\lg_{r_\bk}(2^\bl)
+S_{\bk,\bl}^{\bu}\asymp\bl\).  Thus, we can choose \(c_0,C_0\) so that
\(T\in[c_0\bl,C_0\bl]\) on \(\mathscr S_{\bk,\bl}^{\bu}\).  Let
\(\overline T=(T\vee c_0\bl)\wedge C_0\bl\).  Since
\(T=\tau_\bu\), \(\overline T\) is a \(\mathcal G_t\)-stopping time,
and \(\overline T=T\) on \(\mathscr S_{\bk,\bl}^{\bu}\).

We first complete the proof of the lemma for case~(i), and then indicate how the proof needs to be modified in case~(ii).
We begin by proving part~(a).  Since
\(\mathscr C_{\bk,\bl}^{\bu}=\Omega\) in this case, set
\begin{equation}
\tag{Err-$\mathcal{D}$}
\label{eq:one-block-target-error}
 \mathcal D:=\left|
 \mathbb Q_{\bk,\bl}^{\bu}\left(
  \one\{\operatorname{Tru}_{2^\bl,2^\bu}\}
  \one\{\widetilde{\mathscr E}_{\bk,\bl}^{\bu,L,a,[\bn]}\}
  e^{S_{\bk,\bl}^{\bu}}\Delta_{\bk,\bu}^{[\bn]}
\right)
 -\Delta_{\bk,\bl}^{[\bn]}\right|.
\end{equation}

Put \(Y_t=\Delta_{\bk,\bl}^{[\bn]}-B_t\).  To estimate
\eqref{eq:one-block-target-error}, we replace its discrete barrier by a
Brownian barrier and its weighted endpoint
\(e^{S_{\bk,\bl}^{\bu}}\Delta_{\bk,\bu}^{[\bn]}\) by
\(Y_{\overline T}\).  For \(L'=L\pm Cw_\bl\), define
\begin{equation}\label{eq:clipped-Brownian-barrier}
 \mathscr R_B^{L'}(\overline T)
 :=\left\{Y_t
 \geq-\sqrt2a(q_{\bk,\bl}(2^\bl)+t)^{1/8}-\sqrt2L':
 0\leq t\leq\overline T\right\}.
\end{equation}
Recalling \eqref{eq:continuous-clock-barrier} and \eqref{eq:clock-gap-Brownian}, we have on \(\mathscr S_{\bk,\bl}^{\bu}\) that
\begin{equation}\label{eq:good-clock-Brownian-identification}
 \mathscr R_B^{L'}(\overline T)
 =\widetilde{\mathscr R}_{\bk,\bl}^{\bu,L',a,[\bn]},
 \qquad
 Y_{\overline T}=\widehat\Delta_{\bk}^{[\bn]}(2^\bu).
\end{equation}

We decompose the error in \eqref{eq:one-block-target-error}. First, define the  main Brownian comparison error
\begin{equation}\tag{Err-I}\label{eq:one-block-error-I}
 \mathrm I:=
 \max_{L'=L\pm Cw_\bl}
 \left|\mathbb Q_{\bk,\bl}^{\bu}\left(
  Y_{\overline T}\one\{\mathscr R_B^{L'}(\overline T)\}
\right)
  -\Delta_{\bk,\bl}^{[\bn]}\right|.
\end{equation}
To concisely describe the other errors involved, introduce the following notation.
For an event \(A\), define
\begin{align}
 \operatorname{Err}_{\rm F}(A)
 &:=
 \mathbb Q_{\bk,\bl}^{\bu}\left(
  \one\{A\}\one\{\widetilde{\mathscr E}_{\bk,\bl}^{\bu,L,a,[\bn]}\}
  e^{S_{\bk,\bl}^{\bu}}(1+|\Delta_{\bk,\bu}^{[\bn]}|)
\right),
 \label{eq:one-block-field-error-operator}\\
 \operatorname{Err}_{\rm B}(A;L')
 &:=
 \mathbb Q_{\bk,\bl}^{\bu}\left(
  \one\{A\}\one\{\mathscr R_B^{L'}(\overline T)\}
  (1+|Y_{\overline T}|)
\right).
 \label{eq:one-block-Brownian-error-operator}
\end{align}
Using this notation, define the remaining four errors:
\begin{align}
 \mathrm{II}&:=
 \operatorname{Err}_{\rm F}\!\left(
  \operatorname{Tru}_{2^\bl,2^\bu}
  \cap(\mathscr S_{\bk,\bl}^{\bu})^\complement\right)
 +\max_{L'=L\pm Cw_\bl}
 \operatorname{Err}_{\rm B}\!\left(
  \operatorname{Tru}_{2^\bl,2^\bu}
  \cap(\mathscr S_{\bk,\bl}^{\bu})^\complement;L'\right),
 \tag{Err-II}\label{eq:one-block-error-II}\\
 \mathrm{III}&:=
 \operatorname{Err}_{\rm F}\!\left(
  \operatorname{Tru}_{2^\bl,2^\bu}\cap\mathscr S_{\bk,\bl}^{\bu}
  \cap\mathscr O_{\bl,\bu}^\complement\right)
 +\max_{L'=L\pm Cw_\bl}
 \operatorname{Err}_{\rm B}\!\left(
  \operatorname{Tru}_{2^\bl,2^\bu}\cap\mathscr S_{\bk,\bl}^{\bu}
  \cap\mathscr O_{\bl,\bu}^\complement;L'\right),
 \tag{Err-III}\label{eq:one-block-error-III}\\
 \mathrm{IV}&:=
 \max_{L'=L\pm Cw_\bl}
 \operatorname{Err}_{\rm B}\!\left(
  \operatorname{Tru}_{2^\bl,2^\bu}^\complement;L'\right),
 \tag{Err-IV}\label{eq:one-block-error-IV}\\
 \mathrm V&:=
 \mathbb Q_{\bk,\bl}^{\bu}\left(
  \one\{\operatorname{Tru}_{2^\bl,2^\bu}
         \cap\mathscr S_{\bk,\bl}^{\bu}\cap\mathscr O_{\bl,\bu}\}
  \one\{\widetilde{\mathscr E}_{\bk,\bl}^{\bu,L,a,[\bn]}\}
  \left|e^{S_{\bk,\bl}^{\bu}}\Delta_{\bk,\bu}^{[\bn]}
  -Y_{\overline T}\right|
\right).
 \tag{Err-V}\label{eq:one-block-error-V}
\end{align}

We next claim that
\begin{equation}
\label{eq:one-block-error-reduction}
 \mathcal D\leq\mathrm I+\mathrm{II}+\mathrm{III}+\mathrm{IV}+\mathrm V.
\end{equation}
Indeed, start with the expectation in
\eqref{eq:one-block-target-error}, which already contains
\(\one\{\operatorname{Tru}_{2^\bl,2^\bu}\}\).  The first summands in
\(\mathrm{II}\) and \(\mathrm{III}\) allow us to restrict this expectation
to
\(\operatorname{Tru}_{2^\bl,2^\bu}\cap
\mathscr S_{\bk,\bl}^{\bu}\cap\mathscr O_{\bl,\bu}\).  On this event,
\(\mathrm V\) replaces
\(e^{S_{\bk,\bl}^{\bu}}\Delta_{\bk,\bu}^{[\bn]}\) by
\(Y_{\overline T}\). 
Recalling \eqref{eq:good-clock-Brownian-identification}, we still need to
 replace the barrier
\(\widetilde{\mathscr E}_{\bk,\bl}^{\bu,L,a,[\bn]}\) by the Brownian
barrier $\mathscr R_B^{L'}(\overline T)$. Toward this end, note that
%
during one interval in \eqref{eq:dyadic-bracket-step}, the curved part of the barrier changes
by \(O(\bl^{-7/8})\).  On the dyadic oscillation event
\(\mathscr O_{\bl,\bu}\), the Brownian path changes by at most \(w_\bl\)
on the same interval.  Therefore
\begin{equation}\label{eq:continuous-discrete-sandwich}
 \widetilde{\mathscr R}_{\bk,\bl}^{\bu,L-Cw_\bl,a,[\bn]}
 \subseteq\widetilde{\mathscr E}_{\bk,\bl}^{\bu,L,a,[\bn]}
 \subseteq\widetilde{\mathscr R}_{\bk,\bl}^{\bu,L+Cw_\bl,a,[\bn]}.
\end{equation}
Together with \eqref{eq:good-clock-Brownian-identification}, this gives
on
\(\operatorname{Tru}_{2^\bl,2^\bu}
\cap\mathscr S_{\bk,\bl}^{\bu}\cap\mathscr O_{\bl,\bu}\),
\begin{equation}\label{eq:one-block-good-reduction}
\begin{gathered}
 \mathscr R_B^{L-Cw_\bl}(\overline T)
 \subseteq\widetilde{\mathscr E}_{\bk,\bl}^{\bu,L,a,[\bn]}
 \subseteq\mathscr R_B^{L+Cw_\bl}(\overline T),\\
 \Delta_{\bk,\bu}^{[\bn]}
 =Y_{\overline T}-\sqrt2S_{\bk,\bl}^{\bu},
 \qquad e^{S_{\bk,\bl}^{\bu}}=1+O(e^{-c\bl/2}).
\end{gathered}
\end{equation}
For large \(\bl\), \eqref{eq:clipped-Brownian-barrier} gives, on
\(\mathscr R_B^{L+Cw_\bl}(\overline T)\),
\[
 Y_{\overline T}
 \geq-\sqrt2a(q_{\bk,\bl}(2^\bl)+\overline T)^{1/8}
       -\sqrt2(L+Cw_\bl)\geq0.
\]
Thus \eqref{eq:one-block-good-reduction} and the nonnegativity of
\(Y_{\overline T}\) sandwich the resulting expectation between the
corresponding expectations with
\(\mathscr R_B^{L-Cw_\bl}(\overline T)\) and
\(\mathscr R_B^{L+Cw_\bl}(\overline T)\).  The second summands in
\(\mathrm{II}\) and \(\mathrm{III}\), together with \(\mathrm{IV}\),
remove the restriction to the common event from these two bounds.  The
resulting expectations are precisely those in \(\mathrm I\), which proves
\eqref{eq:one-block-error-reduction}.

It remains to bound the five error terms.  For
\eqref{eq:one-block-error-I}, put
\[
 v=q_{\bk,\bl}(2^\bl),\qquad
 b=-\sqrt2a\,v^{1/8}-\sqrt2L',\qquad
 h(t)=-\sqrt2a\big((v+t)^{1/8}-v^{1/8}\big).
\]
Here \(v\asymp\bl\), \(|b|\leq C_{a,L}\bl^{1/8}\).
Part (a)  of  Lemma
\ref{lem:buffered-brownian}, used with \(-B\), gives for both values of
\(L'\)
\begin{equation}\tag{I}\label{eq:one-block-Brownian-conclusion}
I= \left|\mathbb Q_{\bk,\bl}^{\bu}\left(
  Y_{\overline T}\one\{\mathscr R_B^{L'}(\overline T)\}
\right)-\Delta_{\bk,\bl}^{[\bn]}\right|
\leq C\bl^{1/8}
.
\end{equation}

We next record the endpoint moments used for the exceptional events.
Since \(\bu\leq\bk\), \eqref{eq:S-sign}, positivity of the bracket and
\(\bu-\bl=O(\bl)\) imply
\[
 -C\bl\leq S_{\bk,\bl}^{\bu}\leq 0,
 \qquad e^{S_{\bk,\bl}^{\bu}}\leq 1.
\]
Optional sampling in \eqref{eq:clock-gap-Brownian} gives
\begin{align}
 \mathbb Q_{\bk,\bl}^{\bu}\left(
  (1+|\Delta_{\bk,\bu}^{[\bn]}|)^2
\right)
 &\leq C\left(1+(\Delta_{\bk,\bl}^{[\bn]})^2+\bl^2\right),
 \label{eq:buffered-unconditional-field-L2}\\
 \mathbb  Q_{\bk,\bl}^{\bu}\left(
  (1+|Y_{\overline T}|)^2
\right)
 &\leq C\left(1+(\Delta_{\bk,\bl}^{[\bn]})^2+\bl\right).
 \label{eq:buffered-unconditional-clipped-L2}
\end{align}

We next turn to the control of $\mathrm{II}$.
Using
\eqref{eq-11012026enew}
and
\(R_{2^\bl}\leq C2^{C'\bl}\) we obtain
\[
 \mathbb Q_{\bk,\bl}^{\bu}\left(
  (|S_{\bk,\bl}^{\bu}|-2^{-\bl})_+
  \one\{\operatorname{Tru}_{2^\bl,2^\bu}\}
\right)
 \leq C(2^\bu+2^{C'\bl}\bu)2^{-\bk}.
\]
By \eqref{eq:S-sign},
\(\max_{\bl\leq\bm\leq\bu}|S_{\bk,\bl}^{\bm}|=|S_{\bk,\bl}^{\bu}|\),
and \(2^{-\bl}\leq e^{-c\bl/2}/2\) for \(\bl\) large, so Markov's
inequality at level \(e^{-c\bl/2}/2\), together with
\(\mathfrak c\bu\leq\bk\), \(\mathfrak c>C'+2\), yields
\begin{equation}\label{eq:buffered-good-clock-cost}
 \mathbb Q_{\bk,\bl}^{\bu}\left(
  \one\{\operatorname{Tru}_{2^\bl,2^\bu}\}
  \one\{(\mathscr S_{\bk,\bl}^{\bu})^\complement\}
\right)
 \leq Ce^{-c\bl/2}.
\end{equation}
Cauchy--Schwarz, \eqref{eq:buffered-unconditional-field-L2} and
\eqref{eq:buffered-unconditional-clipped-L2} imply
\begin{equation}\tag{II}\label{eq:one-block-error-II-bound}
 \mathrm{II}\leq Ce^{-c'\bl}
 \Delta_{\bk,\bl}^{[\bn]}.
\end{equation}

For \eqref{eq:one-block-error-III}, \eqref{eq:dyadic-bracket-step}, the
strong Markov property, the reflection principle, and a union bound give
\begin{equation}\label{eq:dyadic-oscillation-cost}
 \mathbb Q_{\bk,\bl}^{\bu}\left(
  \mathscr O_{\bl,\bu}^\complement
\right)
 \leq C\bl e^{-cw_\bl^2}
 \leq Ce^{-c\bl^{1/16}}.
\end{equation}
Thus
\begin{equation}\tag{III}\label{eq:one-block-error-III-bound}
 \mathrm{III}\leq Ce^{-c'\bl^{1/16}}
 \Delta_{\bk,\bl}^{[\bn]}.
\end{equation}
For \eqref{eq:one-block-error-IV}, the block truncation event has exactly
the indices of the coefficients evolved on this block.  Hence
\eqref{eq:tilted-truncation-cost} and Cauchy--Schwarz give
\begin{equation}\tag{IV}\label{eq:one-block-error-IV-bound}
 \mathrm{IV}\leq Ce^{-c\bl}\Delta_{\bk,\bl}^{[\bn]}.
\end{equation}
For \eqref{eq:one-block-error-V}, \eqref{eq:one-block-good-reduction} gives
\(|e^{S_{\bk,\bl}^{\bu}}\Delta_{\bk,\bu}^{[\bn]}-Y_{\overline T}|
\leq Ce^{-c\bl/2}(1+|Y_{\overline T}|)\) on the event in
\eqref{eq:one-block-error-V}, so
\eqref{eq:buffered-unconditional-clipped-L2} gives
\begin{equation}\tag{V}\label{eq:one-block-error-V-bound}
 \mathrm V\leq Ce^{-c'\bl}\Delta_{\bk,\bl}^{[\bn]}.
\end{equation}
Combining \(\mathrm I\)--\(\mathrm V\) proves
\eqref{eq:buffered-one-block-main} in case (i).

For \eqref{eq:buffered-one-block-upper}, repeat the same decomposition
with \(e^{S_{\bk,\bl}^{\bu}}(1+|\Delta_{\bk,\bu}^{[\bn]}|)\) in place
of the signed endpoint.  Part (b) of Lemma~\ref{lem:buffered-brownian},
in place of part (a), gives the main term
\(\Delta_{\bk,\bl}^{[\bn]}+C\bl^{1/8}\), and
\(\mathrm{II}\)--\(\mathrm V\) remove the same exceptional events.
No new estimate is required.

For \eqref{eq:step2-retention}, the comparison with errors
\(\mathrm I\)--\(\mathrm V\) has the following form.  Error
\(\mathrm I\) is replaced by the endpoint-buffer estimate
\eqref{eq:buffered-brownian-buffer-tail}.  Errors \(\mathrm{II}\) and
\(\mathrm{III}\) still remove the clock and dyadic-oscillation
complements.  There is no error \(\mathrm{IV}\): the expectation is
nonnegative, so its truncation indicator can be discarded before the
Brownian estimate is applied.  The comparison corresponding to
\(\mathrm V\) uses the event \(\mathscr O_\bl^*\) defined below to move
the terminal value from the random time \(T\) to the deterministic time
\(T_*\).  The contribution of \((\mathscr W_\bu^B)^\complement\) is not
one of these five errors and is estimated separately.
We now provide the details.

Recall we are in case~(i), where
\(\mathscr C_{\bk,\bl}^{\bu}=\Omega\) and \(r_\bl^*=0\).  Put
\[
 t_0=\lg_{r_\bk}(2^\bu)-\lg_{r_\bk}(2^\bl)\asymp\bl,\qquad
 \delta_\bl=e^{-c\bl/2},\qquad T_*=t_0-\delta_\bl,
\]
and define
\begin{equation}\label{eq:step2-modulus-event}
 \mathscr O_\bl^*
 =\left\{\sup_{s,t\in[t_0-\delta_\bl,t_0+\delta_\bl]}
 |B_s-B_t|\leq1\right\}.
\end{equation}
Thus \(\mathscr O_\bl^*\) controls the Brownian motion throughout the
interval containing the terminal Brownian time.  The reflection principle
gives
\begin{equation}\label{eq:step2-modulus-cost}
 \mathbb Q_{\bk,\bl}^{\bu}\left(
  (\mathscr O_\bl^*)^\complement
\right)
 \leq C\exp(-c/\delta_\bl).
\end{equation}
On \(\mathscr S_{\bk,\bl}^{\bu}\),
\[
 T=t_0+S_{\bk,\bl}^{\bu}
 \in[t_0-\delta_\bl,t_0+\delta_\bl],
 \qquad T_*\leq T.
\]
Recall that \(Y_t=\Delta_{\bk,\bl}^{[\bn]}-B_t\), and put
\[
 b_{L'}(t)=-\sqrt2a(q_{\bk,\bl}(2^\bl)+t)^{1/8}-\sqrt2L'.
\]
For the final indicator in \eqref{eq:step2-retention}, we use
\[
 (\mathscr V_{\bk,\bu}^{\nu,[\bn]}\cap\mathscr W_\bu^B)^\complement
 = (\mathscr V_{\bk,\bu}^{\nu,[\bn]})^\complement
 \mathbin{\dot\cup}
 \bigl(\mathscr V_{\bk,\bu}^{\nu,[\bn]}
       \cap(\mathscr W_\bu^B)^\complement\bigr).
\]
On the event
\[
 \widetilde{\mathscr E}_{\bk,\bl}^{\bu,L,a,[\bn]}
 \cap\mathscr S_{\bk,\bl}^{\bu}
 \cap\mathscr O_{\bl,\bu}
 \cap\mathscr O_\bl^*,
\]
the barrier sandwich \eqref{eq:continuous-discrete-sandwich} implies
\(Y_t\geq b_{L+Cw_\bl}(t)\) for \(0\leq t\leq T_*\).  Moreover,
\[
 |\Delta_{\bk,\bu}^{[\bn]}-Y_{T_*}|
 \leq |B_T-B_{T_*}|+\sqrt2|S_{\bk,\bl}^{\bu}|\leq2.
\]
It follows that
\begin{align}
\label{eq-290826a}
 \{\Delta_{\bk,\bu}^{[\bn]}<\bu^{1/2-\nu}\}
 &\subseteq\{Y_{T_*}<\bu^{1/2-\nu}+2\},
 \{\Delta_{\bk,\bu}^{[\bn]}>\bu^{1/2+\nu}\}
 \subseteq\{Y_{T_*}>\bu^{1/2+\nu}-2\},
 1+|\Delta_{\bk,\bu}^{[\bn]}|\leq3+|Y_{T_*}|.
\end{align}
Part~(c) of Lemma~\ref{lem:buffered-brownian}, applied with \(-B\) at
the deterministic time \(T_*\), therefore bounds the contribution on
\[
 \operatorname{Tru}_{2^\bl,2^\bu}
 \cap\widetilde{\mathscr E}_{\bk,\bl}^{\bu,L,a,[\bn]}
 \cap\mathscr S_{\bk,\bl}^{\bu}
 \cap\mathscr O_{\bl,\bu}
 \cap\mathscr O_\bl^*
 \cap(\mathscr V_{\bk,\bu}^{\nu,[\bn]})^\complement
\]
by \(C\bl^{-3\nu}\Delta_{\bk,\bl}^{[\bn]}\).  In applying the lemma,
we discard the truncation indicator and use
\(e^{S_{\bk,\bl}^{\bu}}\leq e^{\delta_\bl}\).

The contribution of \((\mathscr S_{\bk,\bl}^{\bu})^\complement\) is the
first summand of \eqref{eq:one-block-error-II}, while the contribution of
\(\mathscr S_{\bk,\bl}^{\bu}\cap
(\mathscr O_{\bl,\bu})^\complement\) is the first summand of
\eqref{eq:one-block-error-III}.  They are bounded by
\eqref{eq:one-block-error-II-bound} and
\eqref{eq:one-block-error-III-bound}.  The contribution of
\((\mathscr O_\bl^*)^\complement\) is bounded by conditional
Cauchy--Schwarz, \eqref{eq:step2-modulus-cost}, and
\eqref{eq:buffered-unconditional-field-L2}.  Each of these three bounds
is \(o(\bl^{-3\nu})\Delta_{\bk,\bl}^{[\bn]}\).

It remains to bound the contribution of
\((\mathscr W_\bu^B)^\complement\).  On the terminal barrier
\(\widetilde{\mathscr E}_{\bk,\bl}^{\bu,L,a,[\bn]}\), one has
\(\Delta_{\bk,\bu}^{[\bn]}>0\) for all large \(\bl\).  Thus
\eqref{eq:step2-radial-failure}, conditional Cauchy--Schwarz,
\eqref{eq:buffered-unconditional-field-L2}, and
\(e^{S_{\bk,\bl}^{\bu}}\leq 1\) give
\[
 \mathbb Q_{\bk,\bl}^{\bu}\left(
  \one\{\operatorname{Tru}_{2^\bl,2^\bu}\}
  \one\{\widetilde{\mathscr E}_{\bk,\bl}^{\bu,L,a,[\bn]}\}
  e^{S_{\bk,\bl}^{\bu}}(1+\Delta_{\bk,\bu}^{[\bn]})
  \one\{(\mathscr W_\bu^B)^\complement\}
\right)
 \leq C\bl^{-B/2}(1+\Delta_{\bk,\bl}^{[\bn]}+\bl).
\]
Since \(\Delta_{\bk,\bl}^{[\bn]}\geq\bl^{1/2-\nu}\), choosing \(B\)
large makes the last display at most
\(C\bl^{-3\nu}\Delta_{\bk,\bl}^{[\bn]}\).  This proves
\eqref{eq:step2-retention} in case~(i), and completes the proof of the lemma in that case.

For the conclusions in case~(ii), in the decomposition
\eqref{eq:one-block-error-reduction}, replace
\(\mathscr S_{\bk,\bl}^{\bu}\) by
\(\mathscr C_{\bk,\bl}^{\bu}\).  The first condition in
\eqref{eq:one-block-clock-input} places \(T\) within \(\bl^{-10}\) of
its deterministic value and gives
\[
 e^{S_{\bk,\bl}^{\bu}}=1+O(\bl^{-10}),\qquad
 \Delta_{\bk,\bu}^{[\bn]}
 =\Delta_{\bk,\bl}^{[\bn]}-B_T+O(\bl^{-10}).
\]
Thus \(\mathrm I\), \(\mathrm{III}\), and \(\mathrm{IV}\) are
unchanged, while \(\mathrm V=O(\bl^{-10})
\Delta_{\bk,\bl}^{[\bn]}\).  The second condition in
\eqref{eq:one-block-clock-input}, conditional Cauchy--Schwarz, and
\eqref{eq:buffered-unconditional-field-L2} give
\begin{align}
 &\mathbb Q_{\bk,\bl}^{\bu}\left(
   \one\{\operatorname{Tru}_{2^\bl,2^\bu}\}
   \one\{(\mathscr C_{\bk,\bl}^{\bu})^\complement\}
   e^{S_{\bk,\bl}^{\bu}}
   (1+|\Delta_{\bk,\bu}^{[\bn]}|+r_\bl^*)
\right)
 \leq C\bl^{-20}
  (1+\Delta_{\bk,\bl}^{[\bn]}+\bl+r_\bl^*)
 \leq C\bl^{-10}\Delta_{\bk,\bl}^{[\bn]}.
 \label{eq:one-block-clock-complement-cost}
\end{align}
The last inequality uses
\(\Delta_{\bk,\bl}^{[\bn]}\geq\bl^{1/2-\nu}\) and
\(r_\bl^*\leq\bl^{1/16}\).  This bounds the contribution of
\((\mathscr C_{\bk,\bl}^{\bu})^\complement\) in
\eqref{eq:step2-retention} and gives the new estimate for
\(\mathrm{II}\).  The signed and positive reductions now
give \eqref{eq:buffered-one-block-main},
\eqref{eq:buffered-one-block-upper}.

To complete the proof of \eqref{eq:step2-retention} in case~(ii), split its final complement
as the disjoint union
\begin{align*}
 &(\mathscr C_{\bk,\bl}^{\bu}
   \cap\mathscr V_{\bk,\bu}^{\nu,[\bn]}
   \cap\mathscr W_\bu^B)^\complement
=(\mathscr C_{\bk,\bl}^{\bu})^\complement
 \mathbin{\dot\cup}
 \bigl(\mathscr C_{\bk,\bl}^{\bu}
       \cap(\mathscr V_{\bk,\bu}^{\nu,[\bn]})^\complement\bigr)
 \mathbin{\dot\cup}
 \bigl(\mathscr C_{\bk,\bl}^{\bu}
       \cap\mathscr V_{\bk,\bu}^{\nu,[\bn]}
       \cap(\mathscr W_\bu^B)^\complement\bigr).
\end{align*}
The first term is bounded by
\eqref{eq:one-block-clock-complement-cost}.  For the second term, repeat
the case~(i) Brownian comparison with
\(\mathscr C_{\bk,\bl}^{\bu}\) in place of
\(\mathscr S_{\bk,\bl}^{\bu}\).  Redefine
\(\delta_\bl=\bl^{-10}\) and \(\mathscr O_\bl^*\) by
\eqref{eq:step2-modulus-event}.  Then
\eqref{eq:step2-modulus-cost} bounds the contribution of
\((\mathscr O_\bl^*)^\complement\) using the probability
\(Ce^{-c\bl^{10}}\).  The contribution of
\((\mathscr O_{\bl,\bu})^\complement\) is still bounded by
\eqref{eq:one-block-error-III-bound}, and the truncation indicator may
again be discarded.
The weight \(1+|\Delta_{\bk,\bu}^{[\bn]}|\) is at most
\(3+|Y_{T_*}|\) by \eqref{eq-290826a}, so part~(c) of
Lemma~\ref{lem:buffered-brownian}, applied at the deterministic time
\(T_*\) as in case~(i), bounds its contribution to the second term by
\(C\bl^{-3\nu}\Delta_{\bk,\bl}^{[\bn]}\).  The additional weight
\(r_\bl^*\) is a constant, so for it we need only the probability that
the killed Brownian motion ends outside \([y_-+2,y_+-2]\).  Write
\(x=\Delta_{\bk,\bl}^{[\bn]}\).  For the lower endpoint, the heat-kernel
bound in the proof of part~(c), below \eqref{eq-260826c}, gives the
probability bound
\[
 C\frac{x}{\bl^{3/2}}\int_{-m_\bl}^{y_-+2}(y+m_\bl)\,\dif y
 \leq Cx\,\bl^{-1/2-2\nu}.
\]
Since \(r_\bl^*\leq\bl^{1/16}\), this weight contributes at most
\(C\bl^{-7/16-2\nu}x\leq C\bl^{-3\nu}x\).  For the upper endpoint, the
Gaussian tail gives a probability at most \(e^{-c\bl^{2\nu}}\), since
\(y_+-x\geq c\bl^{1/2+\nu}\) and \(T_*\leq C\bl\).

For the third term in the disjoint union,
\(\mathscr W_\bl^B\) allows us to apply
\eqref{eq:step2-radial-failure} for the fixed block ratio.  Conditional
Cauchy--Schwarz and \eqref{eq:buffered-unconditional-field-L2}, together
with \(e^{S_{\bk,\bl}^{\bu}}\leq e^{\bl^{-10}}\) on
\(\mathscr C_{\bk,\bl}^{\bu}\), bound this term by
\[
 C\bl^{-B/2}
 (1+\Delta_{\bk,\bl}^{[\bn]}+\bl+r_\bl^*).
\]
Choose \(B\) large and use
\(\Delta_{\bk,\bl}^{[\bn]}\geq\bl^{1/2-\nu}\).  The three terms in the
disjoint union are then bounded by
\(C\bl^{-3\nu}\Delta_{\bk,\bl}^{[\bn]}\), which proves
\eqref{eq:step2-retention} in case~(ii).
\end{proof}

\begin{lemma}[One-sided one-block bound]
\label{lem:one-sided-one-block}
Fix \(a<0\), \(L<\infty\), and \(1<\lambda_1<\infty\).  There are
\(\bn_*<\infty\) and a constant \(C\), depending on \(a,L,\lambda_1\),
such that for every \(\bn\geq\bn_*\) and
\[
 \bn\leq\bl<\bu\leq\min\{\lambda_1\bl,\bk\},
\]
the following holds on the initial barrier event
\(\widetilde{\mathscr E}_{\bk,\bl}^{\bl,L,a,[\bn]}\):
\begin{align}
 &\mathbb Q_{\bk,\bl}^{\bu}\left(
   \one\{\operatorname{Tru}_{2^\bl,2^\bu}\}
   \one\{\widetilde{\mathscr E}_{\bk,\bl}^{\bu,L,a,[\bn]}\}
   e^{S_{\bk,\bl}^{\bu}}
   (1+|\Delta_{\bk,\bu}^{[\bn]}|)
\right)
 \leq(1+C\bl^{-1/8})\Delta_{\bk,\bl}^{[\bn]}.
 \label{eq:buffered-one-block-off-buffer}
\end{align}
\end{lemma}

\begin{proof}
We use the conventions from the proof of
Lemma~\ref{lem:buffered-one-block}, including \(B\), \((\mathcal G_t)\),
\(T\), \(\tau_\bm\), 
\(Y_t=\Delta_{\bk,\bl}^{[\bn]}-B_t\), and the event  $\mathscr O_{\bl,\bu}$ of \eqref{eq-300826a}.
Put
\[
 x=\Delta_{\bk,\bl}^{[\bn]},\qquad
 t_0=\lg_{r_\bk}(2^\bu)-\lg_{r_\bk}(2^\bl)\leq C\bl,
 \qquad T=\tau_\bu=t_0+S_{\bk,\bl}^{\bu}.
\]
The initial barrier gives \(x\geq c_a\bl^{1/8}-C_L>0\), and
\(T\leq t_0\) is a bounded \(\mathcal G_t\)-stopping time
since \(S_{\bk,\bl}^{\bu}\leq 0\) by \eqref{eq:S-sign}.  Choose
\(w_\bl=C_*\sqrt{\log(\bl+2)}\), with \(C_*\) large.
The bound \eqref{eq:dyadic-bracket-step}, the reflection
principle, and a union bound give
\begin{equation}\label{eq:off-buffer-oscillation-cost}
 \mathbb Q_{\bk,\bl}^{\bu}\left(
  (\mathscr O_{\bl,\bu})^\complement
\right)\leq C\bl^{-30}.
\end{equation}
On \(\mathscr O_{\bl,\bu}\), if the barrier holds at all the dyadic times, then \(Y_s=x-B_s>b_*(s)\) for all \(s\leq T\), where
\[
 b_*(s)=-\sqrt2a\bigl(\lg_{r_\bk}(2^\bl)+s\bigr)^{1/8}
          -\sqrt2L-Cw_\bl .
\]
This boundary is increasing and positive for all large \(\bl\).  If
\(\tau_*=\inf\{s:Y_s=b_*(s)\}\), optional stopping at the bounded
stopping time \(T\wedge\tau_*\), together with \(b_*>0\), gives
\begin{equation}\label{eq:off-buffer-stopped-endpoint}
 \mathbb Q_{\bk,\bl}^{\bu}\left(
  Y_T\one\{\tau_*>T\}\right)
 =x-\mathbb Q_{\bk,\bl}^{\bu}\left(
  b_*(\tau_*)\one\{\tau_*\leq T\}\right)
 \leq x.
\end{equation}
On \(\widetilde{\mathscr E}_{\bk,\bl}^{\bu,L,a,[\bn]}\) we have
\(\Delta_{\bk,\bu}^{[\bn]}>0\).  Since
\(\Delta_{\bk,\bu}^{[\bn]}=Y_T-\sqrt2S_{\bk,\bl}^{\bu}\) and
\(S_{\bk,\bl}^{\bu}\leq 0\), while \(Y_T>0\) on \(\{\tau_*>T\}\) and
\(\sup_{s\leq 0}e^{s}(1-\sqrt2s)<\infty\),
\begin{equation}\label{eq:off-buffer-clock-penalty}
 e^{S_{\bk,\bl}^{\bu}}(1+\Delta_{\bk,\bu}^{[\bn]})
 =e^{S_{\bk,\bl}^{\bu}}Y_T+e^{S_{\bk,\bl}^{\bu}}(1-\sqrt2S_{\bk,\bl}^{\bu})
 \leq Y_T+C .
\end{equation}
By \eqref{eq:off-buffer-stopped-endpoint}, the contribution of
\(\mathscr O_{\bl,\bu}\) to the left side of
\eqref{eq:buffered-one-block-off-buffer} is therefore at most
\(x+C\leq(1+C\bl^{-1/8})x\), where the last step uses
\(x\geq c_a\bl^{1/8}-C_L\).  Finally, the Cauchy--Schwarz inequality,
\eqref{eq:off-buffer-oscillation-cost}, \(e^{S_{\bk,\bl}^{\bu}}\leq C\),
and the second-moment estimate
\eqref{eq:buffered-unconditional-field-L2} bound the contribution of
\((\mathscr O_{\bl,\bu})^\complement\) by
\(C\bl^{-15}(1+x+\bl)\leq C\bl^{-14}x\).  The truncation indicator can
only decrease this positive expectation.  This proves
\eqref{eq:buffered-one-block-off-buffer}.
\end{proof}

\subsection{Evaluation of the Brownian barrier}
\label{sec:full-barrier}

We now turn to the full-time version of the one-block estimate, the
main result of this section.  The barrier is assembled from the block
barriers \eqref{eq:clock-adapted-barrier} along a partition of the
time axis, and the one-block estimates are iterated along it.

Let \(p\) and \(\lambda_\pm\) be as in \eqref{eq:step2-ratio-choice}.
Put \(\bq_\bk=\lfloor\bk/\mathfrak c^2\rfloor\).  For fixed \(\bn\) and
large \(\bk\), choose a \textbf{partition}
\begin{equation}\label{eq:step2-geometric-grid}
\mathcal P_{\bn,\bk}
=\{\bn=\bn_0<\bn_1<\cdots<\bn_{\bar r}=\bq_\bk\},
\qquad
\lambda_-\leq\frac{\bn_{r+1}}{\bn_r}\leq\lambda_+.
\end{equation}
\begin{minipage}{\textwidth}
The \textbf{Brownian barrier} along the partition is defined by
\begin{equation}\label{eq:brownian-clock-barrier}
\widetilde{\mathscr E}_{\bk,\bn}^{\mathcal P,L,a}
:=
\left(\bigcap_{r=0}^{\bar r-1}
\widetilde{\mathscr E}_{\bk,\bn_r}^{\bn_{r+1},L,a,[\bn]}\right)
\cap
\widetilde{\mathscr E}_{\bk,\bq_\bk}^{\infty,L,a,[\bn]}.
\end{equation}
At a partition endpoint the adjacent local conditions are both retained. 
Throughout, we refer to the intervals between adjacent  $\bn_i$'s as \textbf{blocks}.

\medskip
\begin{center}
\begin{tikzpicture}[x=1.02cm,y=0.67cm,>=stealth,font=\scriptsize]
  \draw[->,gray!70] (0.35,0.35)--(9.25,0.35)
    node[below right,text=black] {logarithmic time \(m\)};
  \draw[->,gray!70] (0.55,0.20)--(0.55,4.95)
    node[above,text=black] {\(X_m\)};

  \foreach \x/\lab in {
    0.85/{\bn_0=\bn},
    3.05/{\bn_1},
    5.25/{\bn_2},
    7.45/{\bn_3=\bq_\bk}}
  {
    \draw[densely dashed,gray!55] (\x,0.35)--(\x,4.65);
    \node[below,text=black] at (\x,0.33) {\(\lab\)};
  }

  \draw[very thick,tabviolet]
    plot[smooth] coordinates {(0.85,2.55) (1.45,2.67) (2.15,2.72) (3.05,2.92)};
  \draw[very thick,tabviolet]
    plot[smooth] coordinates {(3.05,3.30) (3.65,3.38) (4.45,3.45) (5.25,3.60)};
  \draw[very thick,tabviolet]
    plot[smooth] coordinates {(5.25,3.94) (5.90,4.00) (6.65,4.08) (7.45,4.18)};
  \draw[very thick,tabviolet]
    plot[smooth] coordinates {(7.45,4.47) (8.12,4.53) (8.72,4.58)};
  \draw[->,very thick,tabviolet] (8.70,4.58)--(9.10,4.62);

  \fill[tabviolet] (3.05,2.92) circle (1.5pt);
  \draw[tabviolet,fill=white,line width=.6pt] (3.05,3.30) circle (1.8pt);
  \draw[->,tabviolet!80] (3.05,3.00)--(3.05,3.22);
  \node[tabviolet,anchor=west] at (3.10,3.11) {restart};
  \fill[tabviolet] (5.25,3.60) circle (1.5pt);
  \draw[tabviolet,fill=white,line width=.6pt] (5.25,3.94) circle (1.8pt);
  \draw[->,tabviolet!80] (5.25,3.68)--(5.25,3.86);
  \fill[tabviolet] (7.45,4.18) circle (1.5pt);
  \draw[tabviolet,fill=white,line width=.6pt] (7.45,4.47) circle (1.8pt);
  \draw[->,tabviolet!80] (7.45,4.26)--(7.45,4.39);

  \draw[very thick,black]
    plot[smooth] coordinates {
      (0.85,1.55) (1.10,1.83) (1.38,1.48) (1.70,2.05)
      (2.02,1.77) (2.33,2.28) (2.62,1.96) (3.05,2.32)
      (3.35,2.02) (3.67,2.72) (3.96,2.43) (4.28,2.93)
      (4.59,2.58) (4.90,3.08) (5.25,2.82)
      (5.58,3.18) (5.86,2.87) (6.17,3.48) (6.49,3.15)
      (6.80,3.62) (7.10,3.31) (7.45,3.72)
      (7.75,3.45) (8.02,3.92) (8.33,3.64) (8.72,4.02)};

  \node[black,fill=white,inner sep=1pt,anchor=west] at (0.92,1.02)
    {\(X_{\bm}=\sqrt{\beta/8}\,\zeta_{2^{\bm}}(\bn;r_\bk,\theta)\)};
  \node[tabviolet,anchor=south west] at (0.98,2.72)
    {\(b_0(\bm)\)};
  \node[tabviolet,anchor=south west] at (3.55,3.46)
    {\(b_1(\bm)\)};
  \node[tabviolet,anchor=south west] at (5.82,4.06)
    {\(b_2(\bm)\)};
  \node[tabviolet,anchor=west] at (7.70,4.86)
    {terminal block};
\end{tikzpicture}

\parbox{0.93\textwidth}{\small
Here
\[
 b_r(\bm):=q_{\bk,\bn_r}(2^{\bm})
      +a\,q_{\bk,\bn_r}(2^{\bm})^{1/8}+L,
 \qquad \bn_r\leq \bm\leq \bn_{r+1}.
\]
At \(\bn_{r+1}\), the path $X$ is tested against
both the incoming value \(b_r(\bn_{r+1})\) (the filled dot) and the newly
anchored value \(b_{r+1}(\bn_{r+1})\) (the open dot).  At that point
\(q_{\bk,\bn_{r+1}}(2^{\bn_{r+1}})
=\lg_{r_\bk}(2^{\bn_{r+1}})\), and the open-dot barrier then evolves in
Brownian time on the next block. The barrier jumps up at each block end because the incoming clock is \(q_{\bk,\bn_r}(2^{\bn_{r+1}})=\lg_{r_\bk}(2^{\bn_{r+1}})+S_{\bk,\bn_r}^{\bn_{r+1}}\), while \(S_{\bk,\bn_r}^{\bn_{r+1}}\leq 0\) by \eqref{eq:S-sign} and \(u\mapsto u+au^{1/8}\) is increasing for large \(u\).}
\end{center}
\end{minipage}

The quantity evaluated in this subsection is the expectation of the
surviving endpoint under the tilted measure,
\begin{equation}\label{eq:brownian-clock-full-operator}
 \widetilde{\mathcal U}_{\bk,\bn}^{\mathcal P,L,a}
 :=\mathbb Q_{\bk,\bn}^{\infty}\Big(
 \one\{\operatorname{Tru}_{2^\bn}\}
 \one\{\widetilde{\mathscr E}_{\bk,\bn}^{\mathcal P,L,a}\}
 e^{S_{\bk,\bn}^{\infty}}
 (\Delta_{\bk,\infty}^{[\bn]})_+\Big),
\end{equation}
weighted by the clock deviation \(S^{\infty}_{\bk,\bn}\) of \eqref{eq:StLt}.
Let \(\bt_\bk\) be the scale of
\eqref{eq:step3-near-saturation-scales} and \(c_\bk\) the constant of
Lemma~\ref{lem:step3-late-forgetting}.

Throughout this subsection, write \(\Delta_\bm\) for
\(\Delta_{\bk,\bm}^{[\bn]}\), and \(\mathscr V_\bm\), \(\mathscr W_\bm\)
for the events \eqref{eq:buffer-event} and
\eqref{eq:step2-radial-window} at \(\bm\); write \(\bq=\bq_\bk\),
\(\bt=\bt_\bk\); by \eqref{eq:step3-near-saturation-scales},
\(\bk-\bt\asymp(\log\bk)^\alpha\).
All \(o_\bk(1)\) terms are deterministic.

The constant \(c_\bk\) enters through the following lemma, which
evaluates what is left of \eqref{eq:brownian-clock-full-operator}
after time \(2^{\bt}\).

\begin{lemma}[The tail after \(\bt_\bk\)]
\label{lem:bcb-tail}
For \(\bn\leq\bu\leq\bk\), let
\begin{equation}\label{eq:bcb-Hprime}
 H'_\bu:=\mathbb Q_{\bk,\bu}^{\infty}\Big(
  \one\{\operatorname{Tru}_{2^\bu}\}
  \one\{\widetilde{\mathscr E}_{\bk,\bq}^{\infty,L,a,[\bn]}\}
  e^{S_{\bk,\bu}^{\infty}}
  (\Delta_\infty)_+\Big).
\end{equation}
Then
\begin{equation}\label{eq:bcb-tail-moments}
\begin{gathered}
 H'_\bu\leq C\left(1+|\Delta_\bu|+\sqrt{\bk-\bu}\right),
 \\
 \mathbb Q_{\bk,\bu}^{\infty}\Big(
  \one\{\operatorname{Tru}_{2^\bu}\}
  e^{2S_{\bk,\bu}^{\infty}}|\Delta_\infty-\Delta_\bu|^2\Big)^{1/2}
 \leq C\left(1+\sqrt{\bk-\bu}\right).
\end{gathered}
\end{equation}
Moreover, on
\(\mathscr V_\bt\cap\mathscr W_\bt
\cap\widetilde{\mathscr E}_{\bk,\bq}^{\bt,L,a,[\bn]}
\cap\{|S_{\bk,\bq}^{\bt}|\leq1\}\),
\begin{equation}\label{eq:bcb-tail-claim}
 H'_\bt=\bigl(c_\bk+o_\bk(1)\bigr)\Delta_\bt .
\end{equation}
\end{lemma}

\begin{proof}
Let \(Y_s\), \(s\geq2^\bu\), be the martingale
\eqref{eq:clock-gap-Brownian} anchored at \(\bu\), so that
\(Y_{2^\bm}=\Delta_\bm+\sqrt2S_{\bk,\bu}^{\bm}\) for every
\(\bu\leq\bm\leq\infty\), by the definitions \eqref{eq:brownian-clock}
and \eqref{eq:StLt}; in particular \(Y_{2^\bu}=\Delta_\bu\).  By
\eqref{eq:clock-gap-Brownian} with \((\bl,\bu)\) replaced by
\((\bu,\infty)\), \(Y_s-Y_{2^\bu}\) is a continuous martingale under
\(\mathbb Q_{\bk,\bu}^{\infty}\) with bracket
\(\langle Y\rangle_s=q_{\bk,\bu}(s)-q_{\bk,\bu}(2^\bu)\), and
\begin{equation}\label{eq:bcb-tail-bracket}
 \langle Y\rangle_\infty=(\bk-\bu)\log2+S_{\bk,\bu}^{\infty},
 \qquad
 \Delta_\infty-\Delta_\bu=(Y_\infty-Y_{2^\bu})-\sqrt2S_{\bk,\bu}^{\infty}.
\end{equation}
By Lemma~\ref{lem:QVsmallkold}, for every \(p<\infty\),
\begin{equation}\label{eq:bcb-tail-overrun}
 \mathbb Q_{\bk,\bu}^{\infty}\left(\one\{\operatorname{Tru}_{2^\bu}\}
  e^{p(S_{\bk,\bu}^{\infty})_+}\right)\leq C_p .
\end{equation}
The Burkholder--Davis--Gundy inequality applied to \(Y_\cdot-Y_{2^\bu}\),
together with \eqref{eq:bcb-tail-bracket},
\eqref{eq:bcb-tail-overrun}, and Corollary~\ref{cor:fieldincrement},
gives the second bound in \eqref{eq:bcb-tail-moments}; the first
then follows from
\((\Delta_\infty)_+\leq1+|\Delta_\bu|+|\Delta_\infty-\Delta_\bu|\),
the Cauchy--Schwarz inequality and \eqref{eq:bcb-tail-overrun}.

Now let \(\bu=\bt\), and assume the event in the last assertion.  Let
\(F\) be the event that the barrier fails after \(\bt\), that is, that
\(\Delta_\bm+\sqrt2S_{\bk,\bq}^{\bm}<g(q_{\bk,\bq}(2^\bm))\) for some
\(\bm>\bt\), where \(g(u):=-\sqrt2a\,u^{1/8}-\sqrt2L\); since
\(S_{\bk,\bq}^{\bm}=S_{\bk,\bq}^{\bt}+S_{\bk,\bt}^{\bm}\), the left side
is \(Y_{2^\bm}+\sqrt2S_{\bk,\bq}^{\bt}\).  Since the barrier holds up to
\(\bt\), the barrier indicator in \eqref{eq:bcb-Hprime} is
\(\one\{F^\complement\}\), and so
\begin{align}
 H'_\bt
 &=\Delta_\bt\,
  \mathbb Q_{\bk,\bt}^{\infty}\left(\one\{\operatorname{Tru}_{2^\bt}\}e^{S_{\bk,\bt}^{\infty}}\right)
 +\mathbb Q_{\bk,\bt}^{\infty}\left(\one\{\operatorname{Tru}_{2^\bt}\}e^{S_{\bk,\bt}^{\infty}}
   \bigl((\Delta_\infty)_+-\Delta_\bt\bigr)\right)\notag\\
 &\qquad-\mathbb Q_{\bk,\bt}^{\infty}\left(\one\{\operatorname{Tru}_{2^\bt}\}\one\{F\}
   e^{S_{\bk,\bt}^{\infty}}(\Delta_\infty)_+\right).
 \label{eq:bcb-tail-decomp}
\end{align}
We bound the three terms in turn.  By
Lemma~\ref{lem:step3-late-forgetting} on \(\mathscr W_\bt\), the
first term is \((c_\bk+o_\bk(1))\Delta_\bt\).  Since \(\Delta_\bt>0\),
\(|(\Delta_\infty)_+-\Delta_\bt|\leq|\Delta_\infty-\Delta_\bt|\), so
the Cauchy--Schwarz inequality and \eqref{eq:bcb-tail-moments}
bound the second term by \(C(1+\sqrt{\bk-\bt})\), which is
\(o_\bk(1)\Delta_\bt\) because \(\Delta_\bt\geq\bt^{1/2-\nu}\) on
\(\mathscr V_\bt\).

For the third term we apply part~(d) of Lemma~\ref{lem:buffered-brownian},
and we adopt its notation.  Let \(M_s:=Y_s-Y_{2^\bt}\) for
\(s\geq2^\bt\), a continuous martingale with bracket \(\langle M\rangle_s=\langle Y\rangle_s\), and put
\[
 T_0:=(\bk-\bt)\log2,\qquad
 \overline T:=\langle M\rangle_\infty=T_0+S_{\bk,\bt}^{\infty},\qquad
 \mathscr T:=\operatorname{Tru}_{2^\bt},\qquad
 W:=e^{S_{\bk,\bt}^{\infty}}(\Delta_\infty)_+ ,
\]
so that \(\overline T-T_0=S_{\bk,\bt}^{\infty}\).  The two inputs
\eqref{eq:weighted-brownian-clock-inputs} of that lemma are
\eqref{eq:bcb-tail-overrun} and
\begin{equation}\label{eq:bcb-tail-K}
 \mathbb Q_{\bk,\bt}^{\infty}\Big(\one\{\mathscr T\}W^2\Big)^{1/2}
 \leq C\Delta_\bt=:K,
\end{equation}
which follows from \eqref{eq:bcb-tail-moments} and
\(\sqrt{\bk-\bt}\leq\Delta_\bt\).

We now locate the barrier after \(\bt\).  For \(\bm>\bt\), by
\eqref{eq:brownian-clock} and \eqref{eq:StLt},
\begin{equation}\label{eq:bcb-tail-clock-after-t}
 q_{\bk,\bq}(2^\bm)
 =\lg_{r_\bk}(2^\bm)+S_{\bk,\bq}^{\bt}+S_{\bk,\bt}^{\bm}
 \leq\bk\log2+2+(S_{\bk,\bt}^{\infty})_+ .
\end{equation}
Here \(\lg_{r_\bk}(2^\bm)\leq\bk\log2\) by definition,
\(|S_{\bk,\bq}^{\bt}|\leq1\) by assumption, and
\(S_{\bk,\bt}^{\bm}\leq(S_{\bk,\bt}^{\infty})_+\) by the monotonicity
in \eqref{eq:S-sign} and after it.  On
\(\{\overline T\leq2T_0\}=\{S_{\bk,\bt}^{\infty}\leq T_0\}\), the
right side of \eqref{eq:bcb-tail-clock-after-t} is at most
\(2\bk\log2\) for large \(\bk\), and since \(g\) is increasing, the
barrier after \(\bt\), as seen by \(Y\), and the starting point of
\(Y\) satisfy
\begin{equation}\label{eq:bcb-tail-barrier-height}
 g(q_{\bk,\bq}(2^\bm))-\sqrt2S_{\bk,\bq}^{\bt}
 \leq g(2\bk\log2)+2\leq C\bk^{1/8},
 \qquad
 Y_{2^\bt}=\Delta_\bt\geq2C\bk^{1/8},
\end{equation}
the last inequality because \(\Delta_\bt\geq\bt^{1/2-\nu}\) on
\(\mathscr V_\bt\) and \(\nu<3/8\).  Thus, on
\(\{\overline T\leq2T_0\}\), the barrier can fail after \(\bt\) only if
\(Y\) drops by \(\Delta_\bt/2\) from its starting point:
\begin{equation}\label{eq:bcb-tail-hitting}
 F\cap\{\overline T\leq2T_0\}
 \subseteq
 \Big\{\exists s\geq2^\bt:\ \tfrac12\Delta_\bt+M_s\leq0\Big\}
 \cap\{\overline T\leq2T_0\}.
\end{equation}
The right side is the event in \eqref{eq:weighted-brownian-clock-DDS}
with \(h\equiv0\) and \(y=\Delta_\bt/2\).  That estimate gives
\begin{equation}\label{eq:bcb-tail-small-overrun}
 \mathbb Q_{\bk,\bt}^{\infty}\Big(
  \one\{\mathscr T\}\one\{F\}\one\{\overline T\leq2T_0\}W\Big)
 \leq CKe^{-c\Delta_\bt^2/(4T_0)}
 \leq C\Delta_\bt e^{-c\Delta_\bt^2/(\bk-\bt)}.
\end{equation}
On the complementary event, \eqref{eq:bcb-tail-overrun} with
\(p=1\) and Markov's inequality give
\(\mathbb Q_{\bk,\bt}^{\infty}(\mathscr T\cap\{\overline T>2T_0\})\leq C_12^{-(\bk-\bt)}\),
so the Cauchy--Schwarz inequality and the bound on \(W\) give
\begin{equation}\label{eq:bcb-tail-large-overrun}
 \mathbb Q_{\bk,\bt}^{\infty}\Big(
  \one\{\mathscr T\}\one\{F\}\one\{\overline T>2T_0\}W\Big)
 \leq C_1^{1/2}2^{-(\bk-\bt)/2}K
 \leq C\Delta_\bt 2^{-(\bk-\bt)/2}.
\end{equation}
The third term in \eqref{eq:bcb-tail-decomp} is the sum of the
left sides of \eqref{eq:bcb-tail-small-overrun} and
\eqref{eq:bcb-tail-large-overrun}.  Both right sides are
\(o_\bk(1)\Delta_\bt\), because
\(\Delta_\bt^2/(\bk-\bt)\geq c\,\bt^{1-2\nu}/(\log(\bk+2))^\alpha\to\infty\)
and \(\bk-\bt\to\infty\).  Combining the three bounds proves
\eqref{eq:bcb-tail-claim}.
\end{proof}

\begin{proposition}[Evaluation of the Brownian barrier]
\label{prop:brownian-clock-full-tower}
Fix \(a<0\) and \(L<\infty\), and let \(\kappa\) be as in
\eqref{eq:kappa-rho}.  On the initial barrier
\(\widetilde{\mathscr E}_{\bk,\bn}^{\bn,L,a,[\bn]}
=\{\Delta_{\bk,\bn}^{[\bn]}\geq-\sqrt2a\,\lg_{r_\bk}(2^\bn)^{1/8}-\sqrt2L\}\)
one has
\begin{equation}\label{eq:brownian-clock-global-upper}
 0\leq\widetilde{\mathcal U}_{\bk,\bn}^{\mathcal P,L,a}
 \leq C_{a,L}\Delta_{\bk,\bn}^{[\bn]}.
\end{equation}
If in addition
\(R_{2^\bn}\leq C2^{C'\bn}\) and the buffer event
\(\mathscr V_{\bk,\bn}^{\nu,[\bn]}\) hold, then
\begin{equation}\label{eq:brownian-clock-full-asymptotic}
 \left|
 \widetilde{\mathcal U}_{\bk,\bn}^{\mathcal P,L,a}
 -c_\bk\Delta_{\bk,\bn}^{[\bn]}
 \right|
 \leq\left(C\bn^{-\kappa}+o_\bk(1)\right)
 \Delta_{\bk,\bn}^{[\bn]}.
\end{equation}
\end{proposition}

\begin{proof}
\emph{Reduction to blocks.}
Put \(\bn_{\bar r+1}:=\bt_\bk\), so that the last block \([\bq_\bk,\bt_\bk]\)
is \([\bn_{\bar r},\bn_{\bar r+1}]\).  For \(0\leq r\leq\bar r\) and \(\filt_{2^{\bn_{r+1}}}\)-measurable \(X\),
nonnegative or integrable, put
\begin{equation}\label{eq:bcb-T}
 \mathcal T_r[X]
 :=\mathbb Q_{\bk,\bn_r}^{\bn_{r+1}}\Big(
  \one\{\operatorname{Tru}_{2^{\bn_r},2^{\bn_{r+1}}}\}
  \one\{\widetilde{\mathscr E}_{\bk,\bn_r}^{\bn_{r+1},L,a,[\bn]}\}
  e^{S_{\bk,\bn_r}^{\bn_{r+1}}}X\Big),
\end{equation}
the expectation over the block \([\bn_r,\bn_{r+1}]\), weighted by
\(e^{S}\) and restricted to the truncation event and the block barrier.
For \(0\leq r\leq \bar r\), let
\begin{equation}\label{eq:bcb-H}
 H_r:=\mathbb Q_{\bk,\bn_r}^{\infty}\Big(
  \one\{\operatorname{Tru}_{2^{\bn_r}}\}
  \one\Big\{\bigcap_{j=r}^{\bar r-1}
   \widetilde{\mathscr E}_{\bk,\bn_j}^{\bn_{j+1},L,a,[\bn]}
   \cap\widetilde{\mathscr E}_{\bk,\bq_\bk}^{\infty,L,a,[\bn]}\Big\}
  e^{S_{\bk,\bn_r}^{\infty}}
  (\Delta_\infty)_+\Big),
\end{equation}
which is \(\widetilde{\mathcal U}\) with the barrier imposed only from
\(\bn_r\) onward.  Then \(H_0=\widetilde{\mathcal U}_{\bk,\bn}^{\mathcal P,L,a}\),
and since the barrier \eqref{eq:brownian-clock-barrier} is an
intersection over the blocks, the tower property gives
\(H_r=\mathcal T_r[H_{r+1}]\) for \(r<\bar r\).

\emph{The clock on the last block.}
Let \(\mathscr C\) be
the event that \(|S_{\bk,\bq}^{\bm}|\leq\bq^{-10}\) for every
dyadic exponent \(\bm\) with \(\bq\leq\bm\leq\bt\).  We claim that, on
\(\mathscr W_\bq\),
\begin{equation}\label{eq:bcb-clock-event}
 \mathbb Q_{\bk,\bq}^{\bt}\left(
  \mathscr C^\complement\cap\operatorname{Tru}_{2^\bq,2^\bt}\right)
 \leq C\bq^{-40},
\end{equation}
so that \(\mathscr C\), which is \(\filt_{2^\bt}\)-measurable, satisfies
the hypothesis \eqref{eq:one-block-clock-input} of case~(ii) of
Lemma~\ref{lem:buffered-one-block} on the block \([\bq,\bt]\).  On the
grid blocks, where \(\mathfrak c\bn_{r+1}\leq\bk\), case~(i) of that
lemma needs no such input.  Since \(\bt\leq\bk\), \eqref{eq:S-sign}
gives \(\max_{\bq\leq\bm\leq\bt}|S_{\bk,\bq}^{\bm}|=|S_{\bk,\bq}^{\bt}|\),
so the whole block is controlled by the right endpoint, and the right
endpoint has the first-moment bound \eqref{eq-11012026enew} of
Lemma~\ref{lem:QVsmallk}:
\begin{equation}\label{eq:bcb-deviation}
 \mathbb Q_{\bk,\bq}^{\bt}\left(
  \bigl(|S_{\bk,\bq}^{\bt}|-2^{-\bq}\bigr)_+
  \one\{\operatorname{Tru}_{2^\bq,2^\bt}\}\right)
 \leq C2^{-\bk}\left(R_{2^\bq}\,\bt+2^\bt\right).
\end{equation}
On \(\mathscr W_\bq\), that is when \(R_{2^\bq}\leq C_B2^\bq\bq^B\), the
right side is at most \(C2^{\bt-\bk}\bt^{B+1}\), and by
\eqref{eq:step3-near-saturation-scales},
\(2^{\bt-\bk}\leq2^{1-(\log(\bk+2))^\alpha}\) with \(\alpha>1\), so it
is \(o(\bk^{-M})\) for every \(M\).  Markov's inequality applied to
\((|S_{\bk,\bq}^{\bt}|-2^{-\bq})_+\) at level \(\bq^{-10}/2\) proves
\eqref{eq:bcb-clock-event}.

\emph{One block.}
The induction below uses four bounds on \(\mathcal T_r\), all taken
from Section~\ref{sec:one-block}.  For \(0\leq r<\bar r\), case~(i) of
Lemma~\ref{lem:buffered-one-block} applies to the block
\([\bn_r,\bn_{r+1}]\) as soon as \(R_{2^{\bn_r}}\leq C2^{C'\bn_r}\),
which is the hypothesis of the proposition when \(r=0\) and follows
from \(\mathscr W_{\bn_r}\) when \(r\geq1\), since \(C'>1\).  Its
conclusions (a), (b) and (c), in the form recorded after
\eqref{eq:kappa-rho}, say that on \(\mathscr V_{\bn_r}\), which
is contained in the initial barrier by \eqref{eq:buffer-implies-barrier},
\begin{equation}\label{eq:bcb-block}
\begin{gathered}
 \left|\mathcal T_r[\Delta_{\bn_{r+1}}]-\Delta_{\bn_r}\right|
 \leq C\bn_r^{-\kappa}\Delta_{\bn_r},
 \qquad
 \mathcal T_r[1+|\Delta_{\bn_{r+1}}|]
 \leq(1+C\bn_r^{-\kappa})\Delta_{\bn_r},\\
 \mathcal T_r\left[(1+|\Delta_{\bn_{r+1}}|)
  \one\{(\mathscr V_{\bn_{r+1}}\cap\mathscr W_{\bn_{r+1}})^\complement\}\right]
 \leq\rho_{\bn_r}\Delta_{\bn_r},
\end{gathered}
\end{equation}
with \(\rho_{\bn_r}\leq C\bn_r^{-\kappa}\).
Lemma~\ref{lem:one-sided-one-block} needs no event at all: for
\(0\leq r\leq\bar r\), on the initial barrier
\(\widetilde{\mathscr E}_{\bk,\bn_r}^{\bn_r,L,a,[\bn]}\) alone,
\begin{equation}\label{eq:bcb-one-sided}
 \mathcal T_r[1+|\Delta_{\bn_{r+1}}|]
 \leq(1+C\bn_r^{-1/8})\Delta_{\bn_r}.
\end{equation}
On the last block, \(r=\bar r\), case~(ii) of
Lemma~\ref{lem:buffered-one-block} applies on \(\mathscr W_\bq\), with
\(\mathscr C\) as the event in \eqref{eq:one-block-clock-input}, by
\eqref{eq:bcb-clock-event}, and with the weight
\(r_\bq^*=\sqrt{\bk-\bt}\leq\bq^{1/16}\).
Its conclusions are that the first two bounds of
\eqref{eq:bcb-block} hold with \(\one\{\mathscr C\}\) inserted
inside \(\mathcal T_{\bar r}\), and that the third is replaced by
\begin{equation}\label{eq:bcb-last-retention}
 \mathcal T_{\bar r}\left[(1+|\Delta_\bt|+\sqrt{\bk-\bt})
  \one\{(\mathscr C\cap\mathscr V_\bt\cap\mathscr W_\bt)^\complement\}\right]
 \leq\rho_\bq\Delta_\bq,
\end{equation}
both on \(\mathscr V_\bq\cap\mathscr W_\bq\).

\emph{Induction.}
Let \(\mathscr K_r:=\mathscr V_{\bn_r}\cap\mathscr W_{\bn_r}\) for
\(1\leq r\leq\bar r\), and let \(\mathscr K_0\) be the event
\(\{R_{2^\bn}\leq C2^{C'\bn}\}\cap\mathscr V_\bn\) of the hypothesis.
We prove by backward induction on \(r\) that, for
\(0\leq r\leq\bar r\),
\begin{equation}\label{eq:bcb-induction}
 H_r\leq C_r\Delta_{\bn_r}
 \ \text{ on the initial barrier at } \bn_r,
 \qquad
 |H_r-c_\bk\Delta_{\bn_r}|\leq E_r\Delta_{\bn_r}
 \ \text{ on } \mathscr K_r,
\end{equation}
with constants \(C_r\) bounded uniformly in \(r\), and deterministic
\(E_r\) satisfying \(E_{\bar r}=o_\bk(1)\) and
\begin{equation}\label{eq:bcb-recursion}
 E_r\leq(1+C\bn_r^{-\kappa})E_{r+1}+C\bn_r^{-\kappa},
 \qquad 0\leq r<\bar r.
\end{equation}
Since \(H_r\) vanishes off the initial barrier, because its
barrier includes the constraint at time \(2^{\bn_r}\), which is
\(\filt_{2^{\bn_r}}\)-measurable, the first bound in
\eqref{eq:bcb-induction} gives \(H_r\leq C_r(1+|\Delta_{\bn_r}|)\)
everywhere.  At \(r=0\), the two bounds in
\eqref{eq:bcb-induction} are the two statements of the
proposition, because \(\bn_r\geq\lambda_-^r\bn\) makes
\(\prod_r(1+C\bn_r^{-\kappa})\) bounded and
\(\sum_r\bn_r^{-\kappa}\leq C\bn^{-\kappa}\), so that iterating
\eqref{eq:bcb-recursion} gives \(E_0\leq C\bn^{-\kappa}+o_\bk(1)\).

\emph{Base case, \(r=\bar r\).}
Let \(H'_\bt\) and \(H'_\bk\) be as in Lemma~\ref{lem:bcb-tail}, so
that \(H_{\bar r}=\mathcal T_{\bar r}[H'_\bt]\) by the tower property.
For the first bound in \eqref{eq:bcb-induction}, use the tower
property at \(\bk\) instead: \(H'_\bk\leq C(1+|\Delta_\bk|)\) by
\eqref{eq:bcb-tail-moments}, so
Lemma~\ref{lem:one-sided-one-block} on the block \([\bq,\bk]\) gives
\(H_{\bar r}\leq C_{\bar r}\Delta_\bq\) on the initial barrier at \(\bq\).
For the second, split \(H_{\bar r}=\mathcal T_{\bar r}[H'_\bt]\)
according to \(\mathscr C\cap\mathscr V_\bt\cap\mathscr W_\bt\).
Inside \(\mathcal T_{\bar r}\) the barrier up to \(\bt\) is in force,
and \(\mathscr C\subseteq\{|S_{\bk,\bq}^{\bt}|\leq1\}\), so
\eqref{eq:bcb-tail-claim} applies on this event, while
\eqref{eq:bcb-tail-moments} applies off it:
\begin{equation}\label{eq:bcb-base-split}
 H_{\bar r}
 =(c_\bk+o_\bk(1))\,
  \mathcal T_{\bar r}\bigl[\one\{\mathscr C\cap\mathscr V_\bt\cap\mathscr W_\bt\}\Delta_\bt\bigr]
 +O(1)\,\mathcal T_{\bar r}\bigl[(1+|\Delta_\bt|+\sqrt{\bk-\bt})
  \one\{(\mathscr C\cap\mathscr V_\bt\cap\mathscr W_\bt)^\complement\}\bigr].
\end{equation}
The second term is \(o_\bk(1)\Delta_\bq\) by
\eqref{eq:bcb-last-retention}, and the same bound lets us drop
\(\mathscr V_\bt\cap\mathscr W_\bt\) from the first, after which the
first bound of \eqref{eq:bcb-block} with \(\one\{\mathscr C\}\)
inserted gives \((1+o_\bk(1))\Delta_\bq\).  Hence
\(H_{\bar r}=(c_\bk+o_\bk(1))\Delta_\bq\) on
\(\mathscr V_\bq\cap\mathscr W_\bq\), that is, \(E_{\bar r}=o_\bk(1)\).

\emph{Inductive step.}
Suppose \eqref{eq:bcb-induction} holds at \(r+1\).  Since
\(H_r=\mathcal T_r[H_{r+1}]\) and \(H_{r+1}\leq C_{r+1}(1+|\Delta_{\bn_{r+1}}|)\),
the bound \eqref{eq:bcb-one-sided} gives \(H_r\leq C_r\Delta_{\bn_r}\)
on the initial barrier at \(\bn_r\) with
\(C_r:=C_{r+1}(1+C\bn_r^{-1/8})\), and \(\sup_rC_r<\infty\) because
\(\bn_r\geq\lambda_-^r\bn\).  On
\(\mathscr K_r\), by linearity,
\[
 |H_r-c_\bk\Delta_{\bn_r}|
 \leq\mathcal T_r\left[|H_{r+1}-c_\bk\Delta_{\bn_{r+1}}|\right]
  +c_\bk\left|\mathcal T_r[\Delta_{\bn_{r+1}}]-\Delta_{\bn_r}\right| .
\]
The second term is at most \(C\bn_r^{-\kappa}\Delta_{\bn_r}\) by the
first bound in \eqref{eq:bcb-block} and \(c_\bk\leq C\) from
\eqref{eq:step3-ck-bounds}.  Split the first term according to
\(\mathscr K_{r+1}\): on \(\mathscr K_{r+1}\), the induction hypothesis
and the second bound in \eqref{eq:bcb-block} give at most
\(E_{r+1}(1+C\bn_r^{-\kappa})\Delta_{\bn_r}\); off it,
\(|H_{r+1}-c_\bk\Delta_{\bn_{r+1}}|\leq(C_{r+1}+c_\bk)(1+|\Delta_{\bn_{r+1}}|)\) and
the third bound in \eqref{eq:bcb-block} give at most
\(C\bn_r^{-\kappa}\Delta_{\bn_r}\).  This proves
\eqref{eq:bcb-recursion}.
\end{proof}

\section{Comparison of the Brownian and physical barriers}
\label{sec:brownian-physical-comparison}

The barrier of Section~\ref{sec:brownian-clock-approximate-martingale}
was chosen for the Brownian problem: at the fixed radius \(r_\bk\) it
constrains the Gaussian field \(\zeta\) in Brownian time, and at the
boundary \(r=1\), where the clock is deterministic, it is the natural
barrier on the Verblunsky evolution of the random-matrix field.  It is
not the barrier we will ultimately compare with the chaos construction.
The first step toward that comparison, taken in this section, stays at
the radius \(r_\bk\) but replaces the Gaussian field by the field
\(\varphi\) itself and Brownian time by physical time.  For a fixed
radius \(r_\bk\), define the \textbf{physical barrier}
\begin{align}
\mathscr E_{\bk,\bn}^{\bu,L,a}
&:=
\left\{
\sqrt{\frac{\beta}{8}}\,
\varphi_{2^\bm}(r_\bk,\theta)
\leq
\lg_{r_\bk}(2^\bm)
+a\lg_{r_\bk}(2^\bm)^{1/8}+L:
\bn\leq\bm\leq\bu
\right\},\;\label{eq:physical-interface-barrier}
\mathscr E_{\bk,\bn}^{L,a}
:=\bigcap_{\bu\geq\bn}
\mathscr E_{\bk,\bn}^{\bu,L,a}.
\end{align}

Recall the quantity \(\widetilde{\mathcal U}_{\bk,\bn}^{\mathcal P,L,a}\)
of \eqref{eq:brownian-clock-full-operator}, evaluated in
Proposition~\ref{prop:brownian-clock-full-tower}.  We introduce the
\textbf{physical deficits}, the analogues of the deficits
\eqref{eq:buffered-gap-definition} with \(\zeta\) replaced by
\(\varphi\), which measure the distance below the physical barrier
\eqref{eq:physical-interface-barrier}:
\[
 \Delta_{\bk,\bm}^{\varphi}
 :=\sqrt2\lg_{r_\bk}(2^\bm)
 -\sqrt{\frac\beta4}\varphi_{2^\bm}(r_\bk,\theta),
 \qquad
 \Delta_{\bk,\infty}^{\varphi}
 :=\sqrt2\lg_{r_\bk}(\infty)
 -\sqrt{\frac\beta4}\varphi_\infty(r_\bk,\theta).
\]
At the initial time \(2^\bn\),
\(\Delta_{\bk,\bn}^{\varphi}=\Delta_{\bk,\bn}^{[\bn]}\).  The physical
analogue of \(\widetilde{\mathcal U}\) is
\begin{align}\label{eq:physical-global-good-operator}
 \mathcal U_{\bk,\bn}^{\varphi,L,a}
 :=\mathbb Q_{\bk,\bn}^{\infty}\Big(
 &\one\{\mathscr G_\bn\}
 \one\{\mathscr E_{\bk,\bn}^{L,a}\}
 e^{S_{\bk,\bn}^{\infty}+L_{\bk,\bn}^{\infty}}
 (\Delta_{\bk,\infty}^{\varphi})_+
\Big),
\end{align}
where \(L_{\bk,\bn}^{\infty}\) is the linearization error of
\eqref{eq:StLt} and \(\mathscr G_\bn\) is the good event
\eqref{eq:globalgood}.  This is the
quantity that Section~\ref{sec:proof-main-comparison} integrates in
\(\theta\).

The goal of this section is to upgrade
Proposition~\ref{prop:brownian-clock-full-tower} to
\(\mathcal U^{\varphi}\), that is, to show that
\(\mathcal U_{\bk,\bn}^{\varphi,L,a}=(c_\bk+o(1))\Delta_{\bk,\bn}^{[\bn]}\).
We do this by comparison.
Proposition~\ref{prop:brownian-clock-physical-sandwich} sandwiches
\(\mathcal U^{\varphi,L,a}\) between
\(\widetilde{\mathcal U}^{\mathcal P,L-K,a-\eta}\) and
\(\widetilde{\mathcal U}^{\mathcal P,L+K,a+\eta}\), up to an error
\(r_{\bk,\bn}\Delta_{\bk,\bn}^{[\bn]}\) with \(r_{\bk,\bn}\to0\).  The
comparison is made block by block on the partition, where the terms blocks and partitions are similar to
\eqref{eq:step2-geometric-grid}: on a block where
the clock deviation \(S\) is small and the good event holds, the two
barriers differ only by a shift of their parameters, which is
\eqref{eq:block-comparison}. The key to the remaining blocks is
Lemma~\ref{lem:physical-all-good-upper}, which bounds the expectation
under the physical barrier, started from an arbitrary time, by a
multiple of the initial deficit, with no condition on the clock
deviation.  With it, a simple Cauchy--Schwarz bound discharges the
assumption that the clock deviation is small, which fails only with
small probability.  Section~\ref{sec:proof-main-comparison} uses the
lemma as well.

We use the abbreviations of
Section~\ref{sec:brownian-clock-approximate-martingale}:
\(\Delta_\bm=\Delta_{\bk,\bm}^{[\bn]}\), \(\bq=\bq_\bk\), \(\bt=\bt_\bk\),
\(\bn_{\bar r+1}=\bt\), and \(\ell_\bm:=\lg_{r_\bk}(2^\bm)\).  Let
\(\mathscr G_{\bn,\bm}\) be as in \eqref{eq:globalgood-finite}, so that
\(\mathscr G_{\bn,\infty}=\mathscr G_\bn\), and write
\(\mathscr E_{\bk,\bm}^{\infty,L,a}:=\mathscr E_{\bk,\bm}^{L,a}\).

\subsection*{Two consequences of the good event}
Let \(\bn\leq\bl\leq\bv\leq\bu\).  By \eqref{eq:globalgood-finite}, on
\(\mathscr G_{\bn,\bu}\) the event \(\operatorname{Tru}_{2^\bn,2^\bu}\)
holds and
\(|\varphi_{2^\bv}(r_\bk,\theta)-\zeta_{2^\bv}(\bn;r_\bk,\theta)|\leq\epsilon_\bn\);
since
\(\zeta_{2^\bv}(\bl;r_\bk,\theta)-\varphi_{2^\bl}(r_\bk,\theta)
=\zeta_{2^\bv}(\bn;r_\bk,\theta)-\zeta_{2^\bl}(\bn;r_\bk,\theta)\) by
\eqref{eq:Z} and \eqref{eq:zeta-initial-is-field}, the definitions
\eqref{eq:StLt} and \eqref{eq:buffered-gap-definition} give, on
\(\mathscr G_{\bn,\bu}\),
\begin{equation}\label{eq:physical-good-consequences}
 \sup_{\bl\leq\bv\leq\bu}|L_{\bk,\bl}^{\bv}|\leq C\epsilon_\bn,
 \qquad
 \sup_{\bl\leq\bv\leq\bu}
 |\Delta_{\bk,\bv}^{\varphi}-\Delta_{\bk,\bv}^{[\bl]}|
 \leq C\epsilon_\bn,
 \qquad
 \Delta_{\bk,\bl}^{\varphi}=\Delta_{\bk,\bl}^{[\bl]} .
\end{equation}
We take \(\bn_*\) so large that \(C\epsilon_\bn\leq1\) for
\(\bn\geq\bn_*\).  Second, on the initial physical barrier
\(\mathscr E_{\bk,\bl}^{\bl,L,a}
=\{\Delta_{\bk,\bl}^{\varphi}\geq-\sqrt2a\,\ell_\bl^{1/8}-\sqrt2L\}\),
for \(\bl\geq\bn_*\),
\begin{equation}\label{eq:physical-initial-barrier-size}
 \Delta_{\bk,\bl}^{\varphi}\geq1,
 \qquad\text{and}\qquad
 \Delta_{\bk,\bl}^{\varphi}\geq c_a\bl^{1/8}
 \quad\text{if }\bl\leq\bk,
\end{equation}
since \(\ell_\bl=(\bl\wedge\bk)\log2\).

\subsection*{Comparison of the barriers on one block}
The physical barrier \eqref{eq:physical-interface-barrier} and the
Brownian barrier \eqref{eq:clock-adapted-barrier} differ, at a dyadic
time \(2^\bm\) of a block anchored at \(\bl\), in two respects: the
field \(\varphi_{2^\bm}\) replaces \(\zeta_{2^\bm}(\bn)\), and the
deterministic time \(\ell_\bm\) replaces the Brownian time
\(q_{\bk,\bl}(2^\bm)=\ell_\bm+S_{\bk,\bl}^{\bm}\).  The first
discrepancy is at most \(C\epsilon_\bn\) on \(\mathscr G_{\bn,\bu}\).  The
second moves the boundary \(u\mapsto u+au^{1/8}+L\) by at most
\(|S_{\bk,\bl}^{\bm}|(1+C_a\ell_\bm^{-7/8})\).  Both are absorbed by a
small change of the parameters: fix \(\eta>0\) with \(a+\eta<0\), let
\(\eta_0\) be a sufficiently small multiple of \(\eta\), and let
\(K:=2C+2\), with \(C\) the constant of \eqref{eq:physical-good-consequences}.  Then, for \(\bn\leq\bl\leq\bu\leq\infty\), on
\begin{equation}\label{eq:block-clock-event}
 \mathscr G_{\bn,\bu}\cap
 \Big\{|S_{\bk,\bl}^{\bm}|\leq3\eta_0\ell_\bm^{1/8}
 \ \text{ for all }\bl\leq\bm\leq\bu\Big\},
\end{equation}
monotonicity of the boundary in \(L\) and \(a\) gives
\begin{equation}\label{eq:block-comparison}
 \widetilde{\mathscr E}_{\bk,\bl}^{\bu,L-K,a-\eta,[\bn]}
 \subseteq\mathscr E_{\bk,\bl}^{\bu,L,a}
 \subseteq\widetilde{\mathscr E}_{\bk,\bl}^{\bu,L+K,a+\eta,[\bn]}.
\end{equation}
At \(\bm=\bl\) the clock condition is void, since
\(S_{\bk,\bl}^{\bl}=0\), so \eqref{eq:block-comparison} includes the
comparison of the initial barriers.

\begin{lemma}[Upper bound under the physical barrier]
\label{lem:physical-all-good-upper}
Fix \(a<0\) and \(L<\infty\).  There is \(\bn_*<\infty\) such that, for
every \(\bn_*\leq\bn\leq\bm\) and every \(\bu\in\overline{\mathbb N}\)
with \(\bu\geq\bm\), on \(\mathscr G_{\bn,\bm}\cap\mathscr E_{\bk,\bm}^{\bm,L,a}\),
\begin{align}
 &\mathbb Q_{\bk,\bm}^{\bu}\left(
  \one\{\mathscr G_{\bn,\bu}\}
  \one\{\mathscr E_{\bk,\bm}^{\bu,L,a}\}
  e^{S_{\bk,\bm}^{\bu}+L_{\bk,\bm}^{\bu}}
  (1+|\Delta_{\bk,\bu}^{\varphi}|)
\right)
 \leq C_{a,L}(\Delta_{\bk,\bm}^{\varphi})_+ .
 \label{eq:physical-all-good-upper}
\end{align}
\end{lemma}

The proof is the physical analogue of the upper bound in
Proposition~\ref{prop:brownian-clock-full-tower}: a one-block bound,
obtained from Lemma~\ref{lem:one-sided-one-block} through
\eqref{eq:block-comparison}, is iterated along a partition, and the
piece beyond \(\bk\) is a moment bound.

\begin{proof}
By \eqref{eq:physical-good-consequences}, the factor
\(e^{L_{\bk,\bm}^{\bu}}\) is at most \(e\) on \(\mathscr G_{\bn,\bu}\), so
it suffices to prove \eqref{eq:physical-all-good-upper} without it, that is,  to prove that
\begin{align}
 &\mathbb Q_{\bk,\bm}^{\bu}\left(
  \one\{\mathscr G_{\bn,\bu}\}
  \one\{\mathscr E_{\bk,\bm}^{\bu,L,a}\}
  e^{S_{\bk,\bm}^{\bu}}
  (1+|\Delta_{\bk,\bu}^{\varphi}|)
\right)
 \leq C_{a,L}(\Delta_{\bk,\bm}^{\varphi})_+ .
 \label{eq:physical-all-good-upperbis}
\end{align}
Fix \(0<\eta<1/8\); constants depend on \(a\), \(L\) and \(\eta\).

\emph{One block.}
We claim that for \(\bn_*\leq\bl<\bu'\leq\min\{\lambda_+\bl,\bk\}\),
on \(\mathscr E_{\bk,\bl}^{\bl,L,a}\),
\begin{equation}\label{eq:physical-one-block}
 \mathbb Q_{\bk,\bl}^{\bu'}\left(
  \one\{\mathscr G_{\bn,\bu'}\}
  \one\{\mathscr E_{\bk,\bl}^{\bu',L,a}\}
  e^{S_{\bk,\bl}^{\bu'}}
  (1+|\Delta_{\bk,\bu'}^{\varphi}|)
\right)
 \leq\left(1+C\bl^{-1/8}+Ce^{-c\bl^{\eta}}\right)
 \Delta_{\bk,\bl}^{\varphi} .
\end{equation}
Split according to whether \(S_{\bk,\bl}^{\bu'}\leq-\bl^\eta\).  On
that event \(e^{S_{\bk,\bl}^{\bu'}}\leq e^{-\bl^\eta}\), while by
\eqref{eq:physical-good-consequences} and
\eqref{eq:buffered-unconditional-field-L2},
\(\mathbb Q_{\bk,\bl}^{\bu'}(\one\{\mathscr G_{\bn,\bu'}\}(1+|\Delta_{\bk,\bu'}^{\varphi}|))
\leq C(1+\Delta_{\bk,\bl}^{\varphi}+\bl)\); so this part is at most
\(Ce^{-\bl^\eta}(1+\Delta_{\bk,\bl}^{\varphi}+\bl)
\leq Ce^{-c\bl^\eta}\Delta_{\bk,\bl}^{\varphi}\), by
\eqref{eq:physical-initial-barrier-size}.  On the complementary event, that is
\(S_{\bk,\bl}^{\bu'}>-\bl^\eta\), and using that by \eqref{eq:S-sign}
the deviation \(S_{\bk,\bl}^{\bm}\) is nonpositive and nonincreasing
in \(\bm\) on \(\bl\leq\bm\leq\bu'\), we get
\(|S_{\bk,\bl}^{\bm}|\leq\bl^\eta\leq3\eta_0\ell_\bm^{1/8}\)
for all \(\bl\leq\bm\leq\bu'\), because \(\eta<1/8\) and
\(\ell_\bm\geq\bl\log2\).  Hence \eqref{eq:block-comparison} applies on
\([\bl,\bu']\), and this part of \eqref{eq:physical-one-block} is at
most
\[
 \mathbb Q_{\bk,\bl}^{\bu'}\left(
  \one\{\operatorname{Tru}_{2^\bl,2^{\bu'}}\}
  \one\{\widetilde{\mathscr E}_{\bk,\bl}^{\bu',L+K,a+\eta,[\bn]}\}
  e^{S_{\bk,\bl}^{\bu'}}
  (1+|\Delta_{\bk,\bu'}^{[\bn]}|)
\right)
 +C\epsilon_\bn\,
 \mathbb Q_{\bk,\bl}^{\bu'}\left(e^{S_{\bk,\bl}^{\bu'}}\right),
\]
using \(\mathscr G_{\bn,\bu'}\subseteq\operatorname{Tru}_{2^\bl,2^{\bu'}}\)
and \eqref{eq:physical-good-consequences} at time \(\bu'\).  The
initial barrier \(\widetilde{\mathscr E}_{\bk,\bl}^{\bl,L+K,a+\eta,[\bn]}\)
holds by \eqref{eq:block-comparison} at \(\bm=\bl\), so
Lemma~\ref{lem:one-sided-one-block} with parameters \((L+K,a+\eta)\)
bounds the first term by
\((1+C\bl^{-1/8})\Delta_{\bk,\bl}^{[\bn]}
\leq(1+C\bl^{-1/8})\Delta_{\bk,\bl}^{\varphi}+C\epsilon_\bn\).  The
second term is at most \(C\epsilon_\bn\), since \(S_{\bk,\bl}^{\bu'}\leq c2^{-\bl}\).
As \(\Delta_{\bk,\bl}^{\varphi}\geq c_a\bl^{1/8}\) by
\eqref{eq:physical-initial-barrier-size}, both remainders are certainly
\(C\bl^{-1/8}\Delta_{\bk,\bl}^{\varphi}\), which proves
\eqref{eq:physical-one-block}.

\emph{Iteration.}
Put \(\bu_*=\bu\wedge\bk\) if \(\bm<\bk\), and \(\bu_*=\bm\) if
\(\bm\geq\bk\).  Choose \(\bm=\bu_0<\bu_1<\cdots<\bu_{\bar s}=\bu_*\) with
\(\lambda_-\leq \bu_{s+1}/\bu_s\leq\lambda_+\) for \(s<\bar s-1\) and
\(1<\bu_{\bar s}/\bu_{\bar s-1}\leq\lambda_+\), where \(\bar s=0\) if
\(\bu_*=\bm\).  For \(0\leq s\leq\bar s\), put
\[
 H_s:=\mathbb Q_{\bk,\bu_s}^{\bu}\left(
  \one\{\mathscr G_{\bn,\bu}\}
  \one\{\mathscr E_{\bk,\bu_s}^{\bu,L,a}\}
  e^{S_{\bk,\bu_s}^{\bu}}
  (1+|\Delta_{\bk,\bu}^{\varphi}|)\right),
\]
so that the left side of \eqref{eq:physical-all-good-upperbis}
is \(H_0\).  Since the physical barrier is an
intersection over times, \(\mathscr G_{\bn,\bu}\subseteq\mathscr G_{\bn,\bu_{s+1}}\),
and \(S_{\bk,\bu_s}^{\bu}=S_{\bk,\bu_s}^{\bu_{s+1}}+S_{\bk,\bu_{s+1}}^{\bu}\), the
tower property gives, for \(0\leq s<\bar s\),
\begin{equation}\label{eq:physical-one-block-tower}
 H_s=\mathbb Q_{\bk,\bu_s}^{\bu_{s+1}}\left(
  \one\{\mathscr G_{\bn,\bu_{s+1}}\}
  \one\{\mathscr E_{\bk,\bu_s}^{\bu_{s+1},L,a}\}
  e^{S_{\bk,\bu_s}^{\bu_{s+1}}}H_{s+1}\right).
\end{equation}
We claim that, for \(0\leq s\leq\bar s\),
\begin{equation}\label{eq:physical-iteration-claim}
 H_s\leq A_s\Delta_{\bk,\bu_s}^{\varphi}
 \quad\text{on }\mathscr E_{\bk,\bu_s}^{\bu_s,L,a},
 \qquad\text{with }A_s\leq C,
\end{equation}
which for \(s=0\) gives \eqref{eq:physical-all-good-upperbis}.  The case
\(s=\bar s\) is proved in the next paragraph.  The remaining cases follow
by the backward induction (as in \eqref{eq:bcb-induction}).
Namely, the barrier in \eqref{eq:physical-one-block-tower} includes the
constraint at time \(2^{\bu_{s+1}}\), so wherever its integrand is nonzero
the claim \eqref{eq:physical-iteration-claim} at \(s+1\) gives
\(H_{s+1}\leq A_{s+1}(1+|\Delta_{\bk,\bu_{s+1}}^{\varphi}|)\), and then
\eqref{eq:physical-one-block} gives the claim \eqref{eq:physical-iteration-claim} at \(s\) with
\(A_s=(1+C\bu_s^{-1/8}+Ce^{-c\bu_s^\eta})A_{s+1}\).  These factors have a
bounded product because \(\bu_s\geq\lambda_-^s\bm\).

\emph{The piece beyond \(\bk\).}
If \(\bu=\bu_*\), then
\(H_{\bar s}=\one\{\mathscr G_{\bn,\bu}\}\one\{\mathscr E_{\bk,\bu}^{\bu,L,a}\}
(1+|\Delta_{\bk,\bu}^{\varphi}|)\leq2\Delta_{\bk,\bu}^{\varphi}\) on the
initial barrier at \(\bu\), by \eqref{eq:physical-initial-barrier-size},
so \(A_{\bar s}=2\).  Otherwise \(\bu>\bu_*\geq\bk\), and dropping the
barrier indicator by positivity,
\[
 H_{\bar s}\leq\mathbb Q_{\bk,\bu_*}^{\bu}\left(
  \one\{\operatorname{Tru}_{2^{\bu_*},2^\bu}\}
  e^{S_{\bk,\bu_*}^{\bu}}
  (2+|\Delta_{\bk,\bu}^{[\bu_*]}|)\right)
 \leq C(1+|\Delta_{\bk,\bu_*}^{\varphi}|),
\]
by \eqref{eq:physical-good-consequences} and the moment bound
\eqref{eq:bcb-tail-moments}, whose proof applies verbatim to the
starting exponent \(\bu_*\geq\bk\) and the finite horizon \(\bu\): since
\(\lg_{r_\bk}\) is constant beyond \(2^\bk\), the deficit
\(\Delta_{\bk,\bu}^{[\bu_*]}\) differs from
\(\Delta_{\bk,\bu_*}^{[\bu_*]}\) by \(-\sqrt2S_{\bk,\bu_*}^{\bu}\) plus a
martingale with bracket \(S_{\bk,\bu_*}^{\bu}\geq0\), which has bounded
exponential moments by Lemma~\ref{lem:QVsmallkold}.  By
\eqref{eq:physical-initial-barrier-size}, \(A_{\bar s}\leq C\).
\end{proof}

Recall that the event \(\mathscr G_{\bn,\bn}\) of \eqref{eq:globalgood-finite} is \(\filt_{2^\bn}\)-measurable and implies
\(\{\max_{|z|\leq1}|Q_{2^\bn}'(z)|\leq2^{(C_{\rm in}+2)\bn}\}\), by \eqref{eq:initial-derivative}.

\begin{proposition}[Comparison of the Brownian and physical barriers]
\label{prop:brownian-clock-physical-sandwich}
Fix \(a<0\), \(L<\infty\), and \(\eta>0\) with \(a+\eta<0\).  There is
\(K=K(a,L,\eta)<\infty\) such that, on
\(\mathscr G_{\bn,\bn}\cap\mathscr E_{\bk,\bn}^{\bn,L,a}\),
\begin{align}
 \widetilde{\mathcal U}_{\bk,\bn}^{\mathcal P,L-K,a-\eta}
 -r_{\bk,\bn}\Delta_{\bk,\bn}^{[\bn]}
 \leq\mathcal U_{\bk,\bn}^{\varphi,L,a}\leq
 \widetilde{\mathcal U}_{\bk,\bn}^{\mathcal P,L+K,a+\eta}
 +r_{\bk,\bn}\Delta_{\bk,\bn}^{[\bn]},
 \label{eq:brownian-clock-physical-sandwich}
\end{align}
where \(r_{\bk,\bn}\) is deterministic and
\(\lim_{\bn\to\infty}\limsup_{\bk\to\infty}r_{\bk,\bn}=0\).  In
particular, by \eqref{eq:brownian-clock-global-upper},
\(\mathcal U_{\bk,\bn}^{\varphi,L,a}\leq C_{a,L}\Delta_{\bk,\bn}^{[\bn]}\)
on the same event.
\end{proposition}

\begin{proof}
Write \(L_\pm=L\pm K\) and \(a_\pm=a\pm\eta\), and let
\begin{equation}\label{eq:cmp-majorant}
 I:=\one\{\operatorname{Tru}_{2^\bn}\}
 e^{S_{\bk,\bn}^{\infty}}\bigl(1+(\Delta_\infty)_+\bigr)
\end{equation}
be the common majorant of the two integrands: the integrand of
\(\widetilde{\mathcal U}_{\bk,\bn}^{\mathcal P,L',a'}\) is at most
\(\one\{\widetilde{\mathscr E}_{\bk,\bn}^{\mathcal P,L',a'}\}I\), and by
\eqref{eq:physical-good-consequences} the integrand of
\(\mathcal U_{\bk,\bn}^{\varphi,L,a}\) is at most
\(C\one\{\mathscr G_\bn\}\one\{\mathscr E_{\bk,\bn}^{L,a}\}I\).

\emph{The comparison event.}
With \(\eta_0\) and \(K\) as in \eqref{eq:block-comparison}, let
\begin{equation}\label{eq:cmp-event}
 \mathscr A:=\mathscr G_\bn
 \cap\bigcap_{r=0}^{\bar r}
  \Big\{\max_{\bn_r\leq\bm\leq\bn_{r+1}}
  |S_{\bk,\bn_r}^{\bm}|\leq\eta_0\bn_r^{1/8}\Big\}
 \cap\Big\{\sup_{\bm\geq\bt}
  |S_{\bk,\bt}^{\bm}|\leq\eta_0\bk^{1/8}\Big\}.
\end{equation}
On \(\mathscr A\), the clock condition in \eqref{eq:block-clock-event}
holds on every block.  On a grid block,
\(|S_{\bk,\bn_r}^{\bm}|\leq\eta_0\bn_r^{1/8}\leq3\eta_0\ell_\bm^{1/8}\);
on the terminal block, anchored at \(\bq\), for \(\bm\geq\bt\),
\(|S_{\bk,\bq}^{\bm}|\leq|S_{\bk,\bq}^{\bt}|+|S_{\bk,\bt}^{\bm}|
\leq2\eta_0\bk^{1/8}\leq3\eta_0\ell_\bm^{1/8}\), since
\(\ell_\bm\geq\bt\log2\).  Hence \eqref{eq:block-comparison} holds on
each block, and since the Brownian barrier
\eqref{eq:brownian-clock-barrier} is the intersection of the block
barriers,
\begin{equation}\label{eq:cmp-sandwich}
 \mathscr A\cap\widetilde{\mathscr E}_{\bk,\bn}^{\mathcal P,L_-,a_-}
 \subseteq\mathscr A\cap\mathscr E_{\bk,\bn}^{L,a}
 \subseteq\mathscr A\cap\widetilde{\mathscr E}_{\bk,\bn}^{\mathcal P,L_+,a_+}.
\end{equation}
On \(\mathscr A\cap\widetilde{\mathscr E}_{\bk,\bn}^{\mathcal P,L_\pm,a_\pm}\),
letting \(\bm\to\infty\) in the terminal barrier gives
\(\Delta_\infty\geq-\sqrt2a_\pm(\bk\log2/2)^{1/8}-\sqrt2L_\pm-C\eta_0\bk^{1/8}
\geq c\bk^{1/8}\geq1\), so that \(1+(\Delta_\infty)_+\leq2(\Delta_\infty)_+\);
and by \eqref{eq:physical-good-consequences},
\begin{equation}\label{eq:cmp-integrands}
 \left|e^{L_{\bk,\bn}^{\infty}}(\Delta_{\bk,\infty}^{\varphi})_+
 -(\Delta_\infty)_+\right|
 \leq C\epsilon_\bn\bigl(1+(\Delta_\infty)_+\bigr)
 \leq C\epsilon_\bn(\Delta_\infty)_+ .
\end{equation}
Splitting each of the three expectations according to \(\mathscr A\),
and using \eqref{eq:cmp-sandwich} and \eqref{eq:cmp-integrands} on
\(\mathscr A\) together with
\(\widetilde{\mathcal U}_{\bk,\bn}^{\mathcal P,L_\pm,a_\pm}\leq C\Delta_\bn\)
from \eqref{eq:brownian-clock-global-upper}, we obtain
\begin{equation}\label{eq:cmp-reduction}
 \widetilde{\mathcal U}_{\bk,\bn}^{\mathcal P,L_-,a_-}
 -\operatorname{Err}_-
 -C\epsilon_\bn\Delta_\bn
 \leq\mathcal U_{\bk,\bn}^{\varphi,L,a}
 \leq\widetilde{\mathcal U}_{\bk,\bn}^{\mathcal P,L_+,a_+}
 +\operatorname{Err}_{\rm P}
 +C\epsilon_\bn\Delta_\bn,
\end{equation}
where
\begin{equation}
\label{eq-131026a}
\operatorname{Err}_-:=\mathbb Q_{\bk,\bn}^{\infty}(
 \one\{\widetilde{\mathscr E}_{\bk,\bn}^{\mathcal P,L_-,a_-}\}
 \one\{\mathscr A^\complement\}I),\qquad
\operatorname{Err}_{\rm P}:=C\,\mathbb Q_{\bk,\bn}^{\infty}(
 \one\{\mathscr G_\bn\}\one\{\mathscr E_{\bk,\bn}^{L,a}\}
 \one\{\mathscr A^\complement\}I).
\end{equation}
It remains to show that both errors are \(r_{\bk,\bn}\Delta_\bn\).

\emph{A failure is cheap.}
Let \(\bn\leq\bu\leq\bt\) be a block endpoint, that is
\(\bu\in\{\bn_1,\dots,\bn_{\bar r},\bt\}\), and let
\(F\in\filt_{2^\bu}\) satisfy
\(\mathbb Q_{\bk,\bn}^{\bu}(F\cap\operatorname{Tru}_{2^\bn,2^\bu})\leq p\)
on \(\mathscr G_{\bn,\bn}\cap\mathscr E_{\bk,\bn}^{\bn,L,a}\).  We
claim that
\begin{equation}\label{eq:cmp-failure}
 \mathbb Q_{\bk,\bn}^{\infty}\left(
  \one\{\mathscr B\}\one\{F\}I\right)
 \leq C\bu^2p^{1/2}\Delta_\bn
 \qquad\text{for }
 \mathscr B\in\{\widetilde{\mathscr E}_{\bk,\bn}^{\mathcal P,L_-,a_-},\
 \mathscr G_\bn\cap\mathscr E_{\bk,\bn}^{L,a}\}.
\end{equation}
Condition at \(\filt_{2^\bu}\).  In either case, \(\mathscr B\) is an
intersection over times, so it is the intersection of a
\(\filt_{2^\bu}\)-measurable event \(\mathscr B_{\leq\bu}\), the barrier
up to \(\bu\) (together with \(\mathscr G_{\bn,\bu}\) in the physical
case), with the barrier from \(\bu\) onward, \(\mathscr B_{\geq\bu}\),
which contains the initial barrier at \(\bu\).  Since also
\(\operatorname{Tru}_{2^\bn}=\operatorname{Tru}_{2^\bn,2^\bu}
\cap\operatorname{Tru}_{2^\bu}\) and
\(S_{\bk,\bn}^{\infty}=S_{\bk,\bn}^{\bu}+S_{\bk,\bu}^{\infty}\), the
tower property gives
\begin{equation}\label{eq:cmp-failure-tower}
\begin{gathered}
 \mathbb Q_{\bk,\bn}^{\infty}\left(\one\{\mathscr B\}\one\{F\}I\right)
 =\mathbb Q_{\bk,\bn}^{\bu}\Big(
  \one\{\mathscr B_{\leq\bu}\}\one\{F\}
  \one\{\operatorname{Tru}_{2^\bn,2^\bu}\}
  e^{S_{\bk,\bn}^{\bu}}\,\Sigma_\bu\Big),
 \\
 \Sigma_\bu:=\mathbb Q_{\bk,\bu}^{\infty}\Big(
  \one\{\mathscr B_{\geq\bu}\}\one\{\operatorname{Tru}_{2^\bu}\}
  e^{S_{\bk,\bu}^{\infty}}\bigl(1+(\Delta_\infty)_+\bigr)\Big).
\end{gathered}
\end{equation}
We claim that
\begin{equation}\label{eq:cmp-failure-suffix}
 \Sigma_\bu\leq C\bu\,(1+|\Delta_\bu|).
\end{equation}
For the Brownian barrier with \(\bu=\bn_r\) a grid point, the term
\((\Delta_\infty)_+\) in \(\Sigma_\bu\) gives exactly the quantity
\(H_r\) of \eqref{eq:bcb-H}, with parameters \((L_-,a_-)\), which is at
most \(C\Delta_{\bn_r}\) on the initial barrier at \(\bn_r\) by the
first bound in \eqref{eq:bcb-induction}; the term \(1\) gives at most
\(\mathbb Q_{\bk,\bu}^{\infty}(\one\{\operatorname{Tru}_{2^\bu}\}e^{S_{\bk,\bu}^{\infty}})\leq C\)
by Lemma~\ref{lem:QVsmallkold}; and both vanish off the initial barrier
at \(\bu\), which \(\mathscr B_{\geq\bu}\) contains.  For the Brownian
barrier with \(\bu=\bt\), drop \(\one\{\mathscr B_{\geq\bt}\}\) and use
\eqref{eq:bcb-tail-moments}, which gives
\(C(1+|\Delta_\bt|+\sqrt{\bk-\bt})\), and \(\sqrt{\bk-\bt}\leq\bt\).  For
the physical barrier, \eqref{eq:physical-good-consequences} gives
\(e^{S_{\bk,\bu}^{\infty}}(1+(\Delta_\infty)_+)
\leq Ce^{S_{\bk,\bu}^{\infty}+L_{\bk,\bu}^{\infty}}(1+|\Delta_{\bk,\infty}^{\varphi}|)\)
on \(\mathscr G_\bn\), so Lemma~\ref{lem:physical-all-good-upper} with
\(\bm=\bu\) bounds \(\Sigma_\bu\) by \(C(\Delta_\bu^\varphi)_+\leq C(1+|\Delta_\bu|)\),
again by \eqref{eq:physical-good-consequences}.  This proves
\eqref{eq:cmp-failure-suffix}.

Inserting \eqref{eq:cmp-failure-suffix} into
\eqref{eq:cmp-failure-tower}, dropping \(\one\{\mathscr B_{\leq\bu}\}\),
and applying the Cauchy--Schwarz inequality, we get
\begin{equation}\label{eq:cmp-failure-CS}
 \mathbb Q_{\bk,\bn}^{\infty}\left(\one\{\mathscr B\}\one\{F\}I\right)
 \leq C\bu\,
  \mathbb Q_{\bk,\bn}^{\bu}\bigl(F\cap\operatorname{Tru}_{2^\bn,2^\bu}\bigr)^{1/2}
  \mathbb Q_{\bk,\bn}^{\bu}\Big(e^{2S_{\bk,\bn}^{\bu}}(1+|\Delta_\bu|)^2\Big)^{1/2}
 \leq C\bu\,p^{1/2}\bigl(1+\Delta_\bn+\bu\bigr),
\end{equation}
where the second moment is bounded by
\begin{equation}\label{eq:cmp-block-L2}
 \mathbb Q_{\bk,\bn}^{\bu}\Big(e^{2S_{\bk,\bn}^{\bu}}(1+|\Delta_\bu|)^2\Big)
 \leq C\left(1+\Delta_\bn^2+\bu^2\right),
\end{equation}
which follows from optional sampling in \eqref{eq:clock-gap-Brownian}
and \(-C\bu\leq S_{\bk,\bn}^{\bu}\leq C2^{-\bn}\), exactly as
\eqref{eq:buffered-unconditional-field-L2}.  Since \(\Delta_\bn\geq1\)
on the initial barrier by \eqref{eq:physical-initial-barrier-size},
\(1+\Delta_\bn+\bu\leq C\bu\Delta_\bn\), and \eqref{eq:cmp-failure-CS}
proves \eqref{eq:cmp-failure}.

\emph{The failure probabilities.}
The complement of \(\mathscr A\) is covered by the following events.

First, the clock on a block.  On \(\mathscr G_\bn\),
\eqref{eq:julia-derivative} gives \(R_{2^{\bn_r}}\leq2\cdot2^{C'\bn_r}\)
for every \(r\).  By \eqref{eq:S-sign},
\(\max_{\bn_r\leq\bm\leq\bn_{r+1}}|S_{\bk,\bn_r}^{\bm}|=|S_{\bk,\bn_r}^{\bn_{r+1}}|\),
and \(2^{-\bn_r}\leq\eta_0\bn_r^{1/8}/2\) for \(\bn_r\) large, so
Markov's inequality applied to
\((|S_{\bk,\bn_r}^{\bn_{r+1}}|-2^{-\bn_r})_+\) at level
\(\eta_0\bn_r^{1/8}/2\), with \eqref{eq-11012026enew}, gives, on
\(\mathscr G_{\bn,\bn_r}\),
\begin{equation}\label{eq:cmp-clock-failure}
 \mathbb Q_{\bk,\bn_r}^{\bn_{r+1}}\Big(
  \operatorname{Tru}_{2^{\bn_r},2^{\bn_{r+1}}}\cap
  \Big\{\max_{\bn_r\leq\bm\leq\bn_{r+1}}|S_{\bk,\bn_r}^{\bm}|>\eta_0\bn_r^{1/8}\Big\}\Big)
 \leq\frac{C2^{-\bk}\bigl(2^{C'\bn_r}\bn_{r+1}+2^{\bn_{r+1}}\bigr)}{\eta_0\bn_r^{1/8}}
 \leq C\bk\,2^{-(\log(\bk+2))^\alpha},
\end{equation}
because \(\mathfrak c\bn_{r+1}\leq\bk\) with \(\mathfrak c>C'+2\) on the
grid blocks, while \(C'\bq\leq C'\bk/\mathfrak c^2\) and
\(2^{\bt-\bk}\leq2^{1-(\log(\bk+2))^\alpha}\) on the last block.  The
tower property turns this into a bound on
\(\mathbb Q_{\bk,\bn}^{\bn_{r+1}}\) of the same event intersected with
\(\mathscr G_{\bn,\bn_r}\cap\operatorname{Tru}_{2^\bn,2^{\bn_{r+1}}}\).

Second, the tail.  Since \(S_{\bk,\bt}^{\bm}\geq-(\bk-\bt)\log2\) and
\(S_{\bk,\bt}^{\bm}\leq1+(S_{\bk,\bt}^{\infty})_+\), as in the proof of
Lemma~\ref{lem:bcb-tail}, the tail condition in \eqref{eq:cmp-event}
fails only if \(S_{\bk,\bt}^{\infty}>\eta_0\bk^{1/8}-1\), and
Lemma~\ref{lem:QVsmallkold} with Markov's inequality gives
\begin{equation}\label{eq:cmp-tail-failure}
 \mathbb Q_{\bk,\bt}^{\infty}\Big(\operatorname{Tru}_{2^\bt}\cap
  \Big\{\sup_{\bm\geq\bt}|S_{\bk,\bt}^{\bm}|>\eta_0\bk^{1/8}\Big\}\Big)
 \leq Ce^{-\eta_0\bk^{1/8}} .
\end{equation}
This failure is not \(\filt_{2^\bt}\)-measurable, so we apply the
Cauchy--Schwarz inequality under \(\mathbb Q_{\bk,\bt}^{\infty}\) with
\eqref{eq:bcb-tail-moments} first, which bounds the conditional
expectation after \(\bt\) by \(Ce^{-\eta_0\bk^{1/8}/2}(1+|\Delta_\bt|+\sqrt{\bk})\),
and then argue as in \eqref{eq:cmp-failure} with \(\bu=\bt\) and
\(p=1\).

Third, \(\mathscr G_\bn\) itself, which enters \(\operatorname{Err}_-\)
only, since \(\operatorname{Err}_{\rm P}\) carries
\(\one\{\mathscr G_\bn\}\).  By Lemma~\ref{lem:good-event-tilted-tail},
on \(\mathscr G_{\bn,\bn}\) the event
\(\operatorname{Tru}_{2^\bn}\cap\mathscr G_\bn^\complement\) is covered
by the events \(\mathscr N_\bu\), \(\bu\geq\bn\), of
\eqref{eq:good-event-cover}, and
\begin{equation}\label{eq:cmp-good-failure}
 \mathbb Q_{\bk,\bn}^{\bu+1}\bigl(\mathscr N_\bu\cap\operatorname{Tru}_{2^\bn,2^{\bu+1}}\bigr)
 \leq C2^{-9\bu}.
\end{equation}

\emph{Conclusion.}
Write \(\mathscr A_r\), \(0\leq r\leq\bar r\), for the clock event of the
\(r\)-th block in \eqref{eq:cmp-event} and \(\mathscr A_{\rm tail}\)
for the tail event there.  By the union bound,
\begin{equation}\label{eq:cmp-union}
 \one\{\mathscr A^\complement\}
 \leq\one\{\mathscr G_\bn^\complement\}
 +\sum_{r=0}^{\bar r}\one\{\mathscr G_{\bn,\bn_r}\cap\mathscr A_r^\complement\}
 +\one\{\mathscr G_\bn\cap\mathscr A_{\rm tail}^\complement\},
\end{equation}
and we bound the contribution of each term to \(\operatorname{Err}_-\)
and \(\operatorname{Err}_{\rm P}\) of \eqref{eq-131026a}.  For the \(r\)-th clock term, apply
\eqref{eq:cmp-failure} with \(\bu=\bn_{r+1}\),
\(F=\mathscr G_{\bn,\bn_r}\cap\mathscr A_r^\complement\), and \(p\) the
right side of \eqref{eq:cmp-clock-failure}.  For the tail term, use
\eqref{eq:cmp-tail-failure} as described after it.  The term
\(\mathscr G_\bn^\complement\) concerns \(\operatorname{Err}_-\) only, and
is covered by the events \(\mathscr N_{\bu'}\), \(\bu'\geq\bn\), of
\eqref{eq:good-event-cover}: apply \eqref{eq:cmp-failure} to \(\mathscr N_{\bu'}\)
with \(\bu=\bn_{r+1}\leq\lambda_+\bu'\) when \(\bn_r\leq\bu'<\bn_{r+1}\),
and with \(\bu=\bt\) when \(\bq\leq\bu'<\bt\); for \(\bu'\geq\bt\), where
\(\mathscr N_{\bu'}\) is not measurable at any block endpoint, the Cauchy--Schwarz
inequality under \(\mathbb Q_{\bk,\bn}^{\infty}\),
\eqref{eq:cmp-good-failure} and \eqref{eq:bcb-tail-moments} with
\(\bu=\bn\) bound its contribution by
\(C2^{-9\bu'/2}(1+\Delta_\bn+\sqrt\bk)\), which is harmless because
\(\bu'\geq\bt\geq\bk/2\).  Altogether,
\[
 \operatorname{Err}_-+\operatorname{Err}_{\rm P}
 \leq C\Delta_\bn\Big(
  \sum_{\bu'\geq\bn}\bu'^2\,2^{-9\bu'/2}
  +\sqrt\bk\,2^{-9\bk/4}
  +\bk^{4}2^{-(\log(\bk+2))^\alpha/2}
  +\bk^{4}e^{-\eta_0\bk^{1/8}/2}\Big),
\]
where the number of blocks is at most \(C\log\bk\).  The first term is \(o_\bn(1)\) and
the other three are \(o_\bk(1)\), since \(\alpha>1\).  Together with
\eqref{eq:cmp-reduction} this proves
\eqref{eq:brownian-clock-physical-sandwich}, with
\(r_{\bk,\bn}\) the sum of these four terms and \(C\epsilon_\bn\).
\end{proof}

\section{Comparing the radial barrier to physical barriers: the proof of Proposition \ref{prop:main}}
\label{sec:proof-main-comparison}
In the Gaussian multiplicative-chaos construction of
Proposition~\ref{prop:criticalchaos}, the field is evolved through the
radii \(r_{\bl}\uparrow1\) after the Verblunsky evolution has been run
to infinity, and the barrier \(\mathscr A_\bk^{L,a}\) of
\eqref{eq:Abarrier} is a barrier in the radius.  In the random-matrix
construction, the field is instead evolved in Verblunsky time at a
fixed radius, and the barriers of
Sections~\ref{sec:brownian-clock-approximate-martingale}
and~\ref{sec:brownian-physical-comparison} are barriers in time.  The
mixed field \(\varphi_{2^{\bj}}(r_{\bl},\theta)\) interpolates between
these two evolutions: its effective scale is
\(\lg_{r_{\bl}}(2^{\bj})=(\log2)\min\{\bl,\bj\}\), so the time
evolution at radius \(r_\bl\) saturates when \(\bj\) reaches \(\bl\).
This section passes from the physical barrier in time to the radial
barrier of the chaos construction, through the following four barrier
descriptions:
\begin{center}
\begin{tabular}{@{}ccp{0.52\textwidth}@{}}
\toprule
Barrier & Definition & Description\\
\midrule
\(\mathscr A_{\bk,\bj}^{L,a}\)
& \eqref{eq:ABbarrier}
& fixed time \(\bj\), with the radius varying over
  \(r_{\bl}\), \(1\leq\bl\leq\bk\wedge\bj\)\\
\(\mathscr B_{\bk,\bn}^{\bm,L,\underline a}\),
\(\mathscr B_{\bk,\bn}^{L,\underline a}\)
& \eqref{eq:Btruncated}, \eqref{eq:Bbarrier14}
& joint radial/time barrier: for every time $\bj$ with \(\bn\leq\bj\leq\bm\), impose
  all radial constraints \(1\leq\bl\leq\bk\wedge\bj\); without \(\bm\),
  impose them for every \(\bj\geq\bn\)\\
\(\mathscr E_{\bk,\bn}^{\bu,L,a}\),
\(\mathscr E_{\bk,\bn}^{L,a}\)
& \eqref{eq:physical-interface-barrier}
& fixed radius \(r_{\bk}\), with physical time \(2^{\bm}\) over
  \(\bn\leq\bm\leq\bu\), or over all \(\bm\geq\bn\)\\
\(\widetilde{\mathscr E}_{\bk,\bn}^{L,a}\)
& \eqref{eq:brownian-clock-barrier}
& fixed radius \(r_{\bk}\), with Brownian time starting at \(\bn\)\\
\bottomrule
\end{tabular}
\end{center}

For finite exponent \(\bj\), define the \textbf{fixed-time barrier}
\begin{align}
\mathscr A_{\bk,\bj}^{L,a}(\theta)
&:=
\left\{
\sqrt{\frac{\beta}{8}}\,
\varphi_{2^{\bj}}(r_{\bl},\theta)
\leq
\lg_{r_{\bl}}(2^{\bj})
+a\lg_{r_{\bl}}(2^{\bj})^{1/8}+L:
1\leq\bl\leq\bk\wedge\bj
\right\}.
\label{eq:ABbarrier}
\end{align}
At \(\bj=\infty\) this is the Gaussian radial barrier
\(\mathscr A_{\bk}^{L,a}\) from \eqref{eq:Abarrier}.  If
\(\underline a=(a_\bj)\) decreases to \(a_\infty\), define the
\textbf{radial/time barrier}
\begin{equation}\label{eq:Bbarrier14}
\mathscr B_{\bk,\bn}^{L,\underline a}
:=\bigcap_{\bj\geq\bn}\mathscr A_{\bk,\bj}^{L,a_\bj}.
\end{equation}

\begin{lemma}[Eventual finite-time radial barrier]\label{lem-1.20}
Suppose \(a_\bj\downarrow a_\infty\) with
\(a<a_\infty<0\).  For every \(\bk,\bn\),
\begin{equation}\label{eq-13012026a}
\mathscr A_{\bk}^{L,a}
\subseteq
\bigcup_{\bl\geq\bn}\mathscr B_{\bk,\bl}^{L,\underline a},
\qquad
\mathscr B_{\bk,\bl}^{L,\underline a}
\subseteq\mathscr A_{\bk}^{L,a_\infty}.
\end{equation}
\end{lemma}

\begin{proof}
On \(\mathscr A_{\bk}^{L,a}\), for every \(\bl\leq\bk\),
\[
\sqrt{\frac{\beta}{8}}\varphi_\infty(r_\bl,\theta)
\leq \bl\log2+a(\bl\log2)^{1/8}+L.
\]
The convergence
\(\varphi_{2^\bj}(r_\bl,\theta)\to\varphi_\infty(r_\bl,\theta)\)
is simultaneous over the finitely many radii \(\bl\leq\bk\).
Since \(a_\infty>a\), all the inequalities defining
\(\mathscr A_{\bk,\bj}^{L,a_\bj}\) therefore hold once \(\bj\) is
sufficiently large.  This proves the first inclusion.  The second follows
by taking \(\bj\to\infty\) in the defining inequalities.
\end{proof}

We first establish the estimates that transport the barrier information
from the Gaussian construction to the physical process, successively adding
the radial/time barriers and good events.  The proof of
Proposition~\ref{prop:main} at the end of the section gives an overview of
this sequence of approximations.

\subsection{Preliminary derivative-martingale estimates}

The barrier comparisons below are made after integration in
\(\theta\).  We first record two estimates used repeatedly: a
one-sided conditional bound for the remaining physical-time evolution,
followed by derivative-martingale bounds for exceptional radial and
boundary configurations.

\begin{lemma}[$\mathscr{D}$ approximate-martingale bound]
\label{lem:physical-suffix-transfer}
Fix \(a<0\) and \(L<\infty\).  Uniformly for \(\bm\geq\bn\),
\begin{align}
 &\mathbb E_\bm\left(
 \one\{\mathscr G_\bn\}
 \one\{\mathscr E_{\bk,\bm}^{L,a}\}
 \mathscr D_\infty(r_\bk,\theta)\right)\leq C_{a,L}
 \one\{\mathscr G_{\bn,\bm}\}
 \one\{\mathscr E_{\bk,\bm}^{\bm,L,a}\}
 \mathscr D_{2^\bm}(r_\bk,\theta).
 \label{eq:physical-suffix-transfer}
\end{align}
\end{lemma}

\begin{proof}
Both sides of \eqref{eq:physical-suffix-transfer} vanish off
\(\mathscr G_{\bn,\bm}\cap\mathscr E_{\bk,\bm}^{\bm,L,a}\), since this
event is \(\filt_{2^\bm}\)-measurable and contains
\(\mathscr G_\bn\cap\mathscr E_{\bk,\bm}^{L,a}\); on it,
Lemma~\ref{lem:physical-all-good-upper} with \(\bu=\infty\) applies.
The change-of-measure identity \eqref{eq:com} and that lemma give
\begin{align*}
 &\mathbb E_\bm\left(
 \one\{\mathscr G_\bn\}
 \one\{\mathscr E_{\bk,\bm}^{L,a}\}
 \mathscr D_\infty(r_\bk,\theta)\right)\leq
 C\one\{\mathscr G_{\bn,\bm}\}
 \one\{\mathscr E_{\bk,\bm}^{\bm,L,a}\}
 \mathscr M_{\bk,\bm}(\theta)(\Delta_{\bk,\bm}^{\varphi})_+,
\end{align*}
which is the right side of \eqref{eq:physical-suffix-transfer}.
\end{proof}

\emph{The boundary good event.}  For \(\bu\in\mathbb N\) and
\(\theta\in[0,2\pi]\), let
\begin{equation}\label{eq:boundary-deficit}
 d_\bu^\partial(\theta):=\sqrt2\,\bu\log2
 -\sqrt{\tfrac\beta4}\,\varphi_{2^\bu}(1,\theta)
\end{equation}
be the \textbf{boundary deficit}, so that
\(\mathscr D_{2^\bu}(1,\theta)=\mathscr M_{\infty,\bu}(\theta)(d_\bu^\partial(\theta))_+\)
by \eqref{eq:Bk} and \eqref{eq:Mn}.  Fix \(K\in\mathbb N\), put
\(K_0:=\lfloor K/\mathfrak c\rfloor\) and
\(\kappa_K:=\tfrac1{16}\log(K\log2)\), and let \(\mathscr H_K^\partial\)
be the \textbf{boundary good event}, the intersection of the three events
\begin{equation}\label{eq:boundary-good-event}
\begin{gathered}
 \widehat{\mathscr G}_{K_0},\qquad
 \Big\{\sup_{r,\theta}\sup_{\bj\geq K_0}
 |\zeta^{\mathbb C}_{2^\bj}(K_0;r,\theta)-2\log\Phi^*_{2^\bj}(r,\theta)|
 \leq\epsilon_{K_0}\Big\},\\
 \Big\{\sqrt{\tfrac\beta8}\,\varphi_{2^\bj}(1,\theta)
 \leq\bj\log2-\tfrac18\log(\bj\log2)+\kappa_K
 \ \text{ for all }\bj\geq1\text{ and }\theta\in[0,2\pi]\Big\}.
\end{gathered}
\end{equation}
The first two are the derivative and linearization parts of the good
event, see \eqref{eq:globalgooda} and \eqref{eq:globalgood}, and have
probability tending to one as \(K\to\infty\).  The third is the event
\(E_\kappa\), with \(\kappa=\kappa_K\), on which
\cite[Lemma~9.4]{PaquetteZeitouni03} is proved; by
\cite[Lemma~9.3]{PaquetteZeitouni03} the union of these events over
\(\kappa\) has full probability.  The events \(E_\kappa\) increase with \(\kappa\), and \(\kappa_K\to\infty\) as \(K\to\infty\), so \(\mathbb P(E_{\kappa_K})\to1\).
Hence \(\mathbb P((\mathscr H_K^\partial)^\complement)\to0\) as
\(K\to\infty\).  On \(\mathscr H_K^\partial\), for all \(\bj\geq K\) and
all \(\theta\),
\begin{equation}\label{eq:boundary-slack}
 d_\bj^\partial(\theta)
 \geq\sqrt2\Big(\tfrac18\log(\bj\log2)-\kappa_K\Big)
 \geq\tfrac{\sqrt2}{16}\log(\bj\log2)=:s_\bj ,
\end{equation}
and we take \(K\) so large that \(s_K\geq1\); in particular
\(\mathscr D_{2^\bj}(1,\theta)=\mathscr M_{\infty,\bj}(\theta)\,d_\bj^\partial(\theta)\)
there.  For \(\bn\geq K\) we write \(\mathscr H_{K,\bn}^\partial\) for
the event obtained by restricting the three constraints in
\eqref{eq:boundary-good-event} to scales \(\bj\leq\bn\); it is
\(\filt_{2^\bn}\)-measurable, decreasing in \(\bn\), and contains
\(\mathscr H_K^\partial\).

The next lemma is the boundary analogue of
Lemma~\ref{lem:physical-suffix-transfer}: on the boundary good event
the derivative martingale is killed at a barrier with slack
\(s_\bj\), and the slack grows just fast enough that the killed
martingale has bounded conditional expectation, uniformly in the
terminal index.

\begin{lemma}[Boundary suffix transfer]
\label{lem:boundary-suffix-transfer}
There is \(C<\infty\) such that for all large \(K\), all
\(K\leq\bu\leq\bn\), all \(\theta\), and every nonnegative
\(\filt_{2^\bu}\)-measurable \(A\),
\begin{equation}\label{eq:boundary-suffix-transfer}
 \mathbb E\bigl(\one\{\mathscr H_{K,\bn}^\partial\}\mathscr D_{2^\bn}(1,\theta)A\bigr)
 \leq C\,\mathbb E\bigl(\one\{\mathscr H_{K,\bu}^\partial\}\mathscr D_{2^\bu}(1,\theta)A\bigr).
\end{equation}
\end{lemma}

\begin{proof}
By \eqref{eq:com} with \(\bk=\infty\) and base \(\bu\), the left side of
\eqref{eq:boundary-suffix-transfer} equals
\[
 \mathbb E\Big(\mathscr M_{\infty,\bu}(\theta)A\,
 \mathbb Q_{\infty,\bu}^{\bn}\big(e^{S_{\infty,\bu}^{\bn}+L_{\infty,\bu}^{\bn}}
 \one\{\mathscr H_{K,\bn}^\partial\}d_\bn^\partial\big)\Big).
\]
Here \(S^{\bn}_{\infty,\bu}=\sum_{2^\bu\leq j<2^\bn}(j+1)^{-1}-(\bn-\bu)\log2\leq0\)
is deterministic, because \(|Q_j(1,\theta)|=1\), and
\(|L^{\bn}_{\infty,\bu}|\leq C\epsilon_{K_0}\leq C\) on
\(\mathscr H^\partial_{K,\bn}\) by the linearization part of
\eqref{eq:boundary-good-event}.  Since
\(\mathscr H^\partial_{K,\bn}\subseteq\mathscr H^\partial_{K,\bu}\), it
therefore suffices to show that, on \(\mathscr H^\partial_{K,\bu}\),
\begin{equation}\label{eq:boundary-killed-endpoint}
 \mathbb Q_{\infty,\bu}^{\bn}\big(\one\{\mathscr H_{K,\bn}^\partial\}d_\bn^\partial\big)
 \leq d_\bu^\partial+C\leq Cd_\bu^\partial,
\end{equation}
the last step by \eqref{eq:boundary-slack}.

Under \(\mathbb Q^{\bn}_{\infty,\bu}\), by \eqref{eq:Qt} and
\eqref{eq-301024}, the process
\(\sqrt{\beta/4}\,(\zeta_s(\bu;1,\theta)-\varphi_{2^\bu}(1,\theta))-\sqrt2\,c(s)\),
with \(c(s):=\int_{2^\bu}^s\frac{\dif v}{[v]+1}\), is a continuous
martingale whose bracket is the deterministic function \(c(s)\); hence
it equals \(W(c(s))\) for a standard Brownian motion \(W\).  Writing
\(c_\bm:=c(2^\bm)\) and using \eqref{eq:StLt}, for
\(\bu\leq\bm\leq\bn\),
\[
 d_\bm^\partial
 =d_\bu^\partial-W(c_\bm)-\sqrt2\,e_{\bu,\bm}
 -\tfrac1{\sqrt2}L^{\bm}_{\infty,\bu},
\]
with \(e_{\bu,\bm}\in[-2^{-\bu},0]\) the error of \eqref{Qdecompsum}.
Put \(X_\bm:=d_\bu^\partial-W(c_\bm)\), a martingale in \(\bm\) with
\(X_\bu=d_\bu^\partial\).  On \(\mathscr H^\partial_{K,\bn}\),
\eqref{eq:boundary-slack} and the bounds on \(e\) and \(L\) give
\(X_\bm\geq s_\bm-C\geq s_\bm/2\) for all \(\bu\leq\bm\leq\bn\), once
\(K\) is large, and \(d_\bn^\partial\leq X_\bn+C\).  Let \(\tau\) be
the first \(\bm\in(\bu,\bn]\) with \(X_\bm<s_\bm/2\), and
\(\tau=\infty\) if there is none, so that
\(\mathscr H^\partial_{K,\bn}\subseteq\{\tau>\bn\}\).  Optional
stopping gives
\begin{align*}
 \mathbb Q^{\bn}_{\infty,\bu}\big(\one\{\mathscr H^\partial_{K,\bn}\}d^\partial_\bn\big)
 &\leq\mathbb Q^{\bn}_{\infty,\bu}\big(X_\bn\one\{\tau>\bn\}\big)+C
 =d_\bu^\partial-\mathbb Q^{\bn}_{\infty,\bu}\big(X_\tau\one\{\tau\leq\bn\}\big)+C\\
 &\leq d_\bu^\partial
 +\sum_{\bm=\bu+1}^{\bn}\mathbb Q^{\bn}_{\infty,\bu}\big((X_\bm)_-\one\{\tau=\bm\}\big)+C.
\end{align*}
On \(\{\tau=\bm,\ X_\bm<0\}\) we have \(X_{\bm-1}\geq s_{\bm-1}/2\), so
\(W\) moves by more than \(s_{\bm-1}/2\) over \([c_{\bm-1},c_\bm]\),
an interval of length at most \(\log2\), and \((X_\bm)_-\) is at most
the oscillation of \(W\) on that interval.  The Gaussian tail of the
Brownian oscillation therefore gives
\[
 \sum_{\bm>\bu}\mathbb Q^{\bn}_{\infty,\bu}\big((X_\bm)_-\one\{\tau=\bm\}\big)
 \leq C\sum_{\bm>\bu}e^{-cs_{\bm-1}^2}\leq C,
\]
since \(s_\bm\) grows like \(\log\bm\).  This proves
\eqref{eq:boundary-killed-endpoint}.
\end{proof}

Two estimates recur in the rest of this section, and we record them
first.  Several of the arguments below change measure from an earlier
base than \(\bn\).  For the rest of this section we put
\[
 \bn':=\lfloor\bn/\mathfrak c\rfloor .
\]
Changing the base to \(\bn'\) is harmless, for two reasons.  First, the
event \(\mathscr G_{\bn',\bn}\) contains \(\mathscr G_{\bn'}\), so its
probability tends to one.  It is also \(\filt_{2^\bn}\)-measurable.  We may
therefore insert the indicator \(\one\{\mathscr G_{\bn',\bn}\}\) in front of
any of the conditional expectations \(\mathbb E_\bn(\cdot)\) below without
affecting convergence in probability.  Second, the Gaussian approximation
started at \(2^{\bn'}\) remains close to the field on
\(\mathscr G_{\bn',\bn}\cap\mathscr G_{\bn,\bt}\).  Indeed, by \eqref{eq:Z},
for \(\bj\geq\bn\) the approximations started at \(2^{\bn'}\) and at
\(2^\bn\) differ by \(\zeta^\C_{2^\bn}(\bn';\cdot)-2\log\Phi^*_{2^\bn}\).
The two linearization errors therefore add, and on
\(\mathscr G_{\bn',\bn}\cap\mathscr G_{\bn,\bt}\),
\begin{equation}\label{eq:base-change-linearization}
 \bigl|\zeta^\C_{2^\bj}(\bn';r,\theta)-2\log\Phi^*_{2^\bj}(r,\theta)\bigr|
 \leq\epsilon_{\bn'}+\epsilon_\bn,
 \qquad \bn'\leq\bj\leq\bt .
\end{equation}
The event \(\mathscr G_{\bn',\bn}\cap\mathscr G_{\bn,\bt}\) is also contained
in \(\widehat{\mathscr G}_{\bn',\bt}\).

The first estimate bounds the chaos mass of an event \(F\) that is
decided by time \(2^\bt\), when the physical barrier holds from time
\(2^\bt\) on.  The bound is the square root of the tilted probability of
\(F\), times a factor of the size of the deficit.

\begin{lemma}[Rare events before a physical suffix]\label{lem:suffix-CS}
Fix \(b<0\) and \(L<\infty\).  There are \(C\) and \(c>0\) such
that, for \(\bn\) large, \(\bn''\in\{\bn',\bn\}\), \(\bt\geq\bn\),
\(\bk\geq\bn\) and \(F\in\filt_{2^\bt}\),
\begin{align}
 &\mathbb E_{\bn''}\Bigl(\one\{\mathscr G_{\bn'',\bn}\}\one\{\mathscr G_\bn\}\one\{F\}
 \one\{\mathscr E_{\bk,\bt}^{L,b}\}\mathscr D_\infty(r_\bk,\theta)\Bigr)\notag\\
 &\qquad\leq C\mathscr M_{\bk,\bn''}\bigl(1+|\Delta^{\varphi}_{\bk,\bn''}|+\sqrt{\bt-\bn''}\bigr)
 \Bigl(\mathbb Q^{\infty}_{\bk,\bn''}\bigl(F\cap\widehat{\mathscr G}_{\bn''}\bigr)+C2^{-c\bt}\Bigr)^{1/2}.
 \label{eq:suffix-CS}
\end{align}
\end{lemma}

\begin{proof}
The proof follows the argument of
\eqref{eq:cmp-failure-tower}--\eqref{eq:cmp-failure-CS} in the proof of
Proposition~\ref{prop:brownian-clock-physical-sandwich}.  We describe
its three steps.

First, we condition on \(\filt_{2^\bt}\).  By
\eqref{eq:physical-suffix-transfer}, on the physical barrier the
conditional expectation of \(\mathscr D_\infty(r_\bk,\theta)\) is at most
\(C\mathscr D_{2^\bt}(r_\bk,\theta)=C\mathscr M_{\bk,\bt}(\Delta^\varphi_{\bk,\bt})_+\).

Second, we change measure from time \(2^{\bn''}\) by \eqref{eq:com}.  This
produces the factor \(\mathscr M_{\bk,\bn''}\) and the weight
\(e^{S^\bt_{\bk,\bn''}+L^\bt_{\bk,\bn''}}\) under
\(\mathbb Q^\bt_{\bk,\bn''}\).  On the good events, \(L^\bt_{\bk,\bn''}\) is
bounded by \eqref{eq:base-change-linearization}, and so is the difference
between \(\Delta^\varphi_{\bk,\bt}\) and \(\Delta^{[\bn'']}_{\bk,\bt}\).

Third, we apply the Cauchy--Schwarz inequality under
\(\mathbb Q^\bt_{\bk,\bn''}\).  The weighted second moment of the deficit
satisfies
\[
 \mathbb Q^\bt_{\bk,\bn''}\Bigl(e^{2S^\bt_{\bk,\bn''}}\bigl(1+|\Delta^{[\bn'']}_{\bk,\bt}|\bigr)^2\Bigr)
 \leq C\bigl(1+|\Delta^\varphi_{\bk,\bn''}|+\sqrt{\bt-\bn''}\bigr)^2 .
\]
To see this, recall from \eqref{eq:clock-gap-Brownian} that the deficit
is a martingale, minus \(\sqrt2S^\bt_{\bk,\bn''}\).  When \(\bt\leq\bk\), the
clock deviation \(S^\bt_{\bk,\bn''}\) is nonpositive, the bracket of the
martingale is at most \((\bt-\bn'')\log2\), and the weight
\(e^{S}\) absorbs the correction \(\sqrt2S\).  When \(\bt>\bk\), the positive
part of the clock deviation has exponential moments by
Lemma~\ref{lem:QVsmallkold}, and the same bound follows from H\"older's
inequality.

It remains to pass from \(\mathbb Q^\bt_{\bk,\bn''}(F\cap\widehat{\mathscr G}_{\bn'',\bt})\)
to the probability in \eqref{eq:suffix-CS}.  The two measures
\(\mathbb Q^\bt_{\bk,\bn''}\) and \(\mathbb Q^\infty_{\bk,\bn''}\) agree on
\(\filt_{2^\bt}\).  The event \(\widehat{\mathscr G}_{\bn'',\bt}\) differs from
\(\widehat{\mathscr G}_{\bn''}\) only through the truncation and derivative
constraints after time \(2^\bt\).  These fail with probability at most
\(C2^{-c\bt}\), by the proof of \eqref{eq:tilted-truncation-cost} and by
\eqref{eq:qkprimebound}.
\end{proof}

The second estimate applies Lemma~\ref{lem:suffix-CS} to the
increments of the field at a fixed radius after a fixed time.

\begin{corollary}[Weighted increments]\label{cor:weighted-increments}
In the setting of Lemma~\ref{lem:suffix-CS}, let \(\bl\leq\bj\),
\(\bj\geq\bn\) and \(x>0\).  Then
\begin{align*}
 &\mathbb E_\bn\Bigl(\one\{\mathscr G_\bn\}\one\{\mathscr E^{L,b}_{\bk,\bj+1}\}
 \one\Bigl\{\sup_{\bv\geq\bj}|\varphi_{2^\bv}(r_\bl,\theta)-\varphi_{2^\bj}(r_\bl,\theta)|>x\Bigr\}
 \mathscr D_\infty(r_\bk,\theta)\Bigr)\\
 &\qquad\leq C\mathscr M_{\bk,\bn}\bigl(1+|\Delta^{\varphi}_{\bk,\bn}|\bigr)\bigl(1+\sqrt{\bj-\bn}\bigr)
 \Bigl(\min\bigl\{e^{-cx^{1/3}},\,x^{-1/2}2^{-(\bj-\bl)/12}\bigr\}+2^{-c\bj}\Bigr).
\end{align*}
\end{corollary}

\begin{proof}
We decompose the event according to the first dyadic increment that
is large.  Fix \(c_0>0\) with \(\sum_{d\geq0}c_0(10+d)^{-2}\leq1\).  If the
supremum exceeds \(x\), then for some \(\bv\geq\bj\),
\[
 |\varphi_{2^{\bv+1}}(r_\bl,\theta)-\varphi_{2^\bv}(r_\bl,\theta)|>c_0x(10+\bv-\bj)^{-2}.
\]
Let \(\bv\) be the least such index.  The corresponding event is
\(\filt_{2^{\bv+1}}\)-measurable, so Lemma~\ref{lem:suffix-CS} applies with
\(\bn''=\bn\) and \(\bt=\bv+1\).  Corollary~\ref{cor:fieldincrement} bounds
the tilted probability of the event.  The sum over \(\bv\) is then the
computation in the proof of Corollary~\ref{cor:ezfieldincrement}: the
exponential bound is used when \(\bv-\bj\leq x^{1/3}\), and the Markov
bound when \(\bv-\bj>x^{1/3}\).
\end{proof}

Both estimates give bounds of the form
\(C\mathscr M_{\bk,\bn}(1+|\Delta^\varphi_{\bk,\bn}|)\) times a small factor.
The next lemma shows that such bounds imply convergence in probability.

\begin{lemma}[Tightness of the weighted deficit]\label{lem:T-tight}
For \(\bk\geq\bn\) let
\[
 \Upsilon_{\bk,\bn}:=\int_0^{2\pi}\mathscr M_{\bk,\bn}(\theta)\bigl(1+|\Delta^{\varphi}_{\bk,\bn}(\theta)|\bigr)\,d\theta .
\]
\emph{(a)} The family \(\Upsilon_{\bk,\bn}\) is tight in the iterated limit:
\(\lim_{A\to\infty}\limsup_{\bn\to\infty}\limsup_{\bk\to\infty}\Pr(\Upsilon_{\bk,\bn}>A)=0\).

\emph{(b)} Let \(\bn''\in\{\bn',\bn\}\) and \(\varrho\geq0\).  If
\(Z\geq0\) satisfies \(\Exp(Z\mid\filt_{2^{\bn''}})\leq\varrho\,\Upsilon_{\bk,\bn''}\),
then for all \(A,\varepsilon>0\),
\[
 \Pr(Z>\varepsilon)\leq\Pr(\Upsilon_{\bk,\bn''}>A)+A\varrho/\varepsilon .
\]
In particular, if \(\varrho=\varrho_\bn\to0\), then \(Z\to0\) in probability
as \(\bk\to\infty\) and then \(\bn\to\infty\).
\end{lemma}

\begin{proof}
\emph{(a)} We treat separately the integrals of
\(\mathscr M_{\bk,\bn}\), of \(\mathscr M_{\bk,\bn}(\Delta^\varphi_{\bk,\bn})_+\),
and of \(\mathscr M_{\bk,\bn}(\Delta^\varphi_{\bk,\bn})_-\).  The first has
expectation \(2\pi\Exp\mathscr M_{\bk,\bn}\).  This is bounded for \(\bn\)
and \(\bk\) large, by \eqref{eq:identify-step3-ck} and the lower bound in
\eqref{eq:step3-ck-bounds}.

For the other two, fix \(\bn\) and let \(\bk\to\infty\).  The field
\(\varphi_{2^\bn}(r_\bk,\cdot)\) converges uniformly to
\(\varphi_{2^\bn}(1,\cdot)\).  Hence \(\mathscr M_{\bk,\bn}\) and
\(\Delta^\varphi_{\bk,\bn}\) converge uniformly to \(\mathscr M_{\infty,\bn}\)
and \(d^\partial_\bn\).  The limit of the second integral is
\(\int_0^{2\pi}\mathscr M_{\infty,\bn}(d^\partial_\bn)_+\,d\theta=\mathscr B_{2^\bn}\),
which is tight by Theorem~\ref{thm:Bn}.  The limit of the third integral
is the mass of \(\mathscr M_{\infty,\bn}|d^\partial_\bn|\) on the set where
\(d^\partial_\bn\leq0\).  On that set, the indicator \(\chi\) of
Theorem~\ref{thm:Bn} with \(k=2^\bn\) equals one, because the compact set
\(K\) lies in \((0,\infty)\).  So \eqref{eq:Btight} shows that this mass is
small in probability.

\emph{(b)} By the conditional Markov inequality,
\(\Pr(Z>\varepsilon\mid\filt_{2^{\bn''}})\leq\varrho\Upsilon_{\bk,\bn''}/\varepsilon\).
On the event \(\{\Upsilon_{\bk,\bn''}\leq A\}\), the right side is at most
\(A\varrho/\varepsilon\).
\end{proof}

\begin{lemma}[The chaos charges tame fields]
\label{lem:boundary-spine-inputs}
Fix \(b<0\), \(L<\infty\), and \(\eta>0\), and let
\(\mathscr H_K^\partial\) be as in \eqref{eq:boundary-good-event}.

\emph{(i) Shallow boundary deficits.}  For every
\(\alpha\in(1/8,1/6)\) and every \(\varepsilon>0\),
\begin{align}
 \lim_{K\to\infty}\sup_{\bn\geq K}
 \mathbb P\Bigg(&\int_0^{2\pi}\one\{\mathscr H_K^\partial\}
 \one\Big\{\exists\,\bu\in[K,\bn]:
 0<d_\bu^\partial(\theta)\leq C\bu^\alpha\Big\}
 \mathscr D_{2^\bn}(1,\theta)d\theta>\varepsilon\Bigg)=0 .
 \label{eq:shallow-gap-spine}
\end{align}

\emph{(ii) Boundary-to-radius comparison.}  For every
\(\varepsilon>0\),
\begin{align}
 \lim_{K\to\infty}\sup_{\bn\geq K}
 \mathbb P\Bigg(&\int_0^{2\pi}\one\{\mathscr H_K^\partial\}
 \one\Big\{\exists\,\bu\in[K,\bn]:
 |\varphi_{2^\bn}(r_\bu,\theta)
  -\varphi_{2^\bu}(1,\theta)|>\bu^{1/10}\Big\}
 \mathscr D_{2^\bn}(1,\theta)d\theta>\varepsilon\Bigg)=0.
 \label{eq:weighted-radial-spine}
\end{align}

\emph{(iii) Two-radii comparison.}  Put
\[
 \mathscr R_{\bm}^{\eta}
 :=\left\{
 |\varphi_{2^\bm}(r_\bk,\theta)
   -\varphi_{2^\bm}(r_\bm,\theta)|
 >\eta(\bm\log2)^{1/8}\right\}.
\]
With \(\bn'=\lfloor\bn/\mathfrak c\rfloor\), for every measurable
\(0\leq f\leq1\),
\begin{align}
 \int_0^{2\pi}f(\theta)
 \sum_{\bm=\bn}^{\bk-1}
 \mathbb E_\bn\left(
  \one\{\mathscr G_{\bn',\bn}\}
  \one\{\mathscr G_\bn\}
  \one\{\mathscr R_\bm^\eta\}
  \one\{\mathscr E_{\bk,\bm+1}^{L,b}\}
  \mathscr D_\infty(r_\bk,\theta)\right)d\theta
 \Prto[\bn,\bk]0.
 \label{eq:weighted-radial-mismatch}
\end{align}
\end{lemma}

\begin{proof}
\emph{Proof of (i).}
The event \(\{0<d_\bu^\partial\leq C\bu^\alpha\}\) is
\(\filt_{2^\bu}\)-measurable, so Lemma~\ref{lem:boundary-suffix-transfer}
reduces the terminal index from \(\bn\) to \(\bu\):
\begin{align*}
 \mathbb E\int_0^{2\pi}\one\{\mathscr H_K^\partial\}
 \one\{0<d_\bu^\partial(\theta)\leq C\bu^\alpha\}
 \mathscr D_{2^\bn}(1,\theta)d\theta
 &\leq C\,\mathbb E\int_0^{2\pi}\one\{\mathscr H_{K,\bu}^\partial\}
 \one\{0<d_\bu^\partial(\theta)\leq C\bu^\alpha\}
 \mathscr D_{2^\bu}(1,\theta)d\theta\\
 &\leq C\bu^{3\alpha-3/2}.
\end{align*}
The last bound is the Bessel small-gap estimate in the proof of
\cite[Lemma~9.4]{PaquetteZeitouni03}, at the single scale \(2^\bu\); it
is proved by the ballot theorem on the event that the boundary field
stays below the barrier in \eqref{eq:boundary-good-event} up to that
scale, which contains \(\mathscr H^\partial_{K,\bu}\).  Since
\(3\alpha-3/2<-1\), summing over \(\bu\geq K\), Tonelli's theorem and
Markov's inequality prove \eqref{eq:shallow-gap-spine}.

\emph{Proof of (ii).}
Fix \(K\) large and \(K\leq\bu\leq\bn\), and put \(x_\bu:=\bu^{1/10}\).
Since \(\sum_{\bm\geq\bu}(10+\bm-\bu)^{-2}\leq1\), the event
\(\{|\varphi_{2^\bn}(r_\bu,\theta)-\varphi_{2^\bu}(1,\theta)|>x_\bu\}\)
is contained in the union of
\begin{align*}
 B_\bu^{0}&:=\{|\varphi_{2^\bu}(r_\bu,\theta)-\varphi_{2^\bu}(1,\theta)|>x_\bu/2\},\\
 B_{\bu,\bm}&:=\Big\{|\varphi_{2^{\bm+1}}(r_\bu,\theta)-\varphi_{2^\bm}(r_\bu,\theta)|
 >\frac{x_\bu}{2(10+\bm-\bu)^2}\Big\},
 \qquad\bu\leq\bm<\bn .
\end{align*}
Each of these events is measurable at the scale at which it is tested,
\(2^\bu\) for \(B^0_\bu\) and \(2^{\bm+1}\) for \(B_{\bu,\bm}\).  Let
\(B\) be one of them, tested at scale \(2^{\bm'}\), and let \(p_B\) be a
deterministic bound for
\(\mathbb Q^{\bm'}_{\infty,K_0}(\mathscr H^\partial_{K,\bm'}\cap B)\).
Lemma~\ref{lem:boundary-suffix-transfer}, then \eqref{eq:com} at the
base \(K_0\), the bound \(e^{S+L}\leq C\) on
\(\mathscr H^\partial_{K,\bm'}\), Cauchy--Schwarz, and the second-moment
bound
\(\mathbb Q^{\bm'}_{\infty,K_0}(\one\{\mathscr H^\partial_{K,\bm'}\}(d^\partial_{\bm'})^2)
\leq C((d^\partial_{K_0})^2+\bm')\), which follows from the Brownian
representation in the proof of Lemma~\ref{lem:boundary-suffix-transfer},
give
\begin{align*}
 \mathbb E\big(\one\{\mathscr H_K^\partial\}\one\{B\}\mathscr D_{2^\bn}(1,\theta)\big)
 &\leq C\,\mathbb E\Big(\mathscr M_{\infty,K_0}
 \mathbb Q^{\bm'}_{\infty,K_0}\big(\one\{\mathscr H^\partial_{K,\bm'}\}d^\partial_{\bm'}\one\{B\}\big)\Big)\\
 &\leq C(1+\bm')^{1/2}p_B^{1/2}\,
 \mathbb E\big(\mathscr M_{\infty,K_0}(1+|d^\partial_{K_0}|)\big),
\end{align*}
and the last expectation is at most \(C(1+K)\).

For \(B=B^0_\bu\), the linearization part of
\(\mathscr H^\partial_{K,\bu}\) bounds the linearization errors at
time \(2^\bu\) and at the two radii \(r_\bu\) and \(1\) by
\(C\epsilon_{K_0}\), so on \(\mathscr H^\partial_{K,\bu}\) the event
\(B^0_\bu\) forces
\(|\zeta_{2^\bu}(K_0;r_\bu,\theta)-\zeta_{2^\bu}(K_0;1,\theta)|>x_\bu/4\),
and Corollary~\ref{cor:radialincrement} with \(\bl=\bu\geq\mathfrak cK_0\)
gives \(p_{B^0_\bu}\leq Ce^{-c\bu^{1/10}}\).  For \(B=B_{\bu,\bm}\),
Corollary~\ref{cor:fieldincrement} with \(\bl=\bu\) and \(\bj=\bm\)
gives
\[
 p_{B_{\bu,\bm}}\leq C\min\Big\{e^{-x_\bu(10+\bm-\bu)^{-2}/2},\
 \frac{2(10+\bm-\bu)^2\,2^{-(\bm-\bu)/6}}{x_\bu}\Big\}.
\]
Using the first entry for \(\bm-\bu\leq\bu^{1/30}\), where the exponent
is at least \(c\bu^{1/30}\), and the second for
\(\bm-\bu>\bu^{1/30}\),
\[
 (1+\bu)^{1/2}p_{B^0_\bu}^{1/2}
 +\sum_{\bm\geq\bu}(1+\bm)^{1/2}p_{B_{\bu,\bm}}^{1/2}
 \leq C(1+\bu)^{C}e^{-c\bu^{1/30}} .
\]
Summing over \(\bu\in[K,\bn]\) and integrating in \(\theta\),
\begin{align*}
 \mathbb E\int_0^{2\pi}\one\{\mathscr H_K^\partial\}
 \one\Big\{\exists\,\bu\in[K,\bn]:
 |\varphi_{2^\bn}(r_\bu,\theta)-\varphi_{2^\bu}(1,\theta)|>\bu^{1/10}\Big\}&
 \mathscr D_{2^\bn}(1,\theta)d\theta\\
 &\leq C(1+K)\sum_{\bu\geq K}(1+\bu)^{C}e^{-c\bu^{1/30}},
\end{align*}
uniformly in \(\bn\geq K\).  The right side tends to zero as
\(K\to\infty\), and Markov's inequality proves
\eqref{eq:weighted-radial-spine}.

\emph{Proof of (iii).}
We bound each term of the sum separately, and fix \(\bm\) with
\(\bn\leq\bm<\bk\).

The event \(\mathscr R^\eta_\bm\) compares the field at the two radii
\(r_\bk\) and \(r_\bm\).  On \(\mathscr G_{\bn',\bn}\cap\mathscr G_\bn\), we may
replace \(\varphi_{2^\bm}\) by \(\zeta_{2^\bm}(\bn';\cdot)\) at both radii,
at a cost of at most \(1\), by \eqref{eq:base-change-linearization}.  We
then compare each radius with the boundary \(r=1\).  Hence
\(\mathscr R^\eta_\bm\) implies the event \(F_\bm\) that
\[
 |\zeta_{2^\bm}(\bn';r_\bl,\theta)-\zeta_{2^\bm}(\bn';1,\theta)|
 >\tfrac\eta2(\bm\log2)^{1/8}-2
 \qquad\text{for }\bl=\bk\text{ or }\bl=\bm .
\]

The event \(F_\bm\) is controlled by Corollary~\ref{cor:radialincrement}
with base \(\bn'\).  That corollary requires the radius exponent to be
at least \(\mathfrak c\bn'\).  Both radii satisfy this, since
\(\bk\geq\bm\geq\bn\geq\mathfrak c\bn'\).  It gives
\[
 \mathbb Q^\infty_{\bk,\bn'}\bigl(F_\bm\cap\widehat{\mathscr G}_{\bn'}\bigr)\leq Ce^{-c\bm^{1/8}} .
\]
Lemma~\ref{lem:suffix-CS} with \(\bn''=\bn'\), \(\bt=\bm+1\) and
\(F=F_\bm\) now bounds the \(\bm\)-th term.  Summing over \(\bm\), the
conditional expectation of the left side of
\eqref{eq:weighted-radial-mismatch} given \(\filt_{2^{\bn'}}\) is at most
\(\varrho_\bn\Upsilon_{\bk,\bn'}\), where
\[
 \varrho_\bn:=C\sum_{\bm\geq\bn}(1+\bm)\bigl(e^{-c\bm^{1/8}}+2^{-c\bm}\bigr)\longrightarrow0 .
\]
Lemma~\ref{lem:T-tight} then proves \eqref{eq:weighted-radial-mismatch}.

The choice of base \(\bn'\) is what makes this work.  With base \(\bn\),
the event \(\mathscr R^\eta_\bm\) for \(\bm\) close to \(\bn\) would be
essentially \(\filt_{2^\bn}\)-measurable, and its conditional probability
could not be small.

\end{proof}

\subsection{Evaluation of the radial/time barrier}

\begin{lemma}[Equivalence of $\mathscr B$ and $\mathscr E$ barriers]
\label{lem:integrated-barrier-comparison}
Let \(b<0\), \(L<\infty\), and suppose
\(\underline a=(a_\bj)\) decreases to \(a_\infty<0\), with
\(0\leq a_\bj-a_\infty\leq C(1+\bj)^{-1}\).
For every measurable \(0\leq f\leq1\):

\emph{(i) Initial deficit localization.}  Then
\begin{equation}\label{eq:integrated-initial-interface}
 \int_0^{2\pi}f(\theta)
 \one\{\mathscr B_{\bk,\bn}^{\bn,L,\underline a}\}
 \one\{(\mathscr E_{\bk,\bn}^{\bn,L,b})^\complement\}
 \mathscr D_{2^\bn}(r_\bk,\theta)d\theta
 \Prto[\bn,\bk]0 .
\end{equation}

\emph{(ii) Forward comparison.}  If \(a_\infty<2b\), then
\begin{align}
 \int_0^{2\pi}f(\theta)\mathbb E_\bn\left(
 \one\{\mathscr G_\bn\}
 \one\{\mathscr B_{\bk,\bn}^{L,\underline a}\}
 \one\{(\mathscr E_{\bk,\bn}^{L,b})^\complement\}
 \mathscr D_\infty(r_\bk,\theta)\right)d\theta
 \Prto[\bn,\bk]0 .
 \label{eq:integrated-change-of-barrier}
\end{align}

\emph{(iii) Reverse comparison.}  If \(a_\infty>b/2\), then
\begin{align}
 \int_0^{2\pi}f(\theta)\mathbb E_\bn\left(
 \one\{\mathscr G_\bn\}
 \one\{\mathscr B_{\bk,\bn}^{\bn,L,\underline a}\}
 \one\{\mathscr E_{\bk,\bn}^{L,b}\}
 \one\{(\mathscr B_{\bk,\bn}^{L,\underline a})^\complement\}
 \mathscr D_\infty(r_\bk,\theta)\right)d\theta
 \Prto[\bn,\bk]0 .
 \label{eq:integrated-reverse-barrier}
\end{align}
\end{lemma}

\begin{proof}
Write \(c_\beta=\sqrt{\beta/8}\).  We apply
Lemma~\ref{lem:boundary-spine-inputs} with its parameter \(\eta\) chosen
as follows.  For \emph{(i)}, when \(b>a_\infty\), and for \emph{(ii)},
choose \(\eta>0\) so that
\[
 2c_\beta\eta<b-a_\infty.
\]
For \emph{(iii)}, choose \(\eta>0\) so that
\[
 4c_\beta\eta<a_\infty-b;
\]
the factor \(4\) leaves room for the two field differences displayed in
the reverse comparison below and for
\(a_\bj-a_\infty=o(1)\).  No choice of \(\eta\) is needed for
\emph{(i)} when \(b\leq a_\infty\).

We first prove \emph{(i)}.  For fixed \(\bn\) and \(\bk\geq\bn\),
\(\mathscr B_{\bk,\bn}^{\bn,L,\underline a}\) no longer depends on
\(\bk\), while finite-degree continuity gives, uniformly in \(\theta\),
\[
 \varphi_{2^\bn}(r_\bk,\theta)\longrightarrow
 \varphi_{2^\bn}(1,\theta).
\]
Almost surely, \(\theta\mapsto\varphi_{2^\bn}(1,\theta)\) is a
nonconstant trigonometric polynomial, so equality in the limiting
barrier holds only on a Lebesgue-null set.  Dominated convergence
therefore permits \(\bk\to\infty\) in \emph{(i)}, and we may work with
\(\mathscr D_{2^\bn}(1,\theta)\).  If \(b>a_\infty\), the
constraint in \(\mathscr B_{\infty,\bn}^{\bn,L,\underline a}\) at
radius \(r_\bn\) and
\((\mathscr E_{\infty,\bn}^{\bn,L,b})^\complement\) imply
\[
 |\varphi_{2^\bn}(1,\theta)
   -\varphi_{2^\bn}(r_\bn,\theta)|
 >\eta(\bn\log2)^{1/8}
\]
for all large \(\bn\).  This contribution vanishes by
\eqref{eq:weighted-radial-spine}, taking \(\bu=\bn\).  If
\(b\leq a_\infty\), the event
\((\mathscr E_{\infty,\bn}^{\bn,L,b})^\complement\) directly implies
\[
 0<d_\bn^\partial(\theta)\leq C_{a_\infty,b,L}\bn^{1/8}.
\]
By \eqref{eq:shallow-gap-spine}, this contribution to the
derivative-chaos mass vanishes in the limit.  These two alternatives prove
\eqref{eq:integrated-initial-interface} for every \(a_\infty,b<0\).

We next prove \emph{(ii)}.  Because
\(\mathscr E_{\bk,\bn}^{L,b}\) fails, let
\begin{equation*}
 \bm:=\max\left\{\bj\geq\bn:
 \sqrt{\frac\beta8}\varphi_{2^\bj}(r_\bk,\theta)
 >\lg_{r_\bk}(2^\bj)
   +b\lg_{r_\bk}(2^\bj)^{1/8}+L\right\}.
\end{equation*}
For all large \(\bn\), this set contains no \(\bj\geq\bk\): after saturation,
the constraint from \(\mathscr B_{\bk,\bn}^{L,\underline a}\) at radius
\(r_\bk\) has the same leading term and a strictly smaller
\(1/8\)-coefficient.  Thus \(\bn\leq\bm<\bk\), and at this last failure
\[
 |\varphi_{2^\bm}(r_\bk,\theta)
   -\varphi_{2^\bm}(r_\bm,\theta)|
 >\eta(\bm\log2)^{1/8}.
\]
The preceding display is precisely the event \(\mathscr R_\bm^\eta\).
By the definition of \(\bm\), we may also work on
\(\mathscr E_{\bk,\bm+1}^{L,b}\).
We insert the indicator \(\one\{\mathscr G_{\bn',\bn}\}\).  Each
\(\theta\) has at most one last failure \(\bm\).  Hence the left side of
\eqref{eq:integrated-change-of-barrier} is at most the left side of
\eqref{eq:weighted-radial-mismatch}, and \eqref{eq:weighted-radial-mismatch}
proves \eqref{eq:integrated-change-of-barrier}.

We finally prove \emph{(iii)}.  Again we insert
\(\one\{\mathscr G_{\bn',\bn}\}\).  Fix \(\theta\) in the event of (iii).
Since \(\mathscr B_{\bk,\bn}^{\bn,L,\underline a}\) holds but
\(\mathscr B_{\bk,\bn}^{L,\underline a}\) fails, there is a first index
\(\bj>\bn\) at which \(\mathscr A_{\bk,\bj}^{L,a_\bj}\) fails.  Let \(\bu\)
be the smallest radius exponent at which it fails, so that
\begin{equation}\label{eq:reverse-failure}
 \sqrt{\tfrac\beta8}\,\varphi_{2^\bj}(r_\bu,\theta)
 >\bu\log2+a_\bj(\bu\log2)^{1/8}+L .
\end{equation}
All the events below lie in \(\mathscr E^{L,b}_{\bk,\bn}\), since that
indicator appears in (iii).  We fix \(K\) and split according to the size
of \(\bu\).

\emph{The case \(\bu\geq\bn\).}  We compare the failure
\eqref{eq:reverse-failure} with the constraint of
\(\mathscr E^{L,b}_{\bk,\bn}\) at time \(2^\bu\), which concerns the radius
\(r_\bk\).  The two thresholds differ by at least
\((a_\bj-b)(\bu\log2)^{1/8}\), and \(a_\bj-b\geq a_\infty-b>4c_\beta\eta\) by
the choice of \(\eta\).  So at least one of the following two field
differences is large:
\[
 |\varphi_{2^\bu}(r_\bu,\theta)-\varphi_{2^\bu}(r_\bk,\theta)|>\eta(\bu\log2)^{1/8},
 \qquad
 \sup_{\bv\geq\bu}|\varphi_{2^\bv}(r_\bu,\theta)-\varphi_{2^\bu}(r_\bu,\theta)|>\eta(\bu\log2)^{1/8}.
\]
The first event is \(\mathscr R^\eta_\bu\), or is empty if \(\bu=\bk\).  On it,
the physical barrier \(\mathscr E^{L,b}_{\bk,\bu+1}\) holds.  Its contribution
is therefore bounded by the left side of
\eqref{eq:weighted-radial-mismatch}.  The second event concerns the
increments at the radius \(r_\bu\) after time \(2^\bu\).
Corollary~\ref{cor:weighted-increments}, with \(\bj=\bl=\bu\) and
\(x=\eta(\bu\log2)^{1/8}\), bounds its contribution by
\(C\Upsilon_{\bk,\bn}(1+\sqrt{\bu-\bn})e^{-c\bu^{1/24}}\).  Summed over
\(\bu\geq\bn\), this has the form \(\varrho_\bn\Upsilon_{\bk,\bn}\) with
\(\varrho_\bn\to0\).

\emph{The case \(K\leq\bu<\bn\).}  We now compare the failure
\eqref{eq:reverse-failure} with the constraint of
\(\mathscr B_{\bk,\bn}^{\bn,L,\underline a}\) at time \(2^\bn\), at the same
radius \(r_\bu\).  Fix \(\alpha\in(1/8,1/6)\).  For \(K\) large, at least one
of the following three events holds:
\[
 \sup_{\bv\geq\bn}|\varphi_{2^\bv}(r_\bu,\theta)-\varphi_{2^\bn}(r_\bu,\theta)|>\bu^\alpha,\qquad
 |\varphi_{2^\bn}(r_\bu,\theta)-\varphi_{2^\bu}(1,\theta)|>\bu^{1/10},\qquad
 0<d_\bu^\partial(\theta)\leq C\bu^\alpha.
\]
The first event concerns the increments at the radius \(r_\bu\) after
time \(2^\bn\).  Corollary~\ref{cor:weighted-increments}, with \(\bj=\bn\),
\(\bl=\bu\) and \(x=\bu^\alpha\), bounds its contribution by
\(C\Upsilon_{\bk,\bn}e^{-c\bu^{\alpha/3}}\), up to an error that is
exponentially small in \(\bn\).  The sum over \(\bu\geq K\) is small when
\(K\) is large.

The second and third events are \(\filt_{2^\bn}\)-measurable.  For such
events the conditional expectation has a simple bound.  On
\(\mathscr G_{\bn,\bn}\cap\mathscr E^{\bn,L,b}_{\bk,\bn}\), the change of
measure \eqref{eq:com} writes it as \(\mathscr M_{\bk,\bn}\mathcal
U^{\varphi,L,b}_{\bk,\bn}\).  The upper bound in
Proposition~\ref{prop:brownian-clock-physical-sandwich} bounds
\(\mathcal U^{\varphi,L,b}_{\bk,\bn}\) by \(C\Delta^{[\bn]}_{\bk,\bn}\).
Moreover \(\Delta^{[\bn]}_{\bk,\bn}=\Delta^{\varphi}_{\bk,\bn}\geq1\) there, by
\eqref{eq:physical-initial-barrier-size}.  Hence
\begin{equation}\label{eq:reverse-weight}
 \mathbb E_\bn\bigl(\one\{\mathscr G_\bn\}\one\{\mathscr E^{L,b}_{\bk,\bn}\}\mathscr D_\infty(r_\bk,\theta)\bigr)
 \leq C\mathscr D_{2^\bn}(r_\bk,\theta).
\end{equation}
So the second and third events contribute at most
\(C\int_0^{2\pi}\one\{\cdot\}\mathscr D_{2^\bn}(r_\bk,\theta)\,d\theta\).  As
\(\bk\to\infty\), this converges to the same integral with
\(\mathscr D_{2^\bn}(1,\theta)\).  The boundary estimates
\eqref{eq:weighted-radial-spine} and \eqref{eq:shallow-gap-spine} show
that the limit is small in probability for \(K\) large, uniformly in
\(\bn\geq K\).  They are stated on the event \(\mathscr H^\partial_K\), whose
complement has small probability when \(K\) is large.

\emph{The case \(\bu<K\).}  Here the radius exponent \(\bu\) is bounded,
and the field at radius \(r_\bu\) barely moves after time \(2^\bn\).  A
barrier failure after time \(2^\bn\) is then possible only if the field at
time \(2^\bn\) is already very close to the barrier.  We show that the
chaos gives vanishing mass to this.

Define the gap to the barrier at time \(2^\bn\),
\[
 g_\bn(\theta):=\bu\log2+a_\bn(\bu\log2)^{1/8}+L
 -\sqrt{\tfrac\beta8}\,\varphi_{2^\bn}(r_\bu,\theta).
\]
The constraint of \(\mathscr B_{\bk,\bn}^{\bn,L,\underline a}\) at radius
\(r_\bu\) says that \(g_\bn\geq0\).  The failure \eqref{eq:reverse-failure},
together with \(a_\bn-a_\bj\leq C\bn^{-1}\), gives
\[
 g_\bn(\theta)<C_K\bn^{-1}+\sqrt{\tfrac\beta8}\,
 |\varphi_{2^\bj}(r_\bu,\theta)-\varphi_{2^\bn}(r_\bu,\theta)| .
\]
Put \(h:=2^{-\bn/12}\).  Then either the field at radius \(r_\bu\) moves by
more than \(h\) after time \(2^\bn\), or
\(0\leq g_\bn\leq\delta_\bn:=C_K\bn^{-1}+\sqrt{\beta/8}\,h\).
Corollary~\ref{cor:weighted-increments}, with \(\bj=\bn\), \(\bl=\bu\) and
\(x=h\), bounds the contribution of the first alternative by
\(C_K2^{-\bn/24}\Upsilon_{\bk,\bn}\).  The second alternative is
\(\filt_{2^\bn}\)-measurable.  By \eqref{eq:reverse-weight}, it contributes
at most
\(C\int_0^{2\pi}\one\{0\leq g_\bn\leq\delta_\bn\}\mathscr D_{2^\bn}(r_\bk,\theta)\,d\theta\).
As \(\bk\to\infty\), this converges to
\[
 V_\bn:=C\int_0^{2\pi}\one\{0\leq g_\bn(\theta)\leq\delta_\bn\}\,\mathscr D_{2^\bn}(1,\theta)\,d\theta .
\]

We claim that \(V_\bn\to0\) almost surely.  First,
\(g_\bn\to g_\infty:=\bu\log2+a_\infty(\bu\log2)^{1/8}+L-\sqrt{\beta/8}\,\varphi_\infty(r_\bu,\cdot)\)
uniformly, by Theorem~\ref{thm:fieldconverges}.  Also
\(\mathscr D_{2^\bn}(1,\theta)\,d\theta\to\mathscr D_\infty(d\theta)\) weakly,
by Theorem~\ref{thm:Bn}.  Hence, for every \(\delta>0\),
\[
 \limsup_{\bn\to\infty}V_\bn\leq C\mathscr D_\infty\bigl(\{|g_\infty|\leq\delta\}\bigr).
\]
Second, the set \(\{g_\infty=0\}\) has \(\mathscr D_\infty\)-measure zero.
Indeed, the Fourier coefficients of the stationary Gaussian field
\(\varphi_\infty(r_\bu,\cdot)\) are independent centered Gaussians with
variances \(\tfrac4\beta r_\bu^{2m}/m\).  So the field is almost surely a
nonconstant real-analytic function, and its level sets are finite.  The
measure \(\mathscr D_\infty\) has no atoms.  Letting \(\delta\downarrow0\)
proves the claim.

\emph{Conclusion.}  In each case, the contribution is bounded by the left
side of \eqref{eq:weighted-radial-mismatch}, by a quantity of the form
\(\varrho\,\Upsilon_{\bk,\bn}\) with \(\varrho\) small, or by one of the
limits just described.  By Lemma~\ref{lem:T-tight}, letting \(\bk\to\infty\),
then \(\bn\to\infty\), then \(K\to\infty\) proves
\eqref{eq:integrated-reverse-barrier}.
\end{proof}

\begin{lemma}[Integrated main estimate]\label{lem:integrated-mainestimate}
Fix \(b<0\), \(L<\infty\), a nonnegative continuous \(f\), and a
sequence \(\underline a\) as in
Lemma~\ref{lem:integrated-barrier-comparison}.
Then
\begin{align}
 \int_0^{2\pi}
 &\one\{\mathscr B_{\bk,\bn}^{\bn,L,\underline a}(\theta)\}f(\theta)
 \Bigg[
 \mathbb E_\bn\left(
  \one\{\mathscr G_\bn\}
  \one\{\mathscr E_{\bk,\bn}^{L,b}(\theta)\}
  \mathscr D_\infty(r_\bk,\theta)\right)
 -\frac{\mathscr M_{\bk,\bn}(\theta)}
 {2\pi\,\mathbb E\mathscr M_{\bk,\bn}(\theta)}
 \bigl(\Delta_{\bk,\bn}^{[\bn]}(\theta)\bigr)_+
 \Bigg]\,d\theta
 \Prto[\bn,\bk]0 .
 \label{eq:integrated-mainestimate}
\end{align}
\end{lemma}

\begin{proof}
The change-of-measure identity \eqref{eq:com} turns the first
conditional expectation into
\(\mathscr M_{\bk,\bn}\mathcal U_{\bk,\bn}^{\varphi,L,b}\), where
\(\mathcal U_{\bk,\bn}^{\varphi,L,b}\) is defined in
\eqref{eq:physical-global-good-operator}.  On
\(\mathscr G_{\bn,\bn}\) and the initial physical barrier
\(\mathscr E_{\bk,\bn}^{\bn,L,b}\),
Proposition~\ref{prop:brownian-clock-physical-sandwich} sandwiches
\(\mathcal U_{\bk,\bn}^{\varphi,L,b}\) between the two Brownian-barrier
expectations
\(\widetilde{\mathcal U}_{\bk,\bn}^{\mathcal P,L-K,b-\eta}\) and
\(\widetilde{\mathcal U}_{\bk,\bn}^{\mathcal P,L+K,b+\eta}\).  On
\(\mathscr G_{\bn,\bn}\) we have \(R_{2^\bn}\leq2\cdot2^{(C_{\rm in}+2)\bn}\)
by \eqref{eq:julia-derivative} and \eqref{eq:initial-derivative}, so on the initial deficit buffer
\(\mathscr V_{\bk,\bn}^{\nu,[\bn]}\) Proposition
\ref{prop:brownian-clock-full-tower} gives the same leading term
\(c_\bk\Delta_{\bk,\bn}^{[\bn]}\) for both, and
\eqref{eq:identify-step3-ck} identifies
\[
 c_\bk=\frac{1+o_{\bk,\bn}(1)}
 {2\pi\,\mathbb E\mathscr M_{\bk,\bn}}.
\]
Thus the integrand in \eqref{eq:integrated-mainestimate}, restricted to
these events, is bounded in absolute value by
\[
 o_{\bk,\bn}(1)\,
 \mathscr M_{\bk,\bn}
 \bigl(\Delta_{\bk,\bn}^{[\bn]}\bigr)_+ .
\]

It remains to bound the left-hand side of
\eqref{eq:integrated-mainestimate} on the complements of
\(\mathscr G_{\bn,\bn}\),
\(\mathscr V_{\bk,\bn}^{\nu,[\bn]}\), and
\(\mathscr E_{\bk,\bn}^{\bn,L,b}\).  For this purpose, the absolute value
of the expression in square brackets in
\eqref{eq:integrated-mainestimate} is bounded by the common majorant
\begin{align*}
 \mathcal R_{\bk,\bn}(\theta):={}&\mathbb E_\bn\left(
  \one\{\mathscr G_\bn\}
  \one\{\mathscr E_{\bk,\bn}^{L,b}(\theta)\}
  \mathscr D_\infty(r_\bk,\theta)\right)+
 \frac{\mathscr M_{\bk,\bn}(\theta)}
 {2\pi\,\mathbb E\mathscr M_{\bk,\bn}(\theta)}
 \bigl(\Delta_{\bk,\bn}^{[\bn]}(\theta)\bigr)_+.
\end{align*}
Thus each of the three bad-event removals reduces to controlling the two
nonnegative terms defining \(\mathcal R_{\bk,\bn}\).

The event
\(\mathscr G_{\bn,\bn}\) has probability tending to one.  Since
\(\mathscr G_\bn\subset\mathscr G_{\bn,\bn}\), the conditional-chaos
term in \(\mathcal R_{\bk,\bn}\) vanishes on
\((\mathscr G_{\bn,\bn})^\complement\).  By
\eqref{eq:identify-step3-ck} and \eqref{eq:step3-ck-bounds}, the normalized
target term is bounded by a constant multiple of
\(\mathscr D_{2^\bn}(r_\bk,\theta)\).  Its contribution on this complement
is therefore bounded by a constant multiple of
\[
 \one\{(\mathscr G_{\bn,\bn})^\complement\}
 \int_0^{2\pi}\mathscr M_{\bk,\bn}
 (\Delta_{\bk,\bn}^{[\bn]})_+\,d\theta
 =\one\{(\mathscr G_{\bn,\bn})^\complement\}
  \mathscr B_{2^\bn}(r_\bk).
\]
The masses \(\mathscr B_{2^\bn}(r_\bk)\) are tight by
Theorem~\ref{thm:Bn}, since \(\mathscr B_{2^\bn}(r_\bk)\to\mathscr B_{2^\bn}(1)\) as \(\bk\to\infty\) for fixed \(\bn\).  Since
\(\mathbb P((\mathscr G_{\bn,\bn})^\complement)\to0\), the last
display converges to zero in probability.

The initial deficit buffer loses no limiting chaos mass.  For fixed
\(\bn\), \eqref{eq:Btight} and the uniform convergence of the finite-degree
field in \(\theta\) as \(r_\bk\uparrow1\) imply that, for every
\(\varepsilon>0\),
\begin{equation}\label{eq:buffer-removal-integrated}
 \limsup_{\bn\to\infty}\limsup_{\bk\to\infty}
 \mathbb P\left(
 \int_0^{2\pi}
 \one\{(\mathscr V_{\bk,\bn}^{\nu,[\bn]})^\complement\}
 \mathscr D_{2^\bn}(r_\bk,\theta)d\theta>\varepsilon\right)
 \leq\varepsilon .
\end{equation}
The upper bound
\(\mathcal U_{\bk,\bn}^{\varphi,L,b}\leq C\Delta_{\bk,\bn}^{[\bn]}\) of
Proposition~\ref{prop:brownian-clock-physical-sandwich}, which does not require
the buffer, bounds the conditional-chaos term in
\(\mathcal R_{\bk,\bn}\) by
\(C\mathscr D_{2^\bn}(r_\bk,\theta)\).  Its normalized target term satisfies
\[
 \frac{\mathscr M_{\bk,\bn}}
 {2\pi\,\mathbb E\mathscr M_{\bk,\bn}}
 (\Delta_{\bk,\bn}^{[\bn]})_+
 =\frac{\mathscr D_{2^\bn}(r_\bk,\theta)}
 {2\pi\,\mathbb E\mathscr M_{\bk,\bn}}.
\]
By \eqref{eq:step3-ck-bounds}, this is also bounded by
\(C\mathscr D_{2^\bn}(r_\bk,\theta)\).  Thus
\eqref{eq:buffer-removal-integrated} removes both terms from the complement
of the initial deficit buffer.

Finally, on the complement of the initial physical barrier
\(\mathscr E_{\bk,\bn}^{\bn,L,b}\), the conditional-chaos term in
\(\mathcal R_{\bk,\bn}\) vanishes because
\(\mathscr E_{\bk,\bn}^{L,b}\subset
\mathscr E_{\bk,\bn}^{\bn,L,b}\).  The normalized target term is again
bounded by \(C\mathscr D_{2^\bn}(r_\bk,\theta)\).  Under the outer event
\(\mathscr B_{\bk,\bn}^{\bn,L,\underline a}\),
\eqref{eq:integrated-initial-interface} shows that this complement carries no
limiting chaos mass.  Adding the three complements proves
\eqref{eq:integrated-mainestimate}.
\end{proof}

\begin{proposition}[Radial barrier implies time barrier]\label{prop:upperbound}
Fix \(a<0\) and put
\[
 a_\bj=\frac a2+(\bj+10)^{-1}\downarrow\frac a2.
\]
For every measurable \(0\leq f\leq1\),
\begin{align}
 \int_0^{2\pi}f(\theta)\mathbb E_\bn\left(
 \one\{\mathscr A_\bk^{L,a}(\theta)\}
 \one\{(\mathscr B_{\bk,\bn}^{L,\underline a}(\theta))^\complement\}
 \one\{\mathscr G_\bn\}
 \mathscr D_\infty(r_\bk,\theta)\right)d\theta
 \Prto[\bn,\bk]0 .
 \label{eq:upper-radial-time-closure}
\end{align}
\end{proposition}

\begin{proof}
The proof follows the pattern of
Lemma~\ref{lem:integrated-barrier-comparison}.  We locate a barrier
failure, show that a field difference is large there, and bound the
resulting chaos mass by Corollary~\ref{cor:weighted-increments} and by
\eqref{eq:weighted-radial-mismatch}.  Put \(b:=a/4\) and
\(\eta:=|a|/(8\sqrt{\beta/8})\), and insert \(\one\{\mathscr G_{\bn',\bn}\}\).

\emph{The last failure.}  Since \(a<a/2=a_\infty\), Lemma~\ref{lem-1.20}
shows that on \(\mathscr A_\bk^{L,a}\) the constraints
\(\mathscr A^{L,a_\bj}_{\bk,\bj}\) fail for only finitely many \(\bj\).  So on
\(\mathscr A_\bk^{L,a}\cap(\mathscr B^{L,\underline a}_{\bk,\bn})^\complement\)
there is a last index \(J\geq\bn\) at which \(\mathscr A^{L,a_J}_{\bk,J}\)
fails.  Let \(\Lambda\leq\bk\wedge J\) be the smallest radius exponent at
which it fails.  Because \(J\) is the last failure, the radial/time
barrier \(\mathscr B^{L,\underline a}_{\bk,J+1}\) holds.

\emph{A large field difference.}  At the radius \(r_\Lambda\), we compare
the failure at time \(2^J\) with the Gaussian barrier
\(\mathscr A^{L,a}_\bk\), which constrains \(\varphi_\infty(r_\Lambda,\theta)\).
The two thresholds differ by \((a_J-a)(\Lambda\log2)^{1/8}\), and
\(a_J-a\geq|a|/2\).  Hence
\begin{equation}\label{eq:upper-D2}
 |\varphi_\infty(r_\Lambda,\theta)-\varphi_{2^J}(r_\Lambda,\theta)|>c\,\Lambda^{1/8},
\end{equation}
with \(c=c(a,\beta)>0\).  In particular, the field at radius \(r_\Lambda\)
moves by more than \(c\Lambda^{1/8}\) after time \(2^J\).

\emph{On the physical barrier.}  Suppose first that
\(\mathscr E^{L,b}_{\bk,J+1}\) holds.  Corollary~\ref{cor:weighted-increments},
with \(\bj=J\), \(\bl=\Lambda\) and \(x=c\Lambda^{1/8}\), bounds the
contribution of \(\{J=\bj,\Lambda=\bl\}\) by
\[
 C\Upsilon_{\bk,\bn}\bigl(1+\sqrt{\bj-\bn}\bigr)
 \Bigl(\min\bigl\{e^{-c\bl^{1/24}},\,2^{-(\bj-\bl)/12}\bigr\}+2^{-c\bj}\Bigr).
\]
The minimum is small when \(\bl\) is large or when \(\bj-\bl\) is large.
Since \(\bl\leq\bj\), one of the two holds when \(\bj\) is large.
Summing over \(\bl\leq\bj\) and \(\bj\geq\bn\), these contributions are at
most \(\varrho_\bn\Upsilon_{\bk,\bn}\) with \(\varrho_\bn\to0\).  By
Lemma~\ref{lem:T-tight}, they tend to zero in probability.

\emph{Off the physical barrier.}  Suppose now that
\(\mathscr E^{L,b}_{\bk,J+1}\) fails.  Since
\(\mathscr B^{L,\underline a}_{\bk,J+1}\) holds, the argument in the proof of
Lemma~\ref{lem:integrated-barrier-comparison}(ii) applies from time
\(2^{J+1}\).  That argument needs \(b>a_\bm\), and indeed
\(b-a_\bm\geq|a|/8\) for all \(\bm\geq\bn\) once \(\bn\) is large.  It shows
that the event \(\mathscr R^\eta_\bm\cap\mathscr E^{L,b}_{\bk,\bm+1}\) holds
at the last failure \(\bm\) of the physical barrier.  Since \(J\) is
determined by \(\theta\), these contributions are bounded by the left
side of \eqref{eq:weighted-radial-mismatch}, which tends to zero in
probability.
\end{proof}

\begin{lemma}[The chaos charges the early barrier]\label{lem:lowerbound}
Suppose \(\underline a=(a_\bj)\) decreases to \(a_\infty<0\) and
\(0\leq a_\bj-a_\infty\leq C(1+\bj)^{-1}\).  Then
\begin{equation}\label{eq:lower-boundary-closure}
 \int_0^{2\pi}
 \one\{\mathscr G_\bn\}
 \one\{(\mathscr B_{\infty,\bn}^{\bn,L,\underline a})^\complement\}
 \mathscr D_{2^\bn}(1,\theta)d\theta
 \Prto[L,\bn]0.
\end{equation}
\end{lemma}

\begin{proof}
Fix \(\alpha\in(1/8,1/6)\), for instance \(\alpha=1/7\), and use the
events \(\mathscr H_K^\partial\) of \eqref{eq:boundary-good-event}.  Outside the radial-mismatch event
\eqref{eq:weighted-radial-spine}, failure of
\(\mathscr B_{\infty,\bn}^{\bn,L,\underline a}\) at an index
\(K\leq\bu\leq\bn\) implies
\[
0<d_\bu^\partial(\theta)
 \leq C_{\underline a}\bu^{1/8}+C\bu^{1/10}-\sqrt2L
 \leq C_{\underline a}'\bu^\alpha .
\]
Its chaos mass vanishes by
\eqref{eq:shallow-gap-spine}, uniformly in \(\bn\).

For the finitely many \(\bu<K\), failure of the barrier and
\(d_\bu^\partial(\theta)>0\) instead force
\[
 \sqrt{\frac\beta4}\bigl(
 \varphi_{2^\bn}(r_\bu,\theta)-\varphi_{2^\bu}(1,\theta)\bigr)
 \geq \sqrt2L-C_K.
\]
Under the probability measures with density proportional to
\(\mathscr D_{2^\bn}(1,\theta)\,d\mathbb P\,d\theta\), the family
\[
 \left\{\varphi_{2^\bn}(r_\bu,\theta)-\varphi_{2^\bu}(1,\theta):
 \bn\geq K,\ \bu<K\right\}
\]
is tight.  Hence these radii contribute \(o_L(1)\).  Finally,
\(\mathbb P((\mathscr H_K^\partial)^\complement)\to0\).
Let first \(L\to\infty\), then \(K\to\infty\), to prove
\eqref{eq:lower-boundary-closure}.
\end{proof}

\begin{proof}[Proof of Proposition \ref{prop:main}]
Recall that, for every continuous \(f:\mathbb T\to[0,1]\) and every
\(a<0\), we must prove
\begin{equation}\label{eq:main-proof-target}
 \mathscr D_{2^\bn}(1,f)
 -e^{\mathfrak g_\beta}
 \mathbb E_\bn\bigl(\operatorname{CGMC}^{L,a}(f)\bigr)
 \Prto[L,\bn]0.
\end{equation}
For \(\bm\geq\bn\), denote the observed part of the radial/time barrier by
\begin{equation}\label{eq:Btruncated}
 \mathscr B_{\bk,\bn}^{\bm,L,\underline c}
 :=\bigcap_{\bj=\bn}^{\bm}\mathscr A_{\bk,\bj}^{L,c_\bj}.
\end{equation}
Fix \(f\), and let \(\underline c=(c_\bj)\) decrease to \(c_\infty<0\),
with \(c_\bj-c_\infty=O((1+\bj)^{-1})\).  We use the local notation
\begin{align}
 X_{\bk,\bn}^{L,a}(f)
 &:={\int_0^{2\pi}}f(\theta)\,
 \mathbb E_\bn\left(
  \one\{\mathscr G_\bn\}
  \one\{\mathscr A_\bk^{L,a}(\theta)\}
  \mathscr D_\infty(r_\bk,\theta)\right)\dif\theta,
 \label{eq:main-proof-X-definition}\\
 Y_{\bk,\bn}^{L,\underline c}(f)
 &:={\int_0^{2\pi}}f(\theta)\,
 \mathbb E_\bn\left(
  \one\{\mathscr G_\bn\}
  \one\{\mathscr B_{\bk,\bn}^{L,\underline c}(\theta)\}
  \mathscr D_\infty(r_\bk,\theta)\right)\dif\theta,
 \label{eq:main-proof-Y-definition}\\
 T_{\bk,\bn}^{L,\underline c}(f)
 &:={\int_0^{2\pi}}f(\theta)
 \one\{\mathscr B_{\bk,\bn}^{\bn,L,\underline c}(\theta)\}
 \frac{\mathscr M_{\bk,\bn}(\theta)}
 {2\pi\,\mathbb E\mathscr M_{\bk,\bn}(\theta)}
 \bigl(\Delta_{\bk,\bn}^{[\bn]}(\theta)\bigr)_+\dif\theta.
 \label{eq:main-proof-T-definition}
\end{align}
We proceed through four reductions.

\emph{Reduction 1: the Gaussian-chaos approximation.}
It suffices to prove
\begin{equation}\label{eq:main-gaussian-barrier-comparison}
 \mathscr D_{2^\bn}(1,f)
 -e^{\mathfrak g_\beta}X_{\bk,\bn}^{L,a}(f)
 \Prto[L,\bn,\bk]0.
\end{equation}
Indeed, Proposition~\ref{prop:criticalchaos} implies, in \(L^1\) as
\(\bk\to\infty\),
\[
 \int_0^{2\pi}f(\theta)
 \one\{\mathscr A_\bk^{L,a}(\theta)\}
 \mathscr D_\infty(r_\bk,\theta)d\theta
 \longrightarrow \operatorname{CGMC}^{L,a}(f).
\]
Conditional expectation is an \(L^1\)-contraction.  Moreover, the
indicator \(\one\{\mathscr G_\bn\}\) may be removed from
\(X_{\bk,\bn}^{L,a}(f)\): with
\(W_\bk:=\int_0^{2\pi}f(\theta)\one\{\mathscr A_\bk^{L,a}(\theta)\}
\mathscr D_\infty(r_\bk,\theta)\dif\theta\geq0\), the difference is
\(\mathbb E_\bn(\one\{\mathscr G_\bn^\complement\}W_\bk)\), whose
\(L^1\)-norm is \(\mathbb E(\one\{\mathscr G_\bn^\complement\}W_\bk)\).
By \eqref{eq-140924a} with \(p=3/2\), the family \((W_\bk)_\bk\) is
bounded in \(L^{3/2}\), so by H\"older's inequality this is at most
\(\mathbb P(\mathscr G_\bn^\complement)^{1/3}\sup_\bk\|W_\bk\|_{3/2}\),
which tends to zero as \(\bn\to\infty\) uniformly in \(\bk\).
Consequently, \eqref{eq:main-gaussian-barrier-comparison} implies
\eqref{eq:main-proof-target}.

\emph{Reduction 2: evaluation of the radial/time barrier.}
We claim that, for every sequence \(\underline c\) as above,
\begin{equation}\label{eq:main-radial-time-evaluation}
 Y_{\bk,\bn}^{L,\underline c}(f)
 -T_{\bk,\bn}^{L,\underline c}(f)
 \Prto[\bn,\bk]0.
\end{equation}
Choose \(b_+<0\) and \(b_-<0\) so that
\[
 c_\infty<2b_+,
 \qquad
 c_\infty>\frac{b_-}{2}.
\]
For example, one may take \(b_+=c_\infty/3\) and
\(b_-=3c_\infty\).  Since the full radial/time barrier implies its
observed part, the forward and reverse estimates of
Lemma~\ref{lem:integrated-barrier-comparison} give
\begin{align*}
 Y_{\bk,\bn}^{L,\underline c}(f)
 &\leq
 \int_0^{2\pi}f(\theta)\mathbb E_\bn\left(
  \one\{\mathscr G_\bn\}
  \one\{\mathscr B_{\bk,\bn}^{\bn,L,\underline c}\}
  \one\{\mathscr E_{\bk,\bn}^{L,b_+}\}
  \mathscr D_\infty(r_\bk,\theta)\right)d\theta+o_{\bn,\bk}(1),\\
 Y_{\bk,\bn}^{L,\underline c}(f)
 &\geq
 \int_0^{2\pi}f(\theta)\mathbb E_\bn\left(
  \one\{\mathscr G_\bn\}
  \one\{\mathscr B_{\bk,\bn}^{\bn,L,\underline c}\}
  \one\{\mathscr E_{\bk,\bn}^{L,b_-}\}
  \mathscr D_\infty(r_\bk,\theta)\right)d\theta-o_{\bn,\bk}(1),
\end{align*}
where both errors converge to zero in probability in the indicated
iterated limit.  Applying Lemma~\ref{lem:integrated-mainestimate} to
the two conditional expectations yields the common value
\(T_{\bk,\bn}^{L,\underline c}(f)\), proving
\eqref{eq:main-radial-time-evaluation}.

\emph{Reduction 3: radial/time squeeze.}
Set
\[
 c_\bj^+=\frac a2+(\bj+10)^{-1},
 \qquad
 c_\bj^-=a+(\bj+10)^{-1}.
\]
We claim that there are errors
\(\varepsilon_{\bk,\bn}^{\pm}\Prto[\bn,\bk]0\) such that
\begin{equation}\label{eq:main-radial-time-squeeze}
 Y_{\bk,\bn}^{L,\underline c^-}(f)
 -\varepsilon_{\bk,\bn}^-
 \leq X_{\bk,\bn}^{L,a}(f)
 \leq Y_{\bk,\bn}^{L,\underline c^+}(f)
 +\varepsilon_{\bk,\bn}^+.
\end{equation}
Proposition~\ref{prop:upperbound} gives an error
\(\varepsilon_{\bk,\bn}^+\Prto[\bn,\bk]0\) such that
\[
 X_{\bk,\bn}^{L,a}(f)
 \leq Y_{\bk,\bn}^{L,\underline c^+}(f)
 +\varepsilon_{\bk,\bn}^+.
\]
The second inclusion in Lemma~\ref{lem-1.20} gives
\(\mathscr B_{\bk,\bn}^{L,\underline c^-}
\subseteq\mathscr A_\bk^{L,a}\), and hence
\[
 Y_{\bk,\bn}^{L,\underline c^-}(f)
 \leq X_{\bk,\bn}^{L,a}(f).
\]
These inequalities prove \eqref{eq:main-radial-time-squeeze}, with
\(\varepsilon_{\bk,\bn}^-=0\).

\emph{Reduction 4: boundary normalization.}
For every sequence \(\underline c\) as above, we claim that
\begin{equation}\label{eq:main-normalized-spine-limit}
 e^{\mathfrak g_\beta}T_{\bk,\bn}^{L,\underline c}(f)
 -\mathscr D_{2^\bn}(1,f)
 \Prto[L,\bn,\bk]0.
\end{equation}
For fixed \(\bn\), finite-degree continuity and dominated convergence allow
\(\bk\to\infty\) in the definition of \(T\).  By
\cite[Lemma~9.2]{PaquetteZeitouni03}, with the normalization
\eqref{eq:Mn},
\begin{equation}\label{eq:main-normalization-limit}
 \mathbb E\mathscr M_{\infty,\bn}(\theta)
 \longrightarrow \frac{e^{\mathfrak g_\beta}}{2\pi}.
\end{equation}
Moreover,
\(\mathscr M_{\infty,\bn}
(\Delta_{\infty,\bn}^{[\bn]})_+
=\mathscr D_{2^\bn}(1,\theta)\).
Lemma~\ref{lem:lowerbound} removes the complement of the observed
radial/time barrier.  Its factor \(\one\{\mathscr G_\bn\}\) may be
removed because, for the nonnegative omitted mass \(V_{\bn,L}\),
\[
 \mathbb P\left(
 V_{\bn,L}-\one\{\mathscr G_\bn\}V_{\bn,L}>\varepsilon\right)
 \leq\mathbb P(\mathscr G_\bn^\complement)\longrightarrow0.
\]
This proves \eqref{eq:main-normalized-spine-limit}.

Applying Reductions 2 and 4 to both sides of the squeeze in Reduction 3
proves \eqref{eq:main-gaussian-barrier-comparison}.  Reduction 1 now gives
\eqref{eq:main-proof-target}, completing the proof.
\end{proof}

\begin{appendix}
\section{Martingale inequalities}
We recall Freedman's inequality \cite[Theorem 1.6]{Freedman}.
\begin{lemma}[Freedman's inequality]
\label{lem:freedman}
Let $M_n$ be a martingale adapted to a filtration $\filt_n$ satisfying the following. For all $n\geq 1$
\begin{equation}
\label{eq-23012026aba}
|M_{n+1}-M_{n}|\leq 1.
\end{equation}
Let $\langle M\rangle_n \coloneqq \sum_{k=1}^n \mathbb{E}((M_{k+1}-M_{k})^2\,|\,\filt_k)$.
Then for all $x,y \geq 0$,
\[
\Pr[\exists~n : \{|M_n - M_0| > x\} \cap \{\langle M\rangle_n \leq y\} ] \leq 2e^{-x^2/2(x+y)}.
\]
\end{lemma}
The following version will be useful in various arguments in the paper.
\begin{lemma}
\label{lem:appendixmartingale}
Let $M_n$ be a martingale adapted to a filtration $\filt_n$ satisfying the following. There exists a constant $A$ so that, for all $n\geq 1$
\begin{equation}
\label{eq-23012026a}
|M_{n+1}-M_{n}|\leq A.
\end{equation}
and
\begin{equation}
\label{eq-23012026b}
E((M_{n+1}-M_{n})^2\,|\,\filt_n)\leq \frac{A}{n}.
\end{equation}
Then, for any $\alpha,\epsilon>0$, there exists $\delta=\delta(\alpha,\epsilon,A)$ so that, for any $T,u>0$, 
\begin{equation}
\label{eq-23012026c}
\Pr(\exists 0<j<\ell\leq T: M_\ell-M_j\geq \alpha(\log \ell-\log j)+\epsilon(\log T-\log \ell)+u)\leq 2e^{-\delta u}.
\end{equation}
\end{lemma}
\begin{proof}
We may and will assume that $T=2^{\bT}$.
We compare \(M\) at dyadic times and control its oscillation
within each dyadic block.  Comparing \(\log\ell-\log j\) with the
difference of the dyadic indices of \(\ell\) and \(j\), it suffices to
prove the dyadic estimates below with \(\alpha\) replaced by
\(\alpha_0=\alpha\log2/4\); the factor \(1/4\) absorbs the two blocks
containing \(j\) and \(\ell\).  We continue to write \(\alpha\) for
\(\alpha_0\).
We first claim that there exists $\delta_1=\delta_1(\epsilon,A)$ so that
\begin{equation}
\label{eq-230126aa}
\Pr[\max_{s\in [2^{\bl-1},2^{\bl}]} |M_{2^{\bl}}-M_s|\geq \epsilon |\bT-\bl|/8+u/8]\leq 2e^{-\delta_1 [(\bT-\bl)+u]}.
\end{equation}
Indeed, \eqref{eq-230126aa} is just (a weak form of) Freedman's inequality.
It follows that 
\begin{equation}
\label{eq-230126aab}
\Pr[\exists \bl\in \{1,\dots,\bT\}:  \max_{s\in [2^\bl,2^{\bl+1}]} |M_{2^{\bl+1}}-M_s|\geq \epsilon |\bT-\bl|/8+u/8]\leq 2e^{-\delta_1 u}.
\end{equation}
Therefore, \eqref{eq-23012026c}  follows if we show that there exists $\delta_2=\delta_2(\epsilon,A,\alpha)$ so that
\begin{equation}
\label{eq-230126bb}
\Pr[\exists \bj<\bl\in \{1,\dots,\bT\}:   |M_{2^{\bl}}-M_{2^{\bj}}|\geq \alpha(\bl-\bj)+\epsilon |\bT-\bl|/8+u/8]\leq 2e^{-\delta_2 u}.
\end{equation}
Another application of Freedman's inequality gives that for fixed $\bj<\bl\in \{1,\dots,\bT\}$, for some constant $c_1=c_1(A)$,
\begin{equation}
\label{eq-230126cc}
\Pr[ |M_{2^{\bl}}-M_{2^{\bj}}|\geq \alpha(\bl-\bj)+\epsilon |\bT-\bl|/8+u/8]\leq 2e^{-c_1[ \alpha (\bl-\bj)+\epsilon(\bT-\bl)/8+u/8]},
\end{equation}
Summing first over \(\bj<\bl\) and then over \(\bl\) gives
\(Ce^{-c u}\).  After decreasing \(\delta_2\), the bound
\(Ce^{-cu}\wedge1\leq2e^{-\delta_2u}\) yields
\eqref{eq-230126bb}.
\end{proof}

The following elementary consequence of optional stopping shows that a
martingale whose bracket grows must cross below any level that is small
compared with the square root of the bracket.
\begin{lemma}[A martingale with a large bracket crosses a level]
\label{lem:appendix-crossing}
Let \((M_j)_{j_1\leq j\leq j_2}\) be a martingale with respect to \((\filt_j)\) with \(M_{j_1}=0\) and \(\mathbb E M_{j_2}^2<\infty\), and for some \(\rho>0\) put
\[
 \langle M\rangle_j=\sum_{j_1\leq i<j}\mathbb E\bigl((M_{i+1}-M_i)^2\mid\filt_i\bigr),
 \qquad
 \kappa=\sum_{j_1\leq i<j_2}\mathbb E|M_{i+1}-M_i|^{2+\rho}.
\]
Then for all \(x,V>0\) ,
\begin{equation}\label{eq:appendix-crossing}
 \Pr\bigl(M_j>-x\text{ for all }j_1\leq j\leq j_2,\ \langle M\rangle_{j_2}\geq V\bigr) \leq
\frac32 4^{1/3} \left(\frac{x+\kappa^{1/(2+\rho)}}{\sqrt{V}}\right)^{2/3}
+\frac{2\kappa^{2/(2+\rho)}}{V}.
\end{equation}
\end{lemma}
\begin{proof}
We may and will assume that $x+\kappa^{1/(2+\rho)}\leq \sqrt{V}/2$, otherwise there is nothing to prove since the right side of \eqref{eq:appendix-crossing} exceeds one in that case. Fix $\varepsilon\in (0,1)$ so that
 \(y=\varepsilon\sqrt V\geq x\). 
Let \(T\) be the first \(j\in[j_1,j_2]\) with \(M_j\notin(-x,y)\), and \(T=j_2\) if there is no such \(j\).  Let \(E_0\) be the event that \(M_j\in(-x,y)\) for all \(j_1\leq j\leq j_2\), and on \(E_0^{\complement}\) let \(E_+=\{M_T\geq y\}\) and \(E_-=\{M_T\leq-x\}\).  The event on the left side of \eqref{eq:appendix-crossing} is contained in \(E_0\cup E_+\).

Let \(\Delta^*=\max_{j_1\leq i<j_2}|M_{i+1}-M_i|\).  Since \((\Delta^*)^{2+\rho}\leq\sum_{j_1\leq i<j_2}|M_{i+1}-M_i|^{2+\rho}\), Jensen's inequality gives
\[
 \mathbb E\Delta^*\leq\kappa^{1/(2+\rho)},
 \qquad
 \mathbb E(\Delta^*)^2\leq\kappa^{2/(2+\rho)}.
\]
As \(M_{j_1}=0\in(-x,y)\), we have \(T>j_1\) and \(M_{T-1}\in(-x,y)\), hence \(|M_T|\leq y+\Delta^*\) and
\[
 \mathbb E M_T^2\leq2y^2+2\kappa^{2/(2+\rho)}.
\]
Since \(T\) is a bounded stopping time and \(M^2-\langle M\rangle\) is a martingale, \(\mathbb E\langle M\rangle_T=\mathbb E M_T^2\).  On \(E_0\cap\{\langle M\rangle_{j_2}\geq V\}\) we have \(T=j_2\) and \(\langle M\rangle_T\geq V\), so
\[
 V\Pr\bigl(E_0\cap\{\langle M\rangle_{j_2}\geq V\}\bigr)
 \leq\mathbb E\langle M\rangle_T=\mathbb E M_T^2\leq2\varepsilon^2V+2\kappa^{2/(2+\rho)}.
\]
Next, \(M_T\geq y\) on \(E_+\), \(M_T\geq-x-\Delta^*\) on \(E_-\) and \(M_T>-x\) on \(E_0\), so optional stopping gives
\[
 0=\mathbb EM_T\geq y\Pr(E_+)-x\Pr(E_-\cup E_0)-\mathbb E\Delta^*
 \geq\varepsilon\sqrt V\Pr(E_+)-x-\kappa^{1/(2+\rho)}.
\]
Combining the last two displays gives
\[\Pr\bigl(M_j>-x\text{ for all }j_1\leq j\leq j_2,\ \langle M\rangle_{j_2}\geq V\bigr)
 \leq2\varepsilon^2+\frac{x+\kappa^{1/(2+\rho)}}{\varepsilon\sqrt V}
+\frac{2\kappa^{2/(2+\rho)}}{V}.
\]
Optimizing over $\varepsilon$ gives $\varepsilon^3=(x+\kappa^{1/(2+\rho)})/4\sqrt{V}$, which satisfies $y\geq x$ on the assumption 
$x+\kappa^{1/(2+\rho)}\leq\sqrt{V}/2$, since then $(x/\sqrt V)^3\leq(x+\kappa^{1/(2+\rho)})/(4\sqrt V)=\varepsilon^3$, and substituting this value of $\varepsilon$ 
proves \eqref{eq:appendix-crossing}.
\end{proof}

\end{appendix}

\bibliographystyle{alpha}
\bibliography{conv}

\newcommand{\etalchar}[1]{$^{#1}$}
\begin{thebibliography}{DRSV14b}

\bibitem[ABB17]{ABB}
Louis-Pierre Arguin, David Belius, and Paul Bourgade.
\newblock Maximum of the characteristic polynomial of random unitary matrices.
\newblock {\em Communications in Mathematical Physics}, 349(2):703--751, 2017.

\bibitem[ABB{\etalchar{+}}19]{ABBRS}
Louis-Pierre Arguin, David Belius, Paul Bourgade, Maksym Radziwi\l\l, and
  Kannan Soundararajan.
\newblock Maximum of the {R}iemann zeta function on a short interval of the
  critical line.
\newblock {\em Comm. Pure Appl. Math.}, 72(3):500--535, 2019.

\bibitem[ABH17]{ABH17}
Louis-Pierre Arguin, David Belius, and Adam~J. Harper.
\newblock Maxima of a randomized {R}iemann zeta function, and branching random
  walks.
\newblock {\em Ann. Appl. Probab.}, 27(1):178--215, 2017.

\bibitem[ABR20]{ABRI}
Louis-Pierre Arguin, Paul Bourgade, and Maksym Radziwi\l\l.
\newblock The {F}yodorov--{H}iary--{K}eating conjecture. {I}, 2020.
\newblock arXiv:2007.00988.

\bibitem[ABR23]{ABRII}
Louis-Pierre Arguin, Paul Bourgade, and Maksym Radziwi\l\l.
\newblock The {F}yodorov--{H}iary--{K}eating conjecture. {II}, 2023.
\newblock arXiv:2307.00982.

\bibitem[ABZ23]{ABZ}
Fanny Augeri, Rapha{\"e}l Butez, and Ofer Zeitouni.
\newblock A {CLT} for the characteristic polynomial of random {J}acobi
  matrices, and the {G$\beta$E}.
\newblock {\em Probab. Theory Related Fields}, 186(1-2):1--89, 2023.

\bibitem[AH26]{ArguinHamdan}
Louis-Pierre Arguin and Jad Hamdan.
\newblock The {F}yodorov--{H}iary--{K}eating conjecture on mesoscopic
  intervals.
\newblock {\em Ann. Probab.}, 54(5):2173--2215, 2026.

\bibitem[A{\"{\i}}d13]{Aidekon}
Elie A{\"{\i}}d{\'e}kon.
\newblock Convergence in law of the minimum of a branching random walk.
\newblock {\em Ann. Probab.}, 41(3A):1362--1426, 2013.

\bibitem[AMY20]{AglerMcCarthyYoung}
Jim Agler, John~E. McCarthy, and Nicholas~J. Young.
\newblock {\em Operator Analysis: Hilbert Space Methods in Complex Analysis},
  volume 219 of {\em Cambridge Tracts in Mathematics}.
\newblock Cambridge University Press, Cambridge, 2020.

\bibitem[AOR21]{AOR21}
Louis-Pierre Arguin, Fr\'{e}d\'{e}ric Ouimet, and Maksym Radziwi\l\l.
\newblock Moments of the {R}iemann zeta function on short intervals of the
  critical line.
\newblock {\em Ann. Probab.}, 49(6):3106--3141, 2021.

\bibitem[APS19]{APS19}
Juhan Aru, Ellen Powell, and Avelio Sep\'{u}lveda.
\newblock Critical {L}iouville measure as a limit of subcritical measures.
\newblock {\em Electron. Commun. Probab.}, 24:1--16, 2019.

\bibitem[AS14]{AidekonShi14}
Elie A\"{\i}d\'{e}kon and Zhan Shi.
\newblock The {S}eneta-{H}eyde scaling for the branching random walk.
\newblock {\em Ann. Probab.}, 42(3):959--993, 2014.

\bibitem[AZ25]{AugeriZeitouni25}
Fanny Augeri and Ofer Zeitouni.
\newblock Maximum of the characteristic polynomial of random {J}acobi matrices,
  2025.
\newblock arXiv:2512.13289.

\bibitem[BDHK25]{BDHK26}
Paul Bourgade, Guillaume Dubach, Lisa Hartung, and Ahmet Keles.
\newblock Fisher-hartwig asymptotics for non-hermitian random matrices, 2025.
\newblock arXiv:251209123.

\bibitem[BDZ16]{BDZ}
Maury Bramson, Jian Ding, and Ofer Zeitouni.
\newblock Convergence in law of the maximum of the two-dimensional discrete
  {G}aussian free field.
\newblock {\em Comm. Pure Appl. Math.}, 69(1):62--123, 2016.

\bibitem[BF22]{BourgadeFalconet}
Paul Bourgade and Hugo Falconet.
\newblock Liouville quantum gravity from random matrix dynamics, 2022.
\newblock arXiv:2206.03029.

\bibitem[BK97]{BK97}
J.~D. Biggins and A.~E. Kyprianou.
\newblock Seneta-{H}eyde norming in the branching random walk.
\newblock {\em Ann. Probab.}, 25(1):337--360, 1997.

\bibitem[BKN{\etalchar{+}}14]{BKNSW14}
Julien Barral, Antti Kupiainen, Miika Nikula, Eero Saksman, and Christian Webb.
\newblock Critical {M}andelbrot cascades.
\newblock {\em Comm. Math. Phys.}, 325(2):685--711, 2014.

\bibitem[BKN{\etalchar{+}}15]{BKNSW15}
Julien Barral, Antti Kupiainen, Miika Nikula, Eero Saksman, and Christian Webb.
\newblock Basic properties of critical lognormal multiplicative chaos.
\newblock {\em Ann. Probab.}, 43(5):2205--2249, 2015.

\bibitem[BL18]{BL}
Marek Biskup and Oren Louidor.
\newblock Full extremal process, cluster law and freezing for the
  two-dimensional discrete {G}aussian free field.
\newblock {\em Adv. Math.}, 330:589--687, 2018.

\bibitem[BLZ25]{BourgadeLopattoZeitouni}
Paul Bourgade, Patrick Lopatto, and Ofer Zeitouni.
\newblock Optimal rigidity and maximum of the characteristic polynomial of
  {W}igner matrices.
\newblock {\em Geom. Funct. Anal.}, 35(1):161--253, 2025.

\bibitem[BMP22]{BMP}
Paul Bourgade, Krishnan Mody, and Michel Pain.
\newblock Optimal local law and central limit theorem for $\beta$-ensembles.
\newblock {\em Comm. Math. Phys.}, 390:1017--1079, 2022.

\bibitem[BP25]{BerestyckiPowell}
Nathana\"{e}l Berestycki and Ellen Powell.
\newblock {\em Gaussian free field and {L}iouville quantum gravity}.
\newblock Cambridge University Press, Cambridge, 2025.

\bibitem[BWW18]{BWW18}
Nathana{\"e}l Berestycki, Christian Webb, and Mo~Dick Wong.
\newblock Random {H}ermitian matrices and {G}aussian multiplicative chaos.
\newblock {\em Probab. Theory Related Fields}, 172(1-2):103--189, 2018.

\bibitem[CFG26]{CaiFangGuo26}
Yin Cai, Xiang Fang, and Feng Guo.
\newblock Critical characteristic-polynomial chaos for complex iid matrices,
  2026.
\newblock arXiv:2611.01794.

\bibitem[CFLW21]{CFLW}
Tom Claeys, Benjamin Fahs, Gaultier Lambert, and Christian Webb.
\newblock How much can the eigenvalues of a random {H}ermitian matrix
  fluctuate?
\newblock {\em Duke Math. J.}, 170(9):2085--2235, 2021.

\bibitem[CL25]{CL25}
Giorgio Cipolloni and Benjamin Landon.
\newblock Maximum of the characteristic polynomial of i.i.d. matrices.
\newblock {\em Comm. Pure Appl. Math.}, 78(9):1703--1782, 2025.

\bibitem[CMN18]{CMN}
Reda Chhaibi, Thomas Madaule, and Joseph Najnudel.
\newblock On the maximum of the {C$\beta$E} field.
\newblock {\em Duke Mathematical Journal}, 167(12):2243--2345, 2018.

\bibitem[CN26]{CN}
Reda Chhaibi and Joseph Najnudel.
\newblock On the circle, {G}aussian multiplicative chaos and {B}eta ensembles
  match exactly.
\newblock {\em J. Eur. Math. Soc. (JEMS)}, 28(1):1--77, 2026.

\bibitem[CZ20]{CookZeitouni20}
Nicholas Cook and Ofer Zeitouni.
\newblock Maximum of the characteristic polynomial for a random permutation
  matrix.
\newblock {\em Comm. Pure Appl. Math.}, 73(8):1660--1731, 2020.

\bibitem[DRSV14a]{DRSV}
Bertrand Duplantier, R\'{e}mi Rhodes, Scott Sheffield, and Vincent Vargas.
\newblock Critical {G}aussian multiplicative chaos: convergence of the
  derivative martingale.
\newblock {\em Ann. Probab.}, 42(5):1769--1808, 2014.

\bibitem[DRSV14b]{DRSV2}
Bertrand Duplantier, R\'{e}mi Rhodes, Scott Sheffield, and Vincent Vargas.
\newblock Renormalization of critical {G}aussian multiplicative chaos and {KPZ}
  relation.
\newblock {\em Comm. Math. Phys.}, 330(1):283--330, 2014.

\bibitem[DRZ17]{DRZ17}
Jian Ding, Rishideep Roy, and Ofer Zeitouni.
\newblock Convergence of the centered maximum of log-correlated {G}aussian
  fields.
\newblock {\em Ann. Probab.}, 45(6A):3886--3928, 2017.

\bibitem[FB08]{FBo}
Yan~V. Fyodorov and Jean-Philippe Bouchaud.
\newblock Freezing and extreme-value statistics in a random energy model with
  logarithmically correlated potential.
\newblock {\em J. Phys. A}, 41(37):372001, 12, 2008.

\bibitem[FD16]{FLD}
Yan~V. Fyodorov and Pierre~Le Doussal.
\newblock Moments of the position of the maximum for {GUE} characteristic
  polynomials and for log-correlated {G}aussian processes.
\newblock {\em J. Stat. Phys.}, 164(1):190--240, 2016.
\newblock With Appendix I by Alexei Borodin and Vadim Gorin.

\bibitem[FGK18]{FGK18}
Yan~V. Fyodorov, Sven Gnutzmann, and Jonathan~P. Keating.
\newblock Extreme values of {CUE} characteristic polynomials: a numerical
  study.
\newblock {\em J. Phys. A}, 51(46):464001, 2018.

\bibitem[FHK12]{FHK12}
Yan~V. Fyodorov, Ghaith~A. Hiary, and Jonathan~P. Keating.
\newblock Freezing {T}ransition, {C}haracteristic {P}olynomials of {R}andom
  {M}atrices, and the {R}iemann {Z}eta-{F}unction.
\newblock {\em Phys. Rev. Lett.}, 108:170601, 2012.

\bibitem[FK14]{FK}
Yan~V. Fyodorov and Jonathan~P. Keating.
\newblock Freezing transitions and extreme values: random matrix theory, and
  disordered landscapes.
\newblock {\em Philos. Trans. R. Soc. Lond. Ser. A Math. Phys. Eng. Sci.},
  372(2007):20120503, 32, 2014.

\bibitem[FK21]{ForkelKeating21}
Johannes Forkel and Jonathan~P. Keating.
\newblock The classical compact groups and {G}aussian multiplicative chaos.
\newblock {\em Nonlinearity}, 34(9):6050--6119, 2021.

\bibitem[FLDR09]{FLDR09}
Yan~V. Fyodorov, Pierre Le~Doussal, and Alberto Rosso.
\newblock Statistical mechanics of logarithmic {REM}: duality, freezing and
  extreme value statistics of {$1/f$} noises generated by {G}aussian free
  fields.
\newblock {\em J. Stat. Mech. Theory Exp.}, (10):P10005, 2009.

\bibitem[Fre75]{Freedman}
David~A. Freedman.
\newblock On tail probabilities for martingales.
\newblock {\em Ann. Probability}, 3:100--118, 1975.

\bibitem[FS16]{FyodorovSimm16}
Yan~V. Fyodorov and Nick~J. Simm.
\newblock On the distribution of the maximum value of the characteristic
  polynomial of {GUE} random matrices.
\newblock {\em Nonlinearity}, 29(9):2837--2855, 2016.

\bibitem[GW25]{GW25b}
Ofir Gorodetsky and Mo~Dick Wong.
\newblock On the limiting distribution of sums of random multiplicative
  functions, 2025.
\newblock arXiv:2508.12956.

\bibitem[GZ26]{GuZhang26}
Haotian Gu and Zhenyuan Zhang.
\newblock Universality and phase transitions in low moments of secular
  coefficients of critical holomorphic multiplicative chaos.
\newblock {\em Probab. Theory Related Fields}, 194(1-2):613--716, 2026.

\bibitem[Har19]{Harper19}
Adam~J. Harper.
\newblock On the partition function of the {R}iemann zeta function, and the
  {F}yodorov--{H}iary--{K}eating conjecture, 2019.
\newblock arXiv:1906.05783.

\bibitem[Har20a]{Harper20}
Adam~J. Harper.
\newblock Moments of random multiplicative functions, {I}: {L}ow moments,
  better than squareroot cancellation, and critical multiplicative chaos.
\newblock {\em Forum Math. Pi}, 8:Paper No. e1, 2020.

\bibitem[Har20b]{HarperBourbaki}
Adam~J. Harper.
\newblock The {R}iemann zeta function in short intervals [after {N}ajnudel, and
  {A}rguin, {B}elius, {B}ourgade, {R}adziwi\l\l, and {S}oundararajan].
\newblock {\em Ast\'{e}risque}, 422:391--414, 2020.
\newblock S\'{e}minaire Bourbaki, Expos\'{e} 1161.

\bibitem[HKK26]{HKK26}
Haoran Hu, Yujin~H. Kim, and Xaver Kriechbaum.
\newblock Multiplicative chaos from random multiplicative functions, 2026.
\newblock arXiv:2608.13687.

\bibitem[Jeg21]{Jego21}
Antoine Jego.
\newblock Critical {B}rownian multiplicative chaos.
\newblock {\em Probab. Theory Related Fields}, 180(1-2):495--552, 2021.

\bibitem[JM15]{JM15}
Tiefeng Jiang and Sho Matsumoto.
\newblock Moments of traces of circular beta-ensembles.
\newblock {\em Ann. Probab.}, 43(6):3279--3336, 2015.

\bibitem[JS17]{JunnilaSaksman}
Janne Junnila and Eero Saksman.
\newblock Uniqueness of critical {G}aussian chaos.
\newblock {\em Electron. J. Probab.}, 22(11):1--31, 2017.

\bibitem[Jun20]{Junnila20}
Janne Junnila.
\newblock On the multiplicative chaos of non-{G}aussian log-correlated fields.
\newblock {\em Int. Math. Res. Not. IMRN}, (19):6169--6196, 2020.

\bibitem[Kel25]{Kel26}
Ahmet Keles.
\newblock Non-intersecting brownian motions and gaussian multiplicative chaos,
  2025.
\newblock arXiv:2508.11505.

\bibitem[Kiv24]{Kivimae24}
Pax Kivimae.
\newblock Gaussian multiplicative chaos for {G}aussian orthogonal and
  symplectic ensembles.
\newblock {\em Electron. J. Probab.}, 29, 2024.

\bibitem[KN04]{KillipNenciu}
Rowan Killip and Irina Nenciu.
\newblock Matrix models for circular ensembles.
\newblock {\em Int. Math. Res. Not.}, (50):2665--2701, 2004.

\bibitem[KW22]{KeatingWong22}
Jonathan~P. Keating and Mo~Dick Wong.
\newblock On the critical--subcritical moments of moments of random
  characteristic polynomials: a {GMC} perspective.
\newblock {\em Comm. Math. Phys.}, 394(3):1247--1301, 2022.

\bibitem[Lac24]{Lacoin}
Hubert Lacoin.
\newblock Critical {G}aussian multiplicative chaos revisited.
\newblock {\em Ann. Inst. Henri Poincar\'e{} Probab. Stat.}, 60(4):2328--2351,
  2024.

\bibitem[Lam20]{Lambert20Ginibre}
Gaultier Lambert.
\newblock Maximum of the characteristic polynomial of the {G}inibre ensemble.
\newblock {\em Comm. Math. Phys.}, 378(2):943--985, 2020.

\bibitem[LLZ24]{LLZ}
Gaultier Lambert, Thomas Lebl\'e, and Ofer Zeitouni.
\newblock Law of large numbers for the maximum of the two-dimensional {C}oulomb
  gas potential.
\newblock {\em Electron. J. Probab.}, 29:Paper No. 42, 36, 2024.

\bibitem[LN26]{lambert2024subcritical}
Gaultier Lambert and Joseph Najnudel.
\newblock Subcritical multiplicative chaos and the characteristic polynomial of
  the {C$\beta$E}.
\newblock {\em Probab. Theory Related Fields}, 196(3-4):2103--2154, 2026.

\bibitem[LOS18]{LOS18}
Gaultier Lambert, Dmitry Ostrovsky, and Nick Simm.
\newblock Subcritical multiplicative chaos for regularized counting statistics
  from random matrix theory.
\newblock {\em Comm. Math. Phys.}, 360(1):1--54, 2018.

\bibitem[LP19]{LambertPaquette01}
Gaultier Lambert and Elliot Paquette.
\newblock The law of large numbers for the maximum of almost {G}aussian
  log-correlated fields coming from random matrices.
\newblock {\em Probab. Theory Related Fields}, 173(1-2):157--209, 2019.

\bibitem[LP20]{LambertPaquetteEdge}
Gaultier Lambert and Elliot Paquette.
\newblock Strong approximation of {G}aussian {$\beta$}-ensemble characteristic
  polynomials: the edge regime and the stochastic {A}iry function, 2020.
\newblock arXiv:2009.05003.

\bibitem[LP23]{LambertPaquetteHyperbolic}
Gaultier Lambert and Elliot Paquette.
\newblock Strong approximation of {G}aussian {$\beta$} ensemble characteristic
  polynomials: the hyperbolic regime.
\newblock {\em Ann. Appl. Probab.}, 33(1):549--612, 2023.

\bibitem[LS87]{LS}
S.~P. Lalley and T.~Sellke.
\newblock A conditional limit theorem for the frontier of a {Branching
  {B}rownian Motion}.
\newblock {\em Ann. Probab.}, 15(3):1052--1061, 1987.

\bibitem[Mad15]{Madaule}
Thomas Madaule.
\newblock Maximum of a log-correlated {G}aussian field.
\newblock {\em Ann. Inst. Henri Poincar\'e Probab. Stat.}, 51(4):1369--1431,
  2015.

\bibitem[Mad17]{Madaule17}
Thomas Madaule.
\newblock Convergence in law for the branching random walk seen from its tip.
\newblock {\em J. Theoret. Probab.}, 30(1):27--63, 2017.

\bibitem[Naj18]{Najnudel}
Joseph Najnudel.
\newblock On the extreme values of the {R}iemann zeta function on random
  intervals of the critical line.
\newblock {\em Probab. Theory Related Fields}, 172(1-2):387--452, 2018.

\bibitem[NPS23]{NPS23}
Joseph Najnudel, Elliot Paquette, and Nick Simm.
\newblock Secular coefficients and the holomorphic multiplicative chaos.
\newblock {\em Ann. Probab.}, 51(4):1193--1248, 2023.

\bibitem[NSW20]{NSW20}
Miika Nikula, Eero Saksman, and Christian Webb.
\newblock Multiplicative chaos and the characteristic polynomial of the {CUE}:
  the {$L^1$}-phase.
\newblock {\em Trans. Amer. Math. Soc.}, 373(6):3905--3965, 2020.

\bibitem[Oui26a]{OuimetMax26}
Fr\'{e}d\'{e}ric Ouimet.
\newblock Convergence in law of the maximum of the {R}iemann zeta function on
  expanding intervals near the critical line, 2026.
\newblock arXiv:2611.01905.

\bibitem[Oui26b]{OuimetHaar26}
Fr\'{e}d\'{e}ric Ouimet.
\newblock Overlap distributions under the {G}ibbs measure of a {H}aar model of
  $\log|\zeta|$ on and near the critical line, 2026.
\newblock arXiv:2610.14901.

\bibitem[Oui26c]{OuimetPD26}
Fr\'{e}d\'{e}ric Ouimet.
\newblock Poisson--{D}irichlet statistics for the log-modulus of the {R}iemann
  zeta function on expanding intervals near the critical line, 2026.
\newblock arXiv:2610.22821.

\bibitem[Pow18]{Powell18}
Ellen Powell.
\newblock Critical {G}aussian chaos: convergence and uniqueness in the
  derivative normalisation.
\newblock {\em Electron. J. Probab.}, 23:1--26, 2018.

\bibitem[Pow21]{PowellReview}
Ellen Powell.
\newblock Critical {G}aussian multiplicative chaos: a review.
\newblock {\em Markov Process. Related Fields}, 27(4):557--606, 2021.
\newblock arXiv:2006.13767.

\bibitem[PZ18]{PaquetteZeitouni02}
Elliot Paquette and Ofer Zeitouni.
\newblock The maximum of the {CUE} field.
\newblock {\em Int. Math. Res. Not. IMRN}, (16):5028--5119, 2018.

\bibitem[PZ25]{PaquetteZeitouni03}
Elliot Paquette and Ofer Zeitouni.
\newblock The extremal landscape for the {C$\beta$E} ensemble.
\newblock {\em Forum Math. Sigma}, 13:Paper No. e1, 2025.

\bibitem[Rem20]{Remy}
Guillaume Remy.
\newblock The {F}yodorov-{B}ouchaud formula and {L}iouville conformal field
  theory.
\newblock {\em Duke Math. J.}, 169(1):177--211, 2020.

\bibitem[RZ20]{RemyZhu}
Guillaume Remy and Tunan Zhu.
\newblock The distribution of {G}aussian multiplicative chaos on the unit
  interval.
\newblock {\em Ann. Probab.}, 48(2):872--915, 2020.

\bibitem[Sim04]{SimonBook}
Barry Simon.
\newblock {\em Orthogonal Polynomials on the Unit Circle, Parts {I, II}}.
\newblock American Math. Society, Providence, RI, 2004.

\bibitem[SW20]{SaksmanWebb20}
Eero Saksman and Christian Webb.
\newblock The {R}iemann zeta function and {G}aussian multiplicative chaos:
  statistics on the critical line.
\newblock {\em Ann. Probab.}, 48(6):2680--2754, 2020.

\bibitem[SZ15]{SZ}
Eliran Subag and Ofer Zeitouni.
\newblock Freezing and decorated {P}oisson point processes.
\newblock {\em Comm. Math. Phys.}, 337(1):55--92, 2015.

\bibitem[Web15]{Webb15}
Christian Webb.
\newblock The characteristic polynomial of a random unitary matrix and
  {G}aussian multiplicative chaos---the {$L^2$}-phase.
\newblock {\em Electron. J. Probab.}, 20, 2015.

\end{thebibliography}

\end{document}